# Non-convergent sequences of solutions to the massive Vafa-Witten equations with 'interesting' $\mathbb{Z}/2\mathbb{Z}$ self-dual harmonic 2-form limits

Clifford Henry Taubes

Department of Mathematics

Harvard University

Cambridge, MA 02138

chtaubes@math.harvard.edu

ABSTRACT: This paper constructs sequences of solutions to the Vafa-Witten equations with non-zero (but small) mass term on the product of a 2-dimensional torus with a Riemann surface of genus greater than 1. These are divergent sequences (modulo principle bundle automorphisms) that converge after renormalization to define an 'interesting' $\mathbb{Z}/2\mathbb{Z}$ harmonic 2-form data set. This data set consists of a non-empty, codimension 2 submanifold, a real line bundle defined on the complement of that submanifold with no extension across it, and a self-dual, harmonic 2-form with values in that line bundle that does not extend over the submanifold. Even so, the norm of this 2-form does extend as a Hölder continuous function with that submanifold being its zero locus.

## 1. Introduction

The SO(3) Vafa-Witten equations were introduced many years ago in a paper by Cumrun Vafa and Edward Witten [VW].  They constitute a system of equations for a connection on the principle bundle (denoted by A below) and a self-dual 2-form with values in the principle bundle's associated Lie algebra bundle (this 2-form is denoted by ω).  The *massive* Vafa-Witten equations are a slight riff on the version in [VW].  The Vafa-Witten equations with the mass term ask the pair (A, ω) to obey the system of equations below:

- $F_A{}^+ - \frac{1}{2} [\omega; \omega] - m\omega = 0$ ,

- $d_A\omega = 0$ ;

$$(1.1)$$

The notation here is as follow:  What is denoted by $F_A$ signifies A's curvature 2-form; and $F_A{}^+$ signifies the self-dual part of $F_A$.  What is denoted by $[\omega; \omega]$ denotes a self-dual 2-form with values in the principle bundle's associated Lie-algebra bundle whose coefficients are linear combinations of the commutators of the coefficients of ω.  (See the upcoming (1.3) for a precise definition.)  Meanwhile, $m$ is a real number (the mass).  Finally, $d_A\omega$ denotes the exterior, covariant derivative of ω.

Solutions to these equations in the case when $m = 0$ (the original Vafa-Witten equations) can be constructed on Kähler 4-manifolds from algebraic geometric data (see e.g. [Tan], [TT1], [TT2], [Ch]); and solutions on the product of the circle with a compact, oriented 3-manifold can be constructed directly from a (stable) flat $SL(2; \mathbb{C})/\{\pm 1\}$ connection on the 3-manifold (see below in Section 1c for more about this).  Of course, any connection on the principle bundle with anti-self-dual curvature tensor is a solution to the $m = 0$ equations with $\omega \equiv 0$.  Solutions in the case when $m$ is not zero are harder to come by except for what might be termed *reducible* solutions (see e.g. [T1] for some examples).

The $m \neq 0$ modification to the original equations are of interest in part because these $m \neq 0$ solutions are rare.  In particular there was some hope that the moduli space of solutions to the $m \neq 0$ equations are 'less' non-compact than the moduli space for the $m = 0$ equations.  To elaborate:  The equations in (1.1) are invariant under the action of the group of automorphisms of the principle bundle (this is the group of gauge transformations).  This being an infinite dimensional group, there are infinitely many distinct solutions if there is but one.  Even so, the quotient of the solution space by this group action is (formally) finite dimensional–this quotient is called the moduli space of solutions.  The question about compactness (or not) asks whether sequences in this moduli space necessarily have convergent subsequences.

In this regard:  If there are non-convergent sequences in the moduli space of solutions to (1.1), then their limits can be geometrically characterized in the manner described in [T2] (this is true for any value of $m$).  In particular, there are two sorts of non-compactness. The first sort occurs when there is an upper bound along the sequence for the integral over the 4-manifold of $|\omega|^2$. Convergence for a subsequence in this case occurs on the complement of finitely many points in the manifold. The non-compactness here is essentially the same as that for the moduli space of solutions to the anti-self dual equations. See [Mar]. The second sort only occurs when there is no upper bound along the sequence to the integral of $|\omega|^2$. Convergence in this case (for a suitably renormalized subsequence) occurs on the complement of a closed set in the manifold with finite 2-dimensional Hausdorff measure (see [T2], [Zh], [Ch], [P1]). The limiting data for a suitable subsequence in this second case consists of what is called a $\mathbb{Z}/2$ *harmonic, self-dual 2-form*. Such an object consists of a data set $(Z, \mathcal{I}, \nu)$ with Z denoting the closed set in question, with $\mathcal{I}$ denoting a real line bundle (with fiber metric) defined on the complement of Z, and with $\nu$ denoting a self-dual, $\mathcal{I}$-valued 2-form defined on the complement of Z whose norm extends as 0 over Z as a Hölder continuous function on the manifold.

Both sorts of non-compactness are known to occur for solutions to the $m = 0$ version of (1.1). (For example, the $Z \neq 0$ version occurs in the $m = 0$ case on the product of a circle with certain 3-manifolds when the moduli space of $SL(2; \mathbb{C})$ flat connections is suitably non-compact, the product of a circle and a Riemann surface of genus greater than 1 being a relevant example.) Examples of non-convergent sequences of solutions when $m \neq 0$ of the second sort where the corresponding limit line bundle $\mathcal{I}$ does extends over the whole manifold are described in [T1]. To the author's knowledge, there were no examples of the second sort when $m \neq 0$ where the limit line bundle $\mathcal{I}$ has no extension over one or more components of the set Z. The purpose of this paper is to show that these non-extendable limits do occur for some fixed 4-manifolds, and for any given small but still non-zero value of $m$. In particular, the manifolds in question are products of $S^1 \times S^1$ with its product metric and a Riemann surface of genus greater than 1. Below is a formal statement.

**Theorem 1**: *Given a compact Riemann surface of genus greater than* 1 *and a small, non-zero value for* m, *there is a non-convergent sequence of solutions to the corresponding version of (1.1) on the product principle* $SU(2)$ *bundle over the product of* $S^1 \times S^1$ *with that surface such that the integral of* $|\omega|^2$ *diverges along the sequence. However, after renormalization in the manner of* [T2], *that sequence converges to a data set* $(Z, \mathcal{I}, \nu)$ *as described above with* Z *being the product of* $S^1 \times S^1$ *and a non-empty finite set in the Riemann surface, with* $\mathcal{I}$ *denoting a real line bundle defined on the complement of* Z *that does not extend over* Z, *and with* $\nu$ *denoting a nowhere zero,* $\mathcal{I}$-valued self-dual 2-form on *the complement of* Z *that does not extend over* Z *not-with-standing the fact that its norm* *does extends over* $\mathbb{Z}$ *as zero to define a Hölder continuous function (with exponent* $\frac{1}{2}$*) on the product 4-manifold.*

The strategy for constructing these non-convergent sequences is along the lines pioneered by Greg Parker in [P1]-[P4] and in upcoming unpublished work. However, the construction here is far simpler because of the geometric set-up with the manifold in question being the product of a flat torus with a Riemann surface and the set Z being the product of that torus with a finite set in the surface. Because of the product geometry, there is no need for the full suite of technological tools from Parker's work. The construction is none-the-less not straight forward by virtue of the appearance of (potentially) non-trivial obstructions. The argument that these obstructions can be circumvented requires (in part) a detailed analysis of the solutions to the equations that characterize elements in the tangent space to the moduli space of SL(2; $\mathbb{C}$) flat connections on the Riemann surface.

What follows is a list of the various sections of this paper.

1. *Introduction*
2. *The singular limit data* ($\mathbb{Z}$, $\mathfrak{I}$, $\nu$) *and singular non-convergent solution sequences*
3. *Approximate solutions*
4. *The perturbation problem and the associated linear operator.*
5. *The* $m \neq 0$ *equations on each Riemann surface slice*
6. *Putting the slices together to solve the* $m \neq 0$ *Vafa-Witten equations.*
7. *When* $\phi \to \Pi\mathfrak{k}(\phi + \mathfrak{T}_\lambda(\phi))$ *has eigenvalue* $\lambda$ *and* $|\lambda|$ *is small.*
8. *Spectral flow*

*Appendix*

*References.*

What follows elaborates on what is in these sections.

The rest of this first section describes the notation used in this paper and it also specifies various conventions that are used throughout. The section ends with a subsection that elaborates on (1.1) in the case when the 4-manifold is the product of a circle and a 3-dimensional manifold. As explained there, the $m = 0$ equations in (1.1) can be viewed as the equations for a flat SL(2; $\mathbb{C}$) connection on the 3-manifold. The $m \neq 0$ equations are then deformations of the flat SL(2; $\mathbb{C}$) connection equations. This is the relevant case for the body of this paper with the 3-manifold being the product of a circle and a Riemann surface of genus greater than 1. The Riemann surface in question is denoted by $\Sigma$, and its genus is denoted by g. The reduction of the equations in (1.1) to a system of equations on the 3-manifold $S^1 \times \Sigma$ is exploited throughout the body of this paper.

The second section defines the singular limit data $(Z, \mathcal{I}, \nu)$ for the case of interest. This section also defines a sequence of singular solutions to any given $m$ version of (1.1) that converges after renormalization in the manner of [T2] to this singular data. (The solutions in this sequence are smooth outside of Z but singular across this set.) An important point for the future: The singular limit set $(Z, \mathcal{I}, \nu)$ is invariant with respect to the two rotations of the $S^1 \times S^1$ factor of the underlying 4-manifold $S^1 \times S^1 \times \Sigma$, but this is not the case for the members of the singular solution sequence. Each member is invariant only with respect to rotations of the left-most $S^1$-factor (except in the case that $m = 0$). And, each member is far from invariant with respect to the right most factor.

The third section constructs a sequence of smooth, approximate solutions to any given $m$ version of (1.1) which converges in the manner of [T2] to the limit data set $(Z, \mathcal{I}, \nu)$. This is done by replacing the members of the singular sequence with data obtained from a sequence of $SL(2; \mathbb{R})$ flat connections on $\Sigma$ that likewise converges (as described in [T2]) to the limit data set. The use of the $SL(2; \mathbb{R})$ sequence to obtain smooth, but only approximate solutions to (1.1) mimics what is done by Greg Parker in [P1]-[P4] (and that, in turn, is based on work done in [MSWW] which treats some $m = 0$ analogs of (1.1) on Riemann surfaces. See also [F], [FMSW].)

The plan going forward, following the strategy of Greg Parker (see also the aforementioned [MSWW]) is to use perturbation theoretic techniques to turn the sequence of approximate solutions to (1.1) into a sequence of honest solutions to (1.1). This perturbation problem is described at the outset of Section 4. Of relevance here are certain properties of an elliptic, first order, self-adjoint differential operator that is obtained from each approximate solution in the approximation sequence. Section 4 introduces this operator and derives some of its salient properties (Parker's opus [P3] has very similar theorems in a more general context).

Section 5 takes a preliminary step towards turning the approximate solution sequence from Section 3 into a true solution sequence. The analysis in this section 'pretends' that the construction is invariant under all rotations of the $S^1 \times S^1$ factor (not just the left most factor). Doing this in conjunction with the analysis in Section 4 supplies a new sequence of approximate solutions which are (in a quantifiable sense) exponentially closer to being honest solutions to (1.1) when $m$ is small. What is done in Section 5 is equivalent (in part) to honestly solving the equations on $S^1 \times \mathbb{R} \times \Sigma$ for solutions that are invariant with respect to rotations of the left most $S^1$ factor and translations along the $\mathbb{R}$ factor.

By way of an elaboration on the need for Section 5's new approximations: As it turns out, the operator from Section 4 has 12g - 12 eigenvectors (counted with multiplicity over $\mathbb{R}$) with very small normed eigenvalue. In the $m = 0$ case, these eigenvectors have eigenvalue zero and their span is the tangent space to the $SL(2; \mathbb{C})$ flat connection moduli space. When $m$ is non-zero but still small, 6g - 6 of these eigenvalues move away from zero, they have norm on the order of $m$. The other 6g - 6 still have norm much, much less than $m$. When Section 5's more accurate approximation is used, then the 6g - 6 eigenvalues with $\mathcal{O}(m)$ norm are too large to obstruct the deformation to a true solution along each $\Sigma$ slice of $S^1 \times \Sigma$; and the eigenvectors for those very much smaller normed eigenvalues are orthogonal to the deformation (so they automatically don't obstruct). However, this family of $S^1$–parametrized family of solutions cannot be viewed as a solution on $S^1 \times \Sigma$. Even so, certain a priori bounds from Section 4 imply that a slight deformation of this $S^1$ parametrized family of solutions along the slices does define a configuration on $S^1 \times \Sigma$ that is exponentially close to a solution there. This exponential factor has the form $e^{-r/\text{constant}}$ with $r$ being a positive number which is a proxy for the square root of the integral over the manifold of $|\omega|^2$ with $\omega$ being the endomorphism valued, self-dual 2-form component of the given approximate solution. (Remember that those integrals are diverging along the approximation solution sequence.)

Section 6a then takes the exponentially close approximate solutions from Section 5 as input and constructs (again, using the analysis in Section 4) the desired sequence of honest solutions to any small $m$ version of (1.1) given a proposition (Proposition 6.2) which asserts that those eigenvalues which are very much smaller than $o(m)$ do not have exponentially small norm, which is to say that their norms are not less than $e^{-r/\text{constant}}$ for large $r$. Indeed, Proposition 6.2 asserts that the norms of these 6g - 6 very small eigenvalues are all greater than a universal constant times $mr^{-2}$. Results and lemmas in Section 6 subsequent to Section 6a take preliminary steps towards the proof Proposition 6.2 (these subsequent lemmas have antecedents in Greg Parker's paper [P1]).

Section 7 finishes the proof of Proposition 6.2 with some very delicate estimates using the specific representations from Lemma 5.2 of the eigenvectors that are in the kernel of the linear operator that is defined by the $SL(2; \mathbb{C})$ flat connections on the Riemann surface from Section 3. In particular, this analysis shows that certain integrals over the Riemann surface that are defined by these kernel elements do not sum to zero. (Supposing that there is a fatal error in this paper (!), it would likely be in Section 7 because dropping a factor of 2 or -1 would make a difference.)

Section 8 studies the spectral flow for the linearized operator defined by any given pair in the diverging solution sequence that is constructed in the previous sections (thus, Theorem 1's solutions). This spectral flow is relevant by virtue of its interpretation as the relative Morse index between two critical points of a functional of the pairs of connection and endomorphism valued, self-dual 2-form (the equations in (1.1) are the variational equations of this functional). Section 8 states and then proves a proposition to the effect that the norm of this spectral flow is a priori bounded along the sequence of solutions constructed in the previous sections. (See [T1] for more about the spectral flow for the massive Vafa-Witten equations.)

The appendix to this paper gives full proofs of some analytical assertions that can be derived more or less directly from what is said in [MSWW] which concerns various $m = 0$ analogs (1.1) on Riemann surfaces. These analytic assertions constitute Lemmas 3.1-3.3 which are about SL(2; $\mathbb{R}$) flat connections on a Riemann surface, and also Lemma 5.2 which describes in some detail the eigenvectors that define the kernel of the deformation operator on the surface that are defined by those flat connections. I expect that experts will recognize the assertions in these lemmas; and I apologize in advance for not simply referencing specific but more general theorems.

In a similar vein: This paper is self-contained for the most part to provide the intricate details of arguments and proofs which, in some cases, might be found as instances or variations of results in the extant literature. I did this for my own benefit, but also (I hope) for the benefit of the reader who otherwise would have to work through lengthy references and their results written in a much more general and/or somewhat different context (and using different notation) to justify the very special instances of those general theorems that might replace some of the specialized lemmas in this paper.

## a)  Notation

To elaborate on the notation in (1.1): Suppose that X is an oriented manifold with a Riemannian metric. That metric has a corresponding Hodge star operator (denoted by $*$). In the case when X is 4-dimensional, the Hodge star squares to the identity on the bundle of 2-forms and thus it defines orthogonal spliting of this bundle as $\Lambda^+ \oplus \Lambda^-$ with $\Lambda^+$ denoting the bundle of self-dual 2-forms. In particular, if $\{e^0, e^1, e^2, e^3\}$ is an oriented, orthonormal basis for T*X at any given point, then the 2-forms depicted below define an oriented basis for $\Lambda^+$ at that same point:

$$\eta^1 = \tfrac{1}{\sqrt{2}}\,(e^0 \wedge e^1 + e^2 \wedge e^3), \quad \eta^1 = \tfrac{1}{\sqrt{2}}\,(e^0 \wedge e^2 + e^3 \wedge e^1), \ \ \eta^3 = \tfrac{1}{\sqrt{2}}\,(e^0 \wedge e^3 + e^1 \wedge e^2)$$

$$(1.2)$$

Let P $\to$ X denote a principle SO(3) bundle. Supposing that A denotes a connection on this bundle, let $F_A$ denote its curvature 2-form. The self-dual part of $F_A$ can be written using the orthonormal frame in (1.2) as $F_A{}^+{}_k\, \eta^k$. Here and subsequently, repeated indices are implicitly summed over their range of integer values. Now, if $\omega$ is a given self-dual 2-form on X with values in the associated Lie algebra bundle to P, write $\omega$ as $\omega_k\, \eta^k$ also. The top bullet equation in (1.1) can be written now in terms of the components of $F_A$ and $\omega$ as

$$(F_A)^+{}_k - \tfrac{1}{2\sqrt{2}}\,\varepsilon_{ijk}[\omega_i, \omega_j] - m\,\omega_k = 0$$

$$(1.3)$$

with $\{\varepsilon_{ijk}\}_{i,j,k \in \{1,2,3\}}$ denoting the components of the completely antisymmetric 3-tensor with $\varepsilon_{123}$ being equal to 1. More generally, if $\omega$ and $\omega'$ are self-dual 2-forms with values in P's associated Lie algebra bundle, then [ω; ω´] has components with respect to the frame in (1.2) given by the rule

$$[\omega; \omega´]_k = \tfrac{1}{\sqrt{2}} \, \varepsilon_{ijk} \, [\omega_i, \omega´_j] \, .$$

(1.4)

As with the preceding equations, most equations below will be written out using a locally defined, oriented orthonormal frame for T*X (a frame $\{e^0, e^1, e^2, e^3\}$), and then the corresponding frame for $\Lambda^+$. Repeated frame indices in these equations for the T*X frame are to be summed, and likewise for the repeated indices of the $\Lambda^+$ frame. By way of an example, the Hodge $*$ of $d_A\omega$ (which appears in the second bullet of (1.1)) when written using the given local orthonormal frame has components $\{(*d_A\omega)_\alpha\}_{\alpha=0,1,2,3}$ that are given by the rule

$$(*d_A\omega)_\alpha = \eta^k{}_{\alpha\beta}(\nabla_{A\beta}\,\omega)_k$$

(1.5)

where the notation is as follows: Any given $k \in \{1, 2, 3\}$ version of $\{\eta^k{}_{\alpha\beta}\}_{\alpha,\beta=0,1,2,3}$ are the components of the self-dual 2-form $\eta^k$ in (1.2); and where $\{\nabla_{A\beta}\}_{\beta=0,1,2,3}$ are the directional covariant derivatives along the corresponding TX frame vectors that are defined by the connection A and the metric's Levi-Civita connection. The implicit sum in (1.5) is over the indices k (summed from 1 to 3) and the indices labeled by β (summed from 0 to 3).

In general, the covariant derivative that is defined by a given connection A on a given vector bundle is denoted by $\nabla_A$. This same notation is used for the covariant derivative on vector bundle valued differential forms as defined by the connection A and the metric's Levi-Civita acting. (The Levi-Civita connection's covariant derivative acting on $\mathbb{R}$ or $\mathbb{C}$ valued differential forms is denoted by $\nabla$.)

With regards to P's associated Lie algebra bundle: The Lie algebra of SO(3) is isomorphic to the Lie algebra of SU(2) which is the vector space of $\mathbb{C}$-valued, $2 \times 2$, anti-Hermitian matrices. This vector space is denoted by $\mathfrak{su}(2)$ in what follows. The associated Lie algebra bundle to P has $\mathfrak{su}(2)$ as its fiber. This bundle is denoted below by ad(P). The vector space $\mathfrak{su}(2)$ is given its adjoint action invariant inner product which assigns to matrices σ and σ´ the number

$$\langle \sigma\sigma´ \rangle = -\tfrac{1}{2}\, \text{trace}(\sigma\sigma´) \, .$$

(1.5)

This inner product induces an inner product on ad(P) which is denoted by $\langle\,,\,\rangle$ and, with the Riemannian metric, an inner product any bundle of ad(P)-valued tensors. That inner product is also denoted by $\langle\,,\,\rangle$.

## b) Conventions

The 4-manifold of interest in this paper is the product $S^1 \times S^1 \times \Sigma$. In this regard, the torus $S^1 \times S^1$ is defined so that each $S^1$ factor has length 2π–each such factor is identified with $\mathbb{R}/2\pi\mathbb{Z}$; and the coordinate on $\mathbb{R}$ then defines an $\mathbb{R}/2\pi\mathbb{Z}$ valued coordinate on $S^1$. The metric on these $S^1$ factors is chosen so that the differentials of these $\mathbb{R}/2\pi\mathbb{Z}$ coordinates have norm one.

Subsequent sections of this paper will view the 4-manifold $S^1 \times S^1 \times \Sigma$ as the product of the left most $S^1$ with the 3-manifold $S^1 \times \Sigma$ and then focus on $S^1 \times \Sigma$ using the relation that is described in the next subsection between the 4-manifold equations in (1.1) and the 3-manifold equations depicted momentarily in (1.7). With regards to $S^1 \times \Sigma$: The $\mathbb{R}/2\pi\mathbb{Z}$ coordinate for the $S^1$ factor is denoted by t and metric is the product of the metric dt $\otimes$ dt on the $S^1$ factor and a given metric on $\Sigma$.

As for the surface $\Sigma$: This is assumed to have genus greater than 1. It is also assumed to have a Riemannian metric. This given metric then defines a complex structure for $\Sigma$ and hence a splitting of $T^*\Sigma \otimes_{\mathbb{R}} \mathbb{C}$ as $T^{1,0} \oplus T^{0,1}$ with the left most factor denoting the complex line bundle which is spanned near any given point in $\Sigma$ by the exterior derivative of a holomorphic coordinate for a neighborhood of that point. (When z denotes the coordinate, then its differential is denoted by dz). What is denoted in what follows by $T^{2,0}\Sigma$ signifies the holomorphic line bundle $T^{1,0} \otimes_{\mathbb{C}} T^{1,0}$. A holomorphic quadratic differential is a holomorphic section of $T^{2,0}\Sigma$. Note that the vector space over $\mathbb{C}$ of holomorphic quadratic differentials has dimension 3g - 3.

With regards to notational conventions: This paper will use $c_0$ to denote a number greater than 1 that is independent of any relevant parameter in the discussion (for example, $m$, or a given pair of connection on a principal bundle and adjoint bundle valued differential form). However, any given incarnation of $c_0$ is allowed to depend on the chosen metric for $\Sigma$ since that will be fixed throughout. As for these incarnations: The precise value of $c_0$ will increase between success incarnations.

With regards to 'bump' functions: Various constructions in what follows use functions with specified support on $\Sigma$ or $S^1 \times \Sigma$. These are constructed from a 'fiducial' smooth function on $\mathbb{R}$ which is denoted by $\chi$. The key properties of $\chi$ are these: This function is non-increasing, it equals 1 on $(-\infty, \frac{1}{4}]$ and it equals zero on $[\frac{3}{4}, \infty)$. (The use of a fiducial $\chi$ leads directly to a priori bounds for the norms of the derivatives of the resulting bump functions on the manifold in question.)

With regards to the formal adjoint of the covariant derivatives: This paper uses $d^\dagger$ to denote the formal adjoint of the (exterior) derivative operator mapping functions to 1-forms (whether on $S^1 \times \Sigma$ of just on $\Sigma$). Thus, $d^\dagger$ maps 1-forms to functions (when a 1-form, $\nu$, is written with respect to a locally defined orthonormal frame as $\nu_\alpha e^\alpha$, then $d^\dagger\nu$ can be written in that same orthonormal frame as $-\nabla_\alpha \nu^\alpha$). It is defined so that $d^\dagger d$ is the (spectrally) non-negative Laplacian. (Thus, on Euclidean space, $d^\dagger d$ is -1 times the sum of the squares of the partial derivatives in the coordinate directions.) By the same token, the formal adjoint of the covariant derivative defined by a connection A on a vector bundle is denoted by $\nabla_A{}^\dagger$; it maps bundle valued 1-forms to sections of the given vector bundle.

With regards to $d^\dagger d$:  If $f$ denotes a smooth function or section of a vector bundle, then $|f|$ is smooth where $f \neq 0$ but it is only Lipschitz along $f$'s zero locus.  As a consequence, an inequality of the form $d^\dagger d|f| \leq g$ for some given, continuous function g does not make literal sense where $f = 0$.  However, such an inequality does make sense as a distributional inequality.  Equivalently, the inequality $d^\dagger d|f| \leq g$ should be viewed as the $\varepsilon \to 0$ limit of a sequence of inequalities of the form $d^\dagger d(\varepsilon^2 + |f|^2)^{1/2} \leq g_\varepsilon$ where each $g_\varepsilon$ is continuous and these are such that the $\varepsilon \to 0$ limit of the $g_\varepsilon$'s is g.  This interpretation will be implicit (unstated) whenever an inequality of the form $d^\dagger d|f| \leq g$ is written.

### c)  Sl(2;$\mathbb{C}$) flat connections and dimensional reduction

Suppose for the moment that Y denotes a compact, oriented Riemannian 3-manifold.  Let P denote a principle SO(3) bundle over Y.  The group SO(3) is a subgroup of $SL(2; \mathbb{C})/\{\pm 1\}$ and, as such, there is the associated principle $Sl(2; \mathbb{C})/\{\pm 1\}$ bundle.  A connection on the latter bundle can be written as $A + i\mathfrak{a}$ with A being a connection on P and with $\mathfrak{a}$ denoting an ad(P)-valued 1-form.  This connection is flat when $(A, \mathfrak{a})$ obey

$$F_A - \mathfrak{a} \wedge \mathfrak{a} = 0 \quad and \quad d_A\mathfrak{a} = 0 \ .$$

$$(1.6)$$

The connection can be said to be *stable* when $d_A*\mathfrak{a} = 0$ also.  These equations admit an 1-parameter family of deformations, indexed by the real numbers whose $m \in \mathbb{R}$ member requires that $(A, \mathfrak{a})$ obey

- $F_A - \mathfrak{a} \wedge \mathfrak{a} - m*\mathfrak{a} = 0$
- $d_A\mathfrak{a} = 0 \quad and \quad d_A*\mathfrak{a} = 0.$

$$(1.7)$$

As explained directly, a solution to the latter equations on Y defines a solution to (1.1) on $S^1 \times Y$ with it understood that the version of $m$ in (1.1) is $\sqrt{2}$ times the version of $m$ in (1.7).

To start the explanation:  With it understood that Y comes with a Riemannian metric (denoted by $\mathfrak{g}$) and that the circle $S^1$ has its $\mathbb{R}/2\pi\mathbb{Z}$ coordinate s, then $S^1 \times Y$ has its product metric $ds \otimes ds + \mathfrak{g}$.  Supposing that $\{e^1, e^2, e^3\}$ is a local, oriented orthonormal frame for $T^*Y$, then the frame $\{e^0 = ds, e^1, e^2, e^3\}$ is an orthonormal frame for $T^*(S^1 \times Y)$.  Meanwhile, the bundle $\wedge^2 T^*Y$ has the orthonormal frame $\{\frac{1}{2} \varepsilon_{ijk} e^i \wedge e^j\}_{k=1,2,3}$.  This frame with the corresponding $S^1 \times Y$ Hodge star dual 2-forms $\{ds \wedge e^k\}_{k=1,2,3}$ define an orthonormal frame for $\wedge^2 TY$.  Note in this regard that the self-dual projection of any given $k \in \{1, 2, 3\}$ version of either $\frac{1}{2} \varepsilon_{ijk} e^i \wedge e^j$ or its Hodge $*$ dual is the self-dual 2-form $\frac{1}{\sqrt{2}} \eta^k$ with $\eta^k$ as defined in (1.2) using the just noted frame for $T^*(S^1 \times Y)$.

Suppose now that P denotes a principle SO(3) bundle over Y.  Projection to Y from $S^1 \times Y$ pulls the bundle P back to the latter space.  That pull-back is also denoted by P.

Likewise, its adjoint bundle on Y and the pull-back to $S^1 \times Y$ are denoted by ad(P). If A denotes a connection on P defined over Y, then pull-back by the projection defines a connection on the incarnation of P over $S^1 \times Y$ which will also be denoted by A. Having written the connection A's curvature on Y using a given local, oriented orthonormal frame $\{e^1, e^2, e^3\}$ as $F_A = F_{Ak} \frac{1}{2} \varepsilon_{ijk} e^i \wedge e^j$, then the components of the self-dual part of the pull-back of A's curvature to $S^1 \times Y$ when written using the basis in (1.2) are the k = 1, 2, 3 versions of

$$(F_A)^+_k = \frac{1}{\sqrt{2}} F_{Ak} \, .$$

(1.8)

To continue, let $\mathfrak{a} = \mathfrak{a}_k e^k$ denote an ad(P)-valued 1-form on Y. Given this ad(P)-valued 1-form on Y define a corresponding ad(P)-valued self dual 2-form on $S^1 \times Y$ using the same $\{\eta^k\}_{k=1,2,3}$ basis via the rule $\mathfrak{a} \to \mathfrak{a}_k \eta^k$. If (A, $\mathfrak{a}$) obeys (1.7) on Y, then the pair consisting of the pull-back of A to Y and $\omega = \mathfrak{a}_k \eta^k$ obeys

$$(F_A)^+_k - \frac{1}{2\sqrt{2}} \varepsilon_{ijk} [\omega_i, \omega_j] - \frac{1}{\sqrt{2}} m \, \omega_k = 0 \quad and \quad d_A \omega = 0 \, ,$$

(1.9)

which are the equations in (1.1) with the version of $m$ in the latter equation replaced by $\sqrt{2}$ times the version of $m$ in (1.7).

With the preceding understood: The plan for what follows is to focus on the equations in (1.7) for the most part.

## 2. The singular limit data (Z, $\mathcal{I}$, $\nu$) and singular, non-convergent solution sequences

The purpose of this section is to define the singular limit data and then a singular solution to (1.1) that converges to this data in the manner dictated in [T2].

### a) The singular limit data (Z, $\mathcal{I}$, $\nu$)

The singular limit data on a 3-manifold (call it Y for the moment) for a sequence of solutions to (1.7) with the Y-integral of $|\mathfrak{a}|^2$ diverging along the sequence consists of a data set (Z, $\mathcal{I}$, $\varsigma$) with Z being a closed set in Y with finite, 1-dimensional Hausdorff measure, with $\mathcal{I}$ being a real line bundle defined over Y–Z and with $\varsigma$ being an $\mathcal{I}$-valued, harmonic 1-form defined on Y–Z whose norm extends over Z as Hölder continuous function that vanishes on Z (see [T3]). If that sequence is used to define a corresponding sequence of solutions to (1.1) on $S^1 \times Y$ in the manner of Section 1c,, then the latter sequence has a corresponding limit data set as follows: The 4-manifold version of Z is the product of $S^1$ with the 3-manifold version of Z; and the 4-manifold version of $\mathcal{I}$ is the pull-back via the projection to Y of the 3-manifold version. As for the 4-manifold version of $\nu$, it is obtained from $\varsigma$ by first using a local, oriented orthonormal frame for T*Y together with the frame vector $e^0 = ds$ to define the frame $\{\eta^k\}_{k=1,2,3}$ for $\Lambda^+$ in (1.2) and then taking the k'th coefficient of $\nu$ for this frame to be the k'th coefficient of $\varsigma$ for the original 3-manifold frame. (Thus, $\nu = \varsigma_k\eta^k$.)

In the case at hand, the 3-manifold Y is $S^1 \times \Sigma$. The singular data here is defined with the help of a non-trivial, holomorphic quadratic differential on $\Sigma$. (This is a holomorphic section of $T^{2,0}\Sigma$.) This holomorphic, quadratics differential is chosen to have certain properties, one being that it must have only simple zeros (it vanish with degree 1 at each of its zeros. An observation of Bers [B] asserts that there is an open, dense set of holomorphic, quadratic differentials with this property. (This is asserted formally in the upcoming Lemma 2.1). Let $q$ denote one such holomorphic, quadratic differential. Let $\Theta$ denote the zero locus of $q$. Since the zeros of $q$ are first order and since the first Chern number of $T^{2,0}\Sigma$ is 4g - 4, there are 4g-4 zeros, so $\Theta$ is a set of 4g-4 points. The set Z when defined as a subset of $S^1 \times \Sigma$ is $S^1 \times \Theta$.

To define the line bundle $\mathcal{I}$, digress for the moment to define (as done in [DH]) a holomorphic submanifold $\Sigma_q \subset T^{1,0}\Sigma$ as follows: Denote a point in $T^{1,0}\Sigma$ as $(p, \varpi)$ with p being a point in $\Sigma$ and with $\varpi$ being a vector in $T^{1,0}\Sigma|_p$. This point is in $\Sigma_q$ if and only if $\varpi^2 = q$. The projection map from $\Sigma_q$ to $\Sigma$ defines $\Sigma_q$ as a holomorphic, 2-1 branched cover of $\Sigma$ with branch locus $\Theta$. (This projection map is denoted by $\pi_q$ in what follows.) A real line bundle over $\Sigma - \Theta$ is defined by taking the $\mathbb{R}$ span of the two points in $\Sigma_q - \Theta$ over any given point in $\Sigma - \Theta$. This line bundle is denoted by $\mathcal{I}_q$. The line bundle $\mathcal{I}$ over the 3-manifold $S^1 \times (\Sigma - \Theta)$ is the pull-back of $\mathcal{I}_q$ via the projection map to $\Sigma - \Theta$.

The last item is the definition of the $\mathcal{I}$-valued, harmonic, 1-form $\varsigma$. To define the latter, note first that the assignment to any given point $(p, \varpi)$ in $\Sigma_q$ of the element $\varpi$ defines a holomorphic section over $\Sigma_q$ of $T^{1,0}\Sigma_q$ to be denoted by $\varpi_q$. The set $\Theta$ (as a subset in $\Sigma_q$) is precisely the zero locus of $\varpi_q$; and $\varpi_q$ vanishes precisely to second order at each point in $\Theta$. (The Euler characteristic of $\Sigma_q$ is 8g - 8.) To say more about this holomorphic differential $\varpi_q$: The action of -1 on $T^{1,0}\Sigma$ by multiplication defines a holomorphic action of the group $\{1, -1\}$ on $\Sigma_q$ whose fixed point set is the branch locus $\Theta$. The pull-back of $\varpi_q$ by the multiplication by the -1 diffeomorphism of $\Sigma_q$ is $-\varpi_q$. It follows from this that $\varpi_q$ can be viewed from the perspective of $\Sigma - \Theta$ as a section of $T^{1,0}(\Sigma - \Theta) \otimes_{\mathbb{R}} \mathcal{I}_q$. The pull-back of this incarnation of $\varpi_q$ to $S^1 \times (\Sigma - \Theta)$ via the projection map is also denoted by $\varpi_q$. Both the real and imaginary parts of $\varpi_q$ are harmonic, self-dual, $\mathcal{I}$-valued 1-forms on $S^1 \times (\Sigma - \Theta)$. Let $\mathfrak{Re}(\varpi_q)$ denote the real part of $\varpi_q$, and let L denote the integral over $S^1 \times (\Sigma - \Theta)$ of the square of its norm. Then define $\varsigma$ to be the product of $L^{-1/2}\mathfrak{Re}(\varpi_q)$. Because $q$ vanishes only to first order along Z, and because $|\varsigma| \leq c_0 |q|^{1/2}$, the norm of $\varsigma$ extends over Z as a Hölder continuous function with Hölder exponent $\frac{1}{2}$.

By way of notation: That action of -1 as a holomorphic diffeomorphism of $\Sigma_{\hat{q}}$ is denoted by (-1) in what follows.

## b) A non-convergent sequence of singular solutions to the $m = 0$ version of (1.7)

The construction of a non-convergent sequence of singular solutions to (1.7) in the case when $m = 0$ is done here first as a warm-up of sorts for the $m \neq 0$ construction that is described in the next subsection. In this $m = 0$ case, the desired sequence is defined first on $\Sigma$ and then pulled up to $S^1 \times \Sigma$ via the projection map. In particular, the $S^1$ factor plays no role and it is not mentioned until the very end.

The construction has three parts:

*Part 1*: Let $q$ denote again a holomorphic, quadratic differential on $\Sigma$ with only simple zeros. The zero locus of $q$ is denoted again by $\Theta$. Break the set $\Theta$ into two disjoint sets of 2g-2 points. Denote these respective sets by $\Theta_+$ and $\Theta_-$. When viewed as a positive divisor on $\Sigma$, the set $\Theta_+$ defines a pair consisting of a holomorphic line bundle over $\Sigma$ to be denoted by $L_+$, and then a section of this line bundle whose zero locus is $\Theta_+$. (This section is uniquely defined up to multiplication by a non-zero complex number.) A choice for the section is denoted by $\varpi_+$ in what follows. Likewise, $\Theta_-$ defines a holomorphic line bundle $L_-$ and section, $\varpi_-$, whose zero locus is $\Theta_-$. Since $\varpi_+\varpi_-$ has the same zero locus as $q$, this product must be a complex multiple of $q$. As such, no generality is lost by requiring that $\varpi_+$ and $\varpi_-$ to obey $\varpi_+\varpi_- = q$. These sections can be further constrained by requiring that the integrals over $\Sigma$ of the square of their norms are the same.

What follows next are two important points regarding $L_+$ and $L_-$: Both have the same Euler characteristic as $T^{1,0}\Sigma$, and their tensor product is $T^{2,0}\Sigma$. As a consequence of this, these bundles can be written as $L_+ = T^{1,0}\Sigma \otimes L^2$ and $L_- = T^{1,0} \otimes L^{-2}$ with L denoting a holomorphic vector bundle with zero Euler class. The bundle L comes with a flat, Hermitian connection to be denoted by $\theta_L$; it is unique modulo the action of the group of unitary automorphisms of L which is $C^\infty(\Sigma; S^1)$.

Let E denote the $\mathbb{C}^2$ bundle $L \oplus L^{-1}$. Introduce by way of notation $\sigma_3$ to denote the anti-Hermitian endomorphism of this bundle that acts as multiplication by i on the left most summand (thus L) and as multiplication by -i on the right most summand (which is $L^{-1}$). The bundle E has a flat, Hermitian connection to be denoted by $\theta_E$ that is uniquely determined by requiring that its covariant derivative commute with $\sigma_3$ and that it act as $\theta_L$'s covariant derivative on sections of the respective L and $L^{-1}$ summands.

*Part 2*:  By virtue of what was just said about $L_+$ and $L_-$, the pull back by $\pi_q$ of $\varpi_+$ to $\Sigma_q$ can be written as $\mathfrak{f}_\delta\varpi_q$ with $\mathfrak{f}_\delta$ denoting a meromorphic section of $\pi_q{}^*L^2$ with first order zeros on $\Theta_-$ and first order poles on $\Theta_+$.  (Note here that $\varpi_q$ has second order zeros on both $\Theta_-$ and $\Theta_+$ whereas $\pi_q{}^*\varpi_+$ has first order zeros on $\Theta_-$ and third order zeros on $\Theta_+$.)  Meanwhile, $\pi_q{}^*\varpi_-$ can be written as $\mathfrak{f}_\delta{}^{-1}\varpi_q$.  Noting that the diffeomorphism of $\Sigma_q$ given by the action of -1 pulls back the bundle $\pi_q{}^*L^2$ to itself, it follows that the pull-back of $\mathfrak{f}_\delta$ by this diffeomorphism has to be $-\mathfrak{f}_\delta$.  As a consequence, $\mathfrak{f}_\delta$ can be viewed from the perspective of $\Sigma$ as a meromorphic section on $\Sigma-\Theta$ of $L^2 \otimes_\mathbb{R} \mathcal{I}_q$.

The section $\mathfrak{f}_\delta$ will be used now to define an anti-Hermitian endomorphism of $\pi_q{}^*E$ on $\Sigma_q-\Theta$.  To this end, let $u_\delta \equiv |\mathfrak{f}_\delta|^{-1}\mathfrak{f}_\delta$ which is a norm 1 section of $\pi_q{}^*L^2$.  The desired endomorphism is denoted by $\hat{\tau}$ and it is given by the rule below:

$$\hat{\tau} \equiv i \begin{pmatrix} 0 & u_\delta \\ u_\delta{}^{-1} & 0 \end{pmatrix} .$$

(2.1)

This object $\hat{\tau}$ can be viewed from the $\Sigma-\Theta$ vantage point as an isometric homomorphism from either E to $E \otimes_\mathbb{R} \mathcal{I}_q$ or from $E \otimes_\mathbb{R} \mathcal{I}_q$ to E.  (With regards to being isometric, this is because $|\hat{\tau}| = 1$.)

*Part 3*:  An important point now is that $\hat{\tau}$ is covariantly constant when the covariant derivative is defined by the flat connection

$$A = \theta_E - \tfrac{1}{4}\,[\hat{\tau}, d_{\theta_E}\hat{\tau}] .$$

(2.2)

Given the form of $\hat{\tau}$, this connection can also be written as

$$A = \theta_E + \tfrac{i}{2}\,(u_\delta{}^{-1}d_{\theta_L}u_\delta)\,\sigma_3 .$$

(2.3)

Note that A is a smooth connection on $\Sigma-\Theta$ but it has no extension to the whole of $\Sigma$ because $u_\delta$ is singular along $\Theta$.

To continue:  Define a 1-form on $\Sigma$ with values in the bundle of anti-Hermitian endomomorphisms of E by the rule

$$\mathfrak{a} = \mathfrak{Re}(\varpi_q)\,\hat{\tau} .$$

(2.4)

The key points to note are that $d_A\mathfrak{a} = 0$ and $d_A(*_\Sigma\mathfrak{a}) = 0$ and that $\mathfrak{a} \wedge \mathfrak{a} = 0$. Note also that $|\mathfrak{a}|$ extends across $\Theta$ as a Hölder continuous function with $\Theta$ being its zero locus (because this is the case for $|\varpi_q|$).

As a consequence of the preceding identities and the fact that $F_A = 0$, the pair $(A, \mathfrak{a})$ when pulled up to $S^1 \times (\Sigma - \Theta)$ obeys the equations in (1.7). This is also true for the pair $(A, r\,\mathfrak{a})$ when $r$ is any given real number. Taking a sequence of diverging $r$ values gives a sequence of solutions to (1.7) the converges in the manner of [T2] to $(Z = S^1 \times \Theta, \mathcal{I}, \varsigma)$.

## c) Singular solutions to $m \neq 0$ versions of (1.7)

The task for this subsection is to construct a divergent sequence of singular solutions to (1.7) for the $m \neq 0$ cases that converges in the manner of [T2] to Section 2a's singular data set $(Z = S^1 \times \Theta, \mathcal{I}, \varsigma)$. (The singularities are along the set Z.) The construction here has eight parts.

*Part 1*: Return to the milieu of Section 2a and the branched cover $\Sigma_q$ inside $T^{1,0}\Sigma$. As noted, in Section 2a, this complex curve has its canonical holomorphic differential (a section of $T^{1,0}\Sigma_q$) which was denoted by $\varpi_q$. Let $\mathfrak{Im}(\varpi_q)$ denote the imaginary part of $\varpi_q$. This is a closed 1-form on $\Sigma_q$. Integration of $\mathfrak{Im}(\varpi_q)$ around closed, oriented loops in $\Sigma_q$ defines a homomorphism from $H^1(\Sigma_q; \mathbb{Z})$ to $\mathbb{R}$ which is called the *period homomorphism* for $q$. Since $\mathfrak{Im}(\varpi_q)$ pulls back to be $-\mathfrak{Im}(\varpi_q)$ by the action of $(-1)$ as a diffeomorphism of $\Sigma_q$, this period homomorphism is zero on the invariant subspace in $H^1(\Sigma_q; \mathbb{Z})$ (which is the pull-back via the projection map of $H^1(\Sigma; \mathbb{Z})$). The lemma below describes this period homomorphism on the anti-invariant subspace.

**Lemma 2.1**: *Let $\hat{\mathcal{H}}$ denote the vector space of holomorphic quadratic differentials on $\Sigma$. There exists an open, dense set $\hat{\mathcal{H}}_0 \subset \hat{\mathcal{H}}$ which consists of holomorphic quadratic differentials with simple zeros; and there is a dense set in $\hat{\mathcal{H}}_0$ whose period homomorphism maps to the rational numbers.*

If the period homomorphism of a holomorphic, quadratic differential maps to the rational numbers, then multiplying that differential by a large integer will give a holomorphic quadratic differential whose period homomorphism has integer values. Assume henceforth that $q$ is chosen so that this is the case.

***Proof of Lemma 2.1***: Lipman Bers [B] gives a one word proof of the assertion that the set of holomorphic quadratic differential with simple zeros is open and dense in $\hat{\mathcal{H}}$, that word being 'Bertini'. An elaboration based on conversations with Curt McMullen and Joe Harris appears momentarily. The proof of the assertion that there is a dense set in $\mathcal{H}$ with rational period homomorphisms differs only cosmetically from an argument by Duady and Hubbard in [DH] (see Section 3 of [DH] in particular).

As for Bers' one word proof: The set in $\mathcal{H}_0$ is open in $\mathcal{H}$ so it is a matter of showing that the set is dense. This will be the case if, for each point in $\Sigma$, there exists an element in $\mathcal{H}$ that is non-zero at that point. Indeed, if $p \in \Sigma$, and $\hat{h}(p) = 0$, then $\hat{h}$ can be written using a local, holomorphic coordinate centered at p (call it z) as $\hat{h} = (\alpha z^k + \mathcal{O}(z^{k+1}))(dz)^2$ with $\alpha$ a nonzero complex number and with k being an integer. Now suppose that $\hat{k}$ is another element in $\mathcal{H}$ and that $\hat{k}$ at p is non-zero. After multiplying by a non-zero complex number, one can assume that $\hat{k}$ near p has the from $\hat{k} = (1 + \mathcal{O}(z))(dz)^2$. If $\varepsilon$ is any small but non-zero complex number, then $\hat{h} + \varepsilon\hat{k}$ near z has the form $(\alpha z^k + \varepsilon + \mathcal{O}(\varepsilon z + z^{k+1}))(dz)^2$. For $\varepsilon$ sufficiently small (but still not zero), this perturbed quadratic differential will vanish near z at k points, these being very nearly the k different k'th roots of $\alpha^{-1}\varepsilon$.

WIth the preceding understood, suppose for the sake of argument that there exists a point $p \in \Sigma$ such that every holomorphic quadratic differential vanishes at p. To derive nonsense from this assumption: There is a holomorphic line bundle with first Chern number -1 that has a nontrivial section with a single pole, this pole being at the point p. Let $\mathcal{L}$ denote this line bundle and let $\phi$ denote such a section (it is defined up multiplication by a non-zero complex number). Define a map from $\mathcal{H}$ to the vector space of holomorphic sections of $T^{2,0} \otimes_{\mathbb{C}} \mathcal{L}$ by the rule $\hat{h} \rightarrow \phi\hat{h}$. Because this map is injective, the space of holomorphic sections of $T^{2,0} \otimes_{\mathbb{C}} \mathcal{L}$ must be at least 3g - 3 (which is the dimension of $\mathcal{H}$). By Riemann-Roch, the index of the d-bar operator acting on sections of $T^{2,0} \otimes_{\mathbb{C}} \mathcal{L}$ is 3g -4. Thus, the cokernel dimension must be at least 1. This conclusion is nonsensical because the complex conjugate of an element in the cokernel is a holomorphic section of a complex line bundle with negative degree (the degree is -2g + 1) whereas holomorphic sections of line bundles vanish locally only with positive degree.

*Part 2*: By virtue of the integer period constraint on $q$, the 2-form $\frac{1}{2\pi} dt \wedge \mathfrak{Im}(\varpi_q)$ on $S^1 \times \Sigma_q$ defines an integer valued cohomology class which is denoted in this part of the proof by $z_q$. This cohomology class uniquely specifies an isomorphism class of Hermitian, complex line bundle over $S^1 \times \Sigma_q$ because there is no torsion in the 2'nd cohomology. Moreover, every Hermitian, complex line bundle in $z_q$'s isomorphism class has a unitary connection whose curvature 2-form is -2$\pi$ i times the 2-form $\frac{1}{2\pi} dt \wedge \mathfrak{Im}(\varpi_q)$.

Let $L_q$ denote a given complex, Hermitian line bundle whose first Chern class is $z_q$. The diffeomorphism (-1) acting on $\Sigma_q$ lifts to $S^1 \times \Sigma$ as an involution that fixes all points in the $S^1$ factor. Denote this lift diffeomorphism by (-1) also. The pull-back of $L_q$ by the lifted version of (-1) is a complex line bundle with first Chern class $-z_q$; and as a consequence, this pull-back $(-1)^*L_q$ is isomorphic to the complex line bundle $L_q^{-1}$ ($\equiv \mathrm{Hom}(L_q; \mathbb{C})$). Meanwhile, if $A_q$ denotes a unitary connection on $L_q$ with curvature 2-form $-i\, dt \wedge \mathfrak{Im}(\varpi_q)$, then the pull-back connection, $(-1)^*A_q$ is a connection on $(-1)^*L_q$ with curvature 2-form $i\, dt \wedge \mathfrak{Im}(\varpi_q)$.

*Part 3*: Introduce the $\mathbb{C}^2$–bundle $L_q \oplus (-1)^*L_q$. This complex vector bundle has vanishing first Chern class; so it is isomorphic to the product $\mathbb{C}^2$–bundle. An important point for what follows is that there is a fiber-wise $\mathbb{C}$-linear involution from $L_q \oplus (-1)^*L_q$ to itself (to be denoted by $\iota$) that covers the action of the diffeomorphism (-1) on the base space $S^1 \times \Sigma_q$. This map $\iota$ is obtained by first using the pull-back $(-1)^*$ to obtain a fiber-wise $\mathbb{C}$-linear covering map from $L_q \oplus (-1)^*L_q$ to $((-1)^*L_q) \oplus L_q$ and then acting by the endomorphism below

$$\begin{pmatrix} 0 & 1 \\ 1 & 0 \end{pmatrix}$$

(2.5)

which interchanges to two summands.

Let $A_q$ denote a unitary connection on $L_q$ with curvature 2-form $-i\, dt \wedge \mathfrak{Im}(\varpi_q)$. With $A_q$ in hand, define a unitary connection, $\hat{A}$, by requiring that the direct sum splitting of $L_q \oplus (-1)^*L_q$ be $\hat{A}$-covariantly constant, and by requiring that the covariant derivative on the respective $L_q$ and $(-1)^*L_q$ summands are defined by the connections $A_q$ and $(-1)^*A_q$. Note in particular that $\hat{A}$ is invariant with respect to the involution $\iota$.

The question now arises as to whether this $\hat{A}$ connection is an SU(2) connection or a U(2) connection. To rephrase: This question is asking whether the induced connection on determinant line bundle $L_q \otimes (-1)^*L_q$ has trivial holonomy. In this regard, note that the induced connection on this line bundle is necessarily a flat connection because its curvature 2-form is the sum of the curvature 2-forms of $A_q$ and $(-1)^*A_q$.

**Lemma 2.2**: *If the period homomorphism of q maps to $2\mathbb{Z}$, then the connection $A_q$ can be chosen so that its curvature 2-form is $-i\, dt \wedge \mathfrak{Im}(\omega_q)$ and so that the induced flat connection on the determinant line bundle $L_q \otimes (-1)^*L_q$ has trivial holonomies (and thus, there is a covariantly constant section). In general, $A_q$ can be chosen with the desired curvature 2-form so that the holonomies of the induced flat connection on the determinant line bundle are in the subgroup $\{\pm 1\}$.*

The proof of this Lemma is in the upcoming Part 4.

An important remark: Assuming that the initial choice of $q$ has integer period homomorphism, then the period homomorphism for the holomorphic quadratic differential $4q$ will map to $2\mathbb{Z}$. For the purposes of this paper, no generality is lost by assuming that the period homomorphism is even.

A second remark: In the case when the holonomies of the induced connection on the determinant line bundle are in $\{\pm 1\}$, then the connection $\hat{A}$ can be viewed as a connection on a principal bundle with structure group SU(2) $\times_{\{\pm 1\}} \mathbb{Z}/4\mathbb{Z}$ (which is a subgroup of U(2)). In any event, the connection $\hat{A}$ can still be used for the purposes of this paper since the object of interest ultimately is the induced SO(3) connection on the $\mathbb{R}^3$ bundle of skew Hermitian endomorphisms of $L_q \oplus (-1)^*L_q$.

*Part 4*: This part of the subsection is dedicated to proving Lemma 2.2.

**_Proof of Lemma 2.2_**: As noted above, the pull-back of $L_q \otimes (-1)^*L_q$ by the $(-1)$ diffeomorphism is canonically isomorphic to itself. It was also noted that this isomorphism identifies the connection on $L_q \otimes (-1)^*L_q$ induced by $\hat{A}$ with its pull-back. To say this differently, let $\nabla$ denote for the moment the $\hat{A}$ induced covariant derivative on sections of $L_q \otimes (-1)^*L_q$. If $s$ denotes any given section of that bundle, then $\nabla((-1)^*s) = (-1)^*(\nabla s)$.

With regards to sections: This bundle $L_q \otimes (-1)^*L_q$ admits a nowhere vanishing section because it is isomorphic to the product $\mathbb{C}$-bundle (since its first Chern class is zero). Let $s$ denote a section with $|s| = 1$. Because $s$ is nowhere zero, its pull-back $(-1)^*s$ can be written as u·$s$ with u denoting a map from $S^1 \times \Sigma_q$ to $\mathbb{C}$ with norm 1, thus a map to $S^1$. In addition, $(-1)^*u = u^{-1}$. As a map to $S^1$, the map u represents an anti-invariant class in $H^1(S^1 \times \Sigma; \mathbb{Z})$ with respect to the linear action of pull-back by $(-1)$ on the cohomology. If u represents an even element in the anti-invariant part of $H^1(S^1 \times \Sigma_q)$, then $s$ can be modified so that the new version obeys $(-1)^*s = s$. Indeed, if u represents an even element in the anti-invariant first cohomology, then u can be written as $v^2$ with v mapping to $S^1$ and pulling back according to the rule $(-1)^*v = v^{-1}$. (The case $(-1)^*v = -v$ cannot occur because this require v to vanish on the fixed point set $S^1 \times \Theta$.) With v in hand, take the new section, $s'$, to be v$s$. Then $(-1)^*s' = vs$ which is s´.

Suppose now that there is a section, $s$, with $(-1)^*s = s$. Write $\nabla s$ as $\nabla s = \lambda\, s$ with $\lambda$ being an i$\mathbb{R}$-valued, closed 1-form on $S^1 \times \Sigma_q$ (it is closed because $\nabla$ is defined by a flat U(1) connection). Note also that $(-1)^*\lambda = \lambda$ because $\nabla$ commutes with pull-back by $(-1)$. Let

$$A_q' = A_q - \frac{1}{2}\lambda.$$

$$(2.6)$$

Replace $A_q$ by $A_q'$ and construct the connection $\hat{A}$ as instructed but using $A_q'$. This new version of $\hat{A}$ is an SU(2) connection that is invariant with respect to the (-1) covering map ι on $L_q \otimes (-1)^*L_q$. (The section $s$ used above is $\nabla$-covariantly constant when $\nabla$ is defined by the new version of $\hat{A}$.)

As for when u defines an even element in $H^1(S^1 \times \Sigma_q)$, this is guaranteed when the first Chern of $L_q$ is an even element in $H^2(S^1 \times \Sigma; \mathbb{Z})$. To elaborate: If the first Chern class of $L_q$ is an even element, then $L_q$ can be written as $\mathcal{L}^2$ with $\mathcal{L}$ being a complex line bundle. If that can be done, then $L_q \otimes (-1)^*L_q$ is $(\mathcal{L} \otimes (-1)^*\mathcal{L})^2$. Meanwhile, the bundle $(\mathcal{L} \otimes (-1)^*\mathcal{L})$ also has zero first Chern class so it admits a nowhere zero section with norm 1. Let $t$ denote one of the latter. Since the pull-back of $(\mathcal{L} \otimes (-1)^*\mathcal{L})$ by the (-1) diffeomorphism is also canonically isomorphic to itself, the pull-back via (-1) of $t$ can be written $(-1)^*t = vt$ with v mapping to $S^1$. With this done, take $s$ to be $t^2$. It then follows that $(-1)^*s = v^2\,s$.

If the first Chern class of $L_q$ is not even and u does not admit a square root, there is still, locally, a square root of u. Let v denote this locally defined square root. (It can be viewed globally as a map to $S^1$ from a 2-1 cover of $S^1 \times \Sigma_q$.) Set $s' = vs$. Then write $\nabla s'$ as $\alpha\, s'$ just as before. (Note that $\alpha$ is globally defined not-with-standing the fact that $s'$ is not.) Make the replacement in (2.6) to define the new version of $\hat{A}$. It then follows that $s'$ is $\hat{A}$-covariantly constant. Since $s'$ is $\hat{A}$-covariantly constant but defined only modulo $\pm 1$, it follows that this new version of $\hat{A}$ induces a connection on $L_q \otimes (-1)^*L_q$ with holonomy in the group $\{\pm 1\}$.

*Part 5*: Return to the milieu of Parts a) and b) in Section 2b so as to reintroduce the line bundle L over $\Sigma$ and then the norm 1 section $u_\diamond$ of $\pi_q{}^*L^2$. Use the pull-back by the respective projections from $S^1 \times \Sigma$ and $S^1 \times \Sigma_q$ to define L and $u_\diamond$ over these spaces. Define the bundle $\hat{E}_q$ on $S^1 \times \Sigma_q$ by the rule below:

$$\hat{E}_q = (\pi_q{}^*L \otimes L_q) \oplus (\pi_q{}^*L^{-1} \otimes (-1)^*L_q)$$

(2.6)

Pull-back by the diffeomorphism (-1) gives a fiberwise, $\mathbb{C}$-linear map from $\hat{E}_q$ to the bundle

$$(\pi_q{}^*L \otimes (-1)^*L_q) \oplus (\pi_q{}^*L^{-1} \otimes L_q)$$

(2.7)

that covers the map (-1) on the base $S^1 \times \Sigma_q$. (Keep in mind that pull-backs from $S^1 \times \Sigma$ are not affected by the action of (-1).) With the preceding understood, it then follows that the composition of first pull-back by (-1) and then an application of what is denoted by $\hat{\tau}$ in (2.1) defines a fiberwise, $\mathbb{C}$-linear involution of the bundle $\hat{E}_q$ that covers the action of (-1) on $S^1 \times \Sigma_q$. This involution is denoted by $\iota^*$. Thus, $\iota^* = \hat{\tau} \, (-1)^*$.

*Part 6*: The task at hand now is to give $\hat{E}_q$ an $\iota^*$-invariant SU(2) (or SU(2) $\times_{\{\pm 1\}} \mathbb{Z}/4$) connection. To do this, let $A_q$ denote a connection on the line bundle $L_q$ with curvature 2-form $- i \, dt \wedge \mathfrak{Im}(\omega_q)$ of the sort that is supplied by Lemma 2.2. Meanwhile, let $\theta_L$ denote a flat, unitary connection on the line bundle L. This connection induces a unitary connection on $L^{-1}$ which is also denoted by $\theta_L$. Finally, let $\sigma_3$ denote the skew-Hermitian endomorphism of $\hat{E}_q$ that acts as multiplication by i on the left most summand in (2.6) and multiplication by -i on the right most summand. Given this data, a preliminary connection on $\hat{E}_q$ is defined by the following three rules: First, the endomorphism $\sigma_3$ is covariantly constant, which is to say that parallel transport by the connection preserves the respective summands in (2.7). Second, the covariant derivative of the connection restricts to sections of the left most summand to be the covariant derivative that is defined by $\pi_q^*\theta_L$ and $A_q$. Third, the covariant derivative of the connection restricts to sections of the right most summand to be the covariant derivative that is defined by the connection $\pi_q^*\theta_L$ on $L^{-1}$ and by $(-1)^*A_q$. Denote this preliminary connection by $\hat{A}_{\hat{E}}$. The desired connection is

$$\hat{A} = \hat{A}_{\hat{E}} + \frac{i}{2} \, (u_\diamond^{-1}d_{\theta_L}u_\diamond) \, \sigma_3 \ .$$

(2.8)

This connection is $\iota^*$-invariant. To see why, let $\hat{E}_q{}'$ denote the bundle depicted in (2.7) which is the pull-back by (-1) of $\hat{E}_q$. Let $\hat{A}'$ denote the connection on $\hat{E}_q{}'$ that is obtained by pulling back the connection $\hat{A}$ using this same diffeomorphism (-1). What is denoted by $\hat{\tau}$ in (2.1) defines an isomorphism from $\hat{E}_q{}'$ to $\hat{E}_q$. This isomorphism $\hat{\tau}$ intertwines the respective covariant derivatives of $\hat{A}'$ and $\hat{A}$, which is to say that $\nabla_{\hat{A}}\hat{\tau} = \hat{\tau}\nabla_{\hat{A}'}$ (the proof of this is left to the reader, but see (2.2) and (2.3)). (It is important to keep in mind here and subsequently that $\hat{\tau}$ is *not* an endomorphism of either $\hat{E}_q{}'$ or $\hat{E}_q$; it is an isomorphism from the former bundle to the latter.)

For future reference, the $\sigma_3$ endomorphism of $\hat{E}_q$ is $\hat{A}$-covariantly constant. Also for future reference: The curvature 2-form of $\hat{A}$ is

$$F_{\hat{A}} = - \, dt \wedge \mathfrak{Im}(\varpi_q) \, \sigma_3 \, .$$

(2.9)

The $\iota^*$-invariance of $\hat{A}$ implies that $\iota^*F_{\hat{A}} = F_{\hat{A}}$. The fact that $F_{\hat{A}}$ is $\iota^*$ invariant can be seen directly from the formula in (2.9) since $\iota^*\sigma_3 = -\sigma_3 \, \iota^*$ and since the pull-back of $\varpi_q$ by the diffeomorphism (-1) is $- \, \varpi_q$.)

*Part 7*:  Because $\hat{E}_q$ is $\iota^*$-invariant and because $\iota^*$ acts freely on $\hat{E}_q$ over $S^1 \times (\Sigma_q - \Theta)$ covering the action of (-1) on $S^1 \times (\Sigma_q - \Theta)$, the quotient of $\hat{E}_q$ by the action of $\iota^*$ defines a $\mathbb{C}^2$ vector bundle over $S^1 \times (\Sigma - \Theta)$.  This bundle is denoted below by $E_q$.  Because the connection $\hat{A}$ is $\iota^*$–invariant, it defines an unitary connection on $E_q$ (an SU(2) connection if the period homomorphism of $q$ is even).  This connection is denoted by $A_q$.

Not-with-standing the fact that $\iota^*$ does not act freely on the $S^1 \times \Theta$ incarnation inside $S^1 \times \Sigma_q$, the bundle $E_q$ does extend as a smooth $\mathbb{C}^2$ bundle over the $S^1 \times \Theta$ incarnation inside $S^1 \times \Sigma$.  (But the connection A although the norm of its curvature 2-form does extend as a Hölder continuous function over $S^1 \times \Theta$ that is zero on this set.)  To extend $E_q$, let D denote for the moment a small, positive radius disk in $\Sigma$ centered on a point in $\Theta$ with the radius chosen so that the closure of D intersects $\Theta$ only in its center point.  Let p denote that center point.  The restriction of E to $S^1 \times (D-p)$ is isomorphic to a product $\mathbb{C}^2$ bundle and any choice of isomorphism can be used to extend E over $S^1 \times p$.  Moreover, the extension is unique up to bundle isomorphism because the space of isomorphisms on $S^1 \times (D-p)$ between E and the product $\mathbb{C}^2$ bundle is path connected (which is due, in turn, to the fact that any map from a surface to SU(2) is homotopic to a constant map).

With regards to $\sigma_3$:  As noted, this endomorphism of $\hat{E}_q$ anti-commutes with $\iota^*$.  As a consequence, it does not descend to $S^1 \times (\Sigma - \Theta)$ as an endomorphism of $E_q$.  Rather, it descends to give an isometric, $A_q$-covariantly constant homomorphism from the line bundle $\mathcal{I}$ into the bundle of skew-Hermitian endomorphisms of $E_q$.  Alternately, it can be viewed from the vantage of $S^1 \times (\Sigma - \Theta)$ as a section of the tensor product of that bundle with $\mathcal{I}$.  In either guise, it is denoted below by $\hat{\sigma}$.

*Part 8*:  Let $q_1$ denote a holomorphic quadratic differential with simple zeros and even integer period homomorphism.  Supposing that n denotes a positive integer, then $q_n \equiv n^2 q_1$ is also a holomorphic quadratic differential with simple zeros and even integer period homomorphism.  Let $E_n$ denote the $q = q_n$ version of the bundle $E_q$ and let $A_n$ denote the $q = q_n$ version of the connection $A_q$ on $E_n$'s restriction to $S^1 \times (\Sigma - \Theta)$.  Let $m$ denote the constant chosen for (1.7) with it understood that $m > 0$.  Supposing that n denotes a positive integer.  Then the pair

$$\left( A = A_n, \mathfrak{a} = \frac{n}{m}\, \mathfrak{Re}(\varpi_{q_1})\, \hat{\sigma} \right)$$

$$(2.10)$$

is a solution to (1.7) on $S^1 \times (\Sigma - \Theta)$ that does not extend across $S^1 \times \Theta$. Even so, the norm of $\mathfrak{a}$ does extend as zero across this locus to define a Hölder continuous function on $S^1 \times \Sigma$ with Hölder exponent $\frac{1}{2}$.

The various integer n versions of (2.10) defines the promised a sequence of singuluar solutions to (1.7). In particular, this sequence converges in the manner described by [T2] to give the $\mathbb{Z}/2$ data set $(Z = S^1 \times \Theta, \mathcal{I}, \varsigma)$.

## 3. Approximate solutions

This section will use the singular solutions to (1.7) from the previous section to construct sequences of approximate solutions to $m \neq 0$ versions of (1.7) that converge in the manner of [T2] to the data set $(Z, \mathcal{I}, \varsigma)$ from Section 2a. Analogous constructions by Greg Parker [P1]-[P4] in the context of the 2-spinor Seiberg-Witten equations serve as a guiding model. In particular, the constructions here (as in Parker's work) make use of depictions in [MSWW] of flat $Sl(2; \mathbb{R})$ flat connections on $\Sigma$. (Keep in mind that flat $SL(2;\mathbb{R})$ connections define exact solutions to the $m = 0$ version of (1.7)). Because flat $Sl(2; \mathbb{R})$ connections on $\Sigma$ play a central role, the first subsection summarizes their salient properties. Almost all of what is done there will be a review for experts (among others, the authors of [MSWW]). The subsequent subsections use these $m = 0$ solutions with the singular ones from Section 2c to construct the desired sequences of approximate solutions to the $m \neq 0$ versions of (1.7).

## a) Nonconvergent sequences of flat $SL(2; \mathbb{R})$ connections on $\Sigma$

The description that follows of a non-convergent sequence of flat $SL(2; \mathbb{R})$ connections on $\Sigma$ has 7 parts. These parts introduce various constructions and notions that are used ubiquitously in the subsequent parts of this paper.

*Part 1*: Return to the context of Part 1) of Section 2b where $q$ denotes a holomorphic, quadratic differential on $\Sigma$ with only simple zeros. By way of a reminder, the zero locus of $q$ (the set $\Theta$) was divided there into two disjoint sets of 2g-2 points, $\Theta_+$ and $\Theta_-$. These were used to define respective holomorphic line bundles $L_+$ and $L_-$ with holomorphic sections $\varpi_+$ and $\varpi_-$ (the former having zero locus $\Theta_+$ and the latter having zero locus $\Theta_-$). These were chosen so that $\varpi_+\varpi_- = q$ and so that the integrals over $\Sigma$ of $|\varpi_+|^2$ and $|\varpi_-|^2$ are the same. (This section is uniquely defined up to multiplication by a non-zero complex number.) A choice for the section is denoted by $\varpi_+$ in what follows. Likewise, $\Theta_-$ defines a holomorphic line bundle $L_-$ and section, $\varpi_-$, whose zero locus is $\Theta_-$. Since $\varpi_+\varpi_-$ has the same zero locus as $q$, this product must be a complex multiple of $q$. As noted previously, no generality is lost by requiring t$\varpi_+$ and $\varpi_-$ to obey $\varpi_+\varpi_- = q$ and so that the integrals over $\Sigma$ of the square of their norms are the same.

*Part 2*: As noted in that Part 1 of Section 2b, the line bundles $L_+$ and $L_-$ can be holomorphically identified with $T^{1,0}\Sigma \otimes L^2$ and $T^{1,0} \otimes L^{-2}$ with L denoting a holomorphic line bundle with zero Euler class. A $\mathbb{C}^2$ bundle (denoted by E) was then defined to be $L \oplus L^{-1}$. For any given function $\mu$ defined on $\Sigma$, what is written below (denoted by $\mathbb{X}$) defines a section of $\text{End}(E) \otimes T^{1,0}\Sigma$:

$$\mathbb{X} = \begin{pmatrix} 0 & e^{\mu}\varpi_+ \\ e^{-\mu}\varpi_- & 0 \end{pmatrix}.$$

(3.1)

The next equation uses $\mu$ and $\mathbb{X}$ to define a pair consisting of a unitary connection on E (denoted by A) and a 1-form on $S^1 \times \Sigma$ with values in the bundle of anti-Hermitian endomorphisms of E (the latter denoted by $\mathfrak{a}$):

$$A = \theta_{\hat{E}} + \tfrac{1}{2} (*_\Sigma d\mu)\, \sigma_3 \quad and \quad \mathfrak{a} = \tfrac{i}{2} (\mathbb{X} + \mathbb{X}^\dagger).$$

(3.2)

In this instance, $\sigma_3$ denotes the skew-Hermitian endomorphism of E that acts as $+i$ on the L summand in E and as $-i$ on the $L^{-1}$ summand. The definition of A is such that both

$$d_A\mathfrak{a} = 0 \quad and \quad d_A(*_\Sigma\mathfrak{a}) = 0.$$

(3.3)

Granted this last point, then the pair $(A, \mathfrak{a})$ defines a flat $SL(2; \mathbb{C})$ connection on $\Sigma$ if the curvature of A and $\mathfrak{a} \wedge \mathfrak{a}$ are equal which is to say that $F_A - \mathfrak{a} \wedge \mathfrak{a} = 0$. Since the curvature 2-form of A is $\tfrac{1}{2} d(*_\Sigma d\mu) \sigma_3$, flatness is the condition below on $\mu$:

$$d^\dagger d\mu + \tfrac{1}{2} (e^{2\mu}|\varpi_+|^2 - e^{-2\mu}|\varpi_-|^2) = 0$$

(3.4)

where $d^\dagger$ denotes here $-*_\Sigma(d*_\Sigma)$ which is the formal adjoint of the exterior derivative.

By way of a useful but (singular) example: Define the function $\mu_\diamond$ on $\Sigma-\Theta$ by the rule

$$\mu_\diamond = -\tfrac{1}{2}\ln(|\varpi_+|/|\varpi_-|)$$

(3.5)

The $\mu = \mu_\diamond$ version of (3.4) is obeyed on $\Sigma-\Theta$ where $\mu_\diamond$ is defined because the logarithm of the ratio of the norms of the two holomorphic sections is a harmonic function (this is so because the holomorphic structure on L is defined by Hermitian connection with zero curvature). The $\mu_\diamond$ version of (3.3)'s connection A defines a connection on $E|_{\Sigma-\Theta}$ which is denoted still by A. Meanwhile, the $\mu_\diamond$ version of (3.3)'s endomorphism $\mathfrak{a}$ defines a 1-form on $\Sigma - \Theta$ with values in the bundle of skew-Hermitian endomorphisms of E which is denoted still by $\mathfrak{a}$. This $\mu_\diamond$ version of the pair $(A, \mathfrak{a})$ is the pair from (2.3) and (2.4). (To see this, reintroduce from Part 2 of Section 2b the meromorphic section $\mathfrak{f}_\diamond$ of the line bundle $\pi_q^* L^2$ and then note first that $\mu_\diamond$ is $-\ln|\mathfrak{f}_\diamond|$. The identity $-*_\Sigma\, d(\ln|\mathfrak{f}_\diamond|) = i u_\diamond^{-1} d_{\theta_L} u_\diamond$ holds because $\mathfrak{f}_\diamond$ is meromorphic and because $u_\diamond$ is $\mathfrak{f}_\diamond/|\mathfrak{f}_\diamond|$. Meanwhile, the $\mu = \mu_\diamond$ version of $\mathbb{X}$ is $-i\varpi_q\,\hat{\tau}$ by virtue of the definition of $\mu_\diamond$ as $-\ln|\mathfrak{f}_\diamond|$ and the definition of $\mathfrak{f}_\diamond$ as $\omega_+/\omega_q$.)

The upcoming Lemmas 3.1-3.3 say what is needed for now about solutions to (3.4). To set the stage for what is to come: These lemmas use $q_1$ to denote a given holomorphic quadratic differential with simple zeros. Its zero locus is denoted by $\Theta$. Supposing now that $r$ denotes a positive number, then $r^2 q_1$ is a holomorphic quadratic differential with the same zero locus $\Theta$; and the $q = r^2 q_1$ versions of $\varpi_+$ and $\varpi_-$ are $r$ times the respective $q_1$ versions. This implies in particular that the function $\mu_\diamond$ defined by (3.5) when $q = r^2 q_1$ for any given value of $r$ is identical to the version that is defined by $q_1$. It also implies that the $q$ versions of $\varpi_+$ and $\varpi_-$ are $r$ times the $q_1$ versions (which are denoted by $\varpi_{1+}$ and $\varpi_{1-}$).

*Part 3*: Fix a small, positive number to be denoted by $r_0$ which should be chosen so that there is a holomorphic coordinate chart for a neighborhood of any given point in $\Theta$ that identifies the radius $4r_0$ disk in $\mathbb{C}$ centered at the origin with an embedded disk in $\Sigma$ that is well inside a much larger radius Gaussian coordinate chart centered at that point. Any given $p \in \Theta$ version of the embedded $|z| < r_0$ disk in $\Sigma$ is denoted by D in what follows.

To say more about $r_0$: First, $r_0$ should be such that the holomorphic coordinate chart sits inside a radius $1000 r_0$ Gaussian coordinate chart centered at p. In addition, $r_0$ should be such that distance from p to all other points in $\Theta$ is greater than $10,000\ r_0$.

The chosen holomorphic coordinate for a neighborhood of any given $p \in \Theta$ should be chosen so that $\frac{1}{\sqrt{2}} dz$ has length 1 at p. Granted that this is done, then the product of the differential $\frac{1}{\sqrt{2}} dz$ with a function of the form $(1+h)$ with $|h| \leq c_0|z|$ defines a unit length basis section of $T^{1,0}\Sigma$ on the $|z| < 4r_0$ disk in $\Sigma$ and thus isometrically identifies $T^{1,0}\Sigma|_D$ with the product $\mathbb{C}$ bundle over that domain. A particular holomorphic coordinate for any given p from $\Theta$ will be specified in the upcoming Lemma 3.2.

In any event, with a holomorphic coordinate specified for each point in $\Theta$, it proves useful to define a function on $\Sigma$ which gives a measure of the distance to $\Theta$. This function is denoted by $\mathfrak{d}$ and it is defined as follows: For any given $p \in \Theta$, the value of $\mathfrak{d}$ is equal to $|z|$ on the corresponding $|z| \leq 2r_0$ disk. Then, where $|z|$ is between $2r_0$ and $4r_0$, the function $\mathfrak{d}$ is $\chi(\frac{|z|}{2r_0} - 1)z + (1 - \chi(\frac{|z|}{2r_0} - 1))\, 4r_0$. On the complement in $\Sigma$ of the various $p \in \Theta$ versions of the $|z| < 4r_0$ disk, this function $\mathfrak{d}$ has the constant value $4r_0$.

By way of remark regarding $\mathfrak{d}$: This function is commensurate with the metric distance to $\Sigma$ in the sense that $\mathfrak{d}$ can be written as $\mathfrak{d} = (1 + \mathfrak{g})\text{dist}(\cdot, \Theta)$ with $\mathfrak{g}$ being a function on $\Sigma$ whose norm is no greater than $c_0 \, \text{dist}(\cdot, \Theta)$ on any given $p \in \Theta$ version of D.

*Part 4*: This part of the subsection has the first lemma regarding (3.4).

**Lemma 3.1**: *For any $r > 1$, the $q = r^2 q_1$ version of (3.4), which is the equation*

$$d^\dagger d\mu + \tfrac{1}{2}\, r^2 \, (e^{2\mu}|\varpi_{1+}|^2 - e^{-2\mu}|\varpi_{1-}|^2) = 0 \, ,$$

*has a unique, smooth solution on $\Sigma$. Moreover, this solution (denoted by $\mu_r$) depends smoothly on the choice for $r$; and it has the properties listed in the bullets that appear below. The bullets use $\kappa$ to signify a number greater than 1 which is independent of the choice for $r$.*

- $|\mu_r - \mu_\mathfrak{d}| + r^{-2/3}|d(\mu_r - \mu_\mathfrak{d})| + r^{-4/3}|\nabla d(\mu_r - \mu_\mathfrak{d})| \leq \kappa \, e^{-r \mathfrak{d}^{3/2}/\kappa}$  *where $\mathfrak{d} > \kappa r^{-2/3}$.*
- $|\mu_r - \tfrac{1}{3}\ln r| + r^{-2/3}|d\mu_r| + r^{-4/3}|\nabla d\mu_r| \leq \kappa$  *where $\mathfrak{d} < \kappa r^{-2/3}$.*

The bounds that are asserted by the lemma follow directly from the analysis in [MSWW]. Even so, a full proof is given in the Appendix to this paper.

The upcoming Lemma 3.3 says more about the form of $\mu_r$ near the points in $\Theta$. A preliminary Lemma 3.2 comes before Lemma 3.3 to supply some input for Lemma 3.3.

*Part 5*: Given a point $p \in \Theta$, then the chosen holomorphic quadratic differential (this is $q$) can be written using the holomorphic coordinate on p's version of the disk D as $q = \alpha \, z \, (1 + \vartheta) \, (dz)^2$ with $\alpha$ being a non-zero complex number and with $\vartheta$ being a holomorphic function on D which vanishes at $z = 0$ (at p) and such that $1 + \vartheta \neq 0$. In the context of Lemma 3.1, the quadratic differential $q$ has the form $r^2 q_1$ with $r$ denoting a large, positive number. Since $q = r^2 q_1$ with $q_1$ being an $r$-independent, holomorphic quadratic differential $q_1$, the holomorphic function $\vartheta$ is independent of $r$ and the complex number $\alpha$ has the form $r^2 \alpha_1$ with $\alpha_1$ being independent of $r$.

Parallel transport using the flat connection $\theta_L$ along the rays out from $z = 0$ in the coordinate chart (the argument of z is constant on each ray) defines an isometric isomorphism between L on the radius $r_0$ disk centered at p and the produce $\mathbb{C}$-bundle. This isomorphism identifies the connection $\theta_L$ with the product connection. Using this identification of L (and the corresponding identification of $L^{-1}$) with the product structure on $T^{1,0}$ defined by the holomorphic coordinate, and assuming that p is in $\Theta_+$, the sections $\varpi_+$ and $\varpi_-$ appear as

$$\varpi_+ = \sqrt{\alpha}\ \beta\ z\ (1 + \vartheta_+)\ dz \quad \textit{and} \quad \varpi_- = \sqrt{\alpha}\ \beta^{-1}(1 + \vartheta_-)\ dz$$

$$(3.6)$$

with $\sqrt{\alpha}$ denoting a choice of square root for $\alpha$, with $\beta$ denoting a non-zero complex number with $c_0^{-1} \leq |\beta| \leq c_0$, and with $\vartheta_+$ and $\vartheta_-$ denoting holomorphic functions on the $|z| < 4r_0$ disk that vanish at the origin and obey $(1 + \vartheta_+)(1 + \vartheta_-) = (1 + \vartheta)$.

Granted the preceding, define the function $v$ on the $|z| < 4r_0$ disk in $\mathbb{C}$ using the rule

$$v \equiv \mu + \ln|\beta| + \tfrac{1}{2}\ln(|1+\vartheta_-|/|1+\vartheta_+|)\ .$$

$$(3.7)$$

A function $\mu$ obeys (3.4) on the radius $r_0$ disk in $\Sigma$ centered at $p$ if and only the function $v$ obeys the equation

$$-4\partial\bar{\partial}v + \tfrac{1}{2}\alpha\ |1+\vartheta|(e^{2v}\,|z|^2 - e^{-2v}) = 0$$

$$(3.8)$$

on the $|z| < r_0$ disk in $\mathbb{C}$. (The notation here uses $\partial$ to denote differentiation with respect to the complex variable $z$. Meanwhile, $\bar{\partial}$ denotes differentiation with respect to $\bar{z}$.) By way of a singular example of a solution to (3.8): The $\mu = \mu_0$ version of the function $v$ is $-\tfrac{1}{2}\ln|z|$.

There is an analog of $v$ and (3.8) for the case when $p \in \Theta$ which is obtained from (3.8) by changing the $v$ to $-v$ in that equation.

*Part 6*: The next lemma uniquely specifies the holomorphic coordinate for a given point in $\Theta$.

**Lemma 3.2**: *Given* $p \in \Theta$, *there is a unique holomorphic coordinate (the coordinate $z$) for a neighborhood of* $p$ *(the disk* D*) whereby* $\frac{1}{\sqrt{2}}dz$ *has length* 1 *at* $p$ *(which is where* $z = 0$*) and whereby the given holomorphic quadratic differential* $q$ *appears as* $\alpha\ z\ (dz)^2$ *with* $\alpha$ *being a positive real number. The corresponding version of (3.8) has* $\vartheta \equiv 0$ *and thus the corresponding version of* $v$ *obeys*

$$-4\partial\bar{\partial}v + \tfrac{1}{2}\alpha\ (e^{2v}\,|z|^2 - e^{-2v}) = 0\ .$$

Lemma 3.2 is also proved in the appendix.

Lemma 3.2's holomorphic coordinate should be used henceforth. By way of heads up: With this coordinate used, then the functions $(1 + \vartheta_+)$ and $(1 + \vartheta_-)$ that appear in (3.6) and (3.7) are inverses of each other.

A comment is warranted regarding the (apparent) absence of the Riemannian metric in (3.8) and Lemma 3.2. The metric actually does appear, but only implicitly via the complex coordinate since the complex structure on $\Sigma$ is determined by the conformal structure. The implicit dependence of (3.8) only on the underlying conformal structure is a consequence of the fact that the $m = 0$ version of the top equation in (1.7) and the other two equations in (1.7) on a Riemann surface are all conformaly invariant.

By way of a head's up: The existence of a coordinate chart for neighborhoods of the points in $\Theta$ where the equation in (3.4) appears as equation on a disk about the origin in $\mathbb{C}$ with only the constant $\alpha$ as parameter is exploited extensively in much of what follows (see for example Lemma 3.3, Lemma 5.2 and Section 7).

*Part 7*: What follows directly is the promised Lemma 3.3 regarding the solutions to (3.4) near the given point in $\Theta_+$. At the expense of changing $\mu$ to $-\mu$ in (3.4), what is said in Lemma 3.3 holds near the points in $\Theta_-$ as well. The lemma uses s to denote the Euclidean coordinate on $[0, \infty)$.

**Lemma 3.3:** *There exists* $\kappa > 1$ *and a smooth function (denoted by* f*) on* $[0, \infty)$ *with the following significance: Assume that* $r > \kappa$ *and that* $\mu_r$ *is described by Lemma 3.1 (in particular,* $q = r^2 q_1$*). If* p *denotes a given point in* $\Theta_+$, *then the corresponding version of (3.7)'s function* v *on the* $|z| < r_0$ *disk in* $\mathbb{C}$ *can be written as*

$$v = f(\alpha^{1/3}|z|) + \frac{1}{6}\ln(\alpha) + y$$

*with* y *obeying* $|y| + |dy| + |\nabla dy| \leq \kappa e^{-r/\kappa}$ *where* $|z| < \frac{1}{2} r_0$. *As for the function* f:

- *The function* $z \to f(|z|)$ *on* $\mathbb{C}$ *is a radially symmetric solution on the whole of* $\mathbb{C}$ *to the* $\alpha = 1$ *version of the equation in Lemma 3.2,*

$$-\frac{1}{s}\frac{d}{ds}\left(s\frac{d}{ds}f\right) + \frac{1}{2}\left(e^{2f}s^2 - e^{-2f}\right) = 0 \ .$$

- *The derivatives to second order of this function* f *are bounded by* $\kappa$ *on* $[0, 1]$ *and they obey the bounds below on* $[1, \infty)$:

$$|f(s) + \frac{1}{2}\ln(s)| + |s\frac{d}{ds}f + \frac{1}{2}| + |s^2\frac{d^2}{ds^2}f - \frac{1}{2}| \leq \kappa \ e^{-s^{3/2}/\kappa} \ .$$

- *The function* f $+\frac{1}{2}\ln(s)$ *is strictly negative, its derivative is strictly positive and its second derivative is strictly negative. Meanwhile,* $(e^{2f}s^2 - e^{-2f})$ *is strictly negative also, and, in any event, its norm is bounded by* $\kappa \ e^{-s^{3/2}/\kappa}$.

The proof of this lemma can also be deduced from what is said in [MSWW]. Even so, a proof is given in the Appendix

**b) Nonconvergent sequences of approximate solutions to the $m \neq 0$ versions of (1.7)**

The desired sequences of approximate solutions are obtained from the singular sequences of Section 2c by modifying the latter in small radius neighborhoods of the set $S^1 \times \Theta$. The modifications are described explicitly in the six parts of this subsection; but this is done only for neighborhood of the $S^1 \times \Theta_+$ circles because the description of the modification for $S^1 \times \Theta_-$ circles is identical but for some judicious changes of sign.

*Part 1*: To set the stage for what is to come, fix $p \in \Theta_+$ and reintroduce from Lemma 3.2 the corresponding holomorphic coordinate z. The surface $\Sigma_q$ when viewed as a branched 2-fold cover of the $|z| < 4r_0$ part of $\mathbb{C}$ appears as the set of pairs (z, $\varpi$dz) with z in $\mathbb{C}$ obeying $|z| < r_0$ and with $\varpi$ in $T^{1,0}\mathbb{C}|_z$ obeying

$$\varpi^2 = \alpha z.$$

(3.9)

In this regard, it proves convenient to introduce the coordinate $w = \alpha^{-1/2}\varpi$ in what follows and use it as a coordinate for a neighborhood of p in $\Sigma_q$. In terms of the coordinate $w$, the diffeomorphism (-1) sends $w$ to -$w$. The projection map $\pi_q$ when depicted using the coordinate z for the neighborhood of p in $\Sigma$ and $w$ for the neighborhod of p in $\Sigma_q$ is this:

$$w \to z = \pi_q(w) \quad with \quad \pi_q(w) = w^2$$

(3.10)

This map $\pi_q$ is defined on an open set in $\mathbb{C}$ containing the origin that vanishes at the origin. This open set will contain (in particular), a disk in $\mathbb{C}$ centered at the origin with radius greater than $\sqrt{2}\sqrt{r_0}$. This disk where $w$ defines a bonafide holomorphic coordinate for $\Sigma_q$ is denoted by $D_q$ in what follows.

When written using the coordinate $w$ on the $D_q$ part of $\Sigma_q$, the holomorphic objects $\pi_q{}^*q$ and $\varpi_q$ have the form below

$$\pi_q{}^*q = 4\alpha w^4 (\mathrm{d}w)^2 \quad and \quad \varpi_q = \sqrt{\alpha}\, 2w^2 \mathrm{d}w$$

(3.11)

(To see this, use the fact that $q = z(\mathrm{d}z)^2$ and the relation z = $w^2$.)

*Part 2*: The line bundle L on the disk D in $\Sigma$ (where $|z| < 2r_0$) is isomorphic to the product $\mathbb{C}$-bundle via a bundle isomorphism that identifies the flat connection $\theta_L$ with the product connection. This identification of L and its connection $\theta_L$ is used implicitly in what follows.

With regards to the bundle $L_q$ near $S^1 \times p$ in $S^1 \times \Sigma_q$: This bundle is isomorphic to the product $\mathbb{C}$ bundle on $S^1 \times p$. Now, because $S^1 \times p$ is fixed by (-1), the holonomy of the connection $A_q$ on this circle is either $\pm 1$. If $q$'s period homomorphism maps to $4\mathbb{Z}$, then the holonomy has to be $+1$. (Note though that if the period homomorphism maps to $2\mathbb{Z}$, then the holonomy is the same on all of the $S^1 \times \Theta$ circles, either $+1$ on all or $-1$ on all.) For simplicity (but not really losing generality), assume henceforth that the holonomy is $+1$ on all. Granted this assumption, then $L_q$ on $S^1 \times p$ is isomorphic to the product $\mathbb{C}$ bundle via an isomorphism that identifies the connection $A_q$ with the product connection. As a consequence, it is isomorphic to the product $\mathbb{C}$ bundle on the $S^1 \times D_q$ part of $S^1 \times \Sigma_q$ via an isomorphism that identifies $A_q$ with the sum of the product connection (denoted by $\theta_0$) and an $i\mathbb{R}$-valued 1-form obtained from the curvature tensor. Slightly more explicitly, $A_q$ has the schematic form below:

$$A_q = \theta_0 + i \frac{2}{3} \sqrt{\alpha} \, \Im(w^3) \, dt$$

(3.12)

The product structures for L and $L_q$ on the $S^1 \times D_q$ part of $S^1 \times \Sigma_q$ identifies the bundle $\hat{E}_q$ with the product bundle with fiber $\mathbb{C} \oplus \mathbb{C}$ such that the left most $\mathbb{C}$ corresponds to the $L \otimes L_q$ summand in $\hat{E}_q$ and the right most $\mathbb{C}$ corresponds to the $L^{-1} \otimes (-1)^* L_q$ summand in $\hat{E}_q$.

To see about the action of $\iota^*$ on $\hat{E}_q$, note first that $\pi_q^* \varpi_+$ and $\pi_q^* \varpi_-$ will have the following form on $D_q$:

$$\pi_q^* \varpi_+ = 2\sqrt{\alpha} \, \beta \, (1 + \vartheta_+) w^3 dw \quad and \quad \pi_q^* \varpi_- = 2\sqrt{\alpha} \, \beta^{-1} (1 + \vartheta_+)^{-1} \, w \, dw$$

(3.13)

with $\vartheta_+$ viewed here as a holomorphic function of $w^2$ that is defined on the disk in $\mathbb{C}$ corresponding to $D_q$ which vanishes at the origin. It follows as a consequence of these formulae that the isomorphism $u_\Diamond$ that is used in (2.1) is independent of $\alpha$ when written in terms of $w$ since it has the form depicted below:

$$u_\Diamond = \frac{w \, \beta \, (1 + \vartheta_+)}{|w||\beta||1 + \vartheta_+|} \, .$$

(3.14)

To continue with regards to $\iota^*$: Keeping in mind the identification between $\hat{E}_q$ and the product $\mathbb{C} \oplus \mathbb{C}$ bundle, what follows is an $\iota^*$-invariant, orthonormal basis of sections of $\hat{E}_q$ on $S^1 \times (D_q-0)$:

$$s_1 = \tfrac{1}{\sqrt{2}} \left(1, iu_\diamond^{-1}\right) \quad and \quad s_2 = \tfrac{1}{\sqrt{2}} \left(u_\diamond, -i\right) \quad.$$

(3.15)

Being that this basis is $\iota^*$-invariant, it descends to the $S^1 \times \pi_q(D_q-0)$ part of $S^1 \times \Sigma$ as an orthonormal basis for $E_q$. The resulting product structure for $E_q$ on the complement of the point p in its $S^1 \times \pi_q(D_q-0)$ tubular neighborhood is used to extend $E_q$ as a $\mathbb{C}^2$ bundle over $S^1 \times p$. (This product structure on this neighborhood will be referred to in subsequent sections as the $\mathbb{C} \oplus \mathbb{C}$ product structure.)

*Part 3*: Let $q_1$ denote a holomorphic quadratic differential with simple zeros and with period homomorphism mapping to $4\mathbb{Z}$. Supposing that n is a positive integer, set $q = n^2 q_1$. The task in this part of the proof is to depict the pair of singular solutions from (2.10) on the $S^1 \times \pi_q(D_q-0)$ part of $S^1 \times (\Sigma-\Theta)$ using the basis in (3.15) for the $q = n^2 q_1$ version of $E_q$.

Consider first the depiction of $\hat{\sigma}$. In this regard, keep in mind that $\hat{\sigma}$ comes from the skew-Hermitian endomorphism of $\hat{E}_q$ given by $\sigma_3$ (which acts as multiplication by i on the left most summand of $\hat{E}_q$ and multiplication by -i on the right most summand). The latter's action on $s_1$ and $s_2$ is as follows:

$$\sigma_3 s_1 = iu_\diamond^{-1} s_2 \quad and \quad \sigma_3 s_2 = iu_\diamond s_1 \,,$$

(3.16)

This says in effect that the endomorphism $\hat{\sigma}$ when written using the basis $(s_1, s_2)$ to identify $E_q$ with the product $\mathbb{C} \oplus \mathbb{C}$ bundle acts as the $\mathcal{I}$-valued endomorphism $\hat{\tau}$ from (2.1) with $u_\diamond$ given by (3.14). Thus, what is denoted by $\mathfrak{a}$ in (2.15), when written with respect to this $(s_1, s_2)$ identification of $E_q$ with the product $\mathbb{C} \oplus \mathbb{C}$ and for the $q = n^2 q_1$ version of $E_q$ appears as follow:

$$\mathfrak{a} = \tfrac{n}{m} \sqrt{\alpha_1}\, 2\, \mathfrak{Re}(w^2 dw)\, \hat{\tau} \,.$$

(3.17)

To interpret (3.17): Use the relation $w^2 = z$ to write the appearances of $w$ and $dw$ and in $\hat{\tau}$ in terms of the coordinate z. Note in particular that both $w$ and $\hat{\tau}$ are independent of n and also independent of $m$. All dependence of $\mathfrak{a}$ on these two number is in the explicit $\tfrac{n}{m}$ factor in (3.17).

The next order of business is to depict the covariant derivative from the connection in (2.10) when the $(s_1, s_2)$ basis is used to identify $E_q$ with the product $\mathbb{C} \oplus \mathbb{C}$ bundle. To this end, remember that the connection A on $E_q$ comes from the $\iota^*$-invariant connection $\hat{A}$ that is depicted in (2.8). Keeping in mind also the depiction in (3.12) of $A_q$ and what was said above about $\hat{\sigma}$, it follows that the connection $\hat{A}$ from (2.8) on the $S^1 \times (D_q{-}0)$ part of $S^1 \times (\Sigma_q{-}\Theta)$ when written for the $(s_1, s_2)$ identification of $E_q$ with the product $\mathbb{C} \oplus \mathbb{C}$ bundle (denoted by $\theta_{\mathbb{C} \oplus \mathbb{C}}$) is the connection given the upcoming in (3.18). With regards to the notation in (3.18): What is denoted by $\alpha_1$ in (3.18) for the case at hand (when $q = \mathrm{n}^2 q_1$) is $\mathrm{n}^2$ times that $q_1$ version of the positive real number $\alpha$ that appears in (3.9). What is denoted by $\sigma_3$ acts on the product $\mathbb{C} \oplus \mathbb{C}$ bundle as $+i$ on the left most $\mathbb{C}$ summand and as $-i$ on the right most $\mathbb{C}$ summand.

$$\hat{A} = \theta_{\mathbb{C} \oplus \mathbb{C}} + \tfrac{i}{2}(u_\vartheta^{-1}du_\vartheta)\,\sigma_3 + \tfrac{2}{3}\mathrm{n}\sqrt{\alpha_1}\,\mathfrak{Im}(w^3)\,\hat{\tau}\,dt$$

$$(3.18)$$

The connection A on the $S^1 \times \pi_q(D_q{-}0)$ part of $S^1 \times \Sigma$ when written using the $(s_1, s_2)$ isomorphism between $E_q$ and the product $\mathbb{C} \oplus \mathbb{C}$ bundle has the same form as what is written above but for writing appearances of $w$ and its implicit appearance in $u_\vartheta$ using the coordinate z via the formula $w^2 = z$.

*Part 4*: Fix a large positive integer to be denoted by n and set $q$ to be $\mathrm{n}^2 q_1$ as in Part 3 of this subsection. This part of the subsection supplies a preliminary smoothing of the corresponding pair $(A, \mathfrak{a})$ near points in $\Theta$ so that the result is an approximate solution to (1.7) with error that is independent of n (in a suitable norm). The plan for doing this is to modify Section 2c's singular pair only near the set $S^1 \times \Theta$ leaving Section 2c's pair $(A, \mathfrak{a})$ unchanged where $\vartheta > \frac{1}{100} r_0$ on $S^1 \times \Sigma$. Because of this, the discussion that follows will focus solely on the modification where $\vartheta < \frac{1}{100} r_0$ on a given $p \in \Theta_+$ version of $S^1 \times D$. In this regard, the focus will be on a point from $\Theta_+$ with it understood that the modification is identical but for cosmetic changes near points in $\Theta_-$.

To set the stage for the modification, introduce as before, the holomorphic coordinate z for p's version of the disk D. The modification employs a cut-off function to be denoted by $\chi_0$ which is defined by the rule $\chi_0(z) = \chi(200\,\frac{|z|}{r_0} - 1)$. This function is equal to 1 where $|z| < \frac{1}{200} r_0$ and it is equal to zero where $|z| > \frac{1}{100} r_0$.

Let $\mu_\mathrm{n}$ denote the $r = \frac{\mathrm{n}}{m}$ version of the function $\mu_r$ that is described by Lemma 3.1. Meanwhile, define $\mu_\vartheta$ as in (3.5). Then, set

$$\mu_{\lozenge n} \equiv \chi_0 \, \mu_n + (1 - \chi_0) \, \mu_\lozenge \, .$$

(3.19)

Define $\mathbb{X}$ via the $\mu = \mu_{\lozenge n}$ version (3.1) with $\varpi_+$ and $\varpi_-$ defined with the holomorphic, quadratic differential $q$ equal to $r^2 q_1$ (note that this is $r^2 q_1$, not $n^2 q_1$). This $\mu_{\lozenge n}$ version of $\mathbb{X}$ is denoted by $\mathbb{X}_n$ and it is used to replace (3.17) and (3.18) in the $|z| < \frac{1}{100} \, r_0$ part of $S^1 \times D$ by the pair depicted below in (3.20). (Remember that $w^2$ is defined from z by the rule whereby $w^2 = z$, and note that $\hat{\tau}$ is defined via (2.1) using (3.14)'s version of $\mathfrak{u}_\lozenge$.)

- $A = \theta_{\mathbb{C} \oplus \mathbb{C}} + \frac{1}{2} (*_\Sigma d \mu_{\lozenge n}) \, \sigma_3 + \frac{2}{3} \, m r \sqrt{\alpha_1} \, \Im m(w^3) \, \hat{\tau} \, dt$

- $\mathfrak{a} = \frac{i}{2} (\mathbb{X}_n + \mathbb{X}_n^\dagger) \, .$

(3.20)

(This pair is identical to the pair depicted in (3.17) and (3.18) where $|z| > \frac{1}{100} \, r_0$ because $\mu_{\lozenge n} = \mu_\lozenge$ where $|z| > \frac{1}{100} \, r_0$ and because the $\mu = \mu_\lozenge$ version of (3.20) is exactly what is depicted in (3.17) and (3.18).)

With regards to the dt component: The dt component of A is not smooth across the z = 0 locus; it's first derivative is only Hölder continous with exponent $\frac{1}{2}$. This lack of smoothness is fixed in the next part of the subsection.

With regards to solving (1.7): If (3.20)'s version of $(A, \mathfrak{a})$ is used in (1.7), then the the left hand side no longer vanishes. This is depicted schematically below in (3.21).

- $F_A - \mathfrak{a} \wedge \mathfrak{a} - m * \mathfrak{a} = \mathfrak{t}_0$

- $d_A \mathfrak{a} = \mathfrak{t}_1$

- $d_A * \mathfrak{a} = 0.$

(3.21)

It is a consequence of Lemma 3.1 that the pointwise norm of $\mathfrak{t}_0$ and $\mathfrak{t}_1$ are bounded by

$$c_0 \, m \, r^{2/3} \, e^{-r |z|^{3/2} / c_0}$$

(3.22)

with $r$ again denoting $\frac{n}{m}$. (Looking ahead: The bound in (3.22) implies a bound by $c_0 m^2$ on the $S^1 \times \Sigma$ integral of the square of the norm of the left hand side of (3.21).)

To be more explicit about $\mathfrak{t}_0$ and $\mathfrak{t}_1$ in (3.21), first write the dt component of (3.20)'s connection A where $m \, \mathtt{A}^0$ with $\mathtt{A}^0$ being

$$\mathtt{A}^0 = \frac{2}{3} \, r \sqrt{\alpha_1} \, \Im m(w^3) \, \hat{\tau} \, ;$$

(3.23)

and then write (3.20)'s lie algebra valued 1-form $\mathfrak{a}$ as $\mathfrak{a}_1 e^1 + \mathfrak{a}_2 e^2$ where $\{e^1, e^2\}$ constitute an oriented, orthonormal frame for $T^*\Sigma$ where $|z| < r_0$. (These frame 1-forms can be chosen to have the form $e^1 = (1+h)dz_1$ and $e^2 = (1+h)dz_2$ with $z_1$ and $z_2$ denoting the real and imaginary parts of the complex coordinate z, and with h obeying $|h| \leq c_0|z|$.) Then

- $\mathfrak{t}_0 = m \, (d_A A^0 - (\mathfrak{a}_1 e^2 - \mathfrak{a}_2 e^1)) \wedge dt + \eta \, e^1 \wedge e^2$

- $\mathfrak{t}_1 = - m \, [A^0, (\mathfrak{a}_1 e^1 + \mathfrak{a}_2 e^2)] \wedge dt$ .

(3.24)

where $\eta$ is a skew-Hermitian endomorphism of $E_q$ which is supported only where $|z|$ is between $\frac{1}{200} r_0$ and $\frac{1}{100} r_0$. (It is non-zero only on the support of $\chi_0$.) Moreover, $\eta$ is proportional to $\sigma_3$ when it is depicted using Part 2's $\mathbb{C} \oplus \mathbb{C}$ product structure for $E_q$ on that part of $S^1 \times \Sigma$.

*Part 5*: This part of the subsection makes an additional modification to the pair depicted in (3.20) by modifying the dt component of A and adding a dt component of $\mathfrak{a}$. As was the case before, this modification changes (3.20)'s version of $(A, \mathfrak{a})$ only where the $\mathfrak{d}$ is is less than $\frac{1}{100} r_0$. The new version of $(A, \mathfrak{a})$ on the $\mathfrak{d} < \frac{1}{100} r_0$ part of $S^1 \times \Sigma$ has this form:

- $A = \theta_{\mathbb{C} \oplus \mathbb{C}} + \frac{1}{2} (*_\Sigma d\mu_{0n}) \, \sigma_3 + m \, (A^0 + \mathfrak{b}_t) \, dt$

- $\mathfrak{a} = \frac{i}{2} (\mathbb{X}_n + \mathbb{X}_n^\dagger) + m \, \mathfrak{c}_t \, dt$

(3.25)

where $(\mathfrak{b}_t, \mathfrak{c}_t)$ are described below in the upcoming Proposition 3.4.

**Proposition 3.4**: *There exists* $\kappa > 1$ *with the following significance: If* $r \equiv \frac{n}{m}$ *is greater than* $\kappa$*, then* $\mathfrak{b}_t$ *and* $\mathfrak{c}_t$ *can be chosen so that (3.25)'s version of* $(A, \mathfrak{a})$ *are such that*

- $*(F_A - \mathfrak{a} \wedge \mathfrak{a}) - m\mathfrak{a} = m \, \mathfrak{t}_0 + m^2 \mathfrak{s}_0 \, dt$

- $*(d_A \mathfrak{a}) = m \, \mathfrak{k}_1$,

- $*d_A *\mathfrak{a} = m^2 \, \mathfrak{s}_1$ ,

*where* $\mathfrak{k}_0$, $\mathfrak{k}_1$ *and* $\mathfrak{s}_0$ *and* $\mathfrak{s}_1$ *are non-zero only where* $\mathfrak{d} < \frac{1}{100} r_0$ *and their norms are such that*

$$|\mathfrak{k}_0| + |\mathfrak{k}_1| \leq \kappa \, e^{-r/\kappa} \ \ and \ \ |\mathfrak{s}_0| + |\mathfrak{s}_1| \leq \kappa e^{-r \mathfrak{d}^{3/2}/\kappa} \ .$$

*In addition, $\mathfrak{t}_0$ and $\mathfrak{t}_1$ annihilate the vector field $\frac{\partial}{\partial \mathfrak{t}}$. As for $\mathfrak{s}_0$ and $\mathfrak{s}_1$, when they are written on the $\mathfrak{d} < \frac{1}{100} \mathfrak{r}_0$ part of $S^1 \times \Sigma$ using Part 2's $\mathbb{C} \oplus \mathbb{C}$ product structure for $E_q$, the former, $\mathfrak{s}_0$, is proportional to $\sigma_3$ whereas $\mathfrak{s}_1$ is pointwise orthogonal to $\sigma_3$. Finally, each of $\mathfrak{t}_0$, $\mathfrak{t}_1$ and $\mathfrak{s}_0$ and $\mathfrak{s}_1$ are independent of the $\mathfrak{t}$ coordinate for the $S^1$ factor of $S^1 \times \Sigma$ when they are written using that same product structure for $E_q$. With regards to $\mathfrak{b}_\mathfrak{t}$ and $\mathfrak{c}_\mathfrak{t}$: The norms of both are at most $\kappa e^{-r \mathfrak{d}^{3/2}/\kappa}$, and the norms of their covariant derivatives are bounded by $\kappa r^{2/3} e^{-r \mathfrak{d}^{3/2}/\kappa}$. In addition, when $\mathfrak{b}_\mathfrak{t}$ and $\mathfrak{c}_\mathfrak{t}$ are written on the $\mathfrak{d} < \frac{1}{100} \mathfrak{r}_0$ part of $S^1 \times \Sigma$ using Part 2's $\mathbb{C} \oplus \mathbb{C}$ product structure for $E_q$, they appear with $\mathfrak{b}_\mathfrak{t}$ being orthogonal to $\sigma_3$ and $\mathfrak{c}_\mathfrak{t}$ being proportional to $\sigma_3$. Also, both $\mathfrak{b}_\mathfrak{t}$ and $\mathfrak{c}_\mathfrak{t}$ when written with this same product structure, are independent of the $\mathfrak{t}$ coordinate for the $S^1$ factor.*

This proposition is proved in Section 5b.

## 4. The perturbation problem and the associated linear operator

The approximate solutions to (1.7) given by the pairs $(A, \mathfrak{a})$ from Proposition 3.4 are the starting point for a perturbative construction of an honest solution to (1.7). This perturbative approach can lead to a solution if $m$ is small (less than $c_0^{-1}$) and if a certain obstruction vector vanishes. This section first describes the perturbation theory and then starts the analysis of the linear operator that is central to the perturbative approach. (Much of what is done in this section specializes aspects of Parker's [P1] analysis. See also [P3].)

### a) Setting up the perturbation theory

The stage setting for the perturbative approach to (1.7) starts by enlarging that system of equations so as to obtain a first order, elliptic system of differential equations after accounting for the action of the automorphism group of the given principle bundle (the gauge group action). What is said below with regards to (1.7) holds not just for $S^1 \times \Sigma$, but for any compact, oriented, Riemannian 3-manifold. To state things in this more general context, Let Y denote the 3-manifold in question and let $P \to Y$ denote a principle SU(2) or SO(3) bundle. A solution to the enlarged version of (1.7) on Y is a data set $((A, \mathfrak{a}), (\mathfrak{p}, \mathfrak{q}))$ whose constituents are as follows: What is denoted by A signifies a connection on the given principle bundle, and what is denoted by $\mathfrak{a}$ signifies a 1-form on Y with values in the principle bundle's associated Lie algebra vector bundle. Meanwhile, both $\mathfrak{p}$ and $\mathfrak{q}$ denote sections of that same associated Lie algebra bundle. The equations are depicted below in (4.1) when written using an oriented orthonormal frame $\{e^1, e^2, e^3\}$ for the cotangent bundle of Y over a given open set.

- $-\nabla_{Ak}\mathfrak{p} + [\mathfrak{q}, \mathfrak{a}_k] + B_k - \frac{1}{2} \varepsilon_{kij} [\mathfrak{a}_i, \mathfrak{a}_j] - m \, \mathfrak{a}_k = 0$ ,

- $-\nabla_{Ak}q - [\mathfrak{p}, \mathfrak{a}_k] + \varepsilon_{ijk}\nabla_{Ai}\mathfrak{a}_j = 0$ ,
- $\nabla_{Ak}\mathfrak{a}_k = 0$ .

(4.1)

Here, $\{\nabla_{Ak}\}_{k=1,2,3}$ are the directional covariant derivatives defined by A and the metric connection along the dual frame vectors for the tangent space of Y. Meanwhile, what is denoted by $\{B_k\}_{k=1,2,3}$ are the components of the Hodge star of A's curvature tensor, $F_A$. Also, what is denoted by $\{\varepsilon_{ijk}\}_{i,j,k\,=\,1,2,3}$ are the components of the completely antisymmetric 3-tensor on Y with the convention that $\varepsilon_{123} = 1$.

A solution to (1.7) is a $\mathfrak{p} = q = 0$ solution to (4.1). Conversely, the (A, $\mathfrak{a}$) part of a solution to (4.1) obeys (1.7) and, if $\mathfrak{p}$ and/or $q$ are not zero, then [$\mathfrak{a}$, $\mathfrak{p}$] and [$\mathfrak{a}$, $q$] and $\nabla_A\mathfrak{p}$ and $\nabla_A q$ are identically zero.

To prove the preceding claim about the general form of a solutions to (1.7): The left hand side of both the first and second bullets in (4.1) define a 1-form on Y with values in the associated Lie algebra bundle. With this point understood, take the covariant divergence of the expression on the left hand side of the top bullet equation in (4.1) and then use the other bullets to see that $\mathfrak{p}$ obeys the second order equation

$$\nabla_A^\dagger \nabla_A \mathfrak{p} + [\mathfrak{a}_k, [\mathfrak{p}, \mathfrak{a}_k]] = 0 \,.$$

(4.2)

This equation can hold only in the event that $\nabla_A\mathfrak{p} = 0$ and [$\mathfrak{a}$, $\mathfrak{p}$] = 0 (take the inner product of both sides with $\mathfrak{p}$ and then integrate the resulting identity over Y, then integrate by parts). With the preceding about $\mathfrak{p}$ in hand, take the divergence of the expression on the left hand side of the second bullet equation in (4.1) and use the other bullets plus the preceding facts about $\mathfrak{p}$ to see that $q$ also obeys (4.2). Thus, $\nabla_A q$ and [$\mathfrak{a}$, $q$] must both be zero. It follows as a consequence of the preceding observations about $\mathfrak{p}$ and $q$ that (A, $\mathfrak{a}$) must then obey (1.7).

To solve (4.1): Start with a data set of the form (($\hat{A}$, $\hat{\mathfrak{a}}$), ($\mathfrak{p} = 0$, $q = 0$)) with $\hat{A}$ denoting some given connection on the principle bundle P and with $\hat{\mathfrak{a}}$ denoting a given 1-form with values in the associated Lie algebra bundle. With this initial choice in hand, consider adding to that a data set (($\mathfrak{b}$, $\mathfrak{c}$), ($\mathfrak{b}_0$, $\mathfrak{c}_0$)) with $\mathfrak{b}$ and $\mathfrak{c}$ denoting 1-forms with values in the associated Lie algebra bundle and with $\mathfrak{b}_0$ and $\mathfrak{c}_0$ denoting sections of that same bundle. The resulting data set will solve (4.1) and be transveral to the orbit through (A, $\mathfrak{a}$) of the automorphism group of P if $\psi \equiv ((\mathfrak{b}, \mathfrak{c}), (\mathfrak{b}_0, \mathfrak{c}_0))$ solves equation that has the schematic form below (the terms are described subsequently):

$$\mathcal{D}\psi + m\,\psi_\mathfrak{c} + \psi\#\psi + \mathcal{G} = 0 \,.$$

(4.3)

To elaborate: What is denoted above by $\mathcal{D}$ denotes a first order, essentially self-adjoint differential operator whose zero'th order part depends the choice for ($\hat{A}$, $\hat{\mathfrak{a}}$). More is said about $\mathcal{D}$ in the rest of this section since it plays the central role when solving (4.3). What is denoted by $\psi_c$ denotes the image of $\psi$ via the linear endomorphism of $\mathbb{S}$ that sends any given element $\psi \equiv ((\mathfrak{b}, \mathfrak{c}), (\mathfrak{b}_0, \mathfrak{c}_0))$ to $((0, \mathfrak{c}), (0, \mathfrak{c}_0))$. Meanwhile, what is denoted by # in (4.3) signifies a certain covariantly constant, symmetric, bilinear bundle map. To say more about this bundle map: Let $\mathbb{S}$ denote the vector bundle $(\oplus_2 T^*Y \oplus (\oplus_2 \mathbb{R})) \otimes ad(P)$ which is where $\psi$ is from. What is denoted by # comes from a canonical homomorphism from $\mathbb{S} \oplus \mathbb{S}$ to $\mathbb{S}$ whose components (when writtten with respect to any given orthonormal frame for $T^*Y$ are commutators of the respective components of the left and right $\mathbb{S}$ summands in $\mathbb{S} \oplus \mathbb{S}$. Finally, what is denoted by $\mathcal{G}$ in (4.3) is a given section of $\mathbb{S}$. For example, in the context of (4.1), the section $\mathcal{G}$ comes from the left hand side of (4.1) (interpreted as a section of $\mathbb{S}$) when (4.1) has $((A, \mathfrak{a}), (\mathfrak{p}, \mathfrak{q}))$ being the initial data set $((\hat{A}, \hat{\mathfrak{a}}), (\mathfrak{p} = 0, \mathfrak{q} = 0))$.

The perturbative strategy for solving (4.3) starts with the observation that if $\mathcal{G}$ is small in a suitable Hilbert space, and if the linear operator $\mathcal{L} \equiv \mathcal{D} + (\cdot)_c$ is invertible with the norm of the inverse not small relative to the size of $\mathcal{G}$, then $\psi_1 = -\mathcal{L}^{-1}(\mathcal{G})$ will be small and $\psi_1 \# \psi_1$ is (hopefully) much smaller than $\mathcal{G}$. If so, then the map $\psi \to -\mathcal{L}^{-1}(\psi \# \psi + \mathcal{G})$ should be a contraction mapping from a suitably small radius Hilbert space ball to itself. And, if that is so, then this map will have a unique fixed point in that ball, and that fixed point will be (at least in a weak sense) a solution to (4.3).

As it turns out, the versions of $\mathcal{L}$ and $\mathcal{G}$ for the approximate solution from Proposition 3.4 are not suitable for this perturbative approach because, although $\mathcal{G}$ is small, the norm of the corresponding $\mathcal{L}^{-1}$ is very large (it is $\mathcal{O}(mr^2)$). But as explained in Section 5 (see Proposition 5.3), the pair from Proposition 3.4 can be amended slightly so that the resulting version of $\mathcal{G}$ is microscopic when $r$ is large in the sense that the $S^1 \times \Sigma$ integral of $|\mathcal{G}|^2$ is $\mathcal{O}(e^{-r/c_0})$. Sections 6 and 7 then prove that the norm of the corresponding version of $\mathcal{L}$ is still $\mathcal{O}(mr^2)$. As explained in Section 6a, these bounds are more than sufficient to implement the perturbative approached outlined above.

With the preceding as background, the rest of this subsection derives some preliminary facts about the operator $\mathcal{D}$ (which is the dominant part of $\mathcal{L}$).

## b) The operator $\mathcal{D}$

The operator $\mathcal{D}$ plays the starring role in what follows and this is why it is depicted explicitly below in (4.4). To set the notation, suppose that $\psi = ((\mathfrak{b}, \mathfrak{c}), (\mathfrak{b}_0, \mathfrak{c}_0))$ denotes a section of $\mathbb{S}$. Write the $\mathbb{S}$ components of $\mathcal{D}\psi$ as $((\mathfrak{b}', \mathfrak{c}'), (\mathfrak{b}'_0, \mathfrak{c}'_0))$. These are as follows:

- $\mathfrak{b}' = *(d_{\hat{A}}\mathfrak{c} + \hat{\mathfrak{a}} \wedge \mathfrak{b} + \mathfrak{b} \wedge \hat{\mathfrak{a}}) - d_{\hat{A}}\mathfrak{c}_0 + [\hat{\mathfrak{a}}, \mathfrak{b}_0]$ .
- $\mathfrak{c}' = *(d_{\hat{A}}\mathfrak{b} - \hat{\mathfrak{a}} \wedge \mathfrak{c} - \mathfrak{c} \wedge \hat{\mathfrak{a}}) - d_{\hat{A}}\mathfrak{b}_0 - [\hat{\mathfrak{a}}, \mathfrak{c}_0]$ .
- $\mathfrak{b}_0' = *(d_{\hat{A}}*\mathfrak{c} + \mathfrak{b} \wedge *\hat{\mathfrak{a}} - *\hat{\mathfrak{a}} \wedge \mathfrak{b})$ .

- $\mathfrak{c}_0{}' = *(d_{\hat{A}}*\mathfrak{b} - \mathfrak{c} \wedge *\hat{\mathfrak{a}} + *\hat{\mathfrak{a}} \wedge \mathfrak{c})$.

$$(4.4)$$

The operator $\mathcal{D}$ can be written concisely as

$$\mathcal{D} = \gamma_i \nabla_{\hat{A}i} + \rho_i[\hat{\mathfrak{a}}_i, \cdot]$$

$$(4.5)$$

where $\{\gamma_1, \gamma_2, \gamma_3, \rho_1, \rho_2, \rho_3\}$ are 6 of the 7 generators of the action via covariantly constant endomorphisms of a 7-dimensional (real) Clifford action on $(\oplus_2 T^*Y \oplus (\oplus_2 \mathbb{R}))$. In particular, these endomorphisms are antisymmetric (so $\gamma_i{}^T = -\gamma_i$ and $\rho_i{}^T = -\rho_i$) and obey the following algebraic relations:

- $\gamma_i\gamma_j + \gamma_j\gamma_i = -2\delta_{ij}$
- $\rho_i\rho_j + \rho_j\rho_i = -2\delta_{ij}$
- $\rho_i\gamma_j + \gamma_j\rho_i = 0$

$$(4.6)$$

It follows from (4.5) and (4.6) that $\mathcal{D}$ is a symmetric, elliptic operator on the space of sections of $\mathbb{S}$. The next lemma makes a precise statement with the specification of appropriate Hilbert spaces.

**Lemma 4.1**: *Let $\mathbb{L}$ and $\mathbb{H}$ denote the respective Hilbert space completions of the space of smooth sections of $\mathbb{S}$ using the norms whose square are the functionals*

$$\psi \to \int_Y |\psi|^2 \quad and \quad \psi \to \int_Y (|\nabla_{\hat{A}}\psi|^2 + |[\hat{\mathfrak{a}},\psi]|^2 + |\psi|^2)$$

*The operator $\mathcal{D}$ defines an unbounded, self-adjoint operator on $\mathbb{L}$ with discrete spectrum having finite multiplicities and no accumulation points. The space $\mathbb{H}$ is a dense domain for $\mathcal{D}$ and $\mathcal{D}$ defines a bounded, Fredholm map from $\mathbb{H}$ to $\mathbb{L}$ with index zero.*

**Proof of Lemma 4.1**: The assertion is a standard consequence of the first order nature of the operator $\mathcal{D}$ and its ellipticity (which follows from (4.5) and (4.6)). See for example [H], Chapter XIX.

The respective $\mathbb{L}$ and $\mathbb{H}$ norms are denoted below by $\|\cdot\|_{\mathbb{L}}$ and $\|\cdot\|_{\mathbb{H}}$.

By way of a parenthetical remark: The linearized operator that is depicted in (4.4) is not directly the operator that comes from linearizing the equations for a flat $SL(2; \mathbb{C})$ connection and adding a slice condition forcing transversality to the orbit of the automorphism group of the principle bundle. That operator is depicted below:

- $\mathfrak{c}' = *(d_{\hat{A}}\mathfrak{c} + \hat{\mathfrak{a}} \wedge \mathfrak{b} + \mathfrak{b} \wedge \hat{\mathfrak{a}}) - d_{\hat{A}}\mathfrak{c}_0 + [\hat{\mathfrak{a}}, \mathfrak{b}_0]$.

- $\mathfrak{b}' = *(d_{\hat{A}}\mathfrak{b} - \hat{\mathfrak{a}} \wedge \mathfrak{c} - \mathfrak{c} \wedge \hat{\mathfrak{a}}) - d_{\hat{A}}\mathfrak{b}_0 - [\hat{\mathfrak{a}}, \mathfrak{c}_0]$ .
- $\mathfrak{c}_0' = *(d_{\hat{A}}*\mathfrak{c} + \mathfrak{b} \wedge *\hat{\mathfrak{a}} - *\hat{\mathfrak{a}} \wedge \mathfrak{b})$
- $\mathfrak{b}_0' = *(d_{\hat{A}}*\mathfrak{b} - \mathfrak{c} \wedge *\hat{\mathfrak{a}} + *\hat{\mathfrak{a}} \wedge \mathfrak{b})$

$$(4.7)$$

But do note that the operator in (4.7) can be written as $\psi \to \hat{I}\mathcal{D}$ where $\hat{I}$ here denotes the endomorphism of S that sends any given $((\mathfrak{b}, \mathfrak{c}), (\mathfrak{b}_0, \mathfrak{c}_0))$ to $((\mathfrak{c}, \mathfrak{b}), (\mathfrak{c}_0, \mathfrak{b}_0))$. Since $\hat{I}$ is symmetric, and since its square is the identity, what is said about $\mathcal{D}$ holds for (4.7) as well.

## c)  The Bochner-Weitzenboch formula for $\mathcal{D}$

A fundamental tool for probing the spectrum of any given $(\hat{A}, \hat{\mathfrak{a}})$ version of the operator $\mathcal{D}$ is its Bochner-Weitzenboch formula which is this

$$\mathcal{D}^2 = -\nabla_{\hat{A}k}\nabla_{\hat{A}k} + [\hat{\mathfrak{a}}_k, [\cdot, \hat{\mathfrak{a}}_k]] + \frac{1}{2}[\gamma_i, \gamma_j][F_{\hat{A}ij}, \cdot] + \gamma_i \rho_j[\nabla_{\hat{A}i}\hat{\mathfrak{a}}_k, \cdot] + \mathfrak{R}$$

$$(4.8)$$

where the notation has $\{F_{\hat{A}ij}\}_{i,j=1,2,3}$ denoting the components of the curvature tensor of $\hat{A}$, and where $\{\nabla_{\hat{A}i}\hat{\mathfrak{a}}_j\}_{i,j=1,2,3}$ are the components of the covariant derivative of $\hat{\mathfrak{a}}$. Meanwhile, what is denoted by $\mathfrak{R}$ is an endomorphism of $\mathbb{S}$ that comes from the Riemann curvature tensor. Just as useful is an integrated form of (4.16) which says that if $\psi$ is any given element in $\mathbb{H}$, then

$$\int_Y |\mathcal{D}\psi|^2 = \int_Y \left(|\nabla_{\hat{A}}\psi|^2 + |[\hat{\mathfrak{a}},\psi]|^2\right) + \int_Y \left(\langle\psi, \frac{1}{2}[\gamma_i,\gamma_j][F_{\hat{A}ij},\psi] + \gamma_i\rho_j[\nabla_{\hat{A}i}\hat{\mathfrak{a}}_k, \psi]\rangle + \langle\psi, \mathfrak{R}\psi\rangle\right)$$

$$(4.9)$$

The next subsection uses the Bochner-Weitzenboch formula and the $\mathbb{H}$-norm to make some preliminary observations about the spectrum of $\mathcal{D}$ in the case when the pair $(\hat{A}, \hat{\mathfrak{a}})$ is described in Section 3.

## d)  Preliminary observations about the operator $\mathcal{D}$ in the case at hand

Fix a large positive integer (divisible by 4) to be denoted by n and then, if $m \neq 0$, set $r = \frac{n}{m}$. If $m = 0$, take $r$ to be a large positive number. In what follows, the pair $(\hat{A}, \hat{\mathfrak{a}})$ is taken to be a very large $r$ version of either the pair depicted in (3.20) or the pair depicted in (3.25) and Proposition 3.4, or the $\mu = \mu_r$ version of the $m = 0$ pair that is depicted in (3.2) with $\mu_r$ comimg from Lemma 3.1. In this regard, $\hat{A}$ can and will be viewed as a connection on a principal SU(2) bundle over $S^1 \times \Sigma$ whose associated $\mathbb{C}^2$ bundle is either $E_q$ or the bundle $L \oplus L^{-1}$ as the case may be. This bundle is denoted by P in what follows. When viewed in this light, $\hat{\mathfrak{a}}$ is a 1-form on $S^1 \times \Sigma$ with values in ad(P) (which is the associated vector bundle with fiber being the Lie algebra of SU(2).)

The upcoming (4.10) summarizes some of the salient features of each instance of $(\hat{A}, \hat{a})$. By way of notation: The inequalities below are written on open sets of the form $S^1 \times D'$ with $D' \subset \Sigma$ denoting a given open disk; and they refer to an oriented orthonormal frame for $T^*(S^1 \times \Sigma)$ on $S^1 \times D'$ that has the form $\{dt, e^1, e^2\}$ where $\{e^1, e^2\}$ constitute an oriented orthonormal frame for $T^*\Sigma|_{D'}$.

- *The size of $\hat{a}$*:
  - $|\hat{a}_1| + |\hat{a}_2| \geq c_0^{-1} r^{2/3} + r \, \mathfrak{d}^{1/2}$
  - $|\hat{a}_t| \leq c_0 m \, e^{-r \mathfrak{d}^{3/2}/c_0}$
- *The size of the covariant derivatives of $\hat{a}$*:
  - $|\nabla_{\hat{A}} \hat{a}| \leq c_0 r^{4/3} \, e^{-r \mathfrak{d}^{3/2}/c_0} + c_0 r \, \mathfrak{d}^{1/2}$
  - $|\nabla_{\hat{A}t} \hat{a}| + |\nabla_{\hat{A}1} \hat{a}_t| + |\nabla_{\hat{A}2} \hat{a}_t| \leq c_0 m \, r^{2/3} e^{-r \mathfrak{d}^{3/2}/c_0}$ .
  - $|\nabla_{\hat{A}t} \hat{a}_t| \leq c_0 m^2 \, e^{-r \mathfrak{d}^{3/2}/c_0}$ .
- *The size of the curvature tensor of $\hat{A}$*:
  - $|F_{\hat{A}12}| \leq c_0 r^{4/3} \, e^{-r \mathfrak{d}^{3/2}/c_0}$
  - $|F_{\hat{A}t1}| + |F_{\hat{A}t2}| \leq c_0 m \, (r^{2/3} + r \, \mathfrak{d}^{1/2})$ .
  - $|F_{\hat{A}t1}| + |F_{\hat{A}t2}| \geq c_0^{-1} m \, (r^{2/3} + r \, \mathfrak{d}^{1/2})$

$$(4.10)$$

In addition: For each of the instances of $(\hat{A}, \hat{a})$, there exists an isometric homomorphism defined on $S^1 \times (\Sigma - \Theta)$ from the line bundle $\mathcal{I}$ into the bundle ad(P). The image is a unit length, $\mathcal{I}$-valued section of ad(P) which is denoted by $\sigma$ in what follows. The key point about $\sigma$ is that it is nearly $\hat{A}$-covariantly constant and nearly commutes with $\hat{a}$ at distances greater than $c_0 r^{-2/3}$ from $S^1 \times \Theta$. Indeed, the following hold where $\mathfrak{d} > c_0 r^{-2/3}$:

- $|\nabla_{\hat{A}} \sigma| + |[\hat{a}, \sigma]| \leq c_0 \, r^{2/3} \, e^{-r \mathfrak{d}^{3/2}/c_0}$ .
- $|\nabla_{\hat{A}} \nabla_{\hat{A}} \sigma| + |[F_{\hat{A}}, \sigma]| + |[\nabla_{\hat{A}} \hat{a}, \sigma]| \leq c_0 r^{4/3} e^{-r \mathfrak{d}^{3/2}/c_0}$ .

$$(4.11)$$

To be sure: When $m = 0$, what is denoted by $\sigma$ is the section $\hat{\tau}$ from (2.1). When $m \neq 0$, the section $\sigma$ is the section $\hat{\tau}$ from (2.1) when $\sigma$ is written using the $\mathbb{C} \oplus \mathbb{C}$ product structure for $E_q$ near any $p \in \Theta$ version of $S^1 \times p$. On $S^1 \times (\Sigma - \Theta)$ where $E_q$ is the quotient of the bundle $(\pi_q^*L \otimes L_q) \oplus (\pi_q^*L^{-1} \otimes (-1)^*L_q)$ on $S^1 \times (\Sigma_q - \Theta)$, the section $\sigma$ is the image of the endomorphism $\sigma_3$ of $(\pi_q^*L \otimes L_q) \oplus (\pi_q^*L^{-1} \otimes (-1)^*L_q)$ that acts as $+i$ on the left most summand and $-i$ on the right most.)

The lemma the follows (and subsequent lemmas) exploit the preceding bounds.

**Lemma 4.2**: *There exists* $\kappa > 1$ *with the following significance: Suppose that* $r > \kappa$ *and that* $(\hat{A}, \hat{a})$ *is as described above using this value for* $r$. *Then, the inequality below holds for any given section* $\psi$ *of* $\mathbb{S}$.

$$\int_{S^1 \times \Sigma} \left(|\nabla_{\hat{A}} \psi|^2 + |[\hat{a}, \psi]|^2\right) \geq \kappa^{-1} \int_{S^1 \times \Sigma} \left(\frac{r^{4/3}}{(1 + r^{4/3} \mathfrak{d}^2)} |\psi|^2 + r^2 \mathfrak{d} |[\sigma, \psi]|^2\right)$$

The proof is given momentarily.

The lemma that follows implies that any eigenvector of $\mathcal{D}$ with eigenvalue in the range $[-c_0^{-1}, c_0^{-1}]$ can't be too small where the distance to $\Theta$ is on the order of $r^{-2/3}$. It also implies that such an eigenvector can't differ much from a section along $\sigma$ where the distance to $\Theta$ is greater than $\mathcal{O}(r^{-2/3})$.

**Lemma 4.3**: *There exists* $\kappa > 1$ *with the following significance: Fix* $r > \kappa$ *and let* $(\hat{A}, \hat{a})$ *be as described in the first paragraphs of this subsection using this value for* $r$. *Let* $\Sigma_r$ *denote the part of* $\Sigma$ *where* $\mathfrak{d}$ *(the distance to* $\Theta$*) is less than* $\kappa r^{-2/3}$. *Suppose that* $\psi$ *denotes a section of* $\mathbb{S}$ *from the Hilbert space* $\mathbb{H}$ *that obeys the inequality below:*

$$r^{4/3} \int_{S^1 \times \Sigma_r} |\psi|^2 \leq \frac{1}{\kappa} \int_{S^1 \times (\Sigma - \Sigma_r)} \left(\frac{1}{\mathfrak{d}^2} |\psi|^2 + r^2 \mathfrak{d} |[\sigma, \psi]|^2\right).$$

*If this is so, then the* $S^1 \times \Sigma$ *integral of* $|\mathcal{D}\psi|^2$ *is such that*

$$\int_{S^1 \times \Sigma} |\mathcal{D}\psi|^2 \geq \frac{1}{\kappa} \int_{S^1 \times \Sigma} \left(|\nabla_{\hat{A}} \psi|^2 + |[\hat{a}, \psi]|^2\right) + \frac{1}{\kappa} \int_{S^1 \times \Sigma} |\psi|^2 .$$

*In any event,*

$$\int_{S^1 \times \Sigma} |\mathcal{D}\psi|^2 + r^{4/3} \int_{S^1 \times \Sigma_r} |\psi|^2 \geq \frac{1}{\kappa} \int_{S^1 \times \Sigma} \left(|\nabla_{\hat{A}} \psi|^2 + |[\hat{a}, \psi]|^2\right).$$

The proof of this lemma follows directly; then comes the proof of Lemma 4.2.

*Proof of Lemma 4.3*: Use the pointwise bounds in (4.10) and (4.11) with the Bochner-Weitzenboch formula from (4.9) to see that there exists a number $c_\diamond \in (1, c_0)$ with the following significance: Let $\Sigma_r$ denote the $\mathfrak{d} < c_\diamond r^{-2/3}$ part of $\Sigma$. Supposing that $\psi$ signifies a section of $\mathbb{S}$, then

$$\int_{S^1 \times \Sigma} |\mathcal{D}\psi|^2 \geq \int_{S^1 \times \Sigma} \left(|\nabla_{\hat{A}} \psi|^2 + |[\hat{a}, \psi]|^2\right) - c_0 r^{4/3} \int_{S^1 \times \Sigma_r} |\psi|^2$$

(4.12)

With the preceding inequality in hand and supposing Lemma 4.2 is true, invoke Lemma 4.2 to see that

$$\int_{S^1 \times \Sigma} |\mathcal{D}\psi|^2 \geq \tfrac{1}{2} \int_{S^1 \times \Sigma} \left( |\nabla_{\hat{A}}\psi|^2 + |[\hat{a},\psi]|^2 \right) + c_0^{-1} \int_{S^1 \times (\Sigma - \Sigma_r)} \left( \tfrac{1}{\mathfrak{d}^2}|\psi|^2 + r^2 \mathfrak{d} |[\sigma,\psi]|^2 \right)$$
$$- c_0 r^{4/3} \int_{S^1 \times \Sigma_r} |\psi|^2 .$$
$$(4.13)$$

This inequality and the inequality in Lemma 4.2 lead directly to Lemma 4.3's assertion.

***Proof of Lemma 4.2***:  Fix for the moment $c_1 > c_0$ and write $\psi$ where $\mathfrak{d} > c_1 r^{2/3}$ as

$$\psi = \langle \sigma \psi \rangle \sigma + \psi_\perp$$
$$(4.14)$$

with $\langle \psi_\perp \sigma \rangle$ being zero and with $\psi_\perp = \tfrac{1}{4}[\sigma, [\psi, \sigma]]$.  It follows from this and from (4.11) that

- $\nabla_{\hat{A}}\psi = (d\langle \sigma \psi \rangle) \, \sigma + (\nabla_{\hat{A}}\psi)_\perp + \mathfrak{e}_1$
- $[\hat{a},\psi] \geq c_0^{-1} r \, \varsigma \, [\sigma,\psi] + \mathfrak{e}_2$

$$(4.15)$$

where $|\mathfrak{e}_1| + |\mathfrak{e}_2|$ are bounded by $c_0 \, r^{2/3} \, e^{-r\mathfrak{d}^{3/2}/c_0} \, |\psi|$ and where $\varsigma$ denotes an $\mathcal{I}$-valued 1-form with norm obeying $|\varsigma| \geq c_0^{-1} \mathfrak{d}^{1/2}$.  Let $\Sigma_{1r}$ denote the subset in $\Sigma$ where $\mathfrak{d} < c_1 r^{-2/3}$.  Use (4.15) to see that

$$\int_{S^1 \times \Sigma} \left( |\nabla_{\hat{A}}\psi|^2 + |[\hat{a},\psi]|^2 \right) \geq c_0^{-1} \int_{S^1 \times (\Sigma - \Sigma_{1r})} \left( |d\langle \sigma \psi \rangle|^2 + r^2 \mathfrak{d} |[\sigma,\psi]|^2 - c_0 r^{4/3} e^{-r\mathfrak{d}^{3/2}/c_0} \, |\psi|^2 \right) .$$
$$(4.16)$$

To continue:  Reintroduce $r_0$ from Part 3 of Section 3a.  In this regard:  Keep in mind that the embedded disks centered at the points in $\Sigma$ given by the $|z| < 2r_0$ disk in $\mathbb{C}$ using the holomorphic coordinate charts described after Lemma 3.1 and Lemma 3.2 are pairwise disjoint and each such disk is well inside the domain of a Gaussian coordinate chart based at its center point.  Supposing that $r$ is positive but less than $r_0$, let $f$ denote for the moment a section of $\mathcal{I}$ on the radius $r$ circle centered at a point in $\Theta$.  Because the line bundle $\mathcal{I}$ is isomorphic to the Möbius line bundle on this radius $r$ circle, the integral of $|df|^2$ on that circle is no smaller than $c_0^{-1} \, r^{-2}$ times the integral of $|f|^2$ on the same circle.  In particular, this is true for the case when $f = \langle \sigma \psi \rangle$.  With the preceding point understood, let $\Sigma_0$ denote the subset in $\Sigma$ where the distance to $\mathfrak{d}$ is less than $r_0$.  Then (4.16) implies this:

$$\int_{S^1 \times \Sigma} \left( |\nabla_{\hat{A}}\psi|^2 + |[\hat{a},\psi]|^2 \right) \geq c_0^{-1} \int_{S^1 \times (\Sigma_0 - \Sigma_{1r})} \tfrac{1}{\mathfrak{d}^2} |\langle \sigma \psi \rangle|^2$$
$$+ c_0^{-1} \int_{S^1 \times (\Sigma - \Sigma_{1r})} \left( r^2 \mathfrak{d} |[\sigma,\psi]|^2 - c_0 r^{4/3} e^{-r\mathfrak{d}^{3/2}/c_0} \, |\psi|^2 \right) .$$
$$(4.17)$$

Since $r^2\mathfrak{d} > \mathfrak{d}^{-2}$ where $\mathfrak{d} > c_0 r^{-2/3}$, the preceding inequality implies in turn the one below which has $|\langle\sigma\psi\rangle|^2$ in the left most integral on the right hand side of (4.17) replaced by $|\psi|^2$:

$$\int_{S^1\times\Sigma} \left(|\nabla_{\hat{A}}\psi|^2 + |[\hat{a},\psi]|^2\right) \geq c_0^{-1}\int_{S^1\times(\Sigma_0-\Sigma_{1r})} \frac{1}{\mathfrak{d}^2}|\psi|^2$$
$$+ c_0^{-1}\int_{S^1\times(\Sigma-\Sigma_{1r})} \left(r^2\mathfrak{d}|[\sigma,\psi]|^2 - c_0 r^{4/3}e^{-r\mathfrak{d}^{3/2}/c_0}|\psi|^2\right).$$
(4.18)

Now let $\Sigma_{2r}$ denote the subset in $\Sigma$ where the distance to $\Theta$ is less than $2c_1 r^{-2/3}$. Since $|\nabla_{\hat{A}}\psi| > |d|\psi||$, a standard eigenvalue inequality for the Laplacian on the set $\Sigma_{2r}$ can be invoked to see that:

$$\int_{S^1\times\Sigma_{2r}} |\nabla_{\hat{A}}\psi|^2 + r^{4/3}\int_{S^1\times(\Sigma_{2r}-\Sigma_{1r})} |\psi|^2 \geq c_0^{-1}r^{4/3}\int_{S^1\times\Sigma_{1r}} |\psi|^2.$$
(4.19)

Using this with (4.18) leads to the inequality below which has $\Sigma_0-\Sigma_{1r}$ in the left most integral on the right hand side of (4.18) replaced by $\Sigma_0$:

$$\int_{S^1\times\Sigma} \left(|\nabla_{\hat{A}}\psi|^2 + |[\hat{a},\psi]|^2\right) \geq c_0^{-1}\int_{S^1\times\Sigma_0} \frac{r^{4/3}}{(1+r^{4/3}\mathfrak{d}^2)}|\psi|^2$$
$$+ c_0^{-1}\int_{S^1\times\Sigma} r^2\mathfrak{d}|[\sigma,\psi]|^2 - c_0\int_{S^1\times(\Sigma-\Sigma_{1r})} r^{4/3}e^{-r\mathfrak{d}^{3/2}/c_0}|\psi|^2.$$
(4.20)

An analogous eigenvalue inequality for the standard Laplacian on the domain where $\mathfrak{d} > \frac{1}{2}r_0$ leads from (4.20) to the inequality

$$\int_{S^1\times\Sigma} \left(|\nabla_{\hat{A}}\psi|^2 + |[\hat{a},\psi]|^2\right) \geq c_0^{-1}\int_{S^1\times\Sigma} \frac{r^{4/3}}{(1+r^{4/3}\mathfrak{d}^2)}|\psi|^2$$
$$+ c_0^{-1}\int_{S^1\times\Sigma} r^2\mathfrak{d}|[\sigma,\psi]|^2 - c_0\int_{S^1\times(\Sigma-\Sigma_{1r})} r^{4/3}e^{-r\mathfrak{d}^{3/2}/c_0}|\psi|^2$$
(4.21)

which has $\Sigma_0$ in the left-most integral on the right hand side of (4.20) replaced by $\Sigma$.

To finish: It follows now that if $c_1 > c_0$, then the right most integral on the right hand side of (4.21) will be much smaller than the other two terms on the right hand side of (4.21). Thus, if $c_1$ is chosen large (but still $c_1 < c_0$), then (4.21) leads to the inequality stated by Lemma 4.2.

### e) A first look at sections of $\mathbb{S}$ with constrained $\mathcal{D}\psi$

Sections 5 and 6 take a close look at the eigenvectors and eigenvalues of various versions of the operator $\mathcal{D}$. Lemmas 4.2 and 4.3 and the upcoming Lemma 4.4 supply tools that are used in those sections. In this regard: Lemma 4.4 goes beyond what is implied by Lemmas 4.2 and 4.3 by supplying 'semi-local' bounds for norm of $[\sigma, \psi]$ subject to the constraint on $\mathcal{D}\psi$ stated below in (4.22). This constraint involves a number denoted by $c_\ddagger$ which is greater than 1, two positive numbers denoted by $c_1$ and $c_2$ and two non-negative numbers denoted by $C_1$ and $C_2$. With these numbers set, the constraint for an appeal to Lemma 4.4 requires that the inequality below should hold where $\eth > c_0 r^{-2/3}$:

$$|[\sigma, \mathcal{D}\psi]| \le \tfrac{1}{c_\ddagger} r \, \eth^{1/2} \, |[\sigma, \psi]| + C_1 |\psi| e^{-r \eth^{3/2}/c_1} + C_2 \, e^{-r \eth^{3/2}/c_2} \; .$$

(4.22)

An important example: The inequality in (4.22) is obeyed with $C_1 = C_2 = 0$ when $\psi$ is an eigenvector of $\mathcal{D}$ whose eigenvalue norm is at most than $c_\ddagger^{-1} r^{2/3}$.

Lemma 4.4 uses the following notation: Given a number $\rho > 0$, the lemma uses $\Sigma_{\rho, 2\rho}$ to denote the part of $\Sigma$ where $\eth$ is between $\rho$ and $2\rho$.

**Lemma 4.4**: *There exists* $\kappa > 1$ *and* $\kappa_* > \kappa$ *with the following significance: Fix* $r > \kappa$ *and let* $(\hat{A}, \hat{a})$ *denote a pair of connection on* $P$ *and 1-form with values in* $ad(P)$ *that is described in the first paragraphs of Section 4e (and thus obeying (4.10) and (4.11)). Suppose that* $\psi$ *obeys the constraing in (4.22) for all* $\rho > \kappa r^{-2/3}$ *with* $c_\ddagger > \kappa$ *and with any given values for* $c_1$ *and* $c_2$ *and* $C_1$ *and* $C_2$. *Then the inequality below holds for any number* $\rho > \kappa r^{-2/3}$:

$$\int_{S^1 \times \Sigma_{\rho, 2\rho}} |[\sigma, \psi]|^2 \le \kappa_* \big( (1 + C_1^2) \, e^{-r \rho^{3/2}/\kappa} + C_1^2 \, e^{-r \rho^{3/2}/\kappa c_1} \big) \|\psi\|_{\mathbb{L}}^2$$
$$+ \kappa_* \, C_2^2 \, r^{-4/3} \rho^2 (e^{-r \rho^{3/2}/\kappa} + e^{-r \rho^{3/2}/\kappa c_2}).$$

As noted, (4.22) holds with $C_1 = C_2 = 0$ when $\mathcal{D}\psi = \lambda \psi$ with $|\lambda| < c_0^{-1} r^{2/3}$. Thus, the conclusions of Lemma 4.4 hold for these eigenvectors of $\mathcal{D}$.

***Proof of Lemma 4.4*** The argument is a riff on the argument that is used to prove Proposition 4.1 in Section 4e of [T2]. Having fixed a number $c_\diamond > c_0$, let N denote the smallest integer such that $c_\diamond N^{2/3} r^{-2/3}$ is larger than the diameter of $\Sigma$. For each integer k starting at 2 and ending at N, use the standard bump function $\chi$ to construct a non-negative function to be denoted $\chi_k$ that is equal to 1 were $\eth$ is between $c_\diamond \, k^{2/3} r^{-2/3}$ and $c_\diamond (k+1)^{2/3} r^{-2/3}$, and equal to 0 where $\eth$ is either smaller than $c_\diamond \, (k-1)^{2/3} r^{-2/3}$ or larger than $c_\diamond \, (k+2)^{2/3} r^{-2/3}$. This function $\chi_k$ can and should be constructed so that $|d\chi_k|$ is less than $c_0 \, c_\diamond^{-1} \, k^{1/3} \, r^{2/3}$. (An important point for what is to come: The rules for $\chi_k$ have $d\chi_k$ supported inside the union of the respectives sets where $\chi_{k-1}$ and $\chi_{k+1}$ are equal to 1.)

Since the equation

$$\mathcal{D}[\sigma,\psi] = [\sigma, \mathcal{D}\psi] + \gamma_i[\nabla_{\hat{A}i}\sigma,\psi] - \rho_i\,[[\sigma, \hat{a}_i],\psi]$$

(4.23)

always holds, any given integer k version of $\chi_k[\sigma, \psi]$ obeys the following identity:

$$\mathcal{D}(\chi_k[\sigma,\psi]) = \chi_k[\sigma, \mathcal{D}\psi] + \chi_k\gamma_i[\nabla_{\hat{A}}\sigma,\psi] - \chi_k\rho_i\,[[\sigma, \hat{a}_i],\psi] + (\gamma_i\nabla_i\chi_k)(\chi_{k-1}[\sigma,\psi] + \chi_{k+1}[\sigma,\psi])\,.$$

(4.24)

Integrate the square of the norm of both sides in (4.24) over $S^1 \times \Sigma$ and then use (4.20) and the fact that any given $k \geq 2$ version of $d\chi_k$ is supported only on the union of the sets where $\chi_{k-1}$ and $\chi_{k+1}$ are equal to 1 to derive the following integral inequality:

$$\int_{S^1\times\Sigma} |\mathcal{D}(\chi_k[\sigma,\psi])|^2 \leq c_0\int_{S^1\times\Sigma} \chi_k{}^2|[\sigma, \mathcal{D}\psi]|^2 + c_0 r^{4/3}\,e^{-k/c_0}\int_{S^1\times\Sigma} \chi_k{}^2|\psi|^2 +$$
$$c_0\,c_0{}^{-2}\,k^{2/3}r^{4/3}\big(\int_{S^1\times\Sigma} \chi_{k-1}{}^2|[\sigma,\psi]|^2 + \int_{S^1\times\Sigma} \chi_{k+1}{}^2|[\sigma,\psi]|^2\big)$$

(4.25)

Suppose now that $\psi$ obeys a $c_{\ddagger} > c_0$ version of (4.22) and that $c_0$ is large (equal to $100c_0$). If $k \geq c_0$, then the integrated version of the Bochner-Weitzenboch formula in (4.9) (with $\chi_k[\sigma,\psi]$ playing $\psi$'s role in (4.9)) can be used with the bounds in (4.10) to see that (4.25) directly implies the inequality below:

$$10 \int_{S^1\times\Sigma} \chi_k{}^2|[\sigma,\psi])|^2 - \int_{S^1\times\Sigma} \chi_{k-1}{}^2|[\sigma,\psi]|^2 - \int_{S^1\times\Sigma} \chi_{k+1}{}^2|[\sigma,\psi]|^2 \leq$$
$$c_0(e^{-k/c_0} + C_1{}^2e^{-k/c_0c_1})\int_{S^1\times\Sigma} |\psi|^2 + c_0\,C_2{}^2\,r^{-8/3}\,e^{-k/c_0c_2}\,.$$

(4.26)

The preceding inequality has the look of a 2-term, inhomogeneous recursion inequality which can be made manifest by introducing $x_k$ to denote the $S^1 \times \Sigma$ integral of $\chi_k{}^2|[\sigma,\psi])|^2$. With this notation in hand, (4.26)

$$10\,x_k - x_{k-1} - x_{k+1} \leq c_0(e^{-k/c_0} + C_1{}^2e^{-k/c_0c_1})\|\psi\|_{\mathbb{L}}{}^2 + c_0\,C^2\,r^{-8/3}\,e^{-k/c_0c_2}\,.$$

(4.27)

Let $k_1$ denote a fixed $c_0$-size integer. The recursion inequality in (4.27) implies in turn the inequality below for $x_k$:

$$x_k \leq e^{-k/c_0}\,x_{k_1} + c_0(e^{-k/c_0} + C_1{}^2e^{-k/c_0c_1})\,\|\psi\|_{\mathbb{L}}{}^2 + c_0C^2r^{-8/3}(e^{-k/c_0} + e^{-k/c_0c_2}).$$

(4.28)

This inequality leads directly to the lemma's inequality. (With regards to the recursion relation: The left hand side of (4.27) is zero when $x_k$ is a constant multiple of $r^k$ with $r$ being a root of the quadratic equation $r^2 - 10r + 1 = 0$. These roots are $5 \pm \sqrt{24}$. One of the roots is less than 1 (it is approximately 1/10); and the other is greater than one (it is greater than 9). The relevant root is the smaller one since $x_N$ is zero.)

## 5. The $\mathfrak{m} \neq 0$ equations on each Riemann surface slice

Three salient results come from this section: The first result is a proof of Proposition 3.4. The second result is a slightly modified version of any given small $m$ and large $r$ pair from (3.25) and Proposition 3.4 which is far closer to solving (1.7) than the unmodified pair. This second result is the content of Proposition 5.3. The final result (which is actually the first topic of this section) is closely connected to the first two since it is used to prove Proposition 3.4 and Proposition 5.3. This third result also plays a critical role in Section 6. This third result is a detailed analysis of the eigenvectors with small normed eigenvalue for the version of the operator $\mathcal{D}$ that is defined by any large $r$ version of the SL(2;$\mathbb{C}$)-flat connection pair from (3.2) and Lemma 3.1. (The $r = \frac{\mathfrak{n}}{\mathfrak{m}}$ versions of the pairs from (3.2) and Lemma 3.1 are important by virtue of their prominent appearance in the $\mathfrak{d} < r_0$ parts of the pairs that are described by (3.25) and Proposition 3.4 .)

To set the stage for much of what is to come in this section, fix $r > c_0$ and let $(\hat{A}, \hat{a})$ denote a pair that is described by (3.2) and Lemma 3.1. This pair was initially defined on $\Sigma$ as an SU(2) connection on $L \oplus L^{-1}$ and 1-form on $\Sigma$ with values in the bundle of skew-Hermitian endomorphisms of $L \oplus L^{-1}$. But the pair was also viewed via pull-back by the projection map from $S^1 \times \Sigma$ to $\Sigma$ as a pair of connection and skew-Hermitian endomorphism valued 1-form on $S^1 \times \Sigma$. In the latter guise, it defines a version of the operator $\mathcal{D}$ as in (4.5); and an important point is that this version of $\mathcal{D}$ commutes with the corresponding version of $\nabla_{\hat{A}t}$ (more is said about in the next paragraph). As a consequence, any eigenvector of this version of $\mathcal{D}$ with eigenvalue norm less than $c_0^{-1}$ is annihilated by $\nabla_{\hat{A}t}$. Thus, any such eigenvector is the pull-back from $\Sigma$ of an eigenvector of the t-independent operator (denoted by $\mathfrak{D}_0$) which is defined on any constant $t \in S^1$ slice $\{t\} \times \Sigma \subset S^1 \times \Sigma$ by the rule below:

$$\mathfrak{D}_0 \equiv \gamma_1 \nabla_{\hat{A}1} + \gamma_2 \nabla_{\hat{A}2} + \rho_1[\hat{a}_1, \cdot] + \rho_2[\hat{a}_2, \cdot] .$$

(5.1)

(In this and subsequent equations, the subscripts 1 and 2 refer to the components of $\nabla_{\hat{A}}$ and $\hat{a}$ and the corresponding $\gamma$ and $\rho$ endomorphisms when they are written locally using an oriented, orthonormal frame for $T^*\Sigma$.)

With regards to $\mathcal{D}$ and $\mathfrak{D}_0$ and their eigenvalues: Because the bundle $L \oplus L^{-1}$ and the connection $\hat{A}$ are pulled up from $\Sigma$, any given section of $\mathbb{S}$ can be written as a Fourier series with respect to the $S^1$ action, thus $\psi = \sum_{k \in \mathbb{Z}} e^{ikt} \phi_k$ with each $\phi_k$ being an $S^1$-independent section of $\mathbb{S}$ and with $\nabla_{\hat{A}t}$ acting to send $\psi$ to $\sum_{k \in \mathbb{Z}} ik e^{ikt} \phi_k$. Moreover, distinct k terms in the Fourier decomposition are mutually $\mathbb{L}$ orthogonal, and the $\mathbb{L}$-integral of any given $e^{ikt} \phi_k$ term is $2\pi$ times the $\Sigma$ integral of $|\phi_k|^2$. The operator $\mathcal{D}$ preserves this Fourier decomposition in the sense that $\mathcal{D}(e^{ikt} \phi_k) = e^{ikt}(ik\gamma_t\phi_k + \mathcal{D}_0\phi_k)$. Moreover, since

$$\mathcal{D}^2(e^{ikt} \phi_k) = e^{ikt} (k^2 + \mathcal{D}_0{}^2)\phi_k \,,$$

(5.2)

the norm of any eigenvalue of $\mathcal{D}$ can be written as $(k^2 + \lambda^2)^{1/2}$ with k being an integer and with $\lambda$ being an eigenvalue of $\mathcal{D}_0$. In particular, k must be zero if $\mathcal{D}$'s eigenvalue has norm less than 1; and in this case, the corresponding eigenvector is $S^1$-independent and thus the pull-back from $\Sigma$ of an eigenvector of $\mathcal{D}_0$.

## a) Eigenvectors of $\mathcal{D}_0$ with small eigenvalue

The operator $\mathcal{D}_0$ is depicted in (5.1) with $(\hat{A}, \hat{a})$ described by (3.2) with $\mu$ coming from Lemma 3.1 for large $r$ (in the relevant applications, $r = \frac{n}{m}$ with n being a large positive integer that is divisble by 4).

With regards to the domain of $\mathcal{D}_0$: The operator $\mathcal{D}_0$ will be viewed as an operator on $\Sigma$ acting on the space of $S^1$-invariant sections of the bundle $\mathbb{S}$ (which are uniquely determined by their values any one $t \in S^1$ slice of $S^1 \times \Sigma$). In particular, $\mathcal{D}_0$ defines an unbounded, self-adjoint operator on $\Sigma$'s version of the Hilbert space $\mathbb{L}$ from Lemma 4.1 with discrete spectrum having finite multiplicities and no accumulation points. The lemma below says more about the spectrum of $\mathcal{D}_0$.

**Lemma 5.1**: *There exists* $\kappa > 1$ *with the following signficance: Fix* $r > \kappa$ *and use* $r$ *to define* $(\hat{A}, \hat{a})$ *via (3.2) with* $\mu$ *from* $r$*'s version of Lemma 3.1. The norm of any non-zero eigenvalue of* $(\hat{A}, \hat{a})$*'s version of* $\mathcal{D}_0$ *is greater than* $\frac{1}{\kappa}$ *. Also, this same version of* $\mathcal{D}_0$ *has a non-trivial kernel (the eigenspace for the eigenvalue* 0*) which has (real) dimension* 12g - 12.

The proof of this lemma appears momentarily. A description of the kernel of the lemma's version of $\mathcal{D}_0$ appears after the proof (see the upcoming Lemma 5.2). What follows directly is a digression with some observations that play important roles in the subsequent story.

With regards to the space of sections of $\mathbb{S}$ on $\Sigma$: Any given section $\mathbb{S}$ along $\Sigma$ can be written on a disk in $\Sigma$ with respect to an oriented, orthonormal frame for T*$\Sigma$ on the disk as

$$((\mathfrak{b} = (\mathfrak{b}_1, \mathfrak{b}_2, \mathfrak{b}_t), \mathfrak{c} = (\mathfrak{c}_1, \mathfrak{c}_2, \mathfrak{c}_t)) \,, (\mathfrak{b}_0, \mathfrak{c}_0))$$

(5.3)

where $\mathfrak{b}_t$ and $\mathfrak{c}_t$ are the respective dt components of $\mathfrak{b}$ and $\mathfrak{c}$. Use this depiction to introduce, by way of notation,

$$\mathit{b} = \mathfrak{b}_1 + i\mathfrak{b}_2, \quad c = \mathfrak{c}_1 - i\mathfrak{c}_2, \quad u = \mathfrak{b}_0 - i\mathfrak{b}_t, \quad x = \mathfrak{c}_0 + i\mathfrak{c}_t \ .$$

(5.4)

To say more about these objects: Use P to denote the principle bundle of SU(2) frames for $L \oplus L^{-1}$. The pair $u$ and $x$ are sections of $\mathrm{ad}(P) \otimes_{\mathbb{R}} \mathbb{C}$ on $\Sigma$; whereas $\mathit{b}$ and $\mathfrak{c}$ are the restrictions to the given disk of respective sections of the bundles $T^{0,1}\Sigma \otimes_{\mathbb{R}} \mathrm{ad}(P)$ and $T^{1,0} \otimes_{\mathbb{R}} \mathrm{ad}(P)$. This complex ($\mathbb{C}$) depiction of sections of $\mathbb{S}$ is introduce so as to write the operator $\mathfrak{D}_0$ using the complex structure on $\Sigma$. In this regard, the section depicted in (5.3) is an eigenvector of $\Sigma$ with eigenvalue $\lambda$ if and only if $((\mathit{b}, c), (u, x))$ obey the equations that are depicted below in (5.5). These equations introduce $\varphi$ to denote $\hat{\mathfrak{a}}_1 - i\hat{\mathfrak{a}}_2$ which comes from $\mathbb{X}$ as depicted in (3.1) with $\mu$ from Lemma 3.1 by writing $\mathbb{X}$ using an oriented, orthonormal frame for $T^*X$ as $\mathbb{X} = -i\varphi(e^1 + ie^2)$.) What is denoted by $\varphi^*$ signifies $\hat{\mathfrak{a}}_1 + i\hat{\mathfrak{a}}_2$ which is $-\varphi^\dagger$.

- $-(\nabla_{\hat{A}1} + i\nabla_{\hat{A}2})x + [\varphi^*, u] = \lambda\mathit{b}$ ,
- $-(\nabla_{\hat{A}1} - i\nabla_{\hat{A}2})u - [\varphi, x] = \lambda c$ ,
- $(\nabla_{\hat{A}1} + i\nabla_{\hat{A}2})c - [\varphi, \mathit{b}] = \lambda u$ ,
- $(\nabla_{\hat{A}1} - i\nabla_{\hat{A}2})\mathit{b} + [\varphi^*, c] = \lambda x$ .

(5.5)

With regards to $\varphi$: The fact that $(\hat{A}, \hat{\mathfrak{a}})$ satisfies the $m = 0$ version of (1.7) is equivalent to the following assertions about $\hat{A}$ and $\varphi$:

$$(\nabla_{\hat{A}1} + i\nabla_{\hat{A}2})\varphi = 0 \quad \mathit{and} \quad (F_{\hat{A}})_{12} = -\frac{i}{2}[\varphi, \varphi^*] \ .$$

(5.6)

The left most equation in (5.6) asserts that $\varphi$ is a holomorphic section of $T^{1,0}\Sigma \otimes_{\mathbb{R}} \mathrm{ad}(P)$. (Note that $\varphi$ obeys the left most equation in (5.6) if and only if $\varphi^*$ obeys $(\nabla_{\hat{A}1} - i\nabla_{\hat{A}2})\varphi^* = 0$.)

***Proof of Lemma 5.1***: Act on the top equation in (5.5) by $\nabla_{\hat{A}1} - i\nabla_{\hat{A}2}$ and act on the second equation by $\nabla_{\hat{A}1} + i\nabla_{\hat{A}2}$. Then use the other equations with (5.6) to see that if $v$ is either $x$ or $u$, then it obeys the second order equation

$$-(\nabla_{\hat{A}1}^2 + \nabla_{\hat{A}2}^2)v + [\hat{\mathfrak{a}}_1 [v, \hat{\mathfrak{a}}_1]] + [\hat{\mathfrak{a}}_2 [v, \hat{\mathfrak{a}}_2]] = \lambda^2 v \ .$$

(5.7)

This identity implies via $\Sigma$'s version of Lemma 4.2 that either $v \equiv 0$ or $\lambda > c_0^{-1}$. Supposing that $v \equiv 0$, then the third and fourth bullets of (5.5) say that $(\mathit{b}, c)$ are such that

- $(\nabla_{\hat{A}1} + i\nabla_{\hat{A}2})c - [\varphi, \mathit{b}] = 0$ ,

- $(\nabla_{\hat{A}1} - i\nabla_{\hat{A}2})\,\acute{b} + [\varphi^*, c] = 0$ .

$$(5.8)$$

The Atiyah-Singer index theorem can be brought to bear to see that the dimension (over $\mathbb{C}$) of the space of solutions to (5.8) is 6g - 6. Over $\mathbb{R}$, this is 12g - 12.

    To see what the solutions to (5.8) look like, it proves useful to write $\acute{b}$ and $c$ on the complement of $\Theta$ in terms of $\varphi$, $\varphi^*$ and $\sigma_3$. This is done below in (5.9) using $\mathbb{C}$-valued functions $\mu_+$, $\mu_-$, $\varsigma_+$, $\varsigma_-$ on $\Sigma - \Theta$ and respective sections $\mu_3$ and $\varsigma_3$ of $T^{0,1}\Sigma$ and $T^{1,0}\Sigma$ on $\Sigma - \Theta$.

- $\acute{b} = \mu_3\,\sigma_3 + \mu_+\,\varphi^* - \mu_-[\frac{i}{2}\,\sigma_3, \varphi^*]$ ,
- $c = \varsigma_3\,\sigma_3 + \varsigma_+\varphi + \varsigma_-\,[\frac{i}{2}\,\sigma_3, \varphi]$ ,

$$(5.9)$$

Note that $\acute{b}$ and $c$ can be written in this way because $\langle\sigma_3\varphi\rangle$ is identically zero so $\sigma_3$ is orthogonal to $\varphi$ (and thus also to $\varphi^*$); and $[\frac{i}{2}\,\sigma_3, \varphi]$ is orthogonal to $\sigma_3$ and it is orthogonal $\varphi$ except at the points in $\Theta$. The equations in (5.8) when written using the $\mu$'s and $\varsigma$'s devolve into two disjoint systems, one for the set $(\mu_+, \mu_-, \varsigma_3)$ and the other for $(\varsigma_+, \varsigma_-, \mu_3)$. These equations are depicted below.

- *The equations for* $(\mu_+, \mu_-, \varsigma_3)$
  - $(\nabla_1 - i\nabla_2)\mu_+ = 0$
  - $\frac{i}{2}(\nabla_1 - i\nabla_2)\mu_- + \varsigma_3 = 0$
  - $(\nabla_1 + i\nabla_2)\varsigma_3 + \frac{i}{2}\mu_-\langle\sigma_3\,[\varphi^*, [\sigma_3, \varphi]]\rangle + \mu_+\langle\sigma_3\,[\varphi^*, \varphi]\rangle = 0$
- *The equations for* $(\varsigma_+, \varsigma_-, \mu_3)$.
  - $(\nabla_1 + i\nabla_2)\,\varsigma_+ = 0$
  - $\frac{i}{2}(\nabla_1 + i\nabla_2)\varsigma_- + \mu_3 = 0$
  - $(\nabla_1 - i\nabla_2)\mu_3 + \frac{i}{2}\varsigma_-\langle\sigma_3\,[\varphi^*, [\,\sigma_3, \varphi]]\rangle + \varsigma_+\langle\sigma_3\,[\varphi^*, \varphi]\rangle = 0$

$$(5.10)$$

    Taking the second set of equations first (they are relevant first), the top equation of the second set says that $\varsigma_+$ is a holomorphic function on $\Sigma - \Theta$. As will be explained momentarily, an allowable $\varsigma_+$ can extend over $\Sigma$ as a meromorphic function with first order poles at the points in $\Theta$. To see about $\varsigma_-$, use the second equation of the second set of equations to write $\mu_3$ as $-\frac{i}{2}(\nabla_1 + i\nabla_2)\varsigma_-$. Then, substitute this into the third equation of the second set to see that $\varsigma_-$ is a solution to the second order equation below:

$$- (\nabla_1{}^2 + \nabla_2{}^2)\,\varsigma_- + |[\sigma_3, \varphi]|^2\,\varsigma_- - 2i\langle\sigma_3, [\varphi^*, \varphi]\rangle\,\varsigma_+ = 0 .$$

$$(5.11)$$

Since $\langle\sigma_3\varphi\rangle$ is zero, this equation can also be written with out reference to any $T^*\Sigma$ frame as

$$d^{\dagger}d\varsigma_- + 4|\hat{a}|^2\varsigma_- - 4*_\Sigma\langle\sigma_3\,\hat{a}\wedge\hat{a}\rangle\varsigma_+ = 0\;.$$

(5.12)

To consider solutions to the second set in (5.10) and (5.12), the first point to note is that if $\varsigma_+$ extends smoothly over $\Sigma$, then it must be constant. However, $\varsigma_+$ can have a poles at points in $\Theta$ because any first order pole can be canceled by a corresponding first order pole in $\varsigma_-$. To elaborate, use the identity $\varphi(e^1 + ie^2) = i\mathbb{X}$ with (3.1)'s depiction of $\mathbb{X}$ to write

$$(\varsigma_+\varphi + \varsigma_-[\tfrac{i}{2}\sigma_3, \varphi])\;(e^1 + ie^2) = i\begin{pmatrix} 0 & e^{\mu}\varpi_+(\varsigma_+ - \varsigma_-) \\ e^{-\mu}\varpi_-(\varsigma_+ + \varsigma_-) & 0 \end{pmatrix}.$$

(5.13)

Thus, if $p \in \Theta$ is a zero of $\varpi_-$, and if z is the holomorphic coordinate for a disk centered at z, then what is depicted in (5.13) will be smooth across the z = 0 locus if, near z = 0,

$$\varsigma_+ \sim \frac{a}{z} + \mathcal{O}(1)\quad with\; a \neq 0\; as\; long\; as\; \varsigma_- \sim \frac{a}{z} + \mathcal{O}(1)\; also.$$

(5.14)

By the same toke, if $p \in \Theta$ is a zero of $\varpi_+$, then what is depicted in (5.13) will be smooth across the z = 0 locus if

$$\varsigma_+ \sim \frac{a}{z} + \mathcal{O}(1)\quad with\; a \neq 0\; as\; long\; as\; \varsigma_- \sim -\frac{a}{z} + \mathcal{O}(1).$$

(5.15)

Also: If $\varsigma_+$ has a pole of order greater than 1 at any point in $\Theta$, then (5.13) can not define a smooth section section of $T^{1,0}\Sigma \otimes_{\mathbb{R}} ad(P)$ on the whole of $\Sigma$.

The same sort of analysis applies to solutions to the top set of equations in (5.10). Indeed, the top set of equations is obtained from the lower set by complex conjugation. In particular, $\mu_+$ has to be the complex conjugate of a meromorphic function on $\Sigma$ with poles at the points in $\Theta$, and $\mu_-$ must obey the equation in (5.12) with $\mu_-/\mu_+ \to 1$ approaching a zero of $\varpi_-$ where $\mu_+$ has a pole, and $\mu_-/\mu_+ \to$ -1 approaching a zero of $\varpi_+$ where $\mu_+$ has a pole.

The following lemma says in effect that (5.12) has a unique solution for any given meromorphic function $\varsigma_+$ with poles at the points in $\Theta$ subject to the constraints in (5.14) and (5.15). This same lemma says what is needed about $\mu_+$ and $\mu_-$ also because if the latter obey the first set of bullets in (5.10), then $(\varsigma_+, \varsigma_-) \equiv (\bar{\mu}_+, \bar{\mu}_-)$ obey the second set. By way of notation, the lemma reintroduces the function $\mathfrak{f}$ from Lemma 3.3 which is a radially symmetric function on $\mathbb{C}$. Also, given a point from $\Theta$, the lemma reintroduces that point's version of the number $\alpha$ which appears in Lemma 3.2 when using the local holomorphic coordinate as $q = \alpha\,z\,(dz)^2$.

**Lemma 5.2**: *There exists* $\kappa > 1$ *such that if* $r > \kappa$ *and* r *is used to define* $(\hat{A}, \hat{a})$ *using (3.2) with* $\mu$ *from* r*'s version of Lemma 3.1, then the assertions in the bullets that follow* *momentarily hold. Given a point* $\mathrm{p} \in \Theta$, *the upcoming bullets use* $\Xi$ *to denote the function on* $\mathbb{C}$ *given by the rule*

$$\Xi \equiv \mp \frac{2}{3}\left((s\tfrac{d}{ds}f)|_{\alpha^{1/3}z} + \frac{1}{2}\right)$$

*with the* $-$ *sign used for points in* $\Theta_+$ *and the* $+$ *sign used for points in* $\Theta_-$.

- *If* $\varsigma_+ \equiv 1$, *then the corresponding version of* $\varsigma_-$ *and* $\mu_3$ *when written with the holomorphic coordinate on the* $|z| < r_0$ *disk neighborhood of a given point in* $\Theta$ *differ there by at most* $\kappa \, e^{-r/\kappa}$ *from* $\Xi$ *and* $-\frac{i}{2}(\nabla_1 + i\nabla_2)\Xi$ *respectively.*

- *When* $\varsigma_+$ *has the Laurent series*

$$\varsigma_+ = \frac{a_{-1}}{z} + \sum_{k\geq 0} a_k z^k$$

*when written using the holomorphic coordinate on the* $|z| \leq r_0$ *disk neighborhood of a given point in* $\Theta$, *then the corresponding versions of* $\varsigma_-$ *and* $\mu_3$ *respectively differ there by at most* $\kappa \, e^{-r/\kappa}$ *from*

$$\left(\frac{3a_{-1}}{z} + \sum_{k\geq 0} \frac{3}{3+2k}a_k z^k\right)\Xi \quad and \quad -\frac{i}{2}\left(\frac{3a_{-1}}{z} + \sum_{k\geq 0}\frac{3}{3+2k}a_k z^k\right)(\nabla_1 + i\nabla_2)\,\Xi$$

- *In any case,* $|\varsigma_-| + |\mu_3|$ *are bounded by* $\kappa e^{-r/\kappa}\,(\sup_{\mathfrak{d}>r_0/2}|\varsigma_+|)$ *where* $\mathfrak{d} \geq r_0$.

- *With regards to* $\Xi$:
  a) $|\Xi| + r^{-2/3}|d\Xi| \leq \kappa \, e^{-r|z|^{3/2}/\kappa}$.
  b) *The only critical point of* $\Xi$ *is at* $z = 0$ *where it takes its minimum or maximum value (as the case may be). Except at that critical point, the function* $\Xi \frac{d}{d|z|}\Xi$ *is negative. Thus, the norm of* $\Xi$ *decreases monotonically to zero with increasing* $|z|$ *from its maximum at* $z = 0$ *(its maximum is* $\frac{1}{3}$ *).*

The assertions in this lemma are likely well known to the experts, but even so, a proof of this lemma is given in the Appendix to this paper.

The preceding analysis implies that the dimension of the space of solutions to the second set of equations in (5.10) with allowed behavior near the points in $\Theta$ is the dimension of the space of meromorphic functions on $\Sigma$ with at most first order poles at the points in $\Theta$. Since there are g constraints on the positions of the poles of a meromorphic function, and 4g - 4 points in $\Theta$, and one constant meromorphic function, the dimension (over $\mathbb{C}$) of the space of allowed solutions to the second set of equations in (5.10) is equal to 4g - 4 - g + 1 which is 3g - 3. This is also the complex dimension of the space of allowed solutions to the first set of equations in (5.10).

By way of a remark: The vector space of solutions to the second set of equations in (5.10) with allowed poles on $\Theta$ is the same as the vector space of first order deformations of the SL(2; $\mathbb{R}$) flat connection $\hat{A} + i\hat{a}$ in the moduli space of flat SL(2; $\mathbb{R}$) connections on $\Sigma$. The allowed solutions to the first set of equations in (5.10) are the first order deformations of $\hat{A} + i\hat{a}$ to a flat SL(2; $\mathbb{C}$) connection in the moduli space of flat SL(2; $\mathbb{C}$) connections, thus a flat connection whose holonomy is not inside the SL(2; $\mathbb{R}$) subgroup of SL(2; $\mathbb{C}$).

## b) Proposition 3.4

The purpose of this subsection is to elaborate on the pair $(\mathfrak{b}_t, \mathfrak{c}_t)$ in (3.25) and then prove Proposition 3.4. By way of preliminaries: The terms $\mathfrak{k}_0$ and $\mathfrak{k}_1$ and $\mathfrak{s}_0$ and $\mathfrak{s}_1$ in Proposition 3.4 when written explicitly on the $\mathfrak{d} < \frac{1}{100} r_0$ part of $S^1 \times \Sigma$ are depicted below in (5.16). With regards to notation: The pair $(\hat{A}, \hat{a})$ in (5.16) is the $r = \frac{n}{m}$ version of the pair that is used in the previous subsection which is (3.2)'s pair with $\mu$ being $\mu_n$ (which is the $r = \frac{n}{m}$ function $\mu_r$ from Lemma 3.1). The exterior covariant derivatives below in (5.16) (and henceforth) when acting on $\mathfrak{b}_t$ and $\mathfrak{c}_t$ are along $\Sigma$. (This interpretation of the $d_{\hat{A}}$ is implicit any subsequent equation from this section.)

- $\mathfrak{k}_0 = - *_\Sigma \big( d_{\hat{A}} \mathfrak{b}_t - [\hat{a}, \mathfrak{c}_t] + (d_{\hat{A}} \hat{A}_0 - *_\Sigma \hat{a}) \big)$ .
- $\mathfrak{k}_1 = - *_\Sigma \big( d_{\hat{A}} \mathfrak{c}_t + [\hat{a}, \mathfrak{b}_t] + [\hat{a}, \hat{A}_0] \big)$ .
- $\mathfrak{s}_0 = \mathfrak{c}_t$ .
- $\mathfrak{s}_1 = [\hat{A}_0 + \mathfrak{b}_t, \mathfrak{c}_t]$ .

$$(5.16)$$

A 'best' version of $(\mathfrak{b}_t, \mathfrak{c}_t)$ would be a pair for which $\mathfrak{k}_0$ and $\mathfrak{k}_1$ are zero. To see if such a pair exists, it is instructive to write the expressions

$$d_{\hat{A}} \mathfrak{b}_t - [\hat{a}, \mathfrak{c}_t] + (d_{\hat{A}} \hat{A}_0 - *_\Sigma \hat{a}) \quad and \quad d_{\hat{A}} \mathfrak{c}_t + [\hat{a}, \mathfrak{b}_t] + [\hat{a}, \hat{A}_0]$$

$$(5.17)$$

using the complex structure on $\Sigma$ in the manner of the previous subsection. To this end, let $(u, x)$ denote a pair of $\mathbb{C}$-valued sections of ad(P) on $\Sigma$. If the expressions below vanish

- $- (\nabla_{\hat{A}1} + i\nabla_{\hat{A}2})x + [\varphi^*, u] - i[\varphi^*, \hat{A}_0]$ ,
- $- (\nabla_{\hat{A}1} - i\nabla_{\hat{A}2})u - [\varphi, x] + \big(i(\nabla_{\hat{A}1} - i\nabla_{\hat{A}2})A_0 - \varphi\big)$ ,

$$(5.18)$$

and if both $x$ and $u$ are Hermitian endomorphisms (thus $i\mathbb{R}$-valued sections of ad(P)), then (5.17) vanishes when $\mathfrak{b}_t$ is set equal to $iu$ and $\mathfrak{c}_t$ is set equal to $-ix$.

With regards to (5.18): If those expressions vanish, then acting on the top equation in (5.18) with $(\nabla_{\hat{A}1} - i\nabla_{\hat{A}2})$ and using the second equation with the identities in (5.6) lead to the following second order differential equation below for $x$:

$$-(\nabla_{\hat{A}1}{}^2 + \nabla_{\hat{A}2}{}^2)x + [\hat{a}_1, [x, \hat{a}_1]] + [\hat{a}_2, [x, \hat{a}_2]] = [\varphi^*, \varphi] .$$

$$(5.19)$$

Since $[\varphi^*, \varphi]$ is proportional to $\sigma_3$, there is a solution to (5.19) that is also proportional to $\sigma_3$. Moreover, since $[\varphi^*, \varphi]$ is equal to $-2i\,[\mathfrak{a}_1, \mathfrak{a}_2]$, there is a solution having the form $i c\,\sigma_3$ with $c$ denoting the $\mathbb{R}$-valued solution to the equation below:

$$d^{\dagger}dc + 4|\hat{a}|^2 c = -2*_{\Sigma}\langle\sigma_3\,(\hat{a}\wedge\hat{a})\rangle\,.$$

(5.20)

In particular, if (5.19) is viewed as an equation on the whole of $\Sigma$ (not just where $\mathfrak{d} < \frac{1}{100}\,r_0$), then it has a unique solution which has the form $i c\,\sigma_3$ with $c$ obeying (5.20) on the whole of $\Sigma$. Indeed, this global solution $c$ is the function that is denoted by $\varsigma_-$ in the $\varsigma_+ = -\frac{1}{2}$ version of both (5.12) and Lemma 5.2.

With regards to $u$ in (5.18): Fix a point $p \in \Theta$ for the moment and introduce as in Part 3 of Section 3a and Lemma 3.2 the complex, holomorphic coordinate z for the radius $2r_0$ disk centered at p. If the expressions in (5.18) are zero on the $|z| < \frac{1}{100}\,r_0$ disk, then $v \equiv (u - i\hat{A}_0)$ obeys the equation

$$-(\nabla_{\hat{A}1} + i\nabla_{\hat{A}2})(\nabla_{\hat{A}1} - i\nabla_{\hat{A}2})v - [\varphi, [\varphi^*, v]] = 0$$

(5.21)

on that disk subject to the constraint that $v$ approximate $-i\hat{A}_0$ where $|z|$ is nearly $\frac{1}{100}\,r_0$. The first task is to find a solution with this property. To do that, start with the formula for $\hat{A}_0$ in (3.23) to write $-i\hat{A}_0$ as done below:

$$-\frac{1}{3}\,r\,\sqrt{\alpha_1}\,\left(w^3\hat{\tau}\,-\bar{w}^3\,\hat{\tau}\right)\,.$$

(5.22)

(What is denoted by $w$ should be viewed as a $\mathbb{C}$-valued section of the bundle $\mathcal{I}$ which is defined on the complement of the z = 0 point in the $|z| < \frac{1}{100}\,r_0$ disk by the requirement that $w^2 = z$.) Keeping (5.22) in mind, note next that the endomorphism $\hat{\tau}$ as depicted in (3.26) is the $\mu = \mu_\diamond$ version of both of the endomorphisms depicted below (see (3.1) for $\mathbb{X}$ and (3.5) for $\mu_\diamond$):

$$i\mathbb{X}/\varpi_q \quad and \quad i(\mathbb{X}/\varpi_q)^{\dagger}\,.$$

(5.23)

Now let $\mathbb{X}_n$ denote the $\mu = \mu_n$ and $q = r^2 q_1$ version of $\mathbb{X}$ (the function $\mu_n$ is Lemma 3.1's function $\mu_r$ for the case when $r = \frac{n}{m}$.) Since $\mu_n$ differs by at most $c_0\,e^{-r/c_0}$ from $\mu_\diamond$ where $|z|$ is nearly $\frac{1}{100}\,r_0$, what is denoted by $v$ below will be likewise $c_0\,e^{-r/c_0}$ close to $-i\hat{A}_0$ where $|z|$ is nearly $\frac{1}{100}\,r_0$:

$$\nu \equiv -\frac{i}{3} r \sqrt{\alpha_1} \; w^3 \, \mathbb{X}_n / \varpi_q \; + \; \frac{i}{3} r \sqrt{\alpha_1} (w^3 \mathbb{X}_n / \varpi_q)^\dagger$$

(5.24)

What is depicted above in (3.28) extends smoothly over the $z = 0$ locus, a fact which can be seen explicitly by first writing the left most term in (5.24) using (3.11) and (3.13) as

$$-\frac{i}{3} r \sqrt{\alpha_1} \begin{pmatrix} 0 & e^{\mu_n} w^4 \beta (1 + \vartheta_+) \\ e^{-\mu_n} w^2 \beta^{-1} (1 + \vartheta_+)^{-1} & 0 \end{pmatrix},$$

(5.25)

and then, keeping in mind that $\vartheta_+$ is a holomorphic function of $w^2$ since it is pulled up from the disk D in $\Sigma$ by $\pi_q$. Meanwhile, the right most term in (5.24) is -1 times the Hermitian conjugate of what is depicted in (5.25).

With regards to (5.21): Write $\nu$ as $a + a^\dagger$ with $a$ depicted in (5.25). Because both $\nabla_{\hat{A}1} + i\nabla_{\hat{A}2}$ and $[\varphi, \cdot]$ annihilate $a$ (the latter does because both $a$ and $\varphi$ are proportional to $\mathbb{X}_n$), it follows that $\nabla_{\hat{A}1} - i\nabla_{\hat{A}2}$ and $[\varphi^*, \cdot]$ annihilate $a^\dagger$. Thus, the $a^\dagger$ contribution to (5.21) is zero. Meanwhile, the contribution to (5.21) from $a$ is also zero because

$$-(\nabla_{\hat{A}1} + i\nabla_{\hat{A}2})(\nabla_{\hat{A}1} - i\nabla_{\hat{A}2}) - [\varphi, [\varphi^*, \cdot]] = -(\nabla_{\hat{A}1} - i\nabla_{\hat{A}2})(\nabla_{\hat{A}1} + i\nabla_{\hat{A}2})\nu - [\varphi^*, [\varphi, \cdot]],$$

(5.26)

this being an identity that follows from the left most equation in (5.6). (Both sides of (5.26) are the operator $-(\nabla_{\hat{A}1}{}^2 + \nabla_{\hat{A}2}{}^2) + [\hat{a}_1 [ \cdot , \hat{a}_1]] + [\hat{a}_2 [ \cdot , \hat{a}_2]]$.)

With regards to $x$ and $\nu$: Both are Hermitian endomorphisms; they are $c_0 \, e^{-r/c_0}$ close to 0 and $-i\hat{a}_0$ respectively where $|z|$ is nearly $\frac{1}{100} r_0$; and they obey (5.19) and (5.21). Granted all of this, the question now is whether what is written in (5.18) is zero (or nearly zero). To see about this, note that by virtue of (5.19) and (5.21) being zero, the expressions in (5.18) are such that

- $-(\nabla_{\hat{A}1} + i\nabla_{\hat{A}2})x + [\varphi^*, \nu] = \acute{b}$ ,
- $-(\nabla_{\hat{A}1} - i\nabla_{\hat{A}2})\nu - [\varphi, x] - \varphi = c$ ,

(5.27)

where $(\acute{b}, c)$ are solutions to (5.8) where $|z| < \frac{1}{100} r_0$. With regards to $(\acute{b}, c)$, although they need only solve (5.8) on this limited domain, they can none-the-less be written as in (5.9) where $z \neq 0$ and then the coefficients $(\mu_+, \mu_-, \varsigma_3)$ and/or $(\varsigma_+, \varsigma_-, \mu_3)$ must obey the equations that are depicted in (5.10).

With regards to $(\mu_+, \mu_-, \varsigma_3)$: This must be zero because the $c$ component of the corresponding version of $(\acute{b}, c)$ is $\varsigma_3 \sigma_3$ which is proportional to $\sigma_3$ whereas all terms on the left hand side of the lower bullet in (5.27) are pointwise orthogonal to $\sigma_3$. (Keep in mind here that $x$ is proportional to $\sigma_3$ and $\nu$ is orthogonal to $\sigma_3$.)

To study the $(\varsigma_+, \varsigma_-, \mu_3)$ versions of $(\hat{b}, c)$, note first that only the $a$ part of $\nu$ (the (5.25) part) contributes to the left hand side of the lower bullet in (5.27) because of the fact that $(\nabla_{\hat{A}1} - i\nabla_{\hat{A}2})a^\dagger$ is zero. Granted this point, it is useful to consider (5.27) where $|z|$ is between $\frac{1}{1000} r_0$ and $\frac{1}{200} r_0$. In particular, $\mu_n$ is $c_0 e^{-r/c_0}$ close to $\mu_\Diamond$ on this annulus. Likewise $a$ is $c_0 e^{-r/c_0}$ close to $-\frac{1}{3} r w^3 (1 + t_1)\hat{\tau}$. Meanwhile, $\varphi$ is similarly close to $2r w^2 (1 + t_0)\frac{dw}{dz}\hat{\tau}$ (see 3.17 and (3.11)) which is $c_0 e^{-r/c_0}$ close to $-(\nabla_{\hat{A}1} - i\nabla_{\hat{A}2})(\frac{1}{3} r w^3 (1 + t_1)\hat{\tau})$ since $\nabla_{\hat{A}}\hat{\tau}$ is similarly close to zero where $|z|$ is greater than $\frac{1}{1000} r_0$. Therefore, the left hand side of the lower bullet in (5.27) differs from $[\varphi, x]$ by at most $c_0 e^{-r/c_0}$ where $|z|$ is between $\frac{1}{1000} r_0$ and $\frac{1}{200} r_0$. This implies in turn that $|\varsigma_+|$ is bounded by $c_0 e^{-r/c_0}$ on this same annulus. By virtue of $\varsigma_+$ being meromorphic with at worst a first order pole at $z = 0$, this implies that

$$|\varsigma_+| \le c_0 \frac{1}{|z|} e^{-r/c_0}$$

(5.28)

where $|z| \le \frac{1}{200} r_0$.

To see about $\varsigma_-$: Remember that $x = i c \sigma_3$ with $c$ obeying (5.20) on the whole of $\Sigma$. As noted previously, that fact implies that $c$ is the $\varsigma_+ = -\frac{1}{2}$ version of Lemma 5.2's incarnation of $\varsigma_-$. In particular, that same version of Lemma 5.2 implies that $|x| \le c_0 e^{-r/c_0}$ where $|z|$ is between $\frac{1}{1000} r_0$ and $\frac{1}{200} r_0$. Hence, the norm of the current version of $\varsigma_-$ is likewise bounded by $c_0 e^{-r/c_0}$ where $|z|$ is between $\frac{1}{1000} r_0$ and $\frac{1}{200} r_0$. This fact and the fact that the current version obeys (5.11) where $|z| < \frac{1}{200} r_0$ and the bound in (5.28) can be used to see that $|\varsigma_-/\varsigma_+| \le 1 + c_0 e^{-r/c_0}$ where $|z| \le \frac{1}{200} r_0$. The argument is for the most part a replay of the proof of Lemma 5.2 in the Appendix. To elaborate, suppose that the point p is such that $\varpi_-$ at p is zero in which case (5.14) holds. Then (5.12) can be written as

$$d^\dagger d(\varsigma_- - \varsigma_+) + 4|\hat{a}|^2(\varsigma_- - \varsigma_+) = -4(|\hat{a}|^2 - *_\Sigma\langle\sigma_3 \hat{a} \wedge \hat{a}\rangle)\varsigma_+ .$$

(5.29)

Now in the case when $\varpi_-(p)$ is zero, then $|\hat{a}|^2 - *_\Sigma\langle\sigma_3 \hat{a} \wedge \hat{a}\rangle$ is bounded by $c_0 r^2 |z|^2$ where $|z|$ is much less than $c_0 r^{-2/3}$ and by $c_0 r^2 |z|$ otherwise. Thus, (5.29) with (5.28) imply that

$$d^\dagger d|\varsigma_- - \varsigma_+| \le e^{-r/c_0} ;$$

(5.30)

and the latter bound implies in turn that $|\varsigma_- - \varsigma_+| \leq c_0 e^{-r/c_0}$ on the whole of the disk where $|z| < \frac{1}{200} r_0$. (Use the comparison principle with the function $c_0(1 - |z|^2) e^{-r/c_0}$ and the fact that $|\varsigma_-|$ and $|\varsigma_+|$ are at most $c_0 e^{-r/c_0}$ where $|z| > \frac{1}{1000} r_0$.)

To see about $\mu_3$: The Dirichlet Green's function for the Laplacian on the $|z| < \frac{1}{400} r_0$ disk in $\mathbb{C}$ can be used with (5.29) to obtain a depiction of $\nabla(\varsigma_- - \varsigma_+)$ at any given point in the disk where $|z|$ is at most $\frac{1}{400} r_0$ as integral over the $|z| < \frac{1}{200} r_0$ disk of a function that is bounded by $c_0 \frac{1}{|z - (\cdot)|} e^{-r/c_0}$.

Of course, except for some sign changes, the exact same argument can be used for the case when $\varpi_+|_p$ is zero at p.

***Proof of Proposition 3.4***: Use the standard bump function $\chi$ to construct a compactly supported, non-negative function on the $|z| < \frac{1}{700} r_0$ that is equal to 1 where $|z| < \frac{1}{800} r_0$ and equal to zero where $|z| > \frac{1}{700} r_0$. This function can and should be constructed so that the norm of its exterior derivative is at most $c_0 r_0^{-1}$. Denote this function by $\chi_{00}$. For the purposes of Proposition 3.4, set $\mathfrak{b}_t = \chi_{00}(-i\nu - \hat{A}_0)$ and set $\mathfrak{c}_t = i\chi_{00}\chi$. It then follows from what is said previously in this subsection and from what is said in the $\varsigma_+ = -\frac{1}{2}$ version of Lemma 5.2 (whose version of $\varsigma_-$ is the function $c$) that this pair has all of the required properties.

## c) Solving the $m \neq 0$ version of (1.7) along $\Sigma$

This subsection uses the analysis from the previous subsections to modify Proposition 3.4's version of $(A, \mathfrak{a})$ along each $t \in S^1$ version of $S^1 \times \Sigma$ so that the resulting modified pair (henceforth denoted by $(A, \mathfrak{a})$) is such that

- $*(F_A - \mathfrak{a} \wedge \mathfrak{a}) - m\mathfrak{a} = \mathfrak{z}_0$

- $*(d_A\mathfrak{a}) = \mathfrak{z}_1$

- $*d_A*\mathfrak{a} = \mathfrak{z}_2$

(5.31)

with $\mathfrak{z} \equiv (\mathfrak{z}_0, \mathfrak{z}_1, \mathfrak{z}_2)$ obeying $|\mathfrak{z}| \leq c_0 e^{-r/c_0}$ when $r > c_0$ and $m < c_0^{-1}$. To be sure, the Hodge star in (5.31) is the Hodge star from $S^1 \times \Sigma$ and the bounds on $\mathfrak{z}$ hold everywhere on $S^1 \times \Sigma$.

To set the stage for the upcoming modifications, let $(A^1, \mathfrak{a}^1)$ denote Proposition 3.4's pair of $SU(2)$ connection on $E_q$ and 1-form on $S^1 \times \Sigma$ with values in the bundle of skew Hermitian endomorphisms of $E_q$. Use the $\mathbb{C} \oplus \mathbb{C}$ product structure for $E_q$ where $\mathfrak{d} < r_0$ from Part 2 of Section 3b on this part of $S^1 \times \Sigma$ to depict $(A^1, \mathfrak{a}^1)$ there as in (3.25). Now let $(\hat{A}, \hat{\mathfrak{a}})$ denote the version of the pair depicted in (3.2) using $q = r^2 q_1$ and with $\mu$ being the $r = \frac{n}{m}$ version of the function $\mu_r$ from Lemma 3.1. With $(\hat{A}, \hat{\mathfrak{a}})$ defined in this way, then $(A^1, \mathfrak{a}^1)$ can (and will) be written where $\mathfrak{d} < r_0$ as

- $A^1 = \hat{A} + B^1\,\sigma_3\; + m\,A^1\,dt\,,$

- $\mathfrak{a}^1 = \hat{\mathfrak{a}} + D^1 + m\,C^1\,\sigma_3\,dt\,,$

(5.32)

where the notation is as follows: First, $A^1 = A^0 + \mathfrak{b}_t$ with $A^0$ from (3.23) and $\mathfrak{b}_t$ from Proposition 3.4. For future reference, note that $A^1$ is an SU(2)-Lie algebra valued function that is everywhere orthogonal to $\sigma_3$ (this is because both $\hat{A}_0$ and $\mathfrak{b}_t$ have that property.) Meanwhile, $C^1$ is obtained from Proposition 3.4's version of $\mathfrak{c}_t$ by writing the latter as $C^1\sigma_3$. (This $C^1$ is $-\chi_{00}C$ with $\chi_{00}$ as defined in the proof of Proposition 3.4 and with $C$ being the solution to (5.20) on $\Sigma$.) As for $B^1$ and $D^1$: What is denoted by $B^1$ is a section of $T^*\Sigma$ with support only where $\mathfrak{d}$ is between $\frac{1}{4}\,r_0$ and $r_0$; and $D^1$ denotes a Lie algebra valued section of $T^*\Sigma$ with support in this same region which is pointwise orthogonal to $\sigma_3$. Moreover, the norms of both $B^1$ and $D^1$ are bounded by $c_0\,e^{-r/c_0}$.

The modified pair $(A, \mathfrak{a})$ for (5.31) has a form that is similar to (5.32) where $\mathfrak{d} < r_0$ (as before, this is using the $\mathbb{C} \oplus \mathbb{C}$ product structure):

- $A = \hat{A} + (B^1 + B)\sigma_3\; + (mA^1 + A)\,dt\,,$

- $\mathfrak{a} = \hat{\mathfrak{a}} + (D^1 + D) + (mC^1 + C)\,\sigma_3\,dt\,,$

(5.33)

with $A$, $B$, $C$ and $D$ described in the bullets that follows.

- *What are denoted by* $B$ *and* $C$ *signify (respectively) a section of* $T^*\Sigma$ *and a function on* $\Sigma$, *both compactly support where* $\mathfrak{d} < r_0$.
- *What are denoted by* $A$ *and* $D$ *signify (respectively) an* SU(2)*-Lie algebra valued function and section of* $T^*\Sigma$. *Both* $[\hat{\tau}, A]$ *and* $[\hat{\tau}, D]$ *have compact support where* $\mathfrak{d} < r_0$.
- $|[\hat{\tau}, A]| + |B| + |C| + |[\hat{\tau}, D]| \leq c_0 m^2\, r^{-2/3}\, e^{-r\,\mathfrak{d}^{3/2}/c_0}$ *and* $|A| + |D| \leq c_0 m^2 r^{-2/3}$.

(5.34)

Where $\mathfrak{d} > r_0$, the pair $(A, \mathfrak{a})$ is such that

$$A = A^1 + E\,dt\,\sigma \quad and \quad \mathfrak{a} = \mathfrak{a}^1 + F\,\sigma$$

(5.35)

with $E$ and $F$ described momentarily, and with $\sigma$ in (5.35) being the same as its namesake in (4.11); it is a unit normed, $\mathcal{I}$-valued section over $S^1 \times (\Sigma - \Theta)$ of the bundle of skew Hermitian endomorphisms of $E_q$. In particular, $\sigma$ appears as $\hat{\tau}$ when written using the $\mathbb{C} \oplus \mathbb{C}$ product structure of $E_q$ on the $\mathfrak{d} < r_0$ part of $S^1 \times \Sigma$; and $\sigma$ commutes with $\mathfrak{a}^1$ and is $A^1$-covariantly constant on the rest of $S^1 \times \Sigma$. With regards to $E$ and $F$: They are, respectively, an $\mathcal{I}$-valued function on the $\mathfrak{d} \geq r_0$ part of $\Sigma$ and an $\mathcal{I}$-valued section of $T^*\Sigma$ on $\mathfrak{d} > r_0$ part of $\Sigma$ obeying

$$|E| + |F| \leq c_0 m^2\, r^{-2/3}.$$

(5.36)

The proposition that follows momentarily makes a formal statement regarding $(A, \mathfrak{a})$.

With regards to notation: The $\mathbb{H}$ norm in the upcoming proposition is the version in Lemma 4.1 with $(\hat{A}, \hat{\mathfrak{a}})$ being either the pair $(A^1, \mathfrak{a}^1)$ from (5.32) or the pair $(A, \mathfrak{a})$ from the proposition. (It makes no difference which of these pairs is used to define $\|\cdot\|_{\mathbb{H}}$ because of what is said by (5.34) and (5.36) about the norm of the difference between these pairs.)

**Proposition 5.3:** *There exists* $\kappa > 1$ *with the following significance: Suppose that* $m$ *is less than* $\frac{1}{\kappa}$ *and that* $r \equiv \frac{n}{m}$ *. Assuming this, there exists a data set* $\phi = (A, B, C, D, E, F)$ *as described above (with* $c_0 = \kappa$ *in the second bullet of (5.34) and in (5.36)) and skew-Hermitian endomorphisms of* $E_q$ *to be denoted by* $\mathfrak{p}$ *and* $\mathfrak{q}$ *with* $\mathbb{H}$ *norm and pointwise norm bounded by* $\kappa m^2 e^{-r/\kappa}$ *such that*

- $-\nabla_{Ak}\mathfrak{p} + [\mathfrak{q}, \mathfrak{a}_k] + B_k - \frac{1}{2}\varepsilon_{kij}\,[\mathfrak{a}_i, \mathfrak{a}_j] - m\,\mathfrak{a}_k = \mathfrak{w}_{0k}$ ,
- $-\nabla_{Ak}\mathfrak{q} - [\mathfrak{p}, \mathfrak{a}_k] + \varepsilon_{ijk}\nabla_{Ai}\mathfrak{a}_j = \mathfrak{w}_{1k}$ ,
- $\nabla_{Ak}\mathfrak{a}_k = \mathfrak{w}_2$ ,

*with* $\mathfrak{w} \equiv (\mathfrak{w}_0, \mathfrak{w}_1, \mathfrak{w}_2)$ *obeying* $|\mathfrak{w}| \leq \kappa\ m^2\ e^{-r/\kappa}$. *Thus, when* $\mathfrak{z} = (\mathfrak{z}_0, \mathfrak{z}_1, \mathfrak{z}_2)$ *is defined by the* $(A, \mathfrak{a})$ *version of the expression on the left hand side of (5.32), then the* $S^1 \times \Sigma$ *integral of* $|\mathfrak{z}|^2$ *is at most* $\kappa\ m^2\ e^{-r/\kappa}$ *.*

This proposition is a consequence of the upcoming Lemma 5.5 (whose proof occupies the next subsection). By way of a parenthetical remark about $\mathfrak{z}$ in (5.32): One can show (with more work) that $\mathfrak{z}$'s pointwise norm is bounded by $c_0 e^{-r/c_0}$.

To set the stage for Lemma 5.5 and then the proof of Proposition 5.3: The data set $\phi$ is constructed from the solution to a version of (4.3) along $\Sigma$. This version of (4.3) is an equation for a section along $\Sigma$ of the bundle (denoted here by $\mathbb{S}$ ) that is obtained by tensoring the vector bundle $((\oplus_2 T^*(S^1 \times \Sigma)) \oplus (\oplus_2 \mathbb{R}))|_\Sigma$ with the bundle of skew Hermitian endomorphisms of the $\mathbb{C}^2$ bundle $L \oplus L^{-1}$.

With regards to $(\hat{A}, \hat{\mathfrak{a}})$: The pair $(\hat{A}, \hat{\mathfrak{a}})$ that appears in (5.32) and (5.33) extends over $\Sigma$ as a pair of special unitary connection on $L \oplus L^{-1}$ and 1-form on $\Sigma$ with values in the bundle of skew Hermitian endomorphisms of $L \oplus L^{-1}$. As noted before, this pair is depicted in (3.2) using $q = r^2 q_1$ and with $\mu$ being the $r = \frac{n}{m}$ version of the function $\mu_r$ from Lemma 3.1.

The relevant version of (4.3) is depicted in short-hand below.

$$\mathfrak{D}_0 \psi + m(\gamma_t[\chi_{\ddagger} \mathbb{A}^1, \psi] + c^1 \rho_t[\sigma_3, \psi] + \psi_c) + \psi \# \psi + \mathcal{G} = 0 \ .$$

(5.37)

The notation used here is as follows: What is denoted by $\mathfrak{D}_0$ is the operator that is depicted in (5.1) with $(\hat{A}, \hat{a})$ as just described. What is denoted by $\chi_{\ddagger}$ signifies the bump function $\chi(\frac{8 \mathfrak{d}}{r_0} - 1)$ which is equal to 1 where $\mathfrak{d} < \frac{1}{8} r_0$ and equal to zero where $\mathfrak{d} > \frac{1}{4} r_0$. As for $\psi \# \psi$, this is the same as its namesake in (4.3), but more is said about it momentarily. As for $\psi_c$ term: When $\psi$ is written in the form $((\mathfrak{b}, \mathfrak{c}), (\mathfrak{b}_0, \mathfrak{c}_0))$, then $\psi_c$ is the section of $\mathbb{S}$ given by $((0, \mathfrak{c}), (0, \mathfrak{c}_0))$. Finally, $\mathcal{G}$ denotes a section of $\mathbb{S}$ that is defined using what is denoted by $\mathfrak{s}_0$ and $\mathfrak{s}_1$ in Proposition 3.4; it is the section with the $\mathfrak{b}$ and $\mathfrak{c}_0$ components being zero, its $\mathfrak{c}$ component being $\mathfrak{m}^2 \mathfrak{s}_0$ dt and its $\mathfrak{b}_0$ component being $-m^2 \mathfrak{s}_1$.

Regarding (5.37): The $\mathfrak{k}_0$ and $\mathfrak{k}_1$ terms from Proposition 3.4 are not involved, nor are the $\mathbb{B}^1$ and $\mathbb{D}^1$ parts of $(A^1, \mathfrak{a}^1)$. These are ignored for the purposes of this subsection because of their very small, $c_0 e^{-r/c_0}$ norms. Being that they are ignored by (5.37), they contribute to $(\mathfrak{z}_0, \mathfrak{z}_1, \mathfrak{z}_2)$ in (5.32).

With regards to $\psi \# \psi$: As noted after (4.3), $\psi \# \psi$ signifies the image of $\psi$ via a fiber-wise quadratic map from $\mathbb{S}$ to $\mathbb{S}$ whose components are linear combinations of the commutators of the components of $\psi$. Of particular note for what is to come is that $\psi \# \psi$ has the following property with respect to the endomorphism $\sigma_3$:

**Lemma 5.4**: *Write the $\mathfrak{b}$ and $\mathfrak{c}$ components of any given section of $\mathbb{S}$ as $\mathfrak{b} = \mathfrak{b}_\Sigma + \mathfrak{b}_t$ dt and $\mathfrak{c} = \mathfrak{c}_\Sigma + \mathfrak{c}_t$ dt with $\mathfrak{b}_\Sigma$ and $\mathfrak{c}_\Sigma$ annihilating the vector fields tangent to the $\mathbb{S}^1$ factor of $\mathbb{S}^1 \times \Sigma$. Suppose that the conditions below hold for a section $\psi$:*

- $\mathfrak{b}_\Sigma, \mathfrak{c}_t$ *and $\mathfrak{c}_0$ are pointwise proportional to $\sigma_3$.*
- $\mathfrak{c}_\Sigma, \mathfrak{b}_t$ *and $\mathfrak{b}_0$ are pointwise orthogonal to $\sigma_3$.*

*Then, the preceding conditions also hold for the components of $\psi \# \psi$. These conditions also hold for the $\mathfrak{D}_0 \psi$ and $\gamma_t[\chi_{\ddagger} \mathbb{A}^1, \psi]$ and $c^1 \rho_t[\sigma_3, \psi]$ and $\psi_c$ terms in (5.36).*

An important point for the future: The $\mathcal{G}$-term in (5.37) also obeys the conditions of the lemma; this is because it is the commutator of $\hat{\mathbb{A}}^1$ (which is orthogonal to $\sigma_3$) and $\mathfrak{c}_t$ (which is proportional to $\sigma_3$). The implications of this fact about $\mathcal{G}$ and what is stated by the lemma are exploited in the next subsection.

*Proof of Lemma 5.4*: The essential point to keep in mind is that the commutator of $\sigma_3$ with anything orthogonal to $\sigma_3$ is still orthogonal to $\sigma_3$. Meanwhile, the commutator of something orthogonal to $\sigma_3$ with something else orthogonal to $\sigma_3$ is proportional to $\sigma_3$. Keeping the preceding in mind, the two bullets below in (5.37) list the contributions of the commutators of the components of $\psi$ to the respective $\mathfrak{b}_\Sigma$ and $\mathfrak{c}_t$ components of $\psi\#\psi$ (the $\mathfrak{c}_0$ component of $\psi\#\psi$ is zero).

- $[\mathfrak{c}_t, \mathfrak{b}_\Sigma]$ *and* $[\mathfrak{b}_t, \mathfrak{c}_\Sigma]$ *and* $[\mathfrak{b}_0, \mathfrak{c}_\Sigma]$ *and* $[\mathfrak{c}_0, \mathfrak{b}_\Sigma]$.
- $[\mathfrak{b}_\Sigma, \mathfrak{b}_\Sigma]$ *and* $[\mathfrak{c}_\Sigma, \mathfrak{c}_\Sigma]$ *and* $[\mathfrak{b}_t, \mathfrak{b}_0]$ *and* $[\mathfrak{c}_t, \mathfrak{c}_0]$.

$$(5.38)$$

The next three bullets list the the contributions of the commutators of the components of $\psi$ to the respective $\mathfrak{c}_\Sigma$ and $\mathfrak{b}_t$ and $\mathfrak{b}_0$ components of $\psi\#\psi$:

- $[\mathfrak{c}_t, \mathfrak{c}_\Sigma]$ *and* $[\mathfrak{b}_t, \mathfrak{b}_\Sigma]$ *and* $[\mathfrak{b}_\Sigma, \mathfrak{b}_0]$ *and* $[\mathfrak{c}_\Sigma, \mathfrak{c}_0]$ .
- $[\mathfrak{b}_\Sigma, \mathfrak{c}_\Sigma]$ *and* $[\mathfrak{b}_t, \mathfrak{c}_t]$ *and* $[\mathfrak{b}_t, \mathfrak{c}_0]$ *and* $[\mathfrak{c}_t, \mathfrak{b}_0]$ .
- $[\mathfrak{b}_\Sigma, \mathfrak{c}_\Sigma]$ *and* $[\mathfrak{b}_t, \mathfrak{c}_t]$ .

$$(5.39)$$

As for $\mathfrak{D}_0\psi$: Keep in mind that $\sigma_3$ is $\hat{A}$-covariantly constant and that $\hat{a}$ is orthogonal to $\sigma_3$.) With this understood, note that the contributions to the respective $\mathfrak{b}_\Sigma$ and $\mathfrak{c}_t$ and $\mathfrak{c}_0$ components of $\mathfrak{D}_0\psi$ of the $\hat{A}$-covariant derivatives of $\psi$ are those of $\psi$'s version of $\mathfrak{b}_\Sigma \mathfrak{c}_t$ and $\mathfrak{c}_0$. Likewise, those contributions to the respective $\mathfrak{c}_\Sigma$, $\mathfrak{b}_t$ and $\mathfrak{b}_0$ components of $\mathfrak{D}_0\psi$ are those of $\psi$'s version of $\mathfrak{c}_\Sigma \mathfrak{b}_t$ and $\mathfrak{b}_0$.

With regards to $\gamma_t[\chi_{\ddagger}A^1, \psi]$ and $c^1\rho_t[\sigma_3, \psi]$ and $\psi_{\mathfrak{c}}$, the lemma's assertions about these terms are left to the reader to verify.

The promised Lemma 5.5 follows directly..

**Lemma 5.5**: *There exists* $\kappa > 1$ *with the following significance: If* $\mathfrak{m}$ *is less than* $\frac{1}{\kappa}$ *and if* $r \equiv \frac{\mathfrak{n}}{m}$ *is greater than* $\kappa$, *then there is a solution (denoted by* $\psi$) *to (5.37) whose components obey Lemma 5.4's conditions and whose* $\mathbb{H}$ *norm and sup-norm are less than* $\kappa m^2 r^{-2/3}$. *Moreover, this solution* $\psi$ *is such that* $|[\hat{\tau}, \psi]| \leq \kappa m^2 r^{-2/3} e^{-r\mathfrak{d}^{3/2}/\kappa}$ .

The next subsection has the proof of Lemma 5.5 except for the assertions regarding the pointwise norms. These are proved in Section 5e.

***Proof of Proposition 5.3***: With Lemma 5.4 in hand (assume it is true for now), define A, B, C, D, E, and F for Proposition 5.3 as done below in (5.40). By way of notation, (5.40) depicts the section $\psi$ from Lemma 5.5 as $\psi = ((\mathfrak{b}, \mathfrak{c}), (\mathfrak{b}_0, \mathfrak{c}_0))$ and it then writes $\mathfrak{b}$ as $\mathfrak{b} = \mathfrak{b}_\Sigma + \mathfrak{b}_t dt$ and $\mathfrak{c}$ as $\mathfrak{c}_\Sigma + \mathfrak{c}_t dt$ with $\mathfrak{b}_\Sigma$ and $\mathfrak{c}_\Sigma$ being ad(P)-valued sections of T*$\Sigma$. The upcoming (5.40) also reintroduces the bump function $\chi_{00}$ from the proof of Proposition 3.4 (it is equal to 1 where $|z| < \frac{1}{800} r_0$ and equal to 0 where $|z| > \frac{1}{700} r_0$). Here is (5.40):

- A = $\chi_{00}\, \mathfrak{b}_t + (1 - \chi_{00}) \langle \hat{\tau}\mathfrak{b}_t \rangle\, \sigma$ *and* B = $\chi_{00}\, \mathfrak{b}_\Sigma$

- $C = \chi_{00}\, \mathfrak{c}_t$ *and* $D = \chi_{00}\, \mathfrak{c}_\Sigma + (1 - \chi_{00})\, \langle \hat{\tau} \mathfrak{c}_\Sigma \rangle\, \sigma$
- $E = \langle \hat{\tau} \mathfrak{b}_t \rangle$ *and* $F = \langle \hat{\tau} \mathfrak{c}_\Sigma \rangle$

(5.40)

It is a straightforward exercise to check that these choices for $A$, $B$, $C$, $D$, $E$, and $F$ obey the bounds in (5.34) and (5.36) and the conditions set forth in Lemma 5.4.

With regards to the sections $\mathfrak{p}$ and $\mathfrak{q}$: Borrow the function $\chi_\ddagger$ from (5.37) and use it to define skew Hermitian endomorphisms of $E_q$ by the rules below:

$$\mathfrak{p} \equiv \chi_\ddagger \mathfrak{b}_0 + (1 - \chi_\ddagger)\langle \hat{\tau} \mathfrak{b}_0 \rangle\, \sigma \quad \text{and} \quad \mathfrak{q} \equiv \chi_\ddagger \mathfrak{c}_0 + (1 - \chi_\ddagger)\langle \hat{\tau} \mathfrak{c}_0 \rangle\, \sigma \ .$$

(5.41)

Note in particular that $|\mathfrak{p} - \mathfrak{b}_0|$ and $|\mathfrak{q} - \mathfrak{c}_0|$ are both less than $c_0 m^2\, e^{-r/c_0}$. It then follows from Lemma 5.4 that what is denoted by $\mathfrak{w}$ in Proposition 5.3 has the desired $c_0 m^2 e^{-r/c_0}$ pointwise bound. As explained in the next paragraphs, this bound leads to the $c_0 m^2\, e^{-r/c_0}$ bound for the $\mathbb{H}$ and pointwise norms of $\mathfrak{p}$ and $\mathfrak{q}$.

To start the explanation: Use the same manipulations that led to (4.2) starting with the three bullet identities in Proposition 5.3 to derive the equation below for $\mathfrak{p}$ on $S^1 \times \Sigma$:

$$-\nabla_{Ak}\nabla_{Ak}\mathfrak{p} + [\mathfrak{a}_k, [\mathfrak{p}, \mathfrak{a}_k]] + [\mathfrak{q}, \nabla_{Ak}\mathfrak{a}_k] - m\nabla_{Ak}\mathfrak{a}_k = \nabla_k \mathfrak{w}_{0k} - [\mathfrak{a}_k, \mathfrak{w}_{1k}] \ .$$

(5.42)

This last equation can be written (using the three bullet equations in Proposition 5.3) as

$$-\nabla_{Ak}\nabla_{Ak}\mathfrak{p} + [\mathfrak{a}_k, [\mathfrak{p}, \mathfrak{a}_k]] + [\mathfrak{q}, \mathfrak{w}_2] - m\, \mathfrak{w}_2 = \nabla_k \mathfrak{w}_{0k} - [\mathfrak{a}_k, \mathfrak{w}_{1k}] \ .$$

(5.43)

Taking the inner product of both sides of this equation with $\mathfrak{p}$ leads to this equation:

$$\tfrac{1}{2}\, d^\dagger d |\mathfrak{p}|^2 + |\nabla_A \mathfrak{p}|^2 + |[\mathfrak{a}, \mathfrak{p}]|^2 \le c_0 e^{-r/c_0} |\mathfrak{p}||\mathfrak{q}| + \langle \mathfrak{p}(\nabla_{Ak}\mathfrak{w}_{0k} - [\mathfrak{a}_k, \mathfrak{w}_{1k}]) \rangle .$$

(5.44)

Integrate both sides of this equation over $S^1 \times \Sigma$ to see (after some integration by parts on both sides and use of the bound $\|\mathfrak{q}\|_1 < c_0 m^2 r^{-2/3}$) to see that

$$\|\mathfrak{p}\|_{\mathbb{H}}^2 \le c_0 m^2\, e^{-r/c_0} \ .$$

(5.45)

With this bound in hand, return to (5.44). Supposing that $p$ is any given point in $\Sigma$, multiply both sides of (5.44) by the Green's function for the operator $d^\dagger d + 1$ with pole at $p$. After suitable integration by parts on both sides, the result of doing that leads directly to the desired $c_0 m^2 e^{-r/c_0}$ bound for pointwise norm of $\mathfrak{p}$. (This Green's function is bounded by $c_0(1 + |\ln(\text{dist}(p, \cdot)|)$ and its derivative is bounded by $c_0 \frac{1}{\text{dist}(p, \cdot)}$.)

As for $\mathfrak{q}$: Use three bullet equations in Proposition 5.3 to see that $\mathfrak{q}$ obeys

$$-\nabla_{Ak}\nabla_{Ak}q + [\mathfrak{a}_k, [q, \mathfrak{a}_k]] - [\mathfrak{p}, \mathfrak{w}_2] = \nabla_{\hat{A}k}\mathfrak{w}_{1k} - [\mathfrak{w}_{0k}, \mathfrak{a}_k] .$$

(5.46)

But for some minor cosmetic changes, the derivation of the $\mathbb{H}$ and pointwise bounds for $\mathfrak{p}$ can be used to derive a $c_0 m^2 e^{-r/c_0}$ upper bound for the $\mathbb{H}$ and pointwise norm for q.

### d) Proof of Lemma 5.5

The proof of Lemma 5.5 follows the perturbative strategy from Section 4b. But even so, there is preliminary work to be done. These preliminaries and the lemma's proof (except for the pointwise norm bounds) occupy the four parts of this subsection. The pointwise norm bounds are proved in the next subsection (Section 5e).

*Part 1*: To set the stage for the proof, introduce the vector subbundle $\mathbb{S}_{SL} \subset \mathbb{S}|_\Sigma$ whose elements, when written as $((\mathfrak{b}_\Sigma + \mathfrak{b}_t dt, \mathfrak{c}_\Sigma + \mathfrak{c}_t dt), (\mathfrak{b}_0, \mathfrak{c}_0))$, have the properties listed in the two bullets of Lemma 5.4. According to that lemma, if $\psi$ is a section of $\mathbb{S}_{SL}$, then each terms that appears in (5.37) is also sections of $\mathbb{S}_{SL}$. As a consequence, it is sufficient to look for a section of $\mathbb{S}_{SL}$ to solve to (5.37). This restriction is useful by virtue of the fact that the kernel of the operator $\mathfrak{D}_0$ is can be written as a direct sum $K \oplus K^\perp$ with $K \subset C^\infty(\Sigma; \mathbb{S}_{SL})$ and with $K^\perp$ consisiting of elements that are pointwise orthogonal to $\mathbb{S}_{SL}$. In this regard, the elements in $K$ are those in the second bullet of (5.10); whereas those in $K^\perp$ are in those from the top bullet.

The first step in proving Lemma 5.5 constructs an approximate solution to (5.37) this being a particular section of $\mathbb{S}_{SL}$ that obeys the equation

$$\mathfrak{D}_0\psi = -\mathcal{G} .$$

(5.47)

The following lemma says what is needed about the required solution.

**Lemma 5.6**: *There exists $\kappa > 1$ such that if $r > \kappa$, then there is a unique solution to (5.47) in $C^\infty(\Sigma; \mathbb{S}_{SL})$ to be denoted by $\psi_{\mathcal{G}}$ which is $\mathbb{L}$-orthogonal to the kernel of $\mathfrak{D}_0$. This solution has the form $\psi_{\mathcal{G}} = ((\mathfrak{b}_\Sigma, \mathfrak{c}_\Sigma), (0, 0))$ with $\mathfrak{b}_\Sigma$ and $\mathfrak{c}_\Sigma$ being sections of $T^*\Sigma$ with values in the bundle of anti-Hermitian endomorphisms of $L \oplus L^{-1}$. With regards to the size of $\psi_{\mathcal{G}}$: The $\mathbb{L}$ norm of $\psi_{\mathcal{G}}$ is bounded by $\kappa\, m^2\, r^{-4/3}$ and the $\mathbb{H}$ norm is bounded by $\kappa\, m^2\, r^{-2/3}$. As for the pointwise norm, this obeys $|\psi_{\mathcal{G}}| \le \kappa\, m^2 r^{-2/3}$; and also $|[\hat{\tau}, \psi_{\mathcal{G}}]| \le \kappa\, m^2\, r^{-2/3}\, e^{-r\mathfrak{b}^{3/2}/\kappa} .$*

*Proof of Lemma 5.6*: The only non-zero components of $\mathcal{G}$ are its $\mathfrak{b}_0$ and $\mathfrak{c}_t$ components. And, in any event,

$$|\mathcal{G}| < c_0 \; m^2 \; e^{-r \mathfrak{d}^{3/2}/c_0} \; .$$

(5.48)

Since $\mathcal{G}$ has no $\mathfrak{b}_\Sigma$ and $\mathfrak{c}_\Sigma$ components, it can be written as $\mathfrak{D}_0 \psi_\mathcal{G}$ with $\psi_\mathcal{G}$ having only $\mathfrak{b}_\Sigma$ and $\mathfrak{c}_\Sigma$ components (see (5.5)). In this regard, the corresponding $\acute{b}$ and $c$ as defined in (5.4) obey the equations

- $(\nabla_{\hat{A}1} + i\nabla_{\hat{A}2})c - [\varphi, \acute{b}] = u_\mathcal{G}$ ,
- $(\nabla_{\hat{A}1} - i\nabla_{\hat{A}2})\acute{b} + [\varphi^*, c] = x_\mathcal{G}$ ,

(5.49)

where $u_\mathcal{G}$ is $(-\mathcal{G})$'s version of $\mathfrak{b}_0$ and $x_\mathcal{G}$ is i times $(-\mathcal{G})$'s version of $\mathfrak{c}_t$. There is no obstruction to solving this system of equations because the operator defined by the left hand side of (5.49) has only a trivial cokernel. In particular, (5.49) has a unique solution along $\Sigma$ which is $\mathbb{L}$-orthogonal to the kernel of $\mathfrak{D}_0$. Assume this orthogonality for $(\acute{b}, c)$. They can then be written in terms of ad$(P)_\mathbb{C}$ valued functions $u$ and $x$ as done below:

- $-(\nabla_{\hat{A}1} + i\nabla_{\hat{A}2})x + [\varphi^*, u] = \acute{b}$ .
- $-(\nabla_{\hat{A}1} - i\nabla_{\hat{A}2})u - [\varphi, x] = c$ .

(5.50)

The preceding equations imply in turn that $(u, x)$ obey the second order equations

- $-(\nabla_{\hat{A}1}^2 + \nabla_{\hat{A}2}^2)u + [\hat{a}_1, [u, \hat{a}_1]] + [\hat{a}_2, [u, \hat{a}_2]] = u_\mathcal{G}$ .
- $-(\nabla_{\hat{A}1}^2 + \nabla_{\hat{A}2}^2)x + [\hat{a}_1, [x, \hat{a}_1]] + [\hat{a}_2, [x, \hat{a}_2]] = x_\mathcal{G}$ .

(5.51)

Let $v$ denote either $u$ or $x$. It follows from (5.51) and (5.48) that the sequence of inequalities depicted below hold:

- $\int_\Sigma \; (|\nabla_{\hat{A}}v|^2 + |[\hat{a}, v]|^2) \le c_0 m^2 \; (\int_\Sigma \; \frac{1}{\mathfrak{b}^2} |v|^2)^{1/2} \; (\int_\Sigma \; \mathfrak{d}^2 e^{-r \mathfrak{d}^{3/2}/c_0})^{1/2} \; .$
- $\int_\Sigma \; (|\nabla_{\hat{A}}v|^2 + |[\hat{a}, v]|^2) \le c_0 m^2 r^{-4/3} \; (\int_\Sigma \; \frac{1}{\mathfrak{b}^2} |v|^2)^{1/2} \; .$
- $\int_\Sigma \; (|\nabla_{\hat{A}}v|^2 + |[\hat{a}, v]|^2) \le c_0 m^4 r^{-8/3} \; .$

(5.52)

The last inequality implies in turn that

$$\int_\Sigma \; (|\acute{b}|^2 + |c|^2) \le c_0 m^4 r^{-8/3} \; ,$$

(5.53)

which implies (via Lemma 4.3) that

$$\int_\Sigma \ (|\nabla_{\hat{A}}\psi_{\mathcal{G}}|^2 + |[\hat{a}, \psi_{\mathcal{G}}]|^2) \le c_0 m^4 r^{-4/3} \ .$$

(5.54)

It remains only to obtain pointwise bounds for $\psi_{\mathcal{G}}$ and $[\hat{\tau}, \psi_{\mathcal{G}}]$. These are derived in the next subsection, Section 5e.

*Part 2*: If the desired solution $\psi$ to (5.37) is written as $\psi = \psi_{\mathcal{G}} + \phi$, then the section $\phi$ of $\mathbb{S}_{SL}$ must obey an equation that has the schematic form

$$\mathcal{L}\phi + \mathfrak{l}(\phi) + \phi\#\phi + \mathcal{G}' = 0$$

(5.55)

with the notation as follows: What is denoted by $\mathcal{L}$ signifies the operator

$$\phi \to \mathcal{L}\phi \equiv \mathfrak{D}_0\phi + m\big(\gamma_t[\chi_{\ddagger}A^1, \phi] + c^1\rho_t[\sigma_3, \phi] + \phi_c\big) \ .$$

(5.56)

It is viewed in this subsection as an operator acting on sections of $\mathbb{S}_{SL}$. Meanwhile, what is denoted by $\mathfrak{l}(\cdot)$ signifies the endomorphism of $\mathbb{S}_{SL}$ given by the rule

$$\phi \to \mathfrak{l}(\phi) = \phi\#\psi_{\mathcal{G}} + \psi_{\mathcal{G}}\#\phi \ ;$$

(5.57)

and what is denoted by $\mathcal{G}'$ signifies the section below of $\mathbb{S}_{LC}$:

$$\mathcal{G}' \equiv \psi_{\mathcal{G}}\#\psi_{\mathcal{G}} \ - m\big(\gamma_t[\chi_{\ddagger}A^1, \psi_{\mathcal{G}}] + c^1\rho_t[\sigma_3, \psi_{\mathcal{G}}] + (\psi_{\mathcal{G}})_c\big) \ .$$

(5.58)

With regards to the operator $\mathcal{L}$: The key observation is in the upcoming Lemma 5.7. But note first that $\mathcal{L}$ defines an bounded, Fredholm operator from the $\mathbb{S}_{SL}$ version of $\mathbb{H}$ to the corresponding version of $\mathbb{L}$; and it defines an unbounded, self-adjoint operator on this same version of $\mathbb{L}$ with discrete spectrum having no accumulation points and finite multiplicities. This is because it differs from $\mathfrak{D}_0$ by a smooth, Hermitian endomorphism of $\mathbb{S}_{SL}$ and $\mathfrak{D}_0$ has these spectral properties (Lemma 4.1 implies this about $\mathfrak{D}_0$, see also Lemma 5.1).

**Lemma 5.7**: *The exists $\kappa > 1$ such that if $r > \kappa$ and $m < \frac{1}{\kappa^3}$, then $\mathcal{L}$, when viewed as an unbounded, self-adjoint operator on the $\mathbb{L}$ completion of $C^\infty(\Sigma; \mathbb{S}_{SL})$, has no eigenvalues between $-\frac{1}{\kappa} m$ and $\frac{1}{\kappa} m$ (and also none with absolute value between $\kappa m$ and $\frac{1}{\kappa}$ ).*

This lemma is proved in Part 4 of the subsection. (By way of a look ahead to the proof: The lower bound on the eigenvalue norms is explicitly a consequence of the domain for $\mathcal{L}$ being the space of sections of the bundle $\mathbb{S}_{SL}$. The lemma's eigenvalue bound does not hold when $\mathcal{L}$ operates on all sections of $\mathbb{S}$.)

With regards to the endomorphism $\mathfrak{l}(\cdot)$ depicted in (5.57): It is a consequence of the pointwise bounds from Lemma 5.6 that

$$|\mathfrak{l}(\phi)| \leq c_0 \, m^2 r^{-2/3}(|[\hat{\tau}, \phi]| + e^{-r\mathfrak{d}^{3/2}/c_0}|\phi|)$$

(5.59)

This implies in particular that

$$\|\mathfrak{l}(\phi)\|_{\mathbb{L}}^2 \leq c_0 \, m^4 r^{-8/3} \, \|\phi\|_{\mathbb{H}}^2$$

(5.60)

because the $\mathbb{L}$ norm of $[\hat{\tau}, \phi]$ is at most $c_0 r^{-2/3}$ times the $\mathbb{H}$ norm of $\phi$; and because the $\mathbb{L}$ norm of $e^{-r\mathfrak{d}^{3/2}/c_0}|\phi|$ is at most $c_0 r^{-2/3}$ times that of $\frac{r^{2/3}}{(1+ r^{4/3}\mathfrak{d}^2)^{1/2}}\,\phi$ which is, in turn, at most $c_0 r^{-2/3}$ times the $\mathbb{H}$ norm of $\phi$ (see Lemma 4.2).

With regards to $\mathcal{G}'$: Lemmas 3.4 and 5.6 lead directly to the following bound for $\mathcal{G}'$:

$$\|\mathcal{G}'\|_0 \leq c_0 \, m^2 \, r^{-4/3} \ .$$

(5.61)

With regards to the $\phi\#\phi$ term in (5.55): There is a dimension 2 Sobolev inequality that can be invoked which says that if $\mathfrak{f}$ is any given function on $\Sigma$, then

$$\int_\Sigma \mathfrak{f}^4 \leq c_0 \, (\int_\Sigma \, |d\mathfrak{f}|^2)(\int_\Sigma \, \mathfrak{f}^2) + c_0(\int_\Sigma \, \mathfrak{f}^2)^2 \ .$$

(5.62)

(Use (7.26) in Theorem 7.10 from [GP] with u there equal to $\mathfrak{f}^2$.) Here, because $\phi\#\phi$ involves commutators of components of $\phi$, this implies the bound below for $\phi\#\phi$:

$$\int_\Sigma \, |\phi\#\phi|^2 \leq c_0 \, (\int_\Sigma \, |\nabla_{\hat{A}}\phi|^2)(\int_{\Sigma_r} \, |\phi|^2 + (\int_\Sigma \, |\phi|^2)^{1/2}(\int_\Sigma \, |[\hat{\tau}, \phi]|^2)^{1/2}) + c_0(\int_\Sigma \, |\phi|^2)(\int_\Sigma \, |[\hat{\tau}, \phi]|^2) \ .$$

(5.63)

This inequality (with its complicated right hand side) is useful when used with the fact that

$$\int_\Sigma \, |[\hat{\tau}, \phi]|^2 \leq c_0 r^{-4/3} \, \|\phi\|_{\mathbb{H}}^2 \quad and \quad \int_{\Sigma_r} \, |\phi|^2 \leq c_0 r^{-4/3} \, \|\phi\|_{\mathbb{H}}^2.$$

(5.64)

The preceding bounds with what is said by Lemma 5.7 are input for Lemma 5.5's proof which is in the next part of the subsection

*Part 3*: This part of the subsection has the proof of Lemma 5.5 except for the assertions about pointwise norms. Those are proved in Section 5e. With regards to the proof: The assumption in what follows is that Lemma 5.7 is true.

**Proof of Lemma 5.5**: Lemma 5.7's assumptions imply (in part) that $\mathfrak{L}$ defines a bounded, invertible map from the $\mathbb{L}$ completion of $C^\infty(\Sigma; \mathbb{S}_{SL})$ to the $\mathbb{H}$ completion of $C^\infty(\Sigma; \mathbb{S}_{SL})$. Because of this, a map (denoted here by $\mathbb{T}$) from the $\mathbb{H}$ completion of the space $C^\infty(\Sigma; \mathbb{S}_{SL})$ to itself can be defined using the rule

$$\phi \to \mathbb{T}(\phi) = - \mathfrak{L}^{-1}(\phi \# \phi + \mathfrak{l}(\phi) + \mathcal{G}') \ .$$

(5.65)

A fixed point of this map is a solution to (5.55). To find a fixed point with small $\mathbb{H}$ norm, note first that the norm lower bound in Lemma 5.7 for $\mathfrak{L}$'s eigenvalues with (5.59), (5.60) and (5.63), (5.64) lead to the inequality below:

$$\|\mathbb{T}(\phi)\|_{\mathbb{L}}^2 \leq c_0 \frac{1}{m} r^{-2/3} \|\phi\|_{\mathbb{H}^2} \|\phi\|_{\mathbb{H}^2} + c \ m^3 r^{-8/3} \|\phi\|_{\mathbb{H}^2} + c_0 \ m^3 \ r^{-8/3} \ .$$

(5.66)

In turn, this last inequality with the inequality in the third bullet of Lemma 4.3 lead to the bound below for the $\mathbb{H}$ norm of $\mathbb{T}(\phi)$.

$$\|\mathbb{T}(\phi)\|_{\mathbb{H}^2} \leq c_0 \frac{1}{m} r^{2/3} \|\phi\|_{\mathbb{H}^2} \|\phi\|_{\mathbb{H}^2} + c \ m^3 r^{-4/3} \|\phi\|_{\mathbb{H}^2} + c_0 \ m^3 \ r^{-8/3} \ .$$

(5.67)

The preceding inequality implies that $\mathbb{T}$ maps the radius $c_0^{-1} m \ r^{-2/3}$ ball about the origin in $\mathbb{H}$ to itself. Much the same analysis proves that $\mathbb{T}$ maps a larger $c_0$ version of that ball (thus, smaller radius) to itself as a contraction mapping. This implies in particular that $\mathbb{T}$ has a unique fixed point in that (smaller) radius $c_0^{-1} m \ r^{-2/3}$ ball. Moreover the bound in (5.60) implies that the $\mathbb{H}$ norm of this fixed point is at most $c_0 m^{3/2} \ r^{-4/3}$.

The pointwise norm for the fixed point is proved in the final subsection, Section 5e.

*Part 4*: This part of the subsection is dedicated to the proof of Lemma 5.7. The proof first explains why $\mathfrak{L}$ has no eigenvalues with absolute value between $c_0 m$ and $c_0^{-1}$. (This fact holds when $\mathfrak{L}$ acts on the space of sections of $\mathbb{S}$.) The proof then explains why $\mathfrak{L}$ has no eigenvalues with absolute value between 0 and $c_0^{-1} m$ when it is restricted to the space of sections of $\mathbb{S}_{SL}$.

**Proof of Lemma 5.7 (I)**: Supposing that $s \in [0, 1]$, let $\mathfrak{L}_s$ denote the version of $\mathfrak{L}$ in (5.56) that has $m$ replaced by $sm$. Having fixed a value for s in [0, 1], suppose that $\lambda$ is an eigenvalue for $\mathfrak{L}_s$ with $|\lambda| < 1$. Let $\psi$ denote a corresponding eigenvector with $\mathbb{L}$ norm equal to 1. Lemma 4.6 can be brought to bear (when $m < c_0^{-1}$ and $r > c_0$) to see that

$$\| \ [\chi_\ddagger A^1, \psi] \|_{\mathbb{L}} \leq c_0 \ \|\psi\|_{\mathbb{L}} \ .$$

(5.68)

Note that this is true even if $\psi$ is a section of $\mathbb{S}$. It then follows from the definition of $\mathfrak{L}_s$ that

$$\|\mathfrak{D}_0\psi - \lambda\psi\|_{\mathbb{L}}^2 \leq c_0\, s^2 m^2\, \|\psi\|_{\mathbb{L}}^2\,.$$

$$(5.69)$$

Now write $\psi$ as $\phi + \psi_\perp$ with $\phi$ denoting the part of $\phi$ in the kernel of $\mathfrak{D}_0$ and with $\psi_\perp$ being $\mathbb{L}$-orthogonal to the kernel of $\mathfrak{D}_0$. Let $\lambda_1$ denote the smallest of the norms of the non-zero eigenvalues of $\mathfrak{D}_0$. What is written in (5.69) implies that

$$(\lambda^2 - c_0\, s^2 m^2)\|\phi\|_0^2 + \big((\lambda_1 - \lambda)^2 - c_0\, s^2 m^2\big)\|\psi_\perp\|^2 \leq 0\,.$$

$$(5.70)$$

If $sm < c_0^{-1}$, then this inequality can not hold when $|\lambda|$ is between $c_0 sm$ and $c_0^{-1}\lambda_1$ because otherwise, both terms on the left hand side of the inequality would be positive.

Here is the first implication of the preceding observation: If $m < c_0^{-1}$ and $r > c_0$, then there are no eigenvalues for $\mathfrak{L}$ (acting on sections of $\mathbb{S}$) with norm between $c_0 m$ and $c_0^{-1}\lambda_1$. There is also a second implication: Because the eigenvalues are continuous functions of $s$ (see Chapter V, Section 4.3 of [K]), and because no eigenvalue norm crosses the interval between $c_0 m$ and $c_0^{-1}\lambda_1$ as $s$ increases from 0 to 1, the dimension of the span of the eigenvectors of $\mathfrak{L}$ (acting on sections of $\mathbb{S}$) with eigenvalue norm less than $c_0 m$ is equal to that of the kernel of $\mathfrak{D}_0$ which is 12g -12. By the same token, the dimension of the span of the $\mathbb{S}_{SL}$ eigenvectors of $\mathfrak{L}$ with eigenvalue norm less than $c_0 m$ is equal to 6g -6. (The kernel of $\mathfrak{D}_0$ from the space of sections of $\mathbb{S}_{SL}$ are of the sort described by the second bullet in (5.10).)

*Proof of Lemma 5.7* (II): There is one more implication of (5.70) which is needed for what is to come which is this: Suppose again that $m < c_0^{-1}$ and that $r > c_0$. Let $\lambda$ denote an eigenvalue of $\mathfrak{L}$ with norm less than $c_0 m$ and let $\psi$ denote the corresponding eigenvector. Since $\lambda_1$ is no smaller than $c_0^{-1}$, the inequality in (5.70) requires that

$$\|\psi_\perp\|_{\mathbb{L}} \leq c_0\, sm\, \|\phi\|_{\mathbb{L}}\,.$$

$$(5.71)$$

Although the preceding holds when $\psi$ is a section of $\mathbb{S}$, assume henceforth that $\psi$ is a section of $\mathbb{S}_{SL}$. Keeping in mind that $\mathfrak{D}_0\phi = 0$, and keeping in mind (5.71), take the inner $\mathbb{L}$-inner product of both sides of the identity $\mathfrak{L}\psi = \lambda\psi$ with the kernel element $\phi$ to see that the absolute value of the sum of integrals depicted immediately below

$$\int_\Sigma \langle\phi, \phi_\mathrm{c}\rangle + \int_\Sigma \big(\langle\phi, \gamma_t[\chi_{\ddagger}A^1, \phi]\rangle + C^1\langle\phi, \rho_t[\sigma_3, \phi]\rangle\big) - \frac{\lambda}{m}\int_\Sigma |\phi|^2$$

$$(5.72)$$

is no greater than $c_0 m \int_\Sigma |\phi|^2$.

To observe the implication of the preceding observation, note first that what is said by Lemma 5.2 implies that the absolute values of the $\Sigma$ integrals of both $\langle \phi, \gamma_t[\chi_{\ddagger} A^1, \phi] \rangle$ and $c^1 \langle \phi, \rho_t[\sigma_3, \phi] \rangle$ are at most $c_0 \, r^{-2/3}$ times the $\Sigma$ integral of $|\phi|^2$. Meanwhile, supposing that $\phi$ is an $\mathbb{S}_{SL}$ kernel element (thus described by data $(\varsigma_+, \varsigma_-, \mu_3)$ from the second bullet of (5.10)), what is said by Lemma 5.2 about its $\varsigma_-$ and $\mu_3$ parts implies that the $\Sigma$-integral of $\langle \phi, \phi_\varsigma \rangle$ differs from that of $|\phi|^2$ by at most $c_0 r^{-2/3}$ times this same integral of $|\phi|^2$. (This is because the respective $\Sigma$ integrals of $\langle \phi, \phi_\varsigma \rangle$ and $|\phi|^2$ differ from the $\Sigma$-integral of $|\varsigma_+|^2 |\varphi|^2$ by at most $c_0 r^{-2/3}$ times the $\Sigma$ integral of $|\phi|^2$) These last observations when used with the afore-mentioned $c_0 m \int_\Sigma |\phi|^2$ bound for the absolute value of (5.72) need $\lambda$ to obey:

$$|1 - \tfrac{\lambda}{m}| \le c_0 m \,.$$

(5.73)

This bound implies what is asserted by Lemma 5.7 when $m < c_0^{-1}$.

### e) Pointwise bounds for Lemmas 5.5 and 5.6

The subsection ties up some loose ends by proving the asserted pointwise bounds in Lemmas 5.5 and 5.6. These assertions concern solutions on $\Sigma$ to equations that have the schematic form $\mathfrak{D}_0 \psi_{\ddagger} = X$ for a specific version of $X$ in both cases. The lemmas assert $c_0 m^2 r^{-2/3}$ and $c_0 m^2 r^{-2/3} e^{-r \mathfrak{d}^{3/2}/c_0}$ bounds for their versions of $|\psi_{\ddagger}|$ and $|[\hat{\tau}, \psi_{\ddagger}]|$. An important point to keep in mind in the upcoming derivation of these bounds: The versions of $\psi_{\ddagger}$ from these lemmas ($\psi_G$ in the case of Lemma 5.6 and the solution $\psi$ to (5.37) in the case of Lemma 5.5) have $\mathbb{L}$ and $\mathbb{H}$ norm bounds already $c_0 m^2 \, r^{-4/3}$ and $c_0 m^2 \, r^{-2/3}$.

The derivation of the desired bounds has four steps.

*Part 1*: The first step is to derive the sup-norm bound for $|\psi_{\ddagger}|$. To this end, let $\Omega_r$ denote the characteristic function for the $\mathfrak{d} < c_0 r^{-2/3}$ part of $\Sigma$. Use the Bochner-Weitzenboch formula for $\mathfrak{D}_0{}^2$ (this is (4.8) with no $\nabla_{\hat{A}t}$) to derive the inequality

$$\tfrac{1}{2} d^\dagger d |\psi_{\ddagger}|^2 + |\nabla_{\hat{A}} \psi_{\ddagger}|^2 + |[\hat{a}, \psi_{\ddagger}]|^2 \le c_0 r^{4/3} \, \Omega_r \, |\langle \psi_{\ddagger}, [\hat{\tau}, \psi_{\ddagger}] \rangle| + \langle \psi_{\ddagger} \mathfrak{D}_0 X \rangle \,.$$

(5.74)

Fix a point $p \in \Sigma$ where $|\psi_{\ddagger}|$ has its maximum; and then fix a holomorphic coordinate (to be denoted by z) for a neighborhood of p and use this coordinate to identify a neighborhood of p with a radius $c_0^{-1}$ disk about the origin in $\mathbb{C}$. Let r denote for the moment a positive number which is less than the radius of this disk; and let D denote the disk in the coordinate chart where $|z| < r$ (it is assumed here and henceforth that r is small enough so that this disk is well inside the coordinate chart).

The Dirichlet Green's function on D with pole at p for the operator $d^\dagger d$ is the function $\frac{1}{2\pi}|\ln(\frac{|z|}{r})|$. Let $\chi_p$ denote the function $\chi(\frac{2|z|}{r} - 1)$ which is equal to 1 at p and zero near the boundary of D. Multiply both sides of (5.74) by $\chi_p$ times that Green's function and then integrate the resulting inequality over D.) Integration by parts on both sides leads to the inequality in (5.75) below. (The notation in (5.75) and subsequently uses $\|\psi_\ddagger\|_\infty$ to denote the sup-norm of $\psi_\ddagger$. The notation also has $D'$ denoting the $|z| < \frac{1}{4} r$ disk.)

$$\|\psi_\ddagger\|_\infty{}^2 + \int_D \chi_p|\ln(\tfrac{|z|}{r})| |[\hat{a}, \psi_\ddagger]|^2 \le c_0 r^{4/3} \int_D \Omega_r |\ln(\tfrac{|z|}{r})| |\langle \psi_\ddagger, [\hat{\tau}, \psi_\ddagger]\rangle| + c_0 \tfrac{1}{r^2} \int_{D-D'} |\psi_\ddagger|^2$$
$$+ c_0 \int_D \chi_p|\ln(\tfrac{|z|}{r})||X|^2 + c_0(\int_D \tfrac{1}{|z|} |X|)^2 \, .$$
$$(5.75)$$

Noting that the left most integral on the right hand side is restricted to the $\mathfrak{d} < c_0 r^{-2/3}$ part of $\Sigma$, and noting Lemma 4.2, it follows that (5.75) leads to the bound below:

$$\|\psi_\ddagger\|_\infty{}^2 + \int_D \chi_p|\ln(\tfrac{|z|}{r})| |[\hat{a}, \psi_\ddagger]|^2 \le c_0 \|\psi\|_\infty \|\psi_\ddagger\|_\mathbb{H} + c_0 \tfrac{1}{r^2} \int_D |\psi_\ddagger|^2$$
$$+ c_0 \int_D \chi_p|\ln(\tfrac{|z|}{r})||X|^2 + c_0(\int_D \tfrac{1}{|z|} |X|)^2$$
$$. $$
$$(5.76)$$

Also by virtue of Lemma 4.2: If r is the larger of $r^{-2/3}$ and $c_0^{-1}\mathfrak{d}(p)$, which will henceforth be assumed, then $\frac{1}{r^2} \int_D |\psi_\ddagger|^2$ is also bounded by $c_0\|\psi_\ddagger\|_\mathbb{H}{}^2$. Thus with r chosen this way, then (5.76) leads directly to this:

$$\|\psi_\ddagger\|_\infty{}^2 + \int_D \chi_p|\ln(\tfrac{|z|}{r})| |[\hat{a}, \psi_\ddagger]|^2 \le c_0\|\psi_\ddagger\|_\mathbb{H}{}^2 + c_0 \int_D \chi_p|\ln(\tfrac{|z|}{r})||X|^2 + c_0(\int_D \tfrac{1}{|z|} |X|)^2 \, .$$
$$(5.77)$$

Since $\|\psi_\ddagger\|_\mathbb{H} \le c_0 m^4 \, r^{-4/3}$ in the context of both Lemma 5.6 and Lemma 5.5, the desired sup-norm bound will follow if more can be said in each case about the relevant versions of $X$.

*The sup-norm bound for Lemma 5.6*: In this case, the section $\psi_\ddagger$ is $\psi_G$ and what is denoted by $X$ is denoted by $G$ in (5.47). The bound for $|G|$ in (5.48) implies that the $|X|^2$ integral in (5.77) is at most $c_0 m^4 r^2 e^{-r r^{3/2}/c_0}$ which is at most $c_0 m^4 r^{-4/3}$. Meanwhile, that same (5.48) bound implies that the integral with $|X|$ (whose square appears in (5.77)) is at most $c_0 m^2 r \, e^{-r \mathfrak{d}^{3/2}/c_0}$ which is at most $c_0 m^2 r^{-2/3}$. Thus, the right hand side of (5.77) is at most $c_0 m^4 r^{-4/3}$ in the case of Lemma 5.6 which is the desired upper bound for $|\psi_G|$.

*The sup-norm bound for Lemma 5.5*:  In this case, $\psi$ is the solution to (5.37) given by Lemma 5.5.  What is denoted by $X$ is therefore the sum of three terms, the first being $\mathcal{G}$ again, the second being $\psi\#\psi$ and the third depicted below:

$$m\big(\gamma_t[\chi_{\ddagger}A^1, \psi] + c^1\rho_t[\sigma_3, \psi] + (\psi)_\epsilon\big)\,.$$

(5.78)

The right most two terms on the right hand side of (5.77), the terms with integrals of $X$, are no greater than what is obtained by summing $c_0$ times three versions of those terms, the first with $\mathcal{G}$ replacing $X$, the second with $\psi\#\psi$ replacing $X$, and the third with the expression in (5.78) replacing $X$.  In this regard, the two terms on the right hand side of (5.77) with $\mathcal{G}$ replacing $X$ were just discussed in the preceding Lemma 5.6 case.  With this understood, the first case to considered below is that where $X$ is replaced by $\psi\#\psi$.

For the case $X = \psi\#\psi$, the term with the integral of $|X|^2$ in (5.77) is at most

$$c_0\, r\, \big(\textstyle\int_D\ |\psi|^4\big)^{1/2}\, \|\psi\|_\infty^2\,,$$

(5.79)

which is, in turn, no greater than $c_0\, r\, \|\psi\|_{\mathbb{H}}^2 \|\psi\|_\infty^2$.  If $r > c_0$ (as will henceforth be assumed), the term with the $|X|^2$ integral for the case $X = \psi\#\psi$ is at most $\frac{1}{100}\|\psi\|_\infty^2$.  By the same token, the term with the integral of $\frac{1}{|z|}|X|$ when $X$ is $\psi\#\psi$ is at most $c_0\, r^{1/3}\|\psi\|_{\mathbb{H}}^2\|\psi\|_\infty^2$.  This is also at most $\frac{1}{100}\|\psi\|_\infty^2$ when $r > c_0$.

For the case when $X$ is depicted in (5.78):  There are two cases to consider, these being where $\mathfrak{d}(p)$ is less than $c_0 r^{-2/3}$ and where it isn't.  In the first case, $|X| \leq c_0 m^2|\psi_{\ddagger}|^2$ in which case both of the $X$ terms in (5.77) when $X$ is depicted in (5.78) are at most $c_0 m^2\, r^{-4/3}\|\psi\|_\infty^2$ which is at most $\frac{1}{100}\|\psi\|_\infty^2$ when $r > c_0$.  In the second case (which is when $\mathfrak{d}(p) > c_0 r^{-2/3}$), the contributions to the $X$ integrals from the $mc^1\rho_t[\sigma_3, \psi] + (\psi)_\epsilon$ part of $X$ is at most $c_0 m^2\, r^2\, \|\psi\|_\infty^2$ which is at most $\frac{1}{100}\|\psi\|_\infty^2$ if $m < c_0^{-1}$ (which will henceforth be assumed).  Meanwhile, the contribution to the $|X|^2$ integral from the $m\gamma_t[\chi_{\ddagger}A^1, \psi]$ term in $X$ is at most $c_0 m^2$ times the integral with $|[\hat{a}, \psi]|^2$ that appears on the left hand side of (5.77).  This is less than $\frac{1}{100}$ times that left hand side integral if $m < c_0^{-1}$.  (This reference to that left hand integral of $|[\hat{a}, \psi]|^2$ is the reason for keeping it in play.)

As for the contribution of $m\gamma_t[\chi_{\ddagger}A^1, \psi]$ to the term with the square of the integral of $\frac{1}{|z|}|X|$:  Here it proves useful to write $\psi$ as $\psi_{\mathcal{G}} + \phi$ with $\phi$ as in (5.55).  For the $\psi_{\mathcal{G}}$ contribution:  A crucial point in the analysis for this contribution is that (by virtue of (5.48)), the assumptions of Lemma 4.4 are met when $\psi$ in that lemma is $\psi_{\mathbb{G}}$. The conclusions of that lemma imply the following inequality for $\psi_{\mathbb{G}}$:

$$\int_D r^4 \mathfrak{d}^6 |[\hat{\tau}, \psi_{\mathbb{G}}]|^2 \leq c_0 \|\psi_{\mathbb{G}}\|_{\mathbb{L}}^2 \,.$$

(5.80)

Keeping this inequality in mind: Use the description of $A^1$ from (3.23) and Proposition 3.4 to derive the following bound

$$|\gamma_t[\chi_{\ddagger} A^1, \psi_{\mathbb{G}}]| \leq c_0 \|\psi_{\mathbb{G}}\|_\infty^{1/2} (|\psi_{\mathbb{G}}|^2 + r^4 \mathfrak{d}^6 |[\hat{\tau}, \psi_{\mathbb{G}}]|^2)^{1/4}$$

(5.81)

Since $(\frac{1}{|z|})^{4/3}$ is integrable on D (its integral is $c_0 \, r^{2/3}$), the inequality in (5.81) with the integral bound in (5.80) leads to this:

$$\int_D \frac{1}{|z|} |\gamma_t[\chi_{\ddagger} A^1, \psi_{\mathbb{G}}]| \leq c_0 \, r^{1/3} \|\psi_{\mathbb{G}}\|_\infty^{1/2} \|\psi_{\mathbb{G}}\|_{\mathbb{L}}^{1/2} \,,$$

(5.82)

which is an inequality whose right hand side is at most $c_0 m^2 r^{-1}$. Therefore, the contribution to the right most term on the right hand side of (5.77) from $m\gamma_t[\chi_{\ddagger} A^1, \psi]$ is at most $m^5 r^{-2}$.

And, finally, there is the contribution of $m\gamma_t[\chi_{\ddagger} A^1, \phi]$ to the right most term in (5.77) (the square of the integral of $\frac{1}{|z|} |X|$). Breaking the integral over D of the $X = m\gamma_t[\chi_{\ddagger} A^1, \phi]$ version of $\frac{1}{|z|} |X|$ into the part where $|z| > r^{-1}$ and the part where $|z| < r^{-1}$ leads to this:

$$\int_D \frac{1}{|z|} m \, |\gamma_t[\chi_{\ddagger} A^1, \phi]| \leq c_0 \, m \ln(r \, \mathfrak{d}) \|\phi\|_{\mathbb{H}} + c_0 m \|\phi\|_\infty \,.$$

(5.83)

Since $\|\phi\|_{\mathbb{H}} \leq c_0 m^{3/2} r^{-4/3}$ and $\|\phi\|_\infty \leq \|\psi\|_\infty + c_0 m^2 r^{-2/3}$, the right hand side of this no greater than $(c_0 m^2 r^{-2/3} + \frac{1}{100} \|\psi\|_\infty)$.

Putting all of the preceding bounds leads from (5.77) to an inequality that says this:

$$\|\psi\|_\infty^2 \leq \frac{1}{2} \|\psi\|_\infty^2 + c_0 m^2 \, r^{-4/3}$$

(5.84)

when $m < c_0^{-1}$ and $r > c_0$. Lemma 5.5's sup-norm bound follows directly from (5.84).

*Part 2*: This part of the proof and the next part establish the assertion in Lemma 5.6 to the effect that $|[\hat{\tau}, \psi_{\mathbb{G}}]| \leq c_0 m^2 \, r^{-2/3} e^{-r \mathfrak{d}^{3/2}/c_0}$ and the assertion in Lemma 5.5 to the effect that $|[\hat{\tau}, \psi]| \leq c_0 m^2 \, r^{-2/3} e^{-r \mathfrak{d}^{3/2}/c_0}$ assertion in Lemma 5.5. To this end, let $\psi_{\ddagger}$ again denote

either $\psi_G$ or Lemma 5.5's section $\psi$ as the case may be. Supposing that $c > 1$, define a seminorm for $\psi_\ddag$ (or for any other section of $\mathbb{S}$) by the rule

$$\psi_\ddag \to \|\psi_\ddag\|_c \equiv \sup_\Sigma e^{r\mathfrak{d}^{3/2}/c} |[\hat{\tau}, \psi_\ddag]|$$

(5.85)

The task now is to prove that $\|\psi_\ddag\|_c < c_0 m^2 r^{-4/3}$ for some choice of $c \le c_0$ given that $m$ is small ($m < c_0^{-1}$) and $r$ is large ($r > c_0$).

To start the task, reintroduce the equation $\mathfrak{D}_0\psi_\ddag = X$ with $X$ as before. That equation implies that $\mathfrak{D}_0([\hat{\tau}, \psi_\ddag]) = X'$ with $X'$ depicted below:

$$X' = \gamma_k[\nabla_{\hat{A}k}\hat{\tau}, \psi_\ddag] + \rho_k[[\hat{a}_k, \hat{\tau}], \psi_\ddag] + [\hat{\tau}, X].$$

(5.86)

The rest of this part of the proof constitutes a digression to point out some key features of $X'$. The first feature is that the norms of the left most two terms on the right hand side in (5.86) obey the bound

$$|\gamma_k[\nabla_{\hat{A}k}\hat{\tau}, \psi_\ddag]| + |\rho_k[[\hat{a}_k, \hat{\tau}], \psi_\ddag]| < c_0 r^{2/3} e^{-r\mathfrak{d}^{3/2}/c_0} \|\psi_\ddag\|_\infty < c_0 m^2 e^{-r\mathfrak{d}^{3/2}/c_0},$$

(5.87)

the right most inequality being due to the fact that $|\psi_\ddag\|_\infty \le c_0 m^2 r^{-2/3}$.

The second key feature of $X'$ concerns the norm of the right most term in (5.86) which is $[\hat{\tau}, X]$. In the case when $\psi_\ddag = \psi_G$ from Lemma 5.6, this term is a priori bounded by the function on the far right in (5.86) because $X$ is the section $G$ from (5.47) and (5.48). In the case when $\psi_\ddag$ is the section $\psi$ from Lemma 5.5, the section $X$ is the sum of $G$ and $\psi\#\psi$ and what is depicted in (5.78). The important point to note in this regard is that

$$|\psi\#\psi| \le c_0\|\psi\|_\infty|[\hat{\tau}, \psi]| \le c_0 m^2 r^{-2/3}|[\hat{\tau}, \psi]|.$$

(5.88)

Meanwhile, the norm of the commutator of $\hat{\tau}$ with the expression in (5.78) is at most

$$c_0 \, mr \, \mathfrak{d}^{3/2} \, \chi_\ddag |[\hat{\tau}, \psi]| + m\|\psi\|_\infty e^{-r\mathfrak{d}^{3/2}/c_0} + c_0 m|[\hat{\tau}, \psi]|.$$

(5.89)

The preceding bounds imply the following for either case of $\psi_\ddag$: Let $c_1$ denote the version of $\kappa$ from Lemma 4.4. If $m < c_0^{-1}$ and if $r > c_0$, then $|[\hat{\tau}, \mathcal{D}\psi]|$ (which is $X'$) obeys a version of (4.22)'s inequality

$$|[\hat{\tau}, \mathcal{D}\psi_\ddag]| \equiv |X'| \le \min(m, \frac{1}{1000C_1}) \, r \, \mathfrak{d}^{3/2} \, |[\hat{\tau}, \psi_\ddag]| + c_0 m^2 \, e^{-r\mathfrak{d}^{3/2}/c_0}.$$

(5.90)

In particular, if $m < c_0^{-1}$, then the conclusions of Lemma 4.4 can be invoked for either case of $\psi_{\ddagger}$ (this was already noted for the case of $\psi_{\mathbb{G}}$.

*Part 3*: The identity $\mathfrak{D}_0([\hat{\tau}, \psi_{\ddagger}]) = \mathcal{X}'$ leads to an analog of (5.74) which says this where $\mathfrak{d} > c_0 r^{-2/3}$:

$$\tfrac{1}{2} d^{\dagger}d|[\hat{\tau}, \psi_{\ddagger}]|^2 + |\nabla_{\hat{A}}[\hat{\tau}, \psi_{\ddagger}]|^2 + |[\hat{a}, \psi_{\ddagger}]|^2 \leq \langle \psi_{\ddagger} \mathfrak{D}_0 \mathcal{X}' \rangle .$$

(5.91)

To see the implications of this, let p denote a point in $\Sigma$ where the supremum in (5.85) is achieved. No generality is lost by assuming that $\mathfrak{d}(p) > 100 \, c^{2/3} r^{-2/3}$ because if this isn't the case, then $\|\psi_{\ddagger}\|_c \leq e^{1000} \|\psi_{\ddagger}\|_{\infty}$ which is less than $c_0 m^2 r^{-2/3}$. Anyway, assuming $c > c_0$ and this lower bound for $\mathfrak{d}$, fix a holomorphic coordinate (again denoted by z) for the radius $c_0^{-1}$ disk centered at p. Let r denote $c_0^{-1} \mathfrak{d}(p)$ with $c_0$ here chosen so that the $|z| < r$ disk in $\mathbb{C}$ sits well inside the domain of the coordinate z (and thus the $|z| < r$ disk can be viewed at one's discretion as an open set in $\Sigma$). Let D again denote this $|z| < r$ disk, and let D′ denote the concentric disk with radius r/2. There is again a version of (5.75) which says

$$\|\psi_{\ddagger}\|_c^2 + e^{2r\mathfrak{d}(p)^{3/2}/c} \int_D \chi_p |\ln(\tfrac{|z|}{r})| |[\hat{a}, \psi_{\ddagger}]|^2 \leq c_0 \tfrac{1}{r^2} e^{2r\mathfrak{d}(p)^{3/2}/c} \int_{D-D'} |\psi_{\ddagger}|^2$$
$$+ c_0 \, e^{2r\mathfrak{d}(p)^{3/2}/c} \int_D \chi_p |\ln(\tfrac{|z|}{r})| |\mathcal{X}'|^2 + c_0 e^{2r\mathfrak{d}(p)^{3/2}/c} (\int_D \tfrac{1}{|z|} |\mathcal{X}'|)^2 .$$

(5.92)

The task now is to find suitable bounds for each of the terms on the right hand side of (5.92). This is done in the five steps that follow.

Step 1: For the left most term on the right hand side: Let $c_1$ again denote the version of $\kappa$ from Lemma 4.4. An appeal to Lemma 4.4 bounds the left most term on the right hand side of (5.92) by

$$c_0 \tfrac{1}{r^2} e^{2r\mathfrak{d}(p)^{3/2}/c} \, r^{-4/3} \, e^{-r\mathfrak{d}(p)^{3/2}/c_0 c_1} \|\psi\|_{\infty}^2 \ .$$

(5.93)

Note in particular that if $c > c_0 c_1$, then the bound above is no greater than $c_0 m^2 r^{-4/3}$.

Step 2: This step considers the term on the right hand side of (5.92) with the integral of $|\mathcal{X}'|^2$. To this end, let $c_2$ denote the version of $c_0$ that appears in (5.90). If $m < c_0^{-1}$, then the contribution to $|\mathcal{X}'|^2$ from the $r\mathfrak{d}^{1/2} |[\hat{\tau}, \psi_{\ddagger}]|$ term in (5.90) is dominated by the $|[\hat{a}, \psi_{\ddagger}]|^2$ integral on the left hand side of (5.92) at the expense of multiplying the

right most term in (5.90) by $c_0$. With this point taken, then term on the right hand side of (5.92) with the $|\mathcal{X}|^2$ integral can be replaced by

$$c_0\, m^2 r^{-4/3}\, e^{2r\mathfrak{d}(p)^{3/2}/c}\, e^{-r\mathfrak{d}(p)^{3/2}/c_0 c_2}\ .$$

(5.94)

Note in particular that this is no greater than $c_0 m^2 r^{-4/3}$ if $c > c_0 c_2$.

Step 3: The issue comes down now to the right most term in (5.92), the term with the square of the integral of $\frac{1}{|z|}\,|\mathcal{X}|$. The contribution of this term from the right most term on the right hand side of (5.90) is bounded by what is depicted in (5.94) with some larger value of $c_0$.

To bound the contribution of the left most term on the right hand side of (5.90), first bound the contribution to the integrand by

$$\frac{1}{|z|}\, r\, \mathfrak{d}^{3/2}\, |[\hat{\tau}, \psi_{\ddagger}]| \leq c_0\, \|\psi\|_\infty^{1/2}\, \frac{1}{|z|}\, (r^4\, \mathfrak{d}^6 |[\hat{\tau}, \psi_{\ddagger}]|^2)^{1/4}\ .$$

(5.95)

As noted after (5.81), the function $(\frac{1}{|z|})^{4/3}$ is integrable on D with its integral bounded by $c_0 \mathfrak{d}(p)^{2/3}$. Meanwhile, $r^4\, \mathfrak{d}^6 |[\hat{\tau}, \psi_{\ddagger}]|^2$ is also integrable on D (appeal to Lemma 4.4). Granted these facts, then one is lead to bound the contribution of $\frac{1}{|z|}\, r\, \mathfrak{d}^{3/2}\, |[\hat{\tau}, \psi_{\ddagger}]|$ to the right most term in (5.92) by

$$c_0\, e^{2r\mathfrak{d}(p)^{3/2}/c}\, (m^2\, r^{-2/3} e^{-r\mathfrak{d}(p)^{3/2}/c_0 c_2}\, )^2\ .$$

(5.96)

Note in particular that this contribution is no greater than $c_0 m^2 r^{-2/3}$.

Step 4: To summarize: Supposing that $m < c_0^{-1}$, that $r > c_0$ and that $c > c_0 c_1 + c_0 c_2$, then the bounds from Steps 1-3 when used for the right hand side of (5.92) lead to a bound of the form $\|\psi_{\ddagger}\|_\infty \leq c_0 m^2 r^{-4/3}$ bound. And, given the definition of $\|\cdot\|_c$, that bound implies in turn the desired pointwise bound $|[\hat{\tau}, \psi_{\ddagger}]| \leq c_0 m^2 r^{-4/3} e^{-r\mathfrak{d}^{3/2}/c_0 c_2}$.

## 6. Putting the slices together to solve the $m \neq 0$ Vafa-Witten equations

This final section takes as input Proposition 5.3's approximate solutions to the small $m$ versions of (1.7) and uses these with the perturbative construction mentioned in Section 4a to obtain an $|\omega|$-divergent sequence of honest solution to those small $m$ equations.

To see how the construction will go, fix $m < c_0^{-1}$ and then $r > c_0$ (with $r = \frac{n}{m}$ so as to define the bundle $E_q$). Let $(A, \mathfrak{a})$ denote the pair of SU(2) connection for $E_q$ and 1-form on $S^1 \times \Sigma$ with values in the bundle of skew-Hermitian endomorphisms of $E_q$ given by Proposition 5.3. As noted in Section 4a, the equations in (1.7) are equivalent to the equations in (4.3) for a section $\psi$ of $\mathbb{S}$. By way of a reminder, these equations have the form

$$\mathcal{D}\psi + m\psi_c + \psi\#\psi + \mathcal{G} = 0$$

(6.1)

with $\mathcal{D} = \gamma_k \nabla_{Ak} + \rho_k[\mathfrak{a}_k, \cdot]$ and with $-\mathcal{G}$ representing here the left hand side of the three bullets in Proposition 5.3.

## a) Proof of Theorem 1

Theorem 1 follows directly if (6.1) can be solved when $|m|$ is non-zero but small. To see about solving (6.1), let $\mathcal{L}$ denote the operator $\psi \to \mathcal{D}\psi + m\psi_c$. Since $\mathcal{L}$ differs from $\mathcal{D}$ by a symmetric endomorphism of $\mathbb{S}$, it is a self-adjoint operator on $\mathbb{L}$ with $\mathbb{H}$ being a dense domain; and it has discrete spectrum with no accumulation points and finite multiplicities. Assume for the moment that 0 is not an eigenvalue $\mathcal{L}$. In that case, $\mathcal{L}$ has an inverse mapping $\mathbb{L}$ to $\mathbb{H}$; and then a solution to (6.1) is a fixed point of the map from $\mathbb{H}$ to $\mathbb{H}$ depicted below:

$$\psi \to \mathbb{T}(\psi) = -\mathcal{L}^{-1}(\psi\#\psi + \mathcal{G}).$$

(6.2)

In this regard, the bundle map $\psi \to \psi\#\psi$ defines a smooth map from $\mathbb{H}$ to $\mathbb{L}$ by virtue of a Sobolev inequality. (See the upcoming Lemma 6.1.)

To see about a fixed point, let $\mathrm{N}$ denote the norm of $\mathcal{L}^{-1}$ (which is the infimum over the non-zero elements in $\mathbb{L}$ of the ratio $\|\mathcal{L}^{-1}(\eta)\|_{\mathbb{H}}/\|\eta\|_{\mathbb{L}}$). Then, for any section $\psi$ of $\mathbb{S}$,

$$\|\mathbb{T}(\psi)\|_{\mathbb{H}} \leq c_0 \mathrm{N} \left( \|\psi\#\psi\|_{\mathbb{L}} + \|\mathcal{G}\|_{\mathbb{L}} \right).$$

(6.3)

Fix a positive number $\mathrm{R}$ such that the inequality below holds

$$\|\psi_1 \# \psi_2\|_{\mathbb{L}} \leq \mathrm{R}\|\psi_1\|_{\mathbb{H}}\|\psi_2\|_{\mathbb{H}}$$

(6.4)

for any two elements $\psi_1$ and $\psi_2$ from $\mathbb{H}$. Because $\|\mathcal{G}\|_{\mathbb{L}} \leq c_0 m^2 \, e^{-r/c_0}$ (see Proposition 5.3), the inequality in (6.3) says in effect that

$$\|\mathbb{T}(\psi)\|_{\mathbb{H}} \leq c_0 \mathrm{N}(\mathrm{R}\|\psi\|_{\mathbb{H}}^2 + m^2 \, e^{-r/c_0}).$$

(6.5)

Thus, if $c_\ddagger$ is a sufficiently large number ($c_\ddagger > c_0$), then $\mathbb{T}$ maps the radius $(c_\ddagger NR)^{-1}$ ball about the origin in $\mathbb{H}$ to itself if the inequality

$$m^2 \, e^{-r/c_0} < (c_\ddagger N^2 R)^{-1}$$

(6.6)

holds. Moreover, because

$$\|\mathbb{T}(\psi_1) - \mathbb{T}(\psi_2)\|_\mathbb{H} \leq c_0 N \left( \|(\psi_1 - \psi_2)\#\psi_1\|_\mathbb{L} + \|\psi_2\#(\psi_1 - \psi_2)\|_\mathbb{L} \right),$$

(6.7)

the map $\mathbb{T}$ will be a contraction mapping on that radius $(c_\ddagger NR)^{-1}$ ball about the origin in $\mathbb{H}$ if it is the case that $c_\ddagger > c_0$ and (6.6) holds.

To summarize: If (6.6) holds with $c_\ddagger > c_0$, then $\mathbb{T}$ is a contraction mapping on the radius $(c_\ddagger NR)^{-1}$ ball about the origin in $\mathbb{H}$. As a consequence of that, the map $\mathbb{T}$ has a unique fixed point in that ball. This fixed point is the sought after solution to (6.1). (Moreover, the $\mathbb{H}$-norm of that fixed point will be at most $c_0 N m^2 \, e^{-r/c_0}$ ).

By virtue of (6.6), the size of the two numbers R and N are the crucial input. As for the R: The upcoming Lemma 6.1 says that $R < c_0 r^{-1/3}$. As for N: This number can be bounded from above in terms of the smallest of the norms of the eigenvalues of $\mathcal{L}$. To elaborate, let $\lambda_r$ denote that smallest norm. It then follows from what is said by Lemmas 4.2 and 4.3 that N obeys

$$N \leq c_0 r^{2/3} \frac{1}{\lambda_r}.$$

(6.8)

Thus the bound in (6.6) will hold if

$$m^2 \, e^{-r/c_0} < c_\ddagger^{-1} r^{-1} \lambda_r{}^2.$$

(6.9)

Hence, the key question is: How small is $\lambda_r$? (The subscript $r$ is a reminder that this smallest eigenvalue depends on the input $r$.) The upcoming Proposition 6.2 says in effect that $\lambda_r \geq c_0 m^{-1} r^{-2}$. Supposing this is correct, then (6.9) will hold if $r$ is sufficiently large. As a consequence of that, there is a solution to (1.7) for any given small $m$ ($m < c_0^{-1}$) if $r$ is sufficiently large. Thus, taking $r$ ever larger (subject to the constraint $rm \in 4\mathbb{Z}$), there is an $|\omega|$-divergent sequence of solutions to (1.7).

What follows directly are first Lemma 6.1 and then Proposition 6.2.

**Lemma 6.1**: *There exists* $\kappa > 1$ *with the following significance: If $m$ is less than* $\frac{1}{\kappa}$ *and if* $r \equiv \frac{n}{m}$ *is greater than* $\kappa$*, then the bound below holds for any two elements* $\psi_1$ *and* $\psi_2$ *in* $\mathbb{L}_1$:

$$\|\psi_1 \# \psi_2\|_{\mathbb{L}} \leq \kappa r^{-1/3} \|\psi_1\|_{\mathbb{H}} \|\psi_2\|_{\mathbb{H}}.$$

This lemma is proved momentarily.

The promised Proposition 6.2 is next.

**Proposition 6.2**: *There exists* $\kappa > 1$ *such that if* $\mathfrak{m}$ *is less than* $\frac{1}{\kappa}$ *and if* $r \equiv \frac{\mathfrak{n}}{m}$ *is greater than* $\kappa$, *then the following is true: Define the pair* $(A, \mathfrak{a})$ *as done in Proposition 5.3; and use this pair to define the operator* $\mathcal{D}$ *and then set* $\mathcal{L} = \mathcal{D} + \mathfrak{m}(\cdot)_c$. *The smallest of the norms of the eigenvalues of this version of* $\mathcal{L}$ *is no less than* $\frac{1}{\kappa} \mathfrak{m} r^{-2}$.

The proof of this proposition is in the Section 7d. The remaining parts of this section (Section 6) and Sections 7a-7c make preliminary observations that are needed for the proof of the proposition. As noted in the introduction, many of the observations in the later parts of this section regarding Proposition 6.2 have antecedants in the work of Greg Parker (see [P1] and also [P3]). Meanwhile, the lemmas and observations in Sections 7a-7c involve intricate and delicate estimates that are very specific to the problem at hand.

By way of a look ahead to the proof: The proof of Proposition 6.2 first explains why the small normed eigenvalue of $\mathcal{L}$ are nearly identical to small normed eigenvalues of a certain differential operator on the surface $\Sigma$. This operator on $\Sigma$ turns out to be a small perturbation of the operator $\mathfrak{L}$ that is depicted in (5.56). The preceding observation about corresponding small normed eigenvalues of $\mathcal{L}$ and the perturbation of the operator $\mathfrak{L}$ is the jist of what is said in the upcoming Lemmas 6.6-6.8. Section 7 then verifies (a lengthy process indeed) that there are no eigenvalues in a small interval around zero of the relevant perturbation of $\mathfrak{L}$.

*Proof of Lemma 6.1*: There is the dimension 3 Sobolev inequality (see Theorem 7.10 in [GP]) that says this: Supposing that $\mathfrak{f}$ is any function on $S^1 \times \Sigma$ with $|d\mathfrak{f}|^2$ having finite integral, then $\mathfrak{f}^6$ has finite $S^1 \times \Sigma$ integral and, moreover,

$$\left(\int_{S^1 \times \Sigma} \mathfrak{f}^6\right)^{1/3} \leq c_0 \int_{S^1 \times \Sigma} \left(|d\mathfrak{f}|^2 + \mathfrak{f}^2\right).$$

(6.10)

In particular, this inequality holds when $\mathfrak{f} = |\psi|$ for any given element $\psi$ from $\mathbb{H}$ in which case the right hand side is no larger than $c_0 \|\psi\|_{\mathbb{H}}^2$. This Sobolev inequality will be used momentarily.

For the purposes of this proof, let $\Sigma_r$ denote the part of $\Sigma$ where $\mathfrak{d} < c_0 r^{-2/3}$. The $S^1 \times \Sigma_r$ part of the integral of $|\psi_1 \# \psi_2|^2$ is no greater than that of $|\psi_1|^2 |\psi_2|^2$.

$$(\int_{S^1 \times \Sigma_r} |\psi_1|^4)^{1/2}(\int_{S^1 \times \Sigma_r} |\psi_2|^4)^{1/2}.$$

(6.11)

Meanwhile, the integral above is at most $c_0 r^{-2/3} \|\psi_1\|_{\mathbb{H}}^2 \|\psi_2\|_{\mathbb{H}}^2$ because the following inequalities hold for any $\psi$ from $\mathbb{H}$: First,

$$\int_{S^1 \times \Sigma_r} |\psi|^4 \le c_0 \, (\int_{S^1 \times \Sigma_r} |\psi|^6)^{1/2}(\int_{S^1 \times \Sigma_r} |\psi|^2)^{1/2}.$$

(6.12)

Second, that integral of $|\psi|^6$ is at most $c_0 \|\psi\|_{\mathbb{H}}^6$ by virtue of (6.10). Third, the integral $S^1 \times \Sigma_r$ integral of $|\psi|^2$ is at most $c_0 r^{-4/3}$ times the integral of $\mathfrak{d}^{-2}|\psi|^2$ over the same domain which is no greater than $c_0 r^{-4/3} \|\psi\|_{\mathbb{H}}^2$ (see Lemma 4.2).

To bound the integral of $|\psi_1 \# \psi_2|^2$ over the rest of $S^1 \times \Sigma$, note first that

$$|\psi_1 \# \psi_2| \le c_0 |\psi_1| |[\sigma, \psi_2]| + c_0 |\psi_2| |[\sigma, \psi_1]|$$

(6.13)

because the components of $(\cdot)\#(\cdot)$ are linear combinations of commutators. With (6.13) in mind, note next that

$$\int_{S^1 \times (\Sigma - \Sigma_r)} |\psi_1|^2 |[\sigma, \psi_2]|^2 \le c_0 (\int_{S^1 \times (\Sigma - \Sigma_r)} |\psi_1|^6)^{1/3} (\int_{S^1 \times (\Sigma - \Sigma_r)} |[\sigma, \psi_2]|^3)^{2/3},$$

(6.14)

and then note that

$$\int_{S^1 \times (\Sigma - \Sigma_r)} |[\sigma, \psi_2]|^3 \le c_0 (\int_{S^1 \times (\Sigma - \Sigma_r)} |\psi_2|^6)^{1/4} (\int_{S^1 \times (\Sigma - \Sigma_r)} |[\sigma, \psi_2]|^2)^{3/4}.$$

(6.15)

By virtue of Lemma 4.2 and by virtue of (6.10), these inequalities (and the corresponding one with $\psi_1$ and $\psi_2$ switching roles) lead directly to the bound below:

$$\int_{S^1 \times (\Sigma - \Sigma_r)} |\psi_1 \# \psi_2|^2 \le c_0 \, r^{-2/3} \|\psi_1\|_{\mathbb{H}}^2 \|\psi_2\|_{\mathbb{H}}^2.$$

(6.16)

The preceding bound with the $c_0 \, r^{-2/3} \|\psi_1\|_{\mathbb{H}}^2 \|\psi_2\|_{\mathbb{H}}^2$ bound for (6.11) give the bound that is claimed by Lemma 6.1.

### b) The $\nabla_{At}$-covariant derivative of eigenvectors of $\mathcal{L}$ with small normed eigenvalue

Use the pair $(A, \mathfrak{a})$ from Proposition 5.3 (as defined for $m < c_0^{-1}$ and $r > c_0$) to defined the operator $\mathcal{L}$. The central observation in this subsection is a lemma regarding the size of $\nabla_{At}\psi$ when $\psi$ is an eigenvector of $\mathcal{L}$ with its eigenvalue having small norm (the bound also holds for eigenvectors of the $(A, \mathfrak{a})$ version of $\mathcal{D}$). By way of a remark in this regard: The directional covariant derivative along the $S^1$ factor of $S^1 \times \Sigma$ does not commute with either $\mathcal{L}$ or $\mathcal{D}$ when $m \neq 0$. In fact, it is likely that no lift of the $S^1$ action on $S^1 \times \Sigma$ to an action on the bundle $E_q$ with a generator that does commutes with these operators.

**Lemma 6.3**: *There exists* $\kappa > 1$ *with the following significance: Fix* $m < \kappa^{-1}$ *and then* $\frac{n}{m} \equiv r > \kappa$. *and let* $(A, \mathfrak{a})$ *denote the pair of connection on* P *and 1-form with values in* $\mathrm{ad}(P)$ *from Proposition 5.3. Let* $\psi$ *denote an eigenvector of* $\mathcal{L}$ *with* $\mathbb{L}$-*norm equal to 1 and with eigenvalue between* $-\frac{1}{\kappa}$ *and* $\frac{1}{\kappa}$. *Then*

$$\int_{S^1 \times \Sigma} |\nabla_{\mathrm{At}} \psi|^2 \leq m\, \kappa \ .$$

***Proof of Lemma 6.3***: The proof of the lemma has four parts.

*Part 1*: Suppose that $\psi$ is an eigenvector of $\mathcal{L}$ and that $\lambda$ is its eigenvalue. This part of the proof explains why the assumptions of Lemma 4.4 and thus that lemma's conclusions hold for $\psi$ if $|\lambda|$ is at most $c_0^{-1} r^{2/3}$. (The case when $|\lambda| < 1$ is the relevant case.) To see about Lemma 4.4 for $\psi$: Write the eigenvalue equation $\mathcal{L}\psi = \lambda\psi$ as $\mathcal{D}\psi + \psi_{\mathfrak{c}} = \lambda\psi$. Thus, $[\sigma, \mathcal{D}\psi] = \lambda[\sigma, \psi] - [\sigma, \psi]_{\mathfrak{c}}$. But this is not quite good enough because Lemma 4.4 refers to the version of $\mathcal{D}$ defined by a pair that is described by (4.10) and (4.11); for example, this is the case for the pair from (3.25) and Proposition 3.4. Let $(A^1, \mathfrak{a}^1)$ denote that pair, and write $(A, \mathfrak{a})$ in terms of $(A^1, \mathfrak{a}^1)$ as $(A^1 + \mathfrak{e}, \mathfrak{a}^1 + \mathfrak{f})$ with the components of $(\mathfrak{e}, \mathfrak{f})$ coming from (5.33) - (5.36). Let $\mathcal{D}_1$ now denote the $(A^1, \mathfrak{a}^1)$ version of $\mathcal{D}$. The $(A, \mathfrak{a})$ version (still denoted by $\mathcal{D}$) is related to the operator $\mathcal{D}_1$ by the rule below:

$$\mathcal{D} = \mathcal{D}_1 + m\ (\gamma_k[\mathfrak{e}_k, \cdot\ ] + \rho_k[\mathfrak{f}_k, \cdot\ ])$$

(6.17)

What with the bounds from (5.34) and (5.36), this last inequality implies in turn this one:

$$|[\sigma, \mathcal{D}_1\psi]| \leq c_0(|\lambda| + m)|[\sigma, \psi] + c_0 m\ \mathrm{e}^{-r\mathfrak{d}^{3/2}/c_0}|\psi|.$$

(6.18)

Since the pair $(A^1, \mathfrak{a}^1)$ obeys the conditions set forth in (4.10) and (4.11), this last inequality for the norm of $[\sigma, \mathcal{D}_1\psi]$ puts Lemma 4.4 in play for $\psi$.

*Part 2*:  As done in Part 1, write the eigenvalue equation $\mathcal{L}\psi = \lambda\psi$ as $\mathcal{D}\psi + m\psi_c = \lambda\psi$ and then write $\mathcal{D}$ as done in (6.17).  Now separate the covariant derivative in $\mathcal{D}_1$ along the $S^1$ factor from the rest of $\mathcal{D}_1$ by writing it as

$$\mathcal{D}_1 = \gamma_t \nabla_{A^1 t} + \mathfrak{D}_1$$

$$(6.19)$$

Letting h denote some given function on $\Sigma$, use (6.19) and the anti-commutation rules in (4.6) to derive the identity below:

$$\int_{S^1 \times \Sigma} h|\mathcal{D}_1\psi|^2 = \int_{S^1 \times \Sigma} h|\nabla_{A^1 t}\psi|^2 + \int_{S^1 \times \Sigma} h|\mathfrak{D}_1\psi|^2 - \int_{S^1 \times \Sigma} \langle \gamma_a(\nabla_a h)\psi, \gamma_t \nabla_{A^1 t}\psi\rangle$$
$$+ \int_{S^1 \times \Sigma} h\langle\psi, \gamma_t\gamma_a[F_{A^1 ta}, \psi]\rangle + \int_{S^1 \times \Sigma} h\langle\psi, \gamma_t\rho_k[\nabla_{A^1 t}\mathfrak{a}^1{}_k, \psi]\rangle$$

$$(6.20)$$

where the notation is as follows:  The index a is summed over the set $\{1, 2\}$ that label the components of oriented, orthonormal frames for $T^*\Sigma$ over disks in $\Sigma$.  In this regard, $\nabla_1 h$ and $\nabla_2 h$ are the derivatives of h in the directions of the dual vector fields to the frame components.  Finally, $F_{A^1 t1}$ and $F_{A^1 t2}$ are the components of the contraction of $F_{A^1}$ with the vector field $\frac{\partial}{\partial t}$.

With regards to $F_{A^1 t1}$ and $F_{A^1 t2}$:  Their norms are bounded by $c_0 m \, r^{2/3} \, (1 + r^{1/3} \, \mathfrak{d}^{1/2})$, and the norms of their commutator with $\sigma$ are bounded by $c_0 m r^{2/3} e^{-r\mathfrak{d}^{3/2}/c_0}$.  (See (4.10).)  With regards to $\nabla_{A^1 t}\mathfrak{a}^1$:  Its norm is bounded by $c_0 m \, r^{2/3} e^{-r\mathfrak{d}^{3/2}/c_0}$.  (Also see (4.10).)

To exploit the preceding observations, construct a bump function to be denoted by $\chi_1$ which is given by the rule $\chi_1(\cdot) = \chi(\frac{800\,\mathfrak{d}}{r_0} - 1)$.  To be sure, this function is equal to 1 where $\mathfrak{d} < \frac{1}{800} r_0$ and it is equal to zero where $\mathfrak{d} > \frac{1}{400} r_0$.  Now let $\psi$ denote an eigenvector for $\mathcal{L}$ and let $\lambda$ denote its eigenvalue.  Use (6.20) with h = $(1 - \chi_0)$ with (6.17) to see that

$$\int_{S^1 \times \Sigma} (1 - \chi_1)|\nabla_{A^1 t}\psi|^2 \leq (|\lambda| + c_0 m) \int_{S^1 \times \Sigma} (1 - \chi_1)|\psi|^2 - \int_{S^1 \times \Sigma} \langle \gamma_a(\nabla_a \chi_1)\psi, \gamma_t \nabla_{A^1 t}\psi\rangle$$
$$+ c_0 \, m r \int_{S^1 \times \Sigma} (1 - \chi_1)|[\sigma, \psi]|^2$$

$$(6.21)$$

Save the preceding inequality for the moment.  What follows directly in Part 3 derives an analogous inequality on the part of $S^1 \times \Sigma$ where $\mathfrak{d} < \frac{1}{400} r_0$.

*Part 3*:  To start the derivation, fix a point in $\Theta$ and let D $\subset \Sigma$ denote the $\mathfrak{d} < r_0$ disk centered on that point.  Use the holomorphic coordinate system on D that was introduced in Lemma 3.2.  Also, use the vector bundle isomorphism from Section 3b between $E_{q|D}$ and the product $\mathbb{C} \oplus \mathbb{C}$ bundle to depict the connection $A^1$ on $S^1 \times D$ as in (3.25).  Let $\nabla_t$ denote the covariant directional derivative along the $S^1$ factor of $S^1 \times D$ for the connection $\theta_{\mathbb{C} \oplus \mathbb{C}}$.

Then write $\nabla_{A^1 t}$ as $\nabla_t + m\,A^1\,dt$ with $A^1$ denoting an anti-Hermitian endomorphism of the produce $\mathbb{C} \oplus \mathbb{C}$ bundle, and use this depiction of $\nabla_{A^1 t}$ to write

$$\mathcal{D}_1\psi = \gamma_t \nabla_t \psi + \mathfrak{D}_1\psi + m\gamma_t[A^1, \psi]$$

(6.22)

with $\mathfrak{D}_1$ being the same operator here as in (6.19). Three key points with regards to (6.22) are that $\nabla_t$ commutes with $\mathfrak{D}_1$, that $|A^1| \leq c_0 r\, \eth^{3/2}$, and that $|[\sigma, A^1]| \leq c_0\, e^{-r\eth^{3/2}/c_0}$.

There is an analog of (6.21) on $S^1 \times D$ with $\mathcal{D}_1$ replaced by $\gamma_t \nabla_t + \mathfrak{D}_1$ with it understood that the function h has compact support on $S^1 \times D$. This analog has the same form but for the replacement of $\nabla_{A^1 t}$ by $\nabla_t$ and but for there being no $F_{A^1 t 1}$ and $F_{A^1 t 2}$ and $\nabla_{A^1 t}\mathfrak{a}^1$ terms. In particular, this analog holds with h = $\chi_1$. Supposing again that $\mathcal{L}\psi = \lambda\psi$, this analog with (6.22) says this:

$$\int_{S^1 \times \Sigma} \chi_1 |\nabla_t \psi|^2 \leq (|\lambda| + c_0 m) \int_{S^1 \times \Sigma} \chi_1 |\psi|^2 \; + \; \int_{S^1 \times \Sigma} \langle \gamma_a(\nabla_a \chi_1)\psi, \gamma_t \nabla_t \psi \rangle$$
$$+ \; c_0 m \int_{S^1 \times \Sigma} r\eth^{3/2}\chi_1 |[\sigma, \psi]|^2 \,.$$

(6.23)

As will be clear momentarily, it proves useful to write the two occurances of $\nabla_t$ in (6.23) as $\nabla_{\hat{A}t} - m[\hat{A}^1, \cdot\,]$. Do this and then use what is said above about $A^1$ to see that (6.23) leads to the inequality below:

$$\tfrac{1}{2}\int_{S^1 \times \Sigma} \chi_1 |\nabla_{A^1 t}\psi|^2 \leq (|\lambda| + c_0 m) \int_{S^1 \times \Sigma} \chi_1 |\psi|^2 \; + \; \int_{S^1 \times \Sigma} \langle \gamma_a(\nabla_a \chi_1)\psi, \gamma_t \nabla_{A^1 t}\psi \rangle$$
$$+ \; c_0 m \int_{S^1 \times \Sigma} r\eth^{3/2}\chi_1 |[\sigma, \psi]|^2 \,.$$

(6.24)

The next and last part of the proof invokes (6.24) and (6.21).

*Part 4*: Add the preceding inequality to the inequality in (6.20). After noting that the respective integrals in (6.21) and (6.24) with the derivatives of $\chi_1$ are identical but appear with opposite in signs, one is led to this:

$$\tfrac{1}{2}\int_{S^1 \times \Sigma} |\nabla_{A^1 t}\psi|^2 \leq (|\lambda| + c_0 m)\int_{S^1 \times \Sigma} |\psi|^2 + c_0 m \int_{S^1 \times \Sigma} r\eth^{3/2}|[\sigma, \psi]|^2 \,.$$

(6.25)

Because $\nabla_{At} = \nabla_{A^1 t} + \mathfrak{b}_t$ and $|\mathfrak{b}_t| \leq c_0 m^2\, r^{-2/3}$ (see (5.34)), the preceding inequality holds as well with a larger $c_0$ when $\nabla_{A^1 t}$ is replaced by $\nabla_{At}$.

To finish the proof of Lemma 6.3: Use Lemma 4.4 (which can be used if $|\lambda| < 1$ per what is said in Part 1) to bound the right most integral on the right hand side of $\nabla_{At}$'s version of (6.25) by $c_0 m$ times the $S^1 \times \Sigma$ integral of $|\psi|^2$.

### c) The Fourier modes of eigenvectors of $\mathcal{L}$ with small normed eigenvalue

By way of background for the central lemma in this subsection (Lemma 6.4): Assume in this subsection that $(A, \mathfrak{a})$ comes from a version of Proposition 5.3 with $m < c_0^{-1}$ and $r > c_0$. The corresponding version of the operator $\mathcal{L}$ can alway be written as $\mathcal{D} + m(\cdot)_\circ$, and then the $\gamma_t$ part of $\mathcal{D}$ can distinguished by writing $\mathcal{D}$ as $\gamma_t \nabla_{At} + \mathfrak{D}$. The operator $\mathcal{D}$ can also be written as (6.17) in terms of what is denoted there as $\mathcal{D}_1$ which has a similar decomposition as $\gamma_t \nabla_{A^1 t} + \mathfrak{D}_1$. This can be done anywhere on $S^1 \times \Sigma$. The upcoming Lemma 6.4 says in part that an eigenvector of $\mathcal{L}$ with small normed eigenvalue is $\mathcal{O}(e^{-r/c_0})$ close to a section of $\mathbb{S}$ on $S^1 \times \Sigma$ that is annihilated by both $\nabla_{At}$ and $\nabla_{A^1 t}$ where $\mathfrak{d} \geq r_0$.

As for the rest of $S^1 \times \Sigma$: As noted previously, a refined decomposition of $\nabla_{At}$ and $\nabla_{A^1 t}$ is possible on $S^1 \times D$ when D denotes the embedded disk in $\Sigma$ centered at a point in $\Theta$ that is identified with as the $|z| < r_0$ part of $\mathbb{C}$ using the local, holomorphic coordinate that is introduced in Lemma 3.2. This decomposition comes about by writing $\nabla_{A^1 t}$ there using the $\mathbb{C} \oplus \mathbb{C}$ product structure for $E_q$ from Section 3b as in (6.22). This writes $\nabla_{A^1 t}$ as $\nabla_t + m[A^1, \cdot]$ with $\nabla_t$ denoting product connection's directional covariant derivative $\nabla_t$. Lemma 6.4 also says that the afore-mentioned $\mathcal{O}(e^{-r/c_0})$ close section of $\mathbb{S}$ is annihilated by $\nabla_t$ on each of the $p \in \Theta$ versions of $S^1 \times \Sigma$.

**Lemma 6.4**: *There exists* $\kappa > 1$ *with the following significance: Suppose in what follows that* $m < \frac{1}{\kappa}$. *Fix* $r > \kappa$ *and let* $(A, \mathfrak{a})$ *denote the pair from Proposition 5.3. Let* $\psi$ *denote an eigenvector of the corresponding version of* $\mathcal{L}$ *with* $\mathbb{L}$-*norm equal to 1 and with eigenvalue (denoted by* $\lambda$*) between* $-\frac{1}{\kappa}$ *and* $\frac{1}{\kappa}$. *There exists a section of* $\mathbb{S}$ *to be denoted by* $\psi_0$ *with* $\mathbb{L}$-*norm equal to 1 that obeys the constraints listed below.*

- $\mathcal{L}\psi_0 = \lambda\psi_0$ *except where* $\mathfrak{d}$ *is between* $\frac{1}{16} r_0$ *and* $\frac{1}{8} r_0$; *and there,* $|\mathcal{L}\psi_0 - \lambda\psi_0| \leq \kappa e^{-r/\kappa}$.
- $[\sigma, \psi_0] = 0$ *where* $\mathfrak{d} > \frac{1}{8} r_0$.
- $\nabla_t\psi_0 = 0$ *on each* $p \in \Theta$ *version of* D; *and* $\nabla_{At}\psi_0 = 0$ *and* $\nabla_{A^1 t}\psi = 0$ *everywhere else*.
- $\|\psi - \psi_0\|_{\mathbb{H}} < \kappa e^{-r/\kappa}$ *and* $|\psi - \psi_0| \leq \kappa e^{-r/\kappa}$.

With regards to $\mathbb{H}$-norm in the fourth bullet: It doesn't matter whether that norm is defined by the pair $(A^1, \mathfrak{a})$ which defines $\mathcal{D}_1$ or by the pair $(A, \mathfrak{a})$ from Proposition 5.3.)

The rest of this subsection is dedicated to the proof of the lemma.

*Proof of Lemma 6.4* The proof has seven parts.

*Part 1*: Fix $p \in \Theta$ and let D again denote the $|z| < r_0$ part of the holomorphic coordinate chart centered a p as described in Lemma 3.2. Let $\psi$ denote a given section of $\mathbb{S}$ on $S^1 \times D$. This section can be written on $S^1 \times D$ using the Fourier decomposition defined by the product connection on the product $\mathbb{C} \oplus \mathbb{C}$ bundle over $S^1 \times D$ after using the identification from Part 2 of Section 3b between $E_q$ on $S^1 \times D$ and that product bundle. This Fourier decomposition of $\psi$ has the form

$$\psi = \sum_{n \in \mathbb{Z}} e^{int} \psi_n$$

(6.26)

with each $\psi_n$ being independent of t.

As noted previously, $\nabla_t$ commutes with $\mathcal{L}$ on $S^1 \times D$. What that means with regards to (6.26) is this: Writing $\mathcal{L}$ on $S^1 \times D$ as $\mathcal{D} + (\cdot)_{\mathfrak{c}}$ and then $\mathcal{D}$ in terms of $\mathcal{D}_1$ as in (6.17) and then $\mathcal{D}_1$ as in (6.22) ends up writing $\mathcal{L}$ as done below,

$$\mathcal{L} = \gamma_t \nabla_t + \mathcal{D}_1 + m \gamma_t [A^1, \cdot] + m(\gamma_k [\mathfrak{e}_k, \cdot] + \rho_k [\mathfrak{f}_k, \cdot]) + m(\cdot)_{\mathfrak{c}}$$

(6.27)

and thus $\mathcal{L}\psi$ as

$$\mathcal{L}\psi = \sum_{n \in \mathbb{Z}} e^{int} \left( in\gamma_t \psi_n + \mathcal{D}_1 \psi_n + m\gamma_t [A^1, \psi_n] + m(\gamma_k [\mathfrak{e}_k, \psi_n] + \rho_k [\mathfrak{f}_k, \psi_n]) + m(\psi_n)_{\mathfrak{c}} \right).$$

(6.28)

It follows as a consequence that if $\mathcal{L}\psi = \lambda\psi$ on $S^1 \times D$, then for each integer n,

$$\mathcal{L}(e^{int}\psi_n) = \lambda \, e^{int}\psi_n$$

(6.29)

With regards to the size of the $\psi_n$'s: Assume that $r > c_0$ henceforth. If $\psi$ is defined on the whole of $S^1 \times \Sigma$ with $\mathbb{L}$ norm equal to 1, and if $\mathcal{D}\psi = \lambda\psi$ with $|\lambda| < c_0^{-1}$, then it follows from Lemma 6.3 that

$$\sum_{n \in \mathbb{Z}} \left( n^2 \int_{S^1 \times D} |\psi_n|^2 \right) \leq c_0 \, m \, .$$

(6.30)

Save this inequality until Part 3 of the proof.

*Part 2*: Assume in this part of the proof that $\psi$ is defined on the whole of $S^1 \times \Sigma$, that $\mathcal{L}\psi = \lambda\psi$, and that $\psi$ has $\mathbb{L}$ norm equal to 1. For each $p \in \Theta$, let $D'$ denote the disk in p's version of D where $|z|$ is less than $\frac{1}{4} r_0$. Write $\psi$ where $\eth > c_0 r^{-2/3}$ as was done in (4.23),

$$\psi = \langle \sigma\psi \rangle \, \sigma + \psi_\perp \, .$$

(6.31)

In particular, write it this way on the complement of all $p \in \Theta$ versions of $S^1 \times D'$. Since $\nabla_{At}\sigma$ and $[\mathfrak{a}, \sigma]$ are zero on $S^1 \times (\Sigma - \cup_{p \in \Theta} D')$, it follows that $\mathcal{L}(\langle \sigma\psi \rangle \sigma) = 0$ on $S^1 \times D'$.

Equivalently: The $\mathcal{I}$-valued section $\langle\sigma\psi\rangle$ of the bundle $(\oplus_2 T^*(S^1 \times \Sigma) \oplus (\oplus_2 \mathbb{R}))$ obeys the equation

$$\gamma_k \nabla_k \langle\sigma\psi\rangle + \langle\sigma\psi_c\rangle = \lambda\langle\sigma\psi\rangle$$

(6.32)

on the domain $S^1 \times (\Sigma - \cup_{p\in\Theta}D')$. This section $\langle\sigma\psi\rangle$ can also be written as a Fourier series on $S^1 \times (\Sigma - \cup_{p\in\Theta}D')$, thus as

$$\langle\sigma\psi\rangle = \sum_{n\in\mathbb{Z}} e^{in} \langle\sigma\psi\rangle_n$$

(6.33)

with each $n \in \mathbb{Z}$ version of $\langle\sigma\psi\rangle_n$ having no $S^1$-factor dependence. Note in particular that each $n \in \mathbb{Z}$ version of $e^{in} \langle\sigma\psi\rangle_n$ also obeys the eigenvalue equation in (6.32); this is because the line bundle $\mathcal{I}$ is the pull-back of a real line bundle on $(\Sigma - \cup_{p\in\Theta}D')$ and because the metric is $S^1$ invariant. With regards to the size of the $\langle\sigma\psi\rangle_n$'s: As before, if $r > c_0$ and if $|\lambda| < c_0^{-1}$, then Lemma 6.3 can be brought to bear to see that

$$\sum_{n\in\mathbb{Z}} \left(n^2 \int_{\Sigma - \cup_{p\in\Theta}D'} |\langle\sigma\psi\rangle_n|^2\right) \leq c_0 \, m \,.$$

(6.34)

A key observation now concerns the respective Fourier coefficients in (6.26) and (6.33) which is this: Fix $p \in \Theta$. Then $\sigma$ on p's version of $D-D'$ when written with respect to the $\mathbb{C} \oplus \mathbb{C}$ product structure defined in Part 3 of Section 3b is $\hat{\tau}$ from (2.1). With respect to the description of $E_q$ on $S^1 \times (\Sigma - \Theta)$ in Section 2b as the quotient of a bundle on $S^1 \times (\Sigma_q - \Theta)$, $\sigma$ is the endomorphism $\hat{\sigma}$ in Part 7 of Section 2b. Because the transition functions between the quotient depiction of $E_q$ on $S^1 \times (\Sigma - \Theta)$ and the product structure depiction near $S^1 \times \Theta$ are $S^1$ invariant (they are defined by the sections $s_1$ and $s_2$ in (3.15)), it follows that $\langle\hat{\tau}\psi_n\rangle$ which is defined on $S^1 \times (D-p)$ using the $\mathbb{C} \oplus \mathbb{C}$ product structure for $E_q$ is, on $S^1 \times (D-D')$, the same on this domain as the $S^1 \times (\Sigma - \cup_{p\in\Theta}D')$ version of $\langle\sigma\psi\rangle_n$ (which is $\langle\hat{\sigma}\psi\rangle_n$ when written using Section 2b's notation).

*Part 3*: To exploit the observations in Part 2, let $\chi_\Delta$ denote the function on $\Sigma$ given by the rule whereby $\chi_\Delta(\cdot) = \chi(\frac{16\,\mathfrak{d}}{r_0} - 1)$ which is equal to 1 where $\mathfrak{d} < \frac{1}{16} r_0$ and equal to 0 where $\mathfrak{d}$ is greater than twice that (thus $\frac{1}{8} r_0$). Supposing that $\psi$ is an eigenvector of $\mathcal{L}$ with $\mathbb{L}$ norm equal to 1, define a new section of $\mathbb{S}$ to be denoted by $\psi_\Delta$ by setting $\psi_\Delta = \psi$ where $\mathfrak{d} < \frac{1}{16} r_0$ and setting

$$\psi_\Delta \equiv \langle\sigma\psi\rangle\sigma + \chi_\Delta \psi_\perp$$

(6.35)

elsewhere. The important observations in this part of the proof concerns the size of $\psi$ - $\psi_\Delta$ when $r$ is large, and the lemma that follows momentarily says what is needed in this regard. The lemma's notation uses $w_\Delta$ to denote the characteristic function for the $\eth \in [\frac{1}{16} r_0, \frac{1}{8} r_0]$ part of $S^1 \times \Sigma$ (which is where $\chi_\Delta$ is not constant).

**Lemma 6.5**: *There exists* $\kappa > 1$ *with the following significance: Fix* $m < \frac{1}{\kappa}$ *and* $r = \frac{n}{m} > \kappa$; *and let* $(A, \mathfrak{a})$ *denote the pair of connection on* $P$ *and 1-form with values in* $\mathrm{ad}(P)$ *from Proposition 5.3. Let* $\psi$ *denote either an eigenvector of* $\mathcal{L}$ *with* $\mathbb{L}$-*norm equal to 1 and with eigenvalue (denoted by* $\lambda$*) between* $-\frac{1}{\kappa}$ *and* $\frac{1}{\kappa}$. *Then*

- $\|\psi - \psi_\Delta\|_{\mathbb{H}} \le e^{-r/\kappa}$ *and* $|\psi - \psi_\Delta| \le \kappa\, e^{-r/\kappa}$
- $|\mathcal{L}\psi_\Delta - \lambda\psi_\Delta| \le \kappa\, e^{-r/\kappa}\, w_\Delta$

The proof of this lemma is in Part 4 which follows momentarily.

   *Part 4*: This part of the proof is dedicated to proving Lemma 6.5.

***Proof of Lemma 6.5***: Let $\chi_\triangledown$ to denote the function bump function $\chi(\frac{32\,\eth}{r_0} - 1)$. Define the section $\psi_\triangledown$ of $\mathbb{S}$ to be $\langle \sigma\psi \rangle \sigma + \chi_2\psi^\perp$, this being (6.35) but for $\chi_2$ replacing $\chi_\Delta$. The first key observations for the proof are the following ones:

- $\psi - \psi_\triangledown = -\frac{1}{4}(1 - \chi_\triangledown)[\sigma, [\sigma, \psi]]$ .
- $\psi - \psi_\Delta = -\frac{1}{4}(1 - \chi_\Delta)[\sigma, [\sigma, \psi]]$
- $\psi - \psi_\Delta = (1 - \chi_\Delta)(\psi - \psi_\triangledown)$ .

$$(6.36)$$

The third bullet identity implies that a $c_0 e^{-r/c_0}$ bound for the $\mathbb{L}$ or $\mathbb{H}$ or pointwise norm of $\psi - \psi_\triangledown$ leads directly to a $c_0 e^{-r/c_0}$ bound for the corresponding $\mathbb{L}$ or $\mathbb{H}$ or pointwise norm of $\psi - \psi_\Delta$. In addition, the first and second bullets lead from a $c_0 e^{-r/c_0}$ bound for the pointwise norm of $\psi - \psi_\triangledown$ to a $c_0 e^{-r/c_0}$ bound for that of $|\mathcal{L}(\psi - \psi_\Delta) - \lambda(\psi - \psi_\Delta)|$ which is $|\mathcal{L}\psi_\Delta - \psi_\Delta|$. Indeed, because $\sigma$ is A-covariantly constant and commutes with $\mathfrak{a}$ where $\eth > \frac{1}{64} r_0$, the second bullet implies that $\mathcal{L}(\psi - \psi_\Delta) - \lambda(\psi - \psi_\Delta)$ is $\frac{1}{4}(\nabla_k\chi_\Delta)\gamma_k[\sigma, [\sigma, \psi]]$ which is $(\nabla_k\chi_\Delta)\gamma_k(\psi - \psi_\triangledown)$ by virtue of the identity in the first bullet and the fact that $1 - \chi_\triangledown$ is equal to 1 where $\chi_\Delta$ is less than 1.

   With the preceding understood, it is sufficient for the proof of Lemma 6.5 to derive $c_0 e^{-r/c_0}$ bounds for the $\mathbb{H}$ and pointwise norms of $\psi - \psi_\triangledown$. To this end, note first that Lemma 4.4 is in play for the section $\psi$ (see Part 1 of the proof of Lemma 6.3). Granted Lemma 4.4's conclusions, then the top bullet identity in (6.36) leads directly to this bound:

$$\|\psi - \psi_{\bar{v}}\|_{\mathbb{L}} \le c_0 e^{-r/c_0} \ .$$

(6.37)

Meanwhile: The following identity also follows from the top bullet of (6.36) because $\sigma$ is A-covariantly constant and commutes with $\mathfrak{a}$ where $\mathfrak{d} > \frac{1}{64} r_0$:

$$\mathcal{L}(\psi - \psi_{\bar{v}}) = \lambda(\psi - \psi_{\bar{v}}) + \tfrac{1}{4} (\nabla_k \chi_v) \gamma_k [\sigma, [\sigma, \psi]].$$

(6.38)

Given that (6.37) holds, then (3.38) implies this:

$$\|\mathcal{L}(\psi - \psi_{\bar{v}})\|_{\mathbb{L}} \le c_0 e^{-r/c_0} \ .$$

(6.39)

Use this bound with (6.37) and the depiction of $\mathcal{L}$ as $\mathcal{D}_1 + m \, (\gamma_k[\mathfrak{e}_k, \, \cdot \,] + \rho_k[\mathfrak{f}_k, \, \cdot \,]) + m(\cdot)_\mathfrak{c}$ to see that $\|\mathcal{D}_1(\psi - \psi_{\bar{v}})\|_{\mathbb{L}}$ is also bounded by $c_0 e^{-r/c_0}$. Granted the latter bound, then Lemma 4.3 with $(\psi - \psi_{\bar{v}})$ used in lieu of $\psi$ and with $\mathcal{D}_1$ used for $\mathcal{D}$ to obtained the desired $c_0 e^{-r/c_0}$ bound for $\|\psi - \psi_{\bar{v}}\|_{\mathbb{H}}$.

With regards to the pointwise bound for $|\psi - \psi_{\bar{v}}|$: To obtain the desired bound, fix any given point in $S^1 \times \Sigma$ (to be denoted by p) and then fix some small radius ball centered at p that is well inside a Gaussian coordinate chart domain centered at p. Denote this ball by B. The radius of B is assumed to be independent of $r$ (an upper bound for the radius is determined below by the norm of the Riemannian curvature). Fix an isometric isomorphism between $E_q|_p$ with $\mathbb{C}^2$ bundle and then use parallel transport by the connection $A^1$ on the radial geodesic rays from p to obtain an isometric isomorphism between $E_q|_B$ and the product bundle $B \times \mathbb{C}^2$. Use $\theta_B$ to denote the product connection on this bundle. The connection $A^1$ appears now on the product bundle as $A^1 = \theta_B + a^1$ with $a^1$ being a Lie algebra valued 1-form on B. Since the curvature of $A^1$ is bounded by $c_0 r^{4/3}$, the norm of this Lie algebra valued 1-form $a^1$ is bounded by $c_0 r^{4/3}$. As a consequence of Proposition 5.3, the connection A appears with respect to the product structure for $E_q|_B$ as $A = \theta_B + a$ with $|a|$ being at most $c_0 r^{4/3}$ also. Since $|\mathfrak{a}| \le c_0 r$, the operator $\mathcal{L}$ on B when written using the product structure for $E_q|_B$ as

$$\mathcal{L} = \gamma_k \nabla_k + \mathfrak{l}_B(\cdot)$$

(6.40)

with the notation as follows: What is denoted by $\nabla$ signifies the covariant derivative on the product $\mathbb{C}^2$ bundle that is defined by the product connection $\theta_B$. Meanwhile, what is denoted by $\mathfrak{l}_B$ signifies an endomorphism with norm bounded by $c_0 r$.

Now let $\chi_B$ denote a bump function constructed from $\chi$ which is 0 outside of B and equal to 1 on the concentric ball inside B with half of B's radius. Let $\eta = \chi_B(\psi - \psi_\nabla)$, and use (6.40) with (6.38) to see that

$$\gamma_k \nabla_k \eta = -\mathfrak{l}_B(\eta) + \lambda(\eta) + \tfrac{1}{4} \chi_B (\nabla_k \chi_\nabla) \gamma_k [\sigma, [\sigma, \psi]] \,.$$

(6.41)

Since $(\gamma_k \nabla_k)^2 = \nabla^\dagger \nabla + \mathfrak{R}$ with $\mathfrak{R}$ denoting an endomorphism defined by the Riemannian curvature, the operator $(\gamma_k \nabla_k)^2$ has a Dirichlet Green's function in B with pole at p if B's radius is less than $c_0^{-1}$ which will be assumed henceforth. Let $G_p$ denote that Greens function and let $\mathfrak{G}_p$ denote $\gamma_k \nabla_k G_p$. Noting that

$$|G_p| \leq c_0 \, \tfrac{1}{\mathrm{dist}(\cdot,p)} \ \textit{and that} \ |\mathfrak{G}_p| \leq c_0 \, \tfrac{1}{\mathrm{dist}(\cdot,p)^2} \,,$$

(6.42)

it follows from (6.41)

$$|\eta|_p| \leq c_0 \int_B \ \tfrac{1}{\mathrm{dist}(\cdot,p)^2} \ (r \, |\eta| + |\nabla \chi_\nabla| |[\sigma, \psi]|) \,.$$

(6.43)

Since $\tfrac{5}{2}$ is less than 3, this last inequality implies that

$$|\eta|_p| \leq c_0 \, r \, ( \int_{S^1 \times \Sigma} \, |\eta|^5)^{1/5} + c_0 (\int_{S^1 \times \Sigma} \, |\nabla \chi_\nabla| |[\sigma, \psi]|^5)^{1/5} \,.$$

(6.44)

To end the proof: The left most term on the right hand side of (6.44) is no greater than $c_0 r \, \|\eta\|_\mathbb{H}$ (see (6.10)) and thus no greater than $c_0 e^{-r/c_0}$. Meanwhile, the right most term on the right hand side of (6.44) is no greater than

$$( \int_{S^1 \times \Sigma} \, |\psi|^{45/8})^{4/25} \, (\int_{S^1 \times \Sigma} \, |\nabla \chi_\nabla| \, |[\sigma, \psi]|^2)^{1/20}$$

(6.45)

which is less than $c_0 e^{-r/c_0}$ also (use Lemmas 4.3 and 4.4 and (6.10) to see this).

*Part 5*: What follows in this part of Lemma 6.4's proof are observations regarding $\psi_\Delta$ from Lemma 6.5 and the Fourier decompositions from Parts 1 and 2: Given that the Fourier expansions in (6.26) and (6.33) are compatible where their domains overlap, and given that $\chi_\ddagger$ is independent of the coordinate on the $S^1$ factor, it follows that the section $\psi_\ddagger$ has the Fourier series depicted below on the whole of $S^1 \times \Sigma$:

$$\psi_\Delta = \sum_{n \in \mathbb{Z}} \, e^{\mathrm{int}}(\chi_\Delta \psi_n + (1 - \chi_\Delta)\langle \sigma\psi \rangle_n \sigma) \,.$$

(6.46)

Subsequent paragraphs will use $\psi_{\Delta n}$ to denote $\chi_\Delta \psi_n + (1 - \chi_\Delta)\langle\sigma\psi\rangle_n\sigma$ which is the coefficient function on $\Sigma$ for the $e^{int}$ Fourier mode in (6.46).

Because the functions $\{e^{int}\}_{n\in\mathbb{Z}}$ are pairwise orthogonal with respect to the standard $L^2$ Hermitian inner product on $C^\infty(S^1; \mathbb{C})$, the inequality in the second bullet of Lemma 6.5 can hold only in the event that

$$\sum_{n\in\mathbb{Z}} |\mathcal{L}(e^{int}\psi_{\Delta n}) - \lambda(e^{int}\psi_{\Delta n})|^2 \leq c_0 \, e^{-r/c_0} \, w_\Delta \ .$$

(6.47)

(Remember that $w_\Delta$ indicates the characteristic function for the part of $\Sigma$ where $\mathfrak{d}$ is between $[\frac{1}{16} r_0, \frac{1}{8}r_0]$.)

To exploit (6.47): Assume now that $m < c_0^{-1}$ and that $|\lambda| < c_0^{-1}$ and that $r > c_0$. It then follows from (6.30) and (6.34) and the $\mathbb{H}$-norm bound in the the third bullet of Lemma 6.5 that the $\mathbb{L}$ norm of $\psi_{\Delta 0}$ is greater than half that of $\psi$ (which is 1). Hence, from (6.47):

$$|\mathcal{L}\psi_{\Delta 0} - \lambda\psi_{\Delta 0}| \leq c_0 e^{-r/c_0} \, \|\psi_{\Delta 0}\|_{\mathbb{L}} \, w_\Delta \ .$$

(6.48)

This is saying that the $S^1$ independent part of $\psi_\Delta$ is almost an eigenvector of $\mathcal{L}$ with eigenvalue $\lambda$ and that its failure to be one is supported only where $\mathfrak{d} \in [\frac{1}{8} r_0, \frac{1}{4} r_0]$.

With regards to Lemma 6.4's first three bullets: By virtue of (6.48), these first three bullets are satisfied by taking $\psi_0$ in Lemma 6.4 to be $\psi_{\Delta 0}$.

*Part 6*: To say something useful about the non-zero Fourier components of $\psi_\Delta$, let $\psi_\#$ denote $\psi_\Delta - \psi_{\Delta 0}$, which is what is obtained from the right hand side of (6.46) by deleting the $n = 0$ term in the sum. As explained Part 7, the $\mathbb{L}$ norm of $\psi_\#$ must be very small when $m + |\lambda| < c_0^{-1}$ and $r > c_0$. In particular, Part 7 proves that there exists $c_\dagger > 1$ such that if those bounds on $m$, $|\lambda|$ and $r$ hold, then

$$\|\psi_\#\|_{\mathbb{L}} \leq c_\dagger \, e^{-r/c_\dagger} \ .$$

(6.49)

Assume for now that (6.49) holds for some $r$ and $m$ and $\lambda$ independent version of $c_\dagger$ when $m$ and $|\lambda|$ are small and $r$ is large. What follows momentarily is a proof that the $\mathbb{H}$-norm and the pointwise norm of $\psi_\#$ obeys the same sort of bound when $m + |\lambda| < c_0^{-1}$ and $r > c_0$:

$$\|\psi_\#\|_{\mathbb{H}} \leq c_0 e^{-r/c_\dagger c_0} \quad and \quad |\psi_\#| < c_0 e^{-r/c_\dagger c_0} \ .$$

(6.50)

Given the bounds in (6.50), then all of Lemma 6.4's bullets are satified with that Lemma's section $\psi_0$ being $\psi_{\Delta 0}$.

To see how (6.50) follows from (6.59), start with the fact that (6.47) implies that

$$|\mathcal{L}\psi_\# - \lambda\psi_\#| \leq c_0 \; e^{-r/c_0} \; w_\Delta \;.$$

(6.51)

Write $\mathcal{L}$ as $\mathcal{D}_1 + m(\gamma_k[\mathfrak{e}_k, \cdot] + \rho_k[\mathfrak{f}_k, \cdot]) + m(\cdot)_c$. Then (6.51) with (6.49) the bounds on $|\mathfrak{e}|$ and $|\mathfrak{f}|$ from Proposition 5.3 lead to a $c_0 \; e^{-r/c_\dagger c_0}$ bound for $\|\mathcal{D}_1\psi_\#\|_{\mathbb{L}}$. Lemma 4.3 leads from that bound to the desired $c_0 e^{-r/c_\dagger c_0}$ bound for $\|\psi_\#\|_{\mathbb{H}}$.

As for the pointwise bound on $|\psi_\#|$: What with (6.51), the derivation of this bound is a virtual replay of the arguments to obtain the $c_0 \; e^{-r/c_0}$ bound for $|\psi - \psi_\nabla|$ using the depiction in (6.40) of the operator $\mathcal{L}$ and then using the Green's function in (6.42). The only salient difference is that (6.43) is replaced by this:

$$|\psi_\#|_p| \leq c_0 \int_B \; \frac{1}{\text{dist}(\cdot,p)^2} \; (r \; |\psi_\#| + e^{-r/c_0} \; w_\Delta)$$

(6.52)

whose righthand side is no greater than $c_0 r \left( \int_{S^1 \times \Sigma} |\psi_\#|^6 \right)^{1/6} + c_0 e^{-r/c_0}$. The $c_0 e^{-r/c_\dagger c_0}$ bound for $|\psi_\#|$ follows from that since $c_0 r \left( \int_{S^1 \times \Sigma} |\psi_\#|^6 \right)^{1/6}$ is at most $c_0 r \|\psi_\#\|_{\mathbb{H}}$ (see (6.10)).

*Part 7*: This part of the proof derives the bound that is asserted in (6.49). To start the derivation: Let $c_1$ denote the version of $c_0$ that appears on the right hand side of (6.37). Define the positive number $c_2$ by writing $\|\psi_\#\|_{\mathbb{L}}^2$ as $e^{-r/c_2}$. Granted these definitions, then (6.47) implies that

$$|\mathcal{L}\psi_\# - \lambda\psi_\#|^2 \leq c_1 \; e^{-r/c_1} \; e^{r/c_2} \|\psi_\#\|_{\mathbb{L}}^2 \; w_\Delta \;.$$

(6.53)

If $c_2 > 2c_1$, then the preceding inequality implies in turn that

$$|\mathcal{L}\psi_\# - \lambda\psi_\#|^2 \leq c_1 \; e^{-r/2c_1} \|\psi_\#\|_{\mathbb{L}}^2 \; w_\Delta \;,$$

(6.54)

which is saying that $\psi_\#$ is also nearly an eigenvector of $\mathcal{L}$ with eigenvalue $\lambda$ and that its failure to be an eigenvector is supported only where $\mathfrak{d} \in [\frac{1}{16} r_0, \frac{1}{8} r_0]$.

To prove (6.49), it is sufficient to show that (6.54) leads to nonsense if $c_2 > 2c_1 + c_0$ when $m + |\lambda| < c_0^{-1}$ and $r > c_0$. And, the desired nonsense is the contradiction between the following two statements: The first statement is that the $S^1 \times \Sigma$ integral of $|\nabla_{At}\psi_\#|^2$ is at most $c_0 m$ times the $S^1 \times \Sigma$ integral of $|\psi_\#|^2$. The second statement is that the $S^1 \times \Sigma$ integral of $|\nabla_{At}\psi_\#|^2$ is greater than $c_0^{-1}$ times the $S^1 \times \Sigma$ integral of $|\psi_\#|^2$. The next paragraphs justify these two statements about $|\nabla_{At}\psi_\#|^2$.

The starting point for the justifications is the observation that (6.54) allows for an appeal to Lemma 4.4 with $\psi_\#$ used there for $\psi$ and with $\mathcal{D}_1$ used there for $\mathcal{D}$. To prove that Lemma 4.4 is in play, write $\mathcal{L}$ again as $\mathcal{D}_1 + m(\gamma_k[\mathfrak{e}_k, \cdot\,] + \rho_k[\mathfrak{f}_k, \cdot\,]) + m(\cdot)_\mathfrak{c}$ to see from (6.54) and the bounds on $|\mathfrak{e}|$ and $|\mathfrak{f}|$ that

$$|[\sigma, \mathcal{D}_1\psi_\#]| \le c_0(|\lambda| + m)|[\sigma, \psi_\#]| + c_0 m \, |\psi| e^{-r\mathfrak{d}^{3/2}/c_1} + c_0 e^{-r/2c_1}\|\psi_\#\|_{\mathbb{L}} \, w_\Delta \tag{6.55}$$

which is a version of (4.22), and that is the prequisite for Lemma 4.4 with it understood that $m$ and $|\lambda|$ are less than $c_0$ and $r$ is greater than $c_0$.

With regards to the $c_0 m \, \|\psi_\#\|_{\mathbb{L}}^2$ upper bound for the $S^1 \times \Sigma$ integral of $|\nabla_{At}\psi_\#|^2$: With (6.44) as the starting point and with Lemma 4.4 in play, the argument for this differs only cosmetically from the argument that is used to prove Lemma 6.3. Since the modifications of Lemma 6.3's proof are straightforward, no more will be said about the $c_0 m$ lower bound.

With regards to the $c_0^{-1} \|\psi_\#\|_{\mathbb{L}}^2$ lower bound for the $S^1 \times \Sigma$ integral of $|\nabla_{At}\psi_\#|^2$: Because the $\psi_\# = \langle\sigma\psi_\#\rangle\sigma$ and $\nabla_{At}\psi_\# = \nabla_t\langle\sigma\psi_\#\rangle\sigma$ on the $\mathfrak{d} > \frac{1}{4} r_0$ part of $S^1 \times \Sigma$, the contribution from this part of $S^1 \times \Sigma$ to the integral of $|\nabla_{At}\psi_\#|^2$ is the sum of the non-zero Fourier terms in the expression on the left hand side of (6.34). That in turn is no smaller than the integral of $|\psi_\#|^2$ over this part of $S^1 \times \Sigma$.

Meanwhile: Supposing that $p$ is from $\Theta$, let $D$ denote $|z| < r_0$ part of the holomorphic coordinate chart $\Sigma$ centered at $p$ as described in Lemma 3.2. The $S^1 \times D$ integral of $|\nabla_{At}\psi_\#|^2$ can be bounded from below by using the $\mathbb{C} \oplus \mathbb{C}$ product structure for $E_q$ on $S^1 \times D$ from Part 2 of Section 3b to first write $\nabla_{At}$ using the product connection as $\nabla_t + m\,[A^1 + \mathfrak{e}_t, \cdot\,]$. Then, use this depiction of $\nabla_{At}$ to see that

$$\int_{S^1 \times D} |\nabla_{At}\psi_\#|^2 \ge c_0^{-1}\int_{S^1 \times D} |\nabla_t\psi_\#|^2 - c_0 m \int_{S^1 \times D} |\psi_\#|^2 - c_0 m \int_{S^1 \times D} r\mathfrak{d}^{3/2}|[\hat{\tau}, \psi_\#]|^2 \,. \tag{6.56}$$

The left most integral on the right hand side of (6.56) is the contribution to the left hand side of (6.30) from the non-zero Fourier modes. It is therefore larger than $c_0^{-1}$ times the $S^1 \times D$ integral of $|\psi_\#|^2$. Meanwhile, the right most integral on the right hand side of (6.56) is no greater than $c_0 m$ times that same integral of $|\psi_\#|^2$, this by appeal to that $\mathcal{D} = \mathcal{D}_1$ and $\psi = \psi_\#$ version of Lemma 4.4. Thus, supposing that $m$ is less than $c_0^{-1}$, the two negative terms on the right hand side of (6.56) will be less than half of the left most term on that right hand side (which is positive). With this $c_0^{-1}$ bound on $m$ enforced, then the $S^1 \times D$ integral of $|\nabla_{At}\psi_\#|^2$ is greater than $c_0^{-1}$ times the $S^1 \times D$ integral of $|\psi_\#|^2$.

The $c_0^{-1}\|\psi_\#\|_{\mathbb{L}}^2$ lower bound for the $S^1 \times \Sigma$ integral of $|\nabla_{At}\psi_\#|^2$ follows directly from the conclusions of the preceding two paragraphs.

### d) The reappearance of the operator $\mathfrak{L}$

This subsection uses the section $\psi_0$ from Lemma 6.4 to see that the eigenvalue $\lambda$ in Lemma 6.4 is very nearly an eigenvalue of the operator $\mathfrak{L}$ that is depicted in (5.56). To see how this comes about, recall first that $\mathfrak{L}$ acts on sections over $\Sigma$ of the tensor product of the bundle $(\oplus_2(T^*(S^1 \times \Sigma))|_\Sigma \oplus_2 \mathbb{R})$ with the bundle of skew-Hermitian sections of $L \oplus L^{-1}$. Recall also that the respective $\mathbb{C} \oplus \mathbb{C}$ and $\mathbb{C}$ product structures for $E_q$ and $L$ over any $p \in \Theta$ version of $S^1 \times D$ from Part 2 of Section 3b identifies $E_q$ there with $L \oplus L^{-1}$. With this identification understood, then Lemma 6.4's section $\psi_0$ defines an element in the domain of $\mathfrak{L}$ over each $p \in \Theta$ version of $D$ because it is annihilated by directional covariant derivative along the $S^1$ factor of $S^1 \times D$ for the product connection from the $\mathbb{C} \oplus \mathbb{C}$ product structure for $E_q$ on $S^1 \times D$. The section $\psi_0$ is then extended over the rest of $\Sigma$ as an element in the domain of $\mathfrak{L}$ by defining it on $\Sigma - \cup_{p \in \Theta} D$ to be $\langle \sigma \psi_0 \rangle \hat{\tau}$ with $\hat{\tau}$ depicted in (2.1).

Let $\mathbb{S}$ now denote the tensor product of $(\oplus_2(T^*(S^1 \times \Sigma))|_\Sigma \oplus_2 \mathbb{R})$ with the bundle of skew-Hermitian sections of $L \oplus L^{-1}$. With $\psi_0$ identified now as a (not everywhere zero) section of $\mathbb{S}$, it then follows from Lemma 6.4 that this incarnation of $\psi_0$ obeys the inequality

$$|\mathfrak{L}\psi_0 + \mathfrak{l}(\psi_0) - \lambda\psi_0| \leq c_0 e^{-r/c_0} \|\psi_0\|_\mathbb{L} \,,$$

(6.57)

where $\mathfrak{l}(\cdot)$ signifies the symmetric endomorphism of $\mathbb{S}$ with compact support on the $\mathfrak{d} < r_0$ part of $\Sigma$ that is defined below in (6.58) using what are denoted by $\mathbb{A}$, $\mathbb{B}$, $\mathbb{C}$, $\mathbb{D}$ and $\mathbb{B}^1$ and $\mathbb{D}^1$ in (5.33). With regards to the depiction of $\mathfrak{l}(\cdot)$ below: This depiction and (5.33) are written using Section 3's $\mathbb{C} \oplus \mathbb{C}$ product structure for $L \oplus L^{-1}$ on the $\mathfrak{d} < r_0$ part of $\Sigma$. Here is the definition of $\mathfrak{l}(\cdot)$:

$$\mathfrak{l}(\cdot) \equiv \gamma_t[\chi_{\ddagger}\mathbb{A}, \cdot\,] + (\gamma_a(\mathbb{B}^1 + \mathbb{B})_a[\sigma_3, \cdot\,] + \rho_t\mathbb{C}\,[\sigma_3, \cdot\,] + \rho_a[(\mathbb{D}^1 + \chi_{\ddagger}\mathbb{D})_a, \cdot\,] \,.$$

(6.58)

The notation here and subsequently uses the subscript 'a' to label the components of $\mathbb{B}^1$, $\mathbb{B}$, $\mathbb{D}^1$ and $\mathbb{D}$ when they are written on a disk in $\Sigma$ with respect to an oriented orthonormal frame for $T^*\Sigma$; and, when that subsubscript is repeated in a product with the $\gamma$'s (as occurs twice (6.58)), that index is implicitly summed in each product over the set of its values which is the set $\{1, 2\}$. Also by way of notation (a reminder): What is denoted by $\chi_{\ddagger}$ is the function $\chi(\frac{8\mathfrak{d}}{r_0} - 1)$; it is equal to 1 where $\mathfrak{d} < \frac{1}{8} r_0$ and equal to zero where $\mathfrak{d} > \frac{1}{4} r_0$.)

With regards to the $\mathbb{L}$-norm that appears in (6.57): There and henceforth in this section, what is denoted by $\|\psi\|_\mathbb{L}$ signifies the square root of the $\Sigma$-integral of $|\psi|^2$. Likewise, what is denoted subsequently by $\|\psi\|_\mathbb{H}$ signifies the square root of the $\Sigma$-integral of the function $|\nabla_{\hat{A}}\psi|^2 + |[\hat{a}, \psi]|^2 + |\psi|^2$ with $(\hat{A}, \hat{a})$ being the same pair that is used in (5.33)

The next three lemmas summarizes some implications of (6.57) with regards to the spectrum of the operators $\mathcal{L}$, $\mathfrak{L} + \mathfrak{l}(\cdot)$ and $\mathfrak{L}$.

**Lemma 6.6:** *There exists $\kappa > 1$ with the following significance: Fix $\mathfrak{m} < \frac{1}{\kappa}$ and $r = \frac{\mathfrak{n}}{\mathfrak{m}} > \kappa$ and let $(\mathbb{A}, \mathfrak{a})$ denote the corresponding pair from Proposition 5.3. Use this pair to define the operator $\mathfrak{L} + \mathfrak{l}(\cdot)$. Suppose that $\lambda$ is a real number and $\psi_0$ is a not identically zero* *section of* $\mathbb{S}$. *If (6.57) is obeyed by* $\lambda$ *and* $\psi_0$, *then* $\mathcal{L} + \mathfrak{l}(\cdot)$ *has an eigenvalue that differs from* $\lambda$ *by at most* $\kappa e^{-r/\kappa}$.

The second lemma in this series compares the respective spectra of $\mathcal{L} + \mathfrak{l}(\cdot)$ and $\mathcal{L}$.

**Lemma 6.7**: *There exists* $\kappa > 1$ *with the following significance: Fix* $m < \frac{1}{\kappa}$ *and* $r = \frac{n}{m} > \kappa$ *and let* $(A, \mathfrak{a})$ *denote the corresponding pair from Proposition 5.3. Use this pair to define the operator* $\mathcal{L} + \mathfrak{l}(\cdot)$. *A number* $\lambda$ *is an eigenvalue of* $\mathcal{L} + \mathfrak{l}(\cdot)$ *only if it differs by at most* $\kappa m^2 r^{-2/3}$ *from an eigenvalue of* $\mathcal{L}$. *Conversely, a number* $\lambda$ *is an eigenvalue of* $\mathcal{L}$ *only if it differs by at most* $\kappa m^2 r^{-2/3}$ *from an eigenvalue of* $\mathcal{L} + \mathfrak{l}(\cdot)$.

The final lemma in this series circles back to the operator $\mathcal{L}$.

**Lemma 6.8**: *There exists* $\kappa > 1$ *with the following significance: Fix* $m < \frac{1}{\kappa}$ *and* $r = \frac{n}{m} > \kappa$. *Define the operator* $\mathcal{L}$ *as done above using using the pair* $(A, \mathfrak{a})$ *from the* $m$ *and* $r$ *version of Proposition 5.3. There are no eigenvalues of* $\mathcal{L}$ *between* $\frac{1}{\kappa}$ *and* $\kappa m$. *Meanwhile, the span of the space of eigenvectors of* $\mathcal{L}$ *with eigenvalue norm less than* $\kappa m$ *has dimension* 12g - 12.

The rest of this subsection contains the proofs of these lemmas. The next subsection says more about the 12g - 12 dimensional space of eigenvectors of $\mathcal{L}$ with eigenvalue norm below $c_0 m$.

***Proof of Lemma 6.6***: Let $c_1$ denote the value of the version of $c_0$ that appears in (5.57). Integrate the square of both sides of (5.57) over $\Sigma$ to see that

$$\int_\Sigma |\mathcal{L}\psi_0 + \mathfrak{l}(\psi_0) - \lambda\psi_0|^2 \le c_1^2 \, e^{-2r/c_1} \|\psi_0\|_{\mathbb{L}}^2$$

(6.59)

To make something of this, suppose first (for the sake of argument) that $\lambda$ is not an eigenvalue of $\mathcal{L} + \mathfrak{l}(\cdot)$. Letting $\varepsilon$ denote the distance in $\mathbb{R}$ between $\lambda$ and the spectrum of $\mathcal{L} + \mathfrak{l}(\cdot)$, then (6.59) implies that

$$\varepsilon^2 \|\psi_0\|_{\mathbb{L}}^2 \le c_1^2 \, e^{-2r/c_1} \|\psi_0\|_{\mathbb{L}}^2 \, .$$

(6.60)

It is perhaps needless to say that this inequality can hold only in the event that $\varepsilon \le c_1 e^{-r/c_1}$. (To prove that (6.60) holds: Write $\psi_0$ in (6.59) as a linear combination of $\mathbb{L}$-orthonormal eigenvectors for the operator $\mathcal{L} + \mathfrak{l}(\cdot)$.)

**Proof of Lemma 6.7**: Supposing first that $\lambda$ is an eigenvalue of $\mathfrak{L} + \mathfrak{l}(\cdot)$, let $\varepsilon$ now denote the distance between $\lambda$ and the nearest eigenvalue of $\mathfrak{L}$. Since $|\mathfrak{l}(\psi)| \leq c_0 m^2 r^{-2/3}$, the strategy used above for the proof of Lemma 6.6 can be repeated with the result being the that (6.60) is replaced by $\varepsilon^2 \|\psi\|_{\mathbb{L}}^2 \leq c_0 m^4 r^{-4/3} \|\psi\|_{\mathbb{L}}^2$ which implies that $\varepsilon < c_0 m^2 r^{-2/3}$. Reversing the roles of $\mathfrak{L}$ and $\mathfrak{L} + \mathfrak{l}(\cdot)$ with the same argument proves that every eigenvalue of $\mathfrak{L}$ has distance at most $c_0 m^2 r^{-2/3}$ from an eigenvalue of $\mathfrak{L} + \mathfrak{l}(\cdot)$.

**Proof of Lemma 6.8**: The proof of Lemma 5.7 in Part 4 of Section 5d explained why $\mathfrak{L}$ has no eigenvalues between $c_0 m$ and $c_0^{-1}$ if $m$ is small ($< c_0^{-1}$) and $r$ is large ($> c_0$). Therefore, by virtue of Lemma 6.7, the operator $\mathfrak{L} + \mathfrak{l}(\cdot)$ has no eigenvalues between $c_0 m$ and $c_0^{-1}$ if $m$ is small ($< c_0^{-1}$) and $r$ is large ($> c_0$). Therefore, by virtue of Lemma 6.6, the same statement holds for $\mathcal{L}$.

With regards to the dimension span of the set of eigenvectors of $\mathcal{L}$ with eigenvalue norm less than $c_0 m$: To see that this span is at least 12g-12 dimensional, note first that the span of the eigenvectors of $\mathfrak{L}$ with eigenvalue norm less than $c_0 m$ is 12g - 12 dimensional (see the proof of Lemma 5.7). With this understood, let $\psi$ denote for the moment a linear combination of orthonormal eigenvectors of $\mathfrak{L}$ from this span. This section $\psi$ determines a section in the domain of $\mathcal{L}$ (denoted by $\eta_\psi$) as follows: Use the respective $\mathbb{C} \oplus \mathbb{C}$ and $\mathbb{C}$ product structures for $E_q$ and L where $\mathfrak{d} < \frac{1}{8} r_0$ (from Part 2 of Section 3b) to identify $E_q$ there with $L \oplus L^{-1}$. With this identification understood, set $\eta_\psi = \psi$ where $\mathfrak{d} < \frac{1}{8} r_0$ and set $\eta_\psi = \chi_{\divideontimes} \psi + \langle \hat{\tau} \psi \rangle \sigma$ where $\mathfrak{d} > \frac{1}{8} r_0$. Lemma 4.4 can be brought to bear to see that the map $\psi \to \eta_\psi$ is almost isometric in the sense that if $\psi$ and $\psi'$ are any two elements in the span of the eigenvectors of $\mathfrak{L}$ with eigenvalue norm less than $c_0^{-1} m$, then $\|\eta_\psi - \eta_{\psi'}\|_{\mathbb{L}}$ differs from $\|\psi - \psi'\|_{\mathbb{L}}$ by at most $c_0 e^{-r/c_0} \|\psi - \psi'\|_{\mathbb{L}}$. Meanwhile, $\|\mathcal{L}\eta_\psi\|_{\mathbb{L}} \leq c_0 m \|\psi\|_{\mathbb{L}}$ which is less than $c_0 m \|\eta_\psi\|_{\mathbb{L}}$. Thus, the image of $\eta_{(\cdot)}$ is a 12g - 12 dimensional subspace of the domain of $\mathcal{L}$ with the function $\eta \to \|\mathcal{L}(\eta)\|_{\mathbb{L}}^2$ being small ($< c_0 m^2$). This being the case, it follows using the Raleigh-Ritz variational construction of the eigenspaces of $\mathcal{L}^2$ that the span of the eigenvectors of $\mathcal{L}$ with eigenvalue norm less than $c_0 m$ is at least 12g - 12 dimensional.

To see that span of those small normed eigenvectors can't be greater than 12g - 12, one can replay the preceding argument in reverse: Start with an orthonormal basis of small normed eigenvectors of $\mathcal{L}$ and, for each basis eigenvector, construct the corresponding versions of $\psi_0$ using Lemma 6.4. These will be nearly orthonormal. Use Lemmas 6.6 and 6.7 to see that they span a subspace of the domain of $\mathfrak{L}$ where the function $\psi \to \|\mathfrak{L}(\cdot)\|_{\mathbb{L}}^2$ has norm less than $c_0 m$. The latter subspace can't have dimension greater than 12g - 12 because that is the dimension of the span of the eigenvectors of $\mathfrak{L}$ with eigenvalue norm less than a fixed multiple of $m$.

### e) More about the eigenvectors of $\mathfrak{L}$ and $\mathfrak{L} + \mathfrak{l}(\cdot)$ with small normed eigenvalue

This subsection constitutes a digression to say more about the eigenvectors of the operator $\mathfrak{L}$ and $\mathfrak{L} + \mathfrak{l}(\cdot)$ with eigenvalue norm less than $c_0 m$. What is said about these eigenvectors will be used in Section 7d to say more about the eigenvalue of $\mathcal{L}$ with small eigenvalue norm. (The proof of Proposition 6.2 is completed in Section 7d.)

With regards to $\mathfrak{L}$ and $\mathfrak{L} + \mathfrak{l}(\cdot)$: As noted in Lemma 5.4, the operator $\mathfrak{L}$ maps sections of the subbundle $\mathbb{S}_{SL}$ to sections of $\mathbb{S}_{SL}$. Let $\mathbb{S}_{SL}^\perp$ denote the orthogonal complement subbundle. The operator $\mathfrak{L}$ also maps this bundle to itself (see Lemma 5.4). By virtue of what is said by Lemma 5.4, the endomorphism $\mathfrak{l}(\cdot)$ also preserves the orthogonal splitting of $\mathbb{S}$ as $\mathbb{S}_{SL} \oplus \mathbb{S}_{SL}^\perp$. Meanwhile: According to Lemma 5.7, there are no eigenvalues of $\mathfrak{L}$ acting on the $\mathbb{L}$-completion of $C^\infty(\Sigma; \mathbb{S}_{SL})$ with norm less than $c_0^{-1} m$, and so this also the case for $\mathfrak{L} + \mathfrak{l}(\cdot)$ albeit with a larger value of $c_0$; this is because $|\mathfrak{l}|$ is bounded by $c_0 m^2 \, r^{-2/3}$ (see Proposition 5.3). Of course, all of this supposes that $m$ is small ($< c_0^{-1}$) and $r$ is large ($> c_0$) which will be assumed henceforth with updates on the size of $c_0$.

With the preceding understood, the task at hand is to say something useful about the small eigenvalues of $\mathfrak{L}$ and $\mathfrak{L} + \mathfrak{l}(\cdot)$ as operators on the $\mathbb{L}$-completion of $C^\infty(\Sigma; \mathbb{S}_{SL}^\perp)$. To this end, it proves fruitful to view both $\mathfrak{L}$ of $\mathfrak{L} + \mathfrak{l}(\cdot)$ as perturbations of the operator $\mathfrak{D}_0$ with the latter understood to act on $C^\infty(\Sigma; \mathbb{S}_{SL}^\perp)$. This is done by writing either $\mathfrak{L}$ of $\mathfrak{L} + \mathfrak{l}(\cdot)$ as $\mathfrak{D}_0 + \mathfrak{k}(\cdot)$ with $\mathfrak{k}(\cdot)$ denoting the relevant symmetric endomorphism of $\mathbb{S}_{SL}^\perp$.

With regard to $\mathfrak{D}_0$: This is the operator that is depicted in (5.1), and also (in complex form) on the left hand side of (5.5). It is important to remember in what follows that the $\mathbb{S}_{SL}^\perp$ kernel of $\mathfrak{D}_0$ has dimension 6.g -6, these being the kernel elements depicted in (5.9) with $\varsigma_+, \varsigma_-$ and $\mu_3$ being zero; thus they are described by the top bullet in (5.10). In what follows, $\Pi$ is used to denote the $\mathbb{L}$-orthogonal projection to the $\mathbb{S}_{SL}^\perp$ kernel of $\mathfrak{D}_0$.

With regards to $\mathfrak{k}(\cdot)$: In the case of $\mathfrak{L}$, this is the endomorphism that is depicted in (5.78). In the case of $\mathfrak{L} + \mathfrak{l}(\cdot)$, this is the sum of what is depicted in (5.78) and the endomorphism depicted in (6.58).

Now suppose that $\psi$ denotes an eigenvector of one or the other version of $\mathfrak{D}_0 + \mathfrak{k}(\cdot)$ whose eigenvalue has norm less than $c_0^{-1}$. Let $\lambda$ denote this eigenvalue. Write $\psi$ as $\phi + \psi^\perp$ with $\phi$ denoting the $\mathbb{L}$-orthogonal projection on $C^\infty(\Sigma; \mathbb{S}_{SL}^\perp)$ to the $\mathbb{S}_{SL}^\perp$ kernel of $\mathfrak{D}_0$. This is to say that $\phi = \Pi\psi$. Projecting the eigenvalue equation both orthogonal to the $\mathbb{S}_{SL}^\perp$ kernel of $\mathfrak{D}_0$ and to that kernel writes the eigenvalue equation as the following coupled system of equations for $\phi$ and $\psi^\perp$:

- $\mathfrak{D}_0\psi^\perp$ - $\lambda\psi^\perp$ + $(1 - \Pi)(\mathfrak{k}\psi^\perp)$ = -$(1 - \Pi)\mathfrak{k}\phi$ .
- $\Pi\mathfrak{k}(\phi + \psi^\perp) = \lambda\phi$.

$$(6.61)$$

As explained next: The top equation in (6.61) can be used to define a linear map from the $\mathbb{S}_{SL}^{\perp}$ kernel of $\mathfrak{D}_0$ to its $\mathbb{L}$-orthogonal complement in $C^{\infty}(\Sigma; \mathbb{S}_{SL}^{\perp})$. The key observation for doing this is made by the next lemma.

**Lemma 6.9**: *There exists* $\kappa > 1$ *with the following significance: Fix* $m < \frac{1}{\kappa}$ *and* $r = \frac{n}{m} > \kappa$. *For this data, the corresponding version of the operator* $(1 - \Pi)(\mathfrak{D}_0 + \mathfrak{k})(1 - \Pi)$ *when acting on the* $\mathbb{L}$-completion of $(1 - \Pi)C^{\infty}(\Sigma, \mathbb{S}_{SL}^{\perp})$ *has no eigenvalues with norm less than* $\frac{1}{\kappa}$.

This lemma is proved momentarily.

Supposing that Lemma 6.9 is true, then the operator $(1 - \Pi)(\mathfrak{D}_0 + \mathfrak{h} - \lambda)(1 - \Pi)$ is invertible if $|\lambda| < c_0^{-1}$ when acting on $(1 - \Pi)C^{\infty}(\Sigma; \mathbb{S}_{SL}^{\perp})$ with the $\mathbb{L}$-norm of its inverse being at most $c_0$. This inverse can then be used to define a linear map (denoted by $\mathfrak{T}_{\lambda}$) from the $\mathbb{S}_{SL}^{\perp}$ kernel of $\mathfrak{D}_0$ into $(1 - \Pi)C^{\infty}(\Sigma; \mathbb{S}_{SL}^{\perp})$ using the rule below:

$$\phi \to \mathfrak{T}_{\lambda}(\phi) \equiv \text{-}m\left((1 - \Pi)(\mathfrak{D}_0 + \mathfrak{k} - \lambda)(1 - \Pi)\right)^{-1}(\mathfrak{k}\phi) ,$$

$$(6.62)$$

By virtue of the operator norm bound:

$$\|\mathfrak{T}_{\lambda}(\phi)\|_{\mathbb{L}} \leq c_0 m \, \|\mathfrak{k}\phi\|_{\mathbb{L}}$$

$$(6.63)$$

Granted the existence of $\mathfrak{T}_{\lambda}$ when $|\lambda| < c_0^{-1}$, then the top bullet of (6.63) says in effect that $\mathbb{S}_{SL}^{\perp}$ kernel element $\phi$ is an eigenvector with eigenvalue this same $\lambda$ of the ($\lambda$-dependent) linear map below on the $\mathbb{S}_{SL}^{\perp}$ kernel of $\mathfrak{D}_0$:

$$\phi \to \Pi\mathfrak{k}(\phi + \mathfrak{T}_{\lambda}(\phi)) .$$

$$(6.64)$$

Section 7 analyzes this linear map for any given small normed version of $\lambda$. (Section 7d uses that analysis to finish the proof of Proposition 6.2.)

***Proof of Lemma 6.9***: Suppose $\psi^{\perp}$ denotes an eigenvector in $(1 - \Pi)C^{\infty}(\Sigma; \mathbb{S}_{SL}^{\perp})$ of the operator $\mathfrak{D}_0 + (1-\Pi)\mathfrak{k}$ whose eigenvalue (denoted by $\eta$) has norm less than 1. To elaborate, fix an $\mathbb{L}$-orthonormal frame for the $\mathbb{S}_{SL}^{\perp}$ kernel of $\mathfrak{D}_0$; and let $\Lambda$ denote this set of 6g - 6 eigenvectors. The eigenvalue equation for $\psi^{\perp}$ says in effect the following

$$\mathfrak{D}_0\psi^{\perp} + \mathfrak{k}\psi^{\perp} \text{-} \sum_{\eta \in \Lambda} \eta \, \langle\mathfrak{k}\eta, \psi^{\perp}\rangle_{\mathbb{L}} = \eta\psi^{\perp} ,$$

$$(6.65)$$

where $\langle \, , \, \rangle_{\mathbb{L}}$ denotes the Hilbert space inner product on the $\mathbb{L}$-completion of $C^{\infty}(\Sigma; \mathbb{S}_{SL}^{\perp})$.

With regards to $\langle\mathfrak{k}\eta, \psi^{\perp}\rangle_{\mathbb{L}}$: This is bounded by $\|\mathfrak{k}\eta\|_{\mathbb{L}}\|\psi^{\perp}\|_{\mathbb{L}}$; and Lemma 5.2 with what is said about the terms in (5.32) and (5.33) can be brought to bear see that

$$\|\mathfrak{k}\eta\|_{\mathbb{L}} \le c_0 m r^{-1/3}$$

(6.66)

if $m < c_0^{-1}$ and $r > c_0$. This implies in turn that $\|\Pi(\mathfrak{k}\psi^{\perp})\| \le c_0 m r^{-1/3}\|\psi^{\perp}\|_{\mathbb{L}}$.

To bound the norm of $\mathfrak{k}\psi^{\perp}$: Write $\mathfrak{k}$ as

$$\mathfrak{k} = m\gamma_t[\chi_{\ddagger}A^1, \cdot\,] + \mathfrak{k}_1\,.$$

(6.67)

Meanwhile, it follows from what is said in (5.32) and (5.33) that $\|\mathfrak{k}_1\psi^{\perp}\|_{\mathbb{L}} \le c_0 m\|\psi^{\perp}\|_{\mathbb{L}}$. The bound $\|m\gamma_t[A^1, \psi^{\perp}\|_{\mathbb{L}} \le c_0 m\|\psi^{\perp}\|_{\mathbb{L}}$ will also hold provided that Lemma 4.4 can be invoked for $\psi^{\perp}$ so as to bound the $\mathbb{L}$-integral of $r^2\mathfrak{d}^3|[\hat{\tau}, \psi^{\perp}]|^2$ by $c_0\|\psi^{\perp}\|_{\mathbb{L}}^2$. And, indeed it can be invoked when $m < c_0^{-1}$ because (6.65) with what is said about the terms in (5.32) and (5.33) leads to the bound below for $|[\hat{\tau}, \mathfrak{D}_0\psi^{\perp}]|$:

$$|[\hat{\tau}, \mathfrak{D}_0\psi^{\perp}]| \le c_0(m + \lambda + mr\mathfrak{d}^{3/2})|[\hat{\tau}, \psi^{\perp}]| + c_0 m\, e^{-r\mathfrak{d}^{3/2}/c_0}|\psi| + c_0\textstyle\sum_{\eta\in\Lambda}[\hat{\tau}, \eta]\,\langle \mathfrak{k}\eta,\psi^{\perp}\rangle_{\mathbb{L}}$$

(6.68)

which is an inequality whose right most term bounded by $c_0 m\, e^{-r\mathfrak{d}^{3/2}/c_0}\|\psi^{\perp}\|_{\mathbb{L}}$. (This last bound follows by using the $c_0 m r^{-1/3}\|\psi^{\perp}\|_{\mathbb{L}}$ bound for $|\langle \mathfrak{k}\eta, \psi^{\perp}\rangle_{\mathbb{L}}|$ with what is said by Lemma 5.2 about the pointwise behavior of the kernel of $\mathfrak{D}_0$.)

To finish the proof: Granted the preceding, and supposing that $\eta$ in (6.65) has norm less than $c_0^{-1}$ and supposing that $m < c_0^{-1}$ and $r > c_0$, then (6.65) leads to this:

$$\|\mathfrak{D}_0\psi^{\perp} - \eta\psi^{\perp}\|_{\mathbb{L}} \le c_0 m\, \|\psi^{\perp}\|_{\mathbb{L}}$$

(6.69)

which is not tenable when $\eta$ and $m$ are less than $c_0^{-1}$ and $r > c_0$ because the $r > c_0$ versions of $\mathfrak{D}_0$ do not have non-zero eigenvalues with norm less than $c_0^{-1}$.

## 7. When $\phi \to \Pi\mathfrak{k}(\phi + \mathfrak{T}_\lambda(\phi))$ has eigenvalue $\lambda$ and $|\lambda|$ is small

This section studies the $\lambda$-dependent symmetric endomorphism of the $\mathbb{S}_{\mathrm{SL}}{}^{\perp}$ kernel of $\mathfrak{D}_0$ that is defined by the rule $\phi \to \Pi\mathfrak{k}(\phi + \mathfrak{T}_\lambda(\phi))$, this being the endomorphism from (6.64). Sections 7a and 7b analyze the leading order (in $m$) parts of this endomorphism which is the endomorphism denoted in Section 6 by $\phi \to \Pi(\mathfrak{k}\phi)$ with $\mathfrak{k}$ signifying the endomorphism

$$\phi \to \mathfrak{k}\phi = m\left(c^1\rho_t\,[\sigma_3, \phi] + (\phi)_c\right) + \gamma_t[\chi_{\ddagger}A^1, \phi] + \mathfrak{l}(\phi))\,.$$

(7.1)

(What is denoted by $\mathfrak{l}(\cdot)$ is defined in (6.58).) The following proposition summarizes the story regarding the eigenvalues of $\Pi \mathfrak{k} \Pi$.

**Proposition 7.1**: *There exists* $\kappa > 1$ *with the following significance: Fix* $m < \frac{1}{\kappa}$ *and then fix* $r = \frac{n}{m} > \kappa$. *The corresponding* $(m, r)$ *version of the symmetric endomorphism*

$$\phi \to m\, \Pi\big(c^1 \rho_t\, [\sigma_3, \phi] + (\phi)_\epsilon + \gamma_t[\chi_\ddagger A^1, \phi] + \mathfrak{l}(\phi)\big)$$

*of the* $\mathbb{S}_{SL}{}^\perp$ *kernel of* $\mathfrak{D}_0$ *has two eigenvalues (counted with multiplicity) with norm between* $\frac{1}{\kappa}\, mr^{-2}$ *and* $\kappa mr^{-2}$, *and then* $6g - 8$ *eigenvalues (counted with multiplicity) with norm between* $\frac{1}{\kappa}\, mr^{-2/3}$ *and* $\kappa mr^{-2//3}$.

The proof of this proposition is in Section 7b. Section 7c analyzes the effect of adding the $\lambda$-dependent term $\phi \to \Pi \mathfrak{k}(\mathfrak{T}_\lambda(\phi))$ to $\Pi \mathfrak{k} \Pi$. The last subsection, Section 7d, uses the results from the previous subsection and from Section 6 to finish Proposition 6.2's proof.

**a) The eigenvalues of the** $\Pi(mc^1\rho_t\,[\sigma_3, \cdot\,] + (\cdot)_\epsilon)$ **part of** $\Pi \mathfrak{k} \Pi$

The upcoming Proposition 7.2 concerns the largest part of the endomorphism $\Pi \mathfrak{k} \Pi$ which is (looking ahead) the symmetric endomorphism the $\mathbb{S}_{SL}{}^\perp$ kernel of $\mathfrak{D}_0$ given by the rule $\phi \to \Pi(c^1 \rho_t\,[\sigma_3, \phi] + (\phi)_\epsilon)$. The next subsection explains why the $\gamma_t[\chi_\ddagger A^1, \phi]$ and $\mathfrak{l}(\phi)$ contributions to $\Pi(\mathfrak{k}\phi)$ do not substantively change the conclusions of Proposition 7.1.

**Proposition 7.2**: *There exists* $\kappa > 1$ *with the following significance: Fix* $m < \frac{1}{\kappa}$ *and then fix* $r = \frac{n}{m} > \kappa$. *The corresponding version of the symmetric endomorphism*

$$\phi \to m\, \Pi\big(c^1 \rho_t\, [\sigma_3, \phi] + (\phi)_\epsilon\big)$$

*of the* $\mathbb{S}_{SL}{}^\perp$ *kernel of* $\mathfrak{D}_0$ *has two eigenvalues (counted with multiplicity over* $\mathbb{R}$ *) with norm between* $\frac{1}{\kappa}\, mr^{-2}$ *and* $\kappa mr^{-2}$, *and then* $6g - 8$ *eigenvalues (counted with multiplicity over* $\mathbb{R}$ *) with norm between* $\frac{1}{\kappa}\, mr^{-2/3}$ *and* $\kappa mr^{-2//3}$

The rest of this subsection is devoted to the proof of this proposition. (The proof establishes that all of the eigenvalues are, in fact, positive.)

***Proof of Proposition 7.2:*** The proof has 8 parts. By way of notation: The proof writes the proposition's endomorphism as $m \yen$ with $\yen$ denoting the endomorphism of the $\mathbb{S}_{SL}{}^\perp$ kernel of $\mathfrak{D}_0$ given by the following $\phi \to \yen(\phi) \equiv \Pi\big(c^1 \rho_t\, [\sigma_3, \phi] + (\phi)_\epsilon\big)$ .

*Part 1*:  An element in the $\mathbb{S}_{SL}{}^\perp$ kernel of $\mathfrak{D}_0$ is determined by what are denoted by $b$ and $c$ in (5.4) since its $u$ and $x$ components are zero.  The $\Pi((\phi)_c)$ part of ¥ sees only the $c$ component of $\phi$ where as the $\Pi(c^1\rho_t[\sigma_3,\phi])$ part sees only the $b$ component.  To elaborate on this point, suppose that $\phi$ and $\phi'$ are in the $\mathbb{S}_{SL}{}^\perp$ kernel of $\mathfrak{D}_0$.  Then the pointwise inner product between $\phi$ and $¥\phi'$ is the real part of what is depicted below

$$ic^1\,\langle b^*(\,[\sigma_3,\,b'])\rangle + \langle c^*c'\rangle\,.$$

(7.2)

(Here and subsequently, $(\cdot)^*$ denotes $-(\cdot)^\dagger$, the $\mathbb{C}$-anti-linear operation on the bundle of endomorphisms of E that acts as the product of -1 and the Hermitian adjoint operation.)  With regards to (7.2):  If $b$ is written as $\mu_+\varphi^* - \mu_-[\frac{i}{2}\sigma_3,\,\varphi^*]$ and $c$ as $-\frac{i}{2}(\nabla_1 - i\nabla_2)\mu_-$, and if $b'$ and $c'$ are written analogously, then (7.2) is this:

$$c^1\left((\bar{\mu}_+\mu'_+ + \bar{\mu}_-\mu'_-)\langle i\sigma_3[\varphi^*,\varphi]\rangle - 2(\bar{\mu}_+\mu'_- + \bar{\mu}_-\mu'_+)|\varphi|^2\right) + \frac{1}{4}(\nabla_1 + i\nabla_2)\bar{\mu}_-(\nabla_1 - i\nabla_2)\mu'_-$$

(7.3)

This last depiction of ¥ will be exploited momentarily.  Even so, the depiction in (7.2) has an immediate take-away which is that the contribution to the $\mathbb{L}$-inner product between $\phi$ and $¥\phi'$ from the $\mathfrak{d} > \frac{1}{100}\,r_0$ part of $\Sigma$ is bounded by $c_0e^{-r/c_0}\,\|\phi\|_\mathbb{L}\|\phi'\|_\mathbb{L}$.  (This is because $|c|$ is no greater than $c_0\,e^{-r/c_0}$ there and $|b|$ and $|b'|$ are at most $c_0\|\phi\|_\mathbb{L}$ and $\|\phi'\|_\mathbb{L}$ there;and because $|c|$ and $|c'|$ are no greater than $c_0e^{-r/c_0}\,\|\phi\|_\mathbb{L}$ and $c_0e^{-r/c_0}\,\|\phi'\|_\mathbb{L}$ there.)  As a consequence, the focus until until Part 7  is on the contribution of ¥ from the $\mathfrak{d} \le \frac{1}{100}\,r_0$ disk around a given point from $\Theta$.  In this regard, the focus will be on a point in $\Theta_+$ since the analysis for points in $\Theta_-$ is essentially the same.

*Part 2*:  Fix a point in $\Theta_+$ so as to use the holomorphic coordinate from Lemma 3.2 to study the contribution to ¥ from the $|z| < \frac{1}{100}\,r_0$ around that point.  Let $D_1$ denote that disk.  The contribution to ¥ from $D_1$ will be estimated subsequently via an analysis of the symmetric, quadratic form

$$(\phi',\phi) \to \int_{D_1}\,\langle \phi',\,¥(\phi)\rangle$$

(7.4)

on the $\mathbb{S}_{SL}{}^\perp$ kernel of $\mathfrak{D}_0$.  Various approximations and identities will be used to this end.  As a preliminary to their introduction:  Use the holomorphic coordinate to write $\mu_+$ as a Laurent series in the variable $\bar{z}$ on the $|z| \le r_0$ disk in $\mathbb{C}$ as done in the second bullet of Lemma 5.2 for $\varsigma_+$ (which is the complex conjugate of $\mu_+$).  Denote this series by $P(\bar{z})$.  As noted in Lemma 5.2, $\mu_-$ differs by at most $c_0\,e^{-r/c_0}$ from the product of a corresponding

Laurent series in $\bar{z}$ with what is denoted by $\Xi$ in Lemma 5.2. Use $Q(\bar{z})$ to denote that corresponding Laurent series.

The subsequent paragraphs in this part of the proof describe the approximations that are used. The next two parts of the proof derive identities that are used in the proof.

The first approximation replaces $\mu_-$ in (7.3) by $Q\Xi$ which changes the value of (7.4) by at most $c_0 e^{-r/c_0} \|\phi\|_{\mathbb{L}} \|\phi'\|_{\mathbb{L}}$. (See Lemma 5.2.) The same can be said when $\mu'_-$ is replaced by the corresponding primed version of $Q\Xi$ (which is denoted by $Q'\Xi$). By the same token, what is denoted by $c^1$ in (7.3) can be replaced on the $|z| < \frac{1}{50} r_0$ disk by $-\frac{1}{2}\Xi$ which changes (7.4) at most $c_0 e^{-r/c_0} \|\phi\|_{\mathbb{L}} \|\phi'\|_{\mathbb{L}}$. (This is because $c^1$ there is the function $c$ from (5.20) which is $-\frac{1}{2}$ times the $\mu_+ = 1$ version of $\mu_-$. Note also Lemma 3.3.) With these replacements, (7.3) becomes this:

$$-\tfrac{1}{2}\left(\bar{P}P'\Xi + \bar{Q}Q'\Xi^3\right)\langle i\sigma_3[\varphi^*, \varphi]\rangle + \left(\bar{P}Q' + \bar{Q}P'\right)|\varphi|^2\Xi^2 + \tfrac{1}{4}\bar{Q}Q'|d\Xi|^2 \ .$$

$$(7.5)$$

To set the stage for the next set of approximations, reintroduce the function $\mathfrak{f}$ from Lemma 3.3 and use $\mathfrak{f}$ to define a new function on $\mathbb{C}$ to be denoted by $v_{\mathfrak{f}}$ by the rule whereby

$$v_{\mathfrak{f}}(z) = \mathfrak{f}(\alpha^{1/3}|z|) + \tfrac{1}{6}\ln(\alpha) \ .$$

$$(7.6)$$

It follows from what is said by Lemma 3.3 that the functions $|\varphi|^2$ and $\langle i\sigma_3[\varphi^*, \varphi]\rangle$ in (7.5) can be replaced on $D_1$ by the respective functions

- $|\varphi|^2 = \frac{1}{4}\alpha \left(e^{2v_{\mathfrak{f}}}|z|^2 + e^{-2v_{\mathfrak{f}}}\right)$
- $\langle i\sigma_3[\varphi^*, \varphi]\rangle = \frac{1}{2}\alpha \left(e^{2v_{\mathfrak{f}}}|z|^2 - e^{-2v_{\mathfrak{f}}}\right)$

$$(7.7)$$

at the expense of changing the value of (7.4) by at most $c_0 e^{-r/c_0} \|\phi\|_{\mathbb{L}} \|\phi'\|_{\mathbb{L}}$. But keep in mind that these replacements require the use of the Euclidean area 2-form $\frac{i}{2} dz \wedge d\bar{z}$ for the integration in (7.4). This is a manifestation of the conformal invariance that was mentioned with regards to (3.8) and Lemma 3.2 in Part 6 of Section 3a.

When $\phi$, $\phi'$ and $c$ and $|\varphi|^2$ and $\langle i\sigma_3[\varphi^*, \varphi]\rangle$ are changed in (7.4) as described above (and likewise the Euclidean area form is used), then there is an added benefit which is that powers of $\bar{z}$ that appear in $P$, $Q$ and their primed analogs are orthogonal with respect to the quadratic form in (7.4). Indeed, this is because the function $\Xi$ from Lemma 3.3 and $v_{\mathfrak{f}}$ are functions only of $|z|$ and so they are rotationally invariant.

With the preceding points understood, it is sufficient for the analysis of (7.4) (up to an error bounded by $c_0 e^{-r/c_0} \|\phi\|_{\mathbb{L}} \|\phi'\|_{\mathbb{L}}$) to focus on the cases when $P$ and $Q$ and their primed counterparts are homogeneous, which is to say $P$ and $P'$ are the monomial function $\bar{z}^k$ for a given $k \in \{-1, 0, \dots\}$ in which case $Q$ and $Q'$ are the monomial $\frac{3}{3+2k} \bar{z}^k$. In particular, the analysis in Part 4 and 5 of this subsection will focus for the most part on the $k = 0$ case (so $P = 1$ and $Q = 1$ and likewise for their primed counterparts), and the analysis in Parts 6-8 will focus on the case when $k = -1$ (so $P = \frac{1}{\bar{z}}$ and $Q = \frac{3}{\bar{z}}$ and likewise for the primed versions). As it turns out, the focus on these two cases is sufficient for the purposes of proving Proposition 7.2.

*Part 3*: In the event that $P = Q = 1$ and also $P' = Q' = 1$, what is written in (7.5) is

$$-\tfrac{1}{2}\,(\Xi + \Xi^3)\,\langle i\sigma_3[\varphi^*, \varphi]\rangle + 2|\varphi|^2\Xi^2 + \tfrac{1}{4}\,|d\Xi|^2$$

$$(7.8)$$

with it understood that $\langle i\sigma_3[\varphi^*, \varphi]\rangle$ and $|\varphi|^2$ are given by the expressions in (7.7). The task at hand is to determine the value of the integral of this expression over the disk $D_1$. To this end, note first that all of the terms are defined on the whole of $\mathbb{C}$ and that the integral of (7.8) over the whole of $\mathbb{C}$ differs from its integral over $D_1$ by at most $c_0 e^{-r/c_0}$. (This follows from what is said in Lemmas 3.3 and 5.2 to the effect that $|\Xi|$ is bounded by $c_0\, c_0 e^{-r\,|z|^{3/2}/c_0}$ and that $|\varphi|$ where $|z| > r^{-2/3}$ is bounded by $c_0 r\,|z|$.) As explained below, the integral of (7.8) over $\mathbb{C}$ is exactly $\frac{2\pi}{18}\frac{17}{36}$.

Some identities are needed to justify the claim. These identities (and their proofs) use $s$ to denote the radial variable on $\mathbb{C}$, which is to say that $s = |z|$.) The first identity uses (7.7) and the equation obeyed by $v_f$ and the definition of $\Xi$ to write

$$\langle i\sigma_3[\varphi^*, \varphi]\rangle = -\tfrac{3}{2}\,\tfrac{1}{s}\tfrac{d}{ds}\Xi \ .$$

$$(7.9)$$

Indeed, this follows because $\Xi$ is equal to $-\tfrac{2}{3}\,s\tfrac{d}{ds}v_f - \tfrac{1}{3}$ and because $v_f$ is a radially symmetric solution to the equation in Lemma 3.2 (which is to say that $\tfrac{1}{s}\tfrac{d}{ds}(s\tfrac{d}{ds})v_f = \tfrac{1}{2}\,\alpha(e^{2v_f}s^2 - e^{-2v_f})$.)

Using (7.9) in (7.8) leads to this

$$\tfrac{3}{4}\,(\Xi + \Xi^3)\tfrac{1}{s}\tfrac{d}{ds}\Xi \ + 2|\varphi|^2\Xi^2 + \tfrac{1}{4}\,|\tfrac{d}{ds}\Xi|^2$$

$$(7.10)$$

The next lemma asserts three identities that are used to evaluate the $\mathbb{C}$ integral of (7.10).

**Lemma 7.3**: *The following identities hold*

- $\Xi|_{s=0} = -\tfrac{1}{3}\ .$
- $\int_{\mathbb{C}}\ |\tfrac{d}{ds}\Xi|^2 = \ \tfrac{2\pi}{18}\tfrac{17}{18}$
- $\int_{\mathbb{C}}\ |\varphi|^2\Xi^2 = \ \tfrac{2\pi}{36}\tfrac{37}{36}$

(Note that the $-\frac{1}{3}$ value for $\Xi$ at $s = 0$ is predicated on the point chosen from $\Theta$ being in $\Theta_+$. If the point were from $\Theta_-$, then $\Xi|_{s=0}$ would be $\frac{1}{3}$.) This lemma is proved momentarily. Accept it for the time being.

To compute the $\mathbb{C}$ integral of (7.10): Note first that (7.10) is

$$\frac{1}{s}\frac{d}{ds}(\frac{3}{8}\,\Xi^2 + \frac{3}{16}\,\Xi^4) + 2|\varphi|^2\Xi^2 + \frac{1}{4}\,|\frac{d}{ds}\Xi|^2$$

(7.11)

Since the left most term in (7.11) is rotationally invariant (actually, all of the terms are), its integral over $\mathbb{C}$ is $2\pi$ times its integral over $[0, \infty)$ using the measure s ds. With this understood, it follows that the $\mathbb{C}$ integral of the left most term in (7.11) is

$$-2\pi(\frac{3}{8}\,\Xi^2 + \frac{3}{16}\,\Xi^4)|_{s=0}$$

(7.12)

which is (according to Lemma 7.3) equal to $-\frac{2\pi}{18}\frac{19}{24}$. Meanwhile, Lemma 7.2 says that the integral of the middle term in (7.11) is $\frac{2\pi}{18}\frac{37}{36}$, and it says that the integral of the right most term in (7.11) is $\frac{2\pi}{18}\frac{17}{72}$. Adding these values gives $\frac{2\pi}{18}\frac{17}{36}$ for the $\mathbb{C}$ integral of (7.8).

*Part 4*: This part justifies the claims made by Lemma 7.3.

**Proof of Lemma 7.3:** The top bullet's identity follows directly from $\Xi$'s definition in Lemma 5.2 as $-\frac{2}{3}\,s\,\frac{d}{ds}v_f - \frac{1}{3}$. A preliminary step for proving the first integral identity is the integral identity below:

$$\int_{\mathbb{C}} (|\frac{d}{ds}\Xi|^2 + 4|\varphi|^2\Xi^2) = \frac{2\pi}{6}.$$

(7.13)

The first integral identity follows from the one above and the second integral identity.

To prove (7.13), keep in mind that $\Xi$ obeys the equation below on $\mathbb{C}$:

$$-\frac{1}{s}\frac{d}{ds}(s\,\frac{d}{ds}\Xi) + \alpha\,(e^{2v_f}\,s^2 + e^{-2v_f})\,\Xi \;=\; \alpha\,(e^{2v_f}\,s^2 - e^{-2v_f}).$$

(7.14)

With this understood, use (7.9) to rewrite this as

$$-\frac{1}{s}\frac{d}{ds}(s\,\frac{d}{ds}\Xi) + \alpha\,(e^{2v_f}\,s^2 + e^{-2v_f})\,\Xi = -3\,\frac{1}{s}\frac{d}{ds}\Xi \;.$$

(7.15)

Now multiply both sides of this by $\Xi$ to obtain

$$-\frac{1}{2}\frac{1}{s}\frac{d}{ds}(s\frac{d}{ds}\Xi^2) + |\frac{d}{ds}\Xi|^2 + \alpha\,(e^{2v_f}\,s^2 + e^{-2v_f})\,\Xi^2 = -\frac{3}{2}\,\frac{1}{s}\frac{d}{ds}\,\Xi^2$$

$$(7.16)$$

Integrate the latter over $[0,\infty)$ using the measure $sds$ and then integrate by parts. The integration by parts on the left hand side has no $s = 0$ contribution and no $s \to \infty$ contribution (by virtue of the rapid $s \to \infty$ decay to zero of $\Xi$ and $d\Xi$); it integrates to $\frac{1}{2\pi}$ times the integral in (7.13). Meanwhile, the integration by parts on the right hand side has $s = 0$ contribution equal to $\frac{3}{2}\,\Xi^2|_{s=0}$ which is $\frac{1}{6}$ and (likewise) no contribution from $s \to \infty$.

With regards to the identity in the third bullet: The sleight of hand that follows may be more tranparent by writing $v_f = u - \frac{1}{2}\ln(s)$. The function $\Xi$ and (7.9) and (7.14) and the top bullet in (7.7) can be written with the function $u$ as:

- $\Xi = -\frac{2}{3}\,s\frac{d}{ds}u$ .

- $-\frac{1}{s}\frac{d}{ds}(s\frac{d}{ds}\Xi) + \alpha\,s\,(e^{2u} + e^{-2u})\,\Xi\,=\,\alpha\,s\,(e^{2u} - e^{-2u})$ .

- $\frac{1}{2}\,\alpha\,s\,(e^{2u} - e^{-2u}) = -\frac{3}{2}\,\frac{1}{s}\frac{d}{ds}\Xi$ .

- $|\varphi|^2 = \frac{1}{4}\,\alpha\,s\,(e^{2u} + e^{-2u})$

$$(7.17)$$

With regards to $u$: This function is defined on $(0,\infty)$ but not at zero; it is singular as $s \to 0$ (it behaves like $-\frac{1}{2}\ln(s)$ near $s = 0$), but it is bounded where $s > r^{-2/3}$ by $c_0 e^{-r s^{3/2}/c_0}$. With regards to the singularity at $s = 0$: Keep a little bit alert when integrating by parts below.

To derive the value of the second integral in Lemma 7.3: Multiply both sides of the second bullet in (7.17) by the function $s^2\frac{d}{ds}\Xi$ to obtain the identity

$$-\frac{1}{2}\frac{d}{ds}((s\frac{d}{ds}\Xi))^2 + \frac{1}{2}\,\alpha\,s^3\,(e^{2u} + e^{-2u})\,\frac{d}{ds}\Xi^2\,=\,\alpha\,s^3\,(e^{2u} - e^{-2u})\,\frac{d}{ds}\Xi^2\,.$$

$$(7.18)$$

Integrate both sides of this on the domain $[0,\infty)$ and then integrate by parts. The integral of the left most term on the left hand side is zero. The integrals of the other terms give this:

$$-\frac{3}{2}\,\alpha\int_0^\infty\,s^2(e^{2u} + e^{-2u})\Xi^2\,ds\,-\,\alpha\int_0^\infty\,s^3(e^{2u} - e^{-2u})(\frac{d}{ds}u)\,\Xi^2\,ds =$$

$$-\,3\,\alpha\int_0^\infty\,s^2(e^{2u} - e^{-2u})\,\Xi\,ds\,-\,2\alpha\int_0^\infty\,s^3(e^{2u} + e^{-2u})(\frac{d}{ds}u)\,\Xi\,ds\,.$$

$$(7.19)$$

(There is no worry in any case about $s = 0$ or $s \to \infty$ contributions when integrating by parts because of the factors of $s$ being greater than 1 and because of the $c_0 e^{-r s^{3/2}/c_0}$ bounds for $\Xi$ and $u$ at large $s$.) Now use the identity in the top bullet of (7.17) to write (7.19) as

$$-\tfrac{3}{2}\,\alpha \int_0^\infty\, s^2(e^{2u}+e^{-2u})\Xi^2\,ds\; +\tfrac{3}{2}\,\alpha \int_0^\infty\, s^2(e^{2u}-e^{-2u})\,\Xi^3\,ds\;=$$

$$-\,3\,\alpha \int_0^\infty\, s^2(e^{2u}-e^{-2u})\,\Xi\,ds\; +3\alpha \int_0^\infty\, s^2(e^{2u}+e^{-2u})\,\Xi^2\,ds\;.$$

$$(7.20)$$

Moving terms from one side to the other transforms (7.20) to

$$\tfrac{3}{2}\,\alpha \int_0^\infty\, s^2(e^{2u}-e^{-2u})\,\Xi^3\,ds + 3\,\alpha \int_0^\infty\, s^2(e^{2u}-e^{-2u})\,\Xi\,ds\;=\;\tfrac{9}{2}\,\alpha \int_0^\infty\, s^2(e^{2u}+e^{-2u})\,\Xi^2\,ds.$$

$$(7.21)$$

Granted this identity, use the third bullet in (7.17) to rewrite (7.21) as

$$-\tfrac{9}{2}\int_0^\infty\, \Xi^3\tfrac{d}{ds}\Xi\,ds - 9\int_0^\infty\, \Xi\tfrac{d}{ds}\Xi\,ds = \tfrac{9}{2}\,\alpha \int_0^\infty\, s^2(e^{2u}+e^{-2u})\,\Xi^2\,ds\;.$$

$$(7.22)$$

Integation can now be done on the left hand side to obtain

$$\tfrac{1}{8}\,\Xi^4|_{s=0} + \tfrac{1}{2}\,\Xi^2|_{s=0} = \tfrac{1}{2}\,\alpha \int_0^\infty\, s^2(e^{2u}+e^{-2u})\,\Xi^2\,ds$$

$$(7.23)$$

This identity leads directly to the second integral identity in the lemma because $\Xi_{s=0}=-\tfrac{1}{3}$ and because the integral on the right hand side of (7.23) is $\tfrac{1}{\pi}$ times the $\mathbb{C}$-integral of $|\varphi|^2\Xi^2$ (as per the fourth bullet in (7.17)).

   *Part 5:* This part of Proposition 7.2's proof establishes some preliminary facts about the $D_1$ integral of (7.5) in the case when $P = \tfrac{1}{\bar z}$ and $P' = \tfrac{1}{\bar z}$ (and thus $Q$ and $Q'$ are $\tfrac{3}{\bar z}$). In this case, the expression in (7.5) can be written as

$$-\tfrac{1}{2}\tfrac{1}{s^2}\,(\Xi + 9\Xi^3)\,\langle i\sigma_3[\varphi^*,\varphi]\rangle + \tfrac{6}{s^2}\,|\varphi|^2\Xi^2 + \tfrac{9}{4}\tfrac{1}{s^2}\,|\tfrac{d}{ds}\Xi|^2\;.$$

$$(7.24)$$

As before, the integral of (7.24) over $\mathbb{C}$ differs from its $D_1$ integral by at most $c_0e^{-r/c_0}$. In this regard: The $\mathbb{C}$ integral of (7.24) is convergent near $s = 0$ not-with-standing the appearance of negative powers of $s$. Indeed, the integral of the right most term in (7.24) is convergent near $s = 0$ because $|\tfrac{d}{ds}\Xi|$ is $\mathcal{O}(s)$ there (see the second bullet in (7.17)). The integral of the other two terms taken together converges near $s = 0$ as can be seen by writing (7.24) as

$$-\tfrac{1}{2}\tfrac{1}{s^2}\,\Xi(1+3\Xi)^2\,\langle i\sigma_3[\varphi^*,\varphi]\rangle + \tfrac{3}{s^2}\,(2|\varphi|^2 + \langle i\sigma_3[\varphi^*,\varphi]\rangle)\Xi^2 + \tfrac{9}{4}\tfrac{1}{s^2}\,|\tfrac{d}{ds}\Xi|^2\;,$$

$$(7.25)$$

then noting that $1+3\Xi = -\,2s\tfrac{d}{ds}v_f$ which is $\mathcal{O}(s^2)$ as $s \to 0$, and noting (see (7.17)) that

$$(2|\varphi|^2 + \langle i\sigma_3[\varphi^*, \varphi]\rangle) = \alpha \, e^{2v_f} \, s^2 \, .$$

(7.26)

With regards to the respective signs of the terms in (7.25): The middle term and the right most term are manifestly non-negative. Meanwhile, the left most term is non-positive. This is because $\Xi$ is negative (see Lemma 5.2), and so is $\langle i\sigma_3[\varphi^*, \varphi]\rangle$ (see the second bullet of (7.17) and Lemma 3.3).

As for the integral of the expression in (7.25) over $\mathbb{C}$, this is equal to $2\pi$ times the integral of that expression over $[0, \infty)$ using the measure $s \, ds$. The analysis of the latter integral starts with that $[0, \infty)$ integral of the left most term, thus with

$$-\tfrac{1}{2} \int_0^\infty \tfrac{1}{s^2} \, \Xi(1 + 3\Xi)^2 \, \langle i\sigma_3[\varphi^*, \varphi]\rangle \, s ds \, .$$

(7.27)

The first step in analyzing (7.27) invokes (7.17) to write it as

$$\tfrac{3}{4} \int_0^\infty \tfrac{1}{s^2} \, \Xi(1 + 3\Xi)^2 \tfrac{d}{ds}\Xi \, ds \, .$$

(7.28)

The starting point for rewriting (7.28)'s integral is the integral below in (7.29) which has just one power of $(1 + 3\Xi)$ instead of two. (It is convergent because $1 + 3\Xi$ is $\mathcal{O}(s^2)$ as $s \to 0$.)

$$\tfrac{3}{4} \int_0^\infty \tfrac{1}{s^2} \, \Xi(1 + 3\Xi) \tfrac{d}{ds}\Xi \, ds \, .$$

(7.29)

Noting that $\tfrac{1}{s^2}$ is $-\tfrac{d}{ds}(\tfrac{1}{s})$, integration by parts identifies the integral in (7.29) with

$$\tfrac{3}{4} \int_0^\infty \tfrac{1}{s}(1 + 6\Xi)(\tfrac{d}{ds}\Xi)^2 \, ds \; + \; \tfrac{3}{4} \int_0^\infty \tfrac{1}{s}\Xi(1 + 3\Xi) \tfrac{d^2}{ds^2}\Xi \, ds \, ,$$

(7.30)

The next step invokes the second bullet in (7.17) to rewrite the right most integral in (7.30) as follows:

$$\tfrac{3}{4} \int_0^\infty \tfrac{1}{s} \left(\Xi(1 + 3\Xi)\tfrac{d^2}{ds^2}\Xi\right) ds \; = \; -\tfrac{3}{4} \int_0^\infty \tfrac{1}{s^2} \, \Xi(1 + 3\Xi)\tfrac{d}{ds}\Xi \, ds$$
$$+ \; \tfrac{3}{4}\alpha \int_0^\infty \Xi^2(1 + 3\Xi)(e^{2u} + e^{-2u}) ds \; - \; \tfrac{3}{4}\alpha \int_0^\infty \Xi(1 + 3\Xi)(e^{2u} - e^{-2u}) \, ds \, .$$

(7.31)

Now use (7.9) and the third bullet in (7.17) to write the latter identity as

$$\tfrac{3}{4} \int_0^\infty \tfrac{1}{s} \left(\Xi(1 + 3\Xi)\tfrac{d^2}{ds^2}\Xi\right) ds \; = \; \tfrac{3}{2} \int_0^\infty \tfrac{1}{s^2} \, \Xi(1 + 3\Xi)\tfrac{d}{ds}\Xi \, ds \; + \; 3 \int_0^\infty \tfrac{1}{s}\Xi^2(1 + 3\Xi)|\varphi|^2 ds \, .$$

(7.32)

Putting this back into (7.30) leads to the identity below for the integral in (7.29)

$$\frac{3}{4}\int_0^\infty \frac{1}{s^2}\Xi(1+3\Xi)\frac{d}{ds}\Xi\, ds = -3\int_0^\infty \frac{1}{s}\Xi^2(1+3\Xi)|\varphi|^2 ds - \frac{3}{4}\int_0^\infty \frac{1}{s}(1+6\Xi)(\frac{d}{ds}\Xi)^2 ds\,.$$
(7.33)

By adding $\frac{1}{2}\langle i\sigma_3[\varphi^*, \varphi]\rangle$ to $|\varphi|^2$ and then subtracting the same, the left most term on the right of (7.33) can be written as

$$-\frac{3}{2}\int_0^\infty \frac{1}{s}\Xi^2(1+3\Xi)(2|\varphi|^2 + \langle i\sigma_3[\varphi^*, \varphi]\rangle)\, ds + \frac{3}{2}\int_0^\infty \frac{1}{s}\Xi^2(1+3\Xi)\langle i\sigma_3[\varphi^*, \varphi]\rangle)\, ds$$
(7.34)

which is (see (7.9)) the same as the integral below.

$$-\frac{3}{2}\int_0^\infty \frac{1}{s}\Xi^2(1+3\Xi)(2|\varphi|^2 + \langle i\sigma_3[\varphi^*, \varphi]\rangle)\, ds - \frac{9}{4}\int_0^\infty \frac{1}{s^2}\Xi^2(1+3\Xi)\frac{d}{ds}\Xi\, ds\,.$$
(7.35)

The use of (7.35) for the left most term on the right hand side of (7.33) leads directly to the identity below for the integral in (7.29):

$$\frac{3}{4}\int_0^\infty \frac{1}{s^2}\Xi(1+3\Xi)\frac{d}{ds}\Xi\, ds = -\frac{9}{4}\int_0^\infty \frac{1}{s^2}\Xi^2(1+3\Xi)\frac{d}{ds}\Xi\, ds$$
$$-\frac{3}{2}\int_0^\infty \frac{1}{s}\Xi^2(1+3\Xi)(2|\varphi|^2 + \langle i\sigma_3[\varphi^*, \varphi]\rangle)\, ds - \frac{3}{4}\int_0^\infty \frac{1}{s}(1+6\Xi)(\frac{d}{ds}\Xi)^2 ds\,.$$
(7.36)

Moving the left most term on the right hand of (7.36) to the left hand side while changing it sign gives the sought after rewriting of the integral in (7.28):

$$\frac{3}{4}\int_0^\infty \frac{1}{s^2}\Xi(1+3\Xi)^2\frac{d}{ds}\Xi\, ds = -\frac{3}{2}\int_0^\infty \frac{1}{s}\Xi^2(1+3\Xi)(2|\varphi|^2 + \langle i\sigma_3[\varphi^*, \varphi]\rangle)\, ds$$
$$-\frac{3}{4}\int_0^\infty \frac{1}{s}(1+6\Xi)(\frac{d}{ds}\Xi)^2 ds\,.$$
(7.37)

With it understood that (7.27) is the same as (7.28), the preceding identity leads to the expression below for the $[0, \infty)$ integral of (7.25) using the measure s ds:

$$3\int_0^\infty \frac{1}{s}\Xi^2(1 - \frac{1}{2}(1+3\Xi))(2|\varphi|^2 + \langle i\sigma_3[\varphi^*, \varphi]\rangle) ds + \frac{3}{4}\int_0^\infty \frac{1}{s}(2 - 6\Xi)(\frac{d}{ds}\Xi)^2 ds\,.$$
(7.38)

The left most integral in (7.38) is positive and greater than

$$\frac{3}{2}\int_0^\infty \frac{1}{s}\Xi^2(2|\varphi|^2 + \langle i\sigma_3[\varphi^*, \varphi]\rangle) ds$$
(7.39)

because $(1 + 3\Xi)$ is 0 at $s = 0$ and it increases as s increases with limit 1 as $s \to \infty$. (See Lemma 5.2.) Meanwhile, the left most integral in (7.38) is also positive (because $\Xi$ is negative), and it is greater than

$$\frac{3}{2} \int_0^\infty \frac{1}{s} \left(\frac{d}{ds}\Xi\right)^2 ds .$$

(7.40)

This is because the quadratic function $x \to -12x - 27x^2$ on the interval $[-\frac{1}{3}, 0]$ is non-negative. (It can be written as $-3x(4 + 9x)$.) Therefore: The integral of (7.24) (aka (7.25)) over the whole of $\mathbb{C}$ is positive.

Given that the integral in (7.24) is positive, it then follows (by writing s as $\alpha^{-1/3}s'$ and writing the integral in terms of $s'$) that the integral in (7.24) is greater than $c_0^{-1}\alpha^{2/3}$ which is in turn greater than $c_0^{-1}r^{4/3}$.

*Part 6:* Here is where things stand with regards to the $\Sigma$ integral of $\langle \phi', \yen\phi\rangle$ when $\phi'$ and $\phi$ are $\mathbb{S}_{SL}^\perp$ kernel elements: Write the $\mathbb{S}_{SL}^\perp$ kernel of $\mathfrak{D}_0$ as $\mathbb{C} \oplus K$ with the $\mathbb{C}$ summand denoting the kernel elements that correspond to triples $(\mu_+, \mu_-, \varsigma_3)$ with $\mu_+$ being a non-zero complex number. Meanwhile, K here (and below) denotes the $\mathbb{L}$-orthogonal complement to the $\mathbb{C}$-summand in the $\mathbb{S}_{SL}^\perp$ kernel of $\mathfrak{D}_0$. In particular, if $\phi$ is from the K summand, then the corresponding version of $\mu_+$ on each $p \in \Theta$ version of the disk D (where $|z| < r_0$) appears as

$$\mu_+ = \frac{a_{-(p)}}{\bar{z}} + a_{0(p)} + \mathfrak{r}_p(z)$$

(7.41)

where $\mathfrak{r}_p$ is a convergent power series in $\bar{z}$ on the $|z| < 2r_0$ disk in $\mathbb{C}$ with no constant term. (What are denoted by $a_{-(p)}$ and $a_{0(p)}$ in (7.40) are complex numbers.) Moreover, at least one of the $p \in \Theta$ versions $a_{-(p)}$'s is non-zero.

The analysis in the preceding Parts 3-5 of Proposition 7.2's proof implies directly what is written below in (7.42) about the $\Sigma$ integral of $\langle\phi', \yen\phi\rangle$ when $\phi$ and $\phi'$ are from the $\mathbb{S}_{SL}^\perp$ kernel of $\mathfrak{D}_0$.

- *The symmetric bilinear form* $(\phi', \phi) \to \int_\Sigma \langle\phi', \yen\phi\rangle$ *is positive definite on the $\mathbb{C}$ summand. Moreover, supposing that $\phi$ is from this summand and that $a \in \mathbb{C}$ denotes its corresponding version of* $\mu_+$ *then,* $\int_\Sigma \langle\phi, \yen\phi\rangle$ *is within* $c_0 e^{-r/c_0}$ *of the number*

$$\frac{2\pi}{18} \frac{17}{36} (4g - 4) |a|^2 .$$

- *When $\phi$ is from the K summand and $\phi'$ is from the $\mathbb{C}$-summand: Write each $p \in \Theta$ version of $\mu_+$ as in (7.40) and let $a'$ denote the constant $\phi'$ version of $\mu_+$. Then* $\int_\Sigma \langle\phi', \yen\phi\rangle$ *is within* $c_0 e^{-r/c_0}$ *of the number*

$$\frac{2\pi}{18}\frac{17}{36}\sum_{p\in\Theta}\mathfrak{Re}(\bar{a}_{0(p)}a')\,.$$

- *The symmetric bilinear form* $(\phi',\phi)\rightarrow\int_{\Sigma}\langle\phi',\maltese\phi\rangle$ *is positive definite on the* K-*summand,. Moreover, there is a number* $z>c_0^{-1}$ *with the following significance: Supposing that* $\phi$ *is from this summand, and supposing that each* $p\in\Theta$ *version of its corresponding version* $\mu_+$ *is written as done in (7.40), then* $\int_{\Sigma}\langle\phi,\maltese\phi\rangle$ *differs from*

$$z\,r^{4/3}\sum_{p\in\Theta}\alpha_{1p}^{2/3}\,|a_{-(p)}|^2$$

*by at most* $c_0\sum_{p\in\Theta}(|a_{0(p)}|^2+\sup_{|z|<r_0}|\mathfrak{r}_p|^2))$

(7.42)

By way of a reminder: Each $p\in\Theta$ version of $\alpha_{1p}$ in the third bullet of (7.41) is a positive number that is defined from the holomorphic, quadratic differential $q_1$ by using Lemma 3.2 to write $q_1$ near p as $\alpha_{1p}z\,(dz)^2$ with $\alpha_{1p}$ being positive. (Remember that $q$ is $r^2q_1$). As for the number $z$ in that third bullet: This number is obtained from the $\mathbb{C}$ integral of the expression in (7.24) by writing that integral as $z\,\alpha_p^{2/3}$ with $\alpha_p$ being $r^2\alpha_{1p}$.

To elaborate on the third bullet's error bound: For any given $p\in\Theta$, the lower order terms in (7.41)'s expression for $\mu_+$ when $\phi$ is from K will certainly contribute to the expression in (7.5) when computing the $D_1$ integral of $\langle\phi,\maltese\phi\rangle$. None-the-less, the size of their contribution to (7.5) is at most

$$c_0\,(|a_{0(p)}|^2+\sup_{|z|<r_0}|\mathfrak{r}_p|^2)\,(|\Xi||\langle i\sigma_3[\varphi^*,\varphi]\rangle|+|\varphi|^2|\Xi|^2+\tfrac{1}{4}\,|d\Xi|^2)$$

(7.43)

Since neither $\Xi$ nor $\frac{d}{ds}\Xi$ change sign and their signs are opposite (see Lemma 5.2), the identity in (7.9) can be used to bound the expression in (7.43) by

$$c_0\,(|a_{0(p)}|^2+\sup_{|z|<r_0}|\mathfrak{r}_p|^2)\,(-\Xi\tfrac{d}{ds}\Xi\,+|\varphi|^2|\Xi|^2+\tfrac{1}{4}\,|d\Xi|^2)$$

(7.44)

The integral of this expression over $D_1$ is no greater than its integral over $\mathbb{C}$ which is, according to Lemma 7.3, at most $c_0\sup_{|z|<r_0}(|a_{0(p)}|^2+|\mathfrak{r}_p|^2)$.

As explained directly, the size of the numbers a, $a_{-(p)}$, $a_{0(p)}$ and $\sup_{|z|<r_0}|\mathfrak{r}_p|$ that appear in (7.42) relate as follows to the $\mathbb{L}$-norm of $\phi$

- *If* $\phi$ *is from the* $\mathbb{C}$ *summand and* $a\in\mathbb{C}$ *is the value of* $\phi$*'s version of* $\mu_+$*, then*

$$c_0^{-1}|a|r\leq\,\|\phi\|_{\mathbb{L}}\,\leq\,c_0|a|r.$$

- *If* $\phi$ *is from the* K *summand and if* $\phi$*'s version of* $\mu_+$ *is depicted near each* $p\in\Theta$ *as* $\mu_+=\frac{a_{-(p)}}{\bar{z}}+a_0+\mathfrak{w}_p$ *with* $a_{-(p)}$ *and* $a_{0(p)}$ *being complex numbers and with* $\mathfrak{w}_p$ *denoting an anti-holomorphic function on the* $|z|<2r_0$ *disk in* $\mathbb{C}$ *that vanishes at* $z=0$*, then*

- ○ $c_0^{-1}(\sum_{p \in \Theta}|a_{-(p)}|^2)^{1/2}\, r \leq \|\phi\|_{\mathbb{L}} \leq c_0\,(\sum_{p \in \Theta}|a_{-(p)}|^2)^{1/2}\, r$.

- ○ $\sum_{p \in \Theta}\,(|a_{0(p)}| + \sup_{|z|<r_0}|\mathfrak{r}_p|) \leq c_0\|\phi\|_{\mathbb{L}}$.

$$(7.45)$$

To explain: The first bullet follows from what is said by Lemma 5.2 using the fact that the $\Sigma$–integral of $|\phi|^2$ is between $c_0^{-1}r^2$ and $c_0 r^2$. The second bullet's assertions follow from the preceding and the fact that the sum $\sum_{p \in \Theta}|a_{-(p)}|^2$ when viewed as a function on the $\|\phi\|_{\mathbb{L}} = r$ ball in K has a $c_0^{-1}$ lower bound.

The upcoming Part 8 of the proof explains why Proposition 7.2 follows from what is said by (7.42) with regards to $\int_{\Sigma} \langle \phi, \yen\phi \rangle$ using the observations from (7.44) about $\mathbb{L}$-norms and some of the finite dimensional, linear perturbation theory in the next part of the proof.

*Part 7:* What follows are the finite dimensional perturbation theory lemmas that are needed to finish the proof of Proposition 7.2. (Infinite dimensional analogs of these lemmas lie behind the analysis that led to (6.64).)

**Lemma 7.4**: *Having fixed* $k \in \{1, 2, ...\}$, *there exists* $\kappa > 1$ *with the following significance: Let* $\mathfrak{S}$ *and* $\mathfrak{s}$ *denote a pair of symmetric (or Hermitian)* $k \times k$ *matrices. Then the eigenvalues of their sum,* $\mathfrak{S} + \mathfrak{s}$, *differ from those of* $\mathfrak{S}$ *by at most* $\kappa|\mathfrak{s}|$.

The proof of this lemma is given momentarily.

**Lemma 7.5**: *Having fixed* $k \in \{2, 3, ...\}$, *there exists* $\kappa > 1$ *with the following significance: Let* M *denote a symmetric (or Hermitian),* $k \times k$ *matrix that can be written in block diagonal form as*

$$M = \begin{pmatrix} \mathfrak{s} & \boldsymbol{\ell} \\ \boldsymbol{\ell}^\dagger & S \end{pmatrix}$$

*with* $\mathfrak{s}$ *denoting a symmetric (or Hermitian),* d *dimensional matrix (for* $d \in \{1, ..., k\text{-}1\}$) *with the norm of its eigenvectors being at most* $\frac{1}{\kappa}$; *with* S *denoting a symmetric (or Hermitian),* (k-d) × (k-d) *matrix whose eigenvalues have norm greater than 1; and with* $\boldsymbol{\ell}$ *denoting a* d × (k-d) *dimensional matrix with norm at most* $\frac{1}{\kappa^2}$. *Then* M *has* d *eigenvalues (counted with multiplicity) that have distance at most* $\kappa|\ell|^2$ *from the set of eigenvalues of* $\mathfrak{s}$, *and also* (n - d) *eigenvalues (counted with multiplicity) that have distance at most* $\kappa|\ell|^2$ *from the set of eigenvalues of* S.

These lemmas will be used again in the next two subsections. (Likely you noticed that infinite dimensional analogs of these lemmas lie behind the analysis that led to (6.64).)

***Proof of Lemma 7.4***: Let $\mathbf{x}$ denote an eigenvector of $\mathfrak{S} + \mathfrak{s}$ and let $\lambda$ denote its eigenvalue. Suppose that $\varepsilon > 0$ and that there is no eigenvector of $\mathfrak{S}$ with distance less than $\varepsilon$ from $\lambda$. Let I denote the identity matrix. The eigenvalue equation $(\mathfrak{S} + \mathfrak{s})\mathbf{x} = \lambda\mathbf{x}$ implies directly that $|(\mathfrak{S} - \lambda I)\mathbf{x}| \leq c_0 |\mathfrak{S}||\mathbf{x}|$. Meanwhile, $|(\mathfrak{S} - \lambda I)\mathbf{x}| \geq \varepsilon|\mathbf{x}|$ because $\mathbf{x}$ can be written as a linear combination of eigenvectors of $\mathfrak{S}$ and none have eigenvalue within $\varepsilon$ of $\lambda$. These two inequalities are inconsistent if $\varepsilon > c_0|\mathfrak{S}|$.

***Proof of Lemma 7.5***: Suppose first that v is an eigenvalue of M with eigenvalue $\lambda$. Write the vector v as

$$v = \begin{pmatrix} x \\ \mathbf{X} \end{pmatrix}$$

(7.46)

with x the first d components and with $\mathbf{X}$ denoting the last (k-d) components. The assertion that v is an eigenvector of M with eigenvalue $\lambda$ is equivalent to the two conditions

- $\mathfrak{s}x + \boldsymbol{\ell}\mathbf{X} = \lambda x$ ,
- $S\mathbf{X} + \boldsymbol{\ell}^\dagger x = \lambda\mathbf{X}$ .

(7.47)

Supposing that $\lambda$ is not equal to an eigenvalue of $\mathfrak{s}$, then the first bullet in (7.47) can be used to write x as

$$x = -\frac{1}{\mathfrak{s} - \lambda I} \, \boldsymbol{\ell}\mathbf{X}$$

(7.48)

with I denoting (here and below) the identity matrix of the appropriate dimension (d × d in this case). The second bullet of (7.47) can be written using (7.48)

$$(S - \lambda I)\mathbf{X} = \boldsymbol{\ell}^\dagger \frac{1}{\mathfrak{s} - \lambda I}(\boldsymbol{\ell}\mathbf{X}) \, .$$

(7.49)

With (7.49) in hand, fix R >1 for the moment, and suppose for the sake of continuing the argument that $|\mathfrak{s} - \lambda I| > R|\boldsymbol{\ell}|^2$ If this is the case, and if $R > c_0$, then Lemma 7.4 can be invoked with $\mathfrak{S} = S$ and $\mathfrak{s} = \boldsymbol{\ell}^\dagger \frac{1}{\mathfrak{s} - \lambda I} \boldsymbol{\ell}$ to see that $|\lambda|$ has to be greater than $\frac{1}{2}$ (remember that the norms of the eigenvalues of S are greater than 1). But if $|\lambda| > \frac{1}{2}$, then the norm of $\boldsymbol{\ell}^\dagger \frac{1}{\mathfrak{s} - \lambda I} \boldsymbol{\ell}$ is bounded by $c_0|\boldsymbol{\ell}|^2$ and so Lemma 7.4 requires that $\lambda$ is within $c_0|\boldsymbol{\ell}|^2$ of an eigenvalue of S.

With regards to eigenspace dimensions: Suppose for the sake of argument that M has fewer than d linearly independent eigenvectors whose eigenvalue is within $c_0|\boldsymbol{\ell}|^2$ of an eigenvalue of $\mathfrak{s}$. Nonsense will be derived from this assumption by using the fact that if v is a non-zero, linear combination of eigenvectors of M whose eigenvalues have norm greater than $\frac{1}{2}$, then $|Mv| > \frac{1}{2}|v|$. To derive the desired nonsense if M has less than d eigenvalues (counted with multiplicity) with small norm, note first that if $|\boldsymbol{\ell}|$ and $|\mathfrak{s}|$ are less than $c_0^{-1}$, then M must have more than n - d linearly independent eigenvectors whose eigenvalues have norm greater than $\frac{1}{2}$. As a consequence, there is a non-zero, linear combination of these eigenvectors (call it v) whose **X** component in (7.46) is zero. Then $|Mv|$ for this version of v obeys the bound $|Mv| \leq c_0(|\mathfrak{s}| + |\boldsymbol{\ell}|)|v|$ which is less than $\frac{1}{4}|v|$ if both $|\mathfrak{s}|$ and $|\boldsymbol{\ell}|$ are less than $c_0^{-1}$. This is nonsense because it violates the $|Mv| > \frac{1}{2}|v|$ condition.

Much the same sort of argument (which is left to the reader) will prove that M must have exactly d eigenvalues with distance at most $c_0|\boldsymbol{\ell}|^2$ from the set of eigenvalues of $\mathfrak{s}$ if both $|\mathfrak{s}|$ and $|\boldsymbol{\ell}|$ are less than $c_0^{-1}$.

*Part 8:* To start this last part of the proof, write $\Pi ¥ \Pi$ as a block diagonal matrix with respect to the $\mathbb{C} \oplus K$ decomposition of the $\mathbb{S}_{SL}^{\perp}$ kernel of $\mathfrak{D}_0$ as below:

$$\begin{pmatrix} \mathfrak{h} & \ell \\ \ell^{\dagger} & H \end{pmatrix}$$

$$(7.50)$$

with $\mathfrak{h}$ denoting a symmetric matrix mapping the $\mathbb{C}$ summand to itself, with H denoting a symmetric matrix mapping the K summand to itself and with $\ell$ denoting a matrix that maps the K summand to the $\mathbb{C}$ summand. According to the top bullet in (7.42) and what is said subsequently in (7.45) about norms, the matrix $\mathfrak{h}$ has two nearly identical eigenvalues which are between $c_0^{-1}r^{-2}$ and $c_0 r^{-2}$ (see Lemma 7.4 with regards to the $c_0 e^{-r/c_0}$ ambiguity that is asserted by the top bullet in (7.42).) Meanwhile: According to the third bullet in (7.42) and what is said in (7.44) about norms, the matrix H is positive definite with all of its 6g - 8 eigenvalues in the interval between $c_0^{-1}r^{-2/3}$ and $c_0 r^{-2/3}$.

As for the matrix $\ell$: According to (7.42) and (7.44), the norm of this matrix is bounded by $c_0 r^{-2}$. This being the case, it then follows from Lemma 7.5 that the eigenvalues of (7.50) differ by at most $c_0 r^{-4}$ from those of the diagonal matrix

$$\begin{pmatrix} \mathfrak{h} & 0 \\ 0 & H \end{pmatrix}.$$

$$(7.51)$$

More to the point: Lemma 7.5 implies that the matrix in (7.50) has 2 eigenvalues between $c_0^{-1}r^{-2}$ and $c_0 r^{-2}$, and then 6g - 8 eigenvalues between $c_0^{-1}r^{-2/3}$ and $c_0 r^{-2/3}$. And, this is the content of Proposition 7.2.

### b) The proof of Proposition 7.1

This section contains the proof of Proposition 7.1. The proof has two parts. The first part analyzes the effect of adding the endomorphism $\phi \to m\Pi(\gamma_t[\chi_{\ddagger}A^1, \phi])$ to Proposition 7.2's endomorphism. The second part of the proof analyzes the effect of the adding the endomorphism $\phi \to \Pi(\mathfrak{l}(\phi))$ to the endomorphism $m\Pi\big(c^1\rho_t[\sigma_3, \phi] + (\phi)_c + \gamma_t[\chi_{\ddagger}A^1, \phi]\big)$.

*Part 1*:  Although the endomorphism

$$\phi \to \Pi(\gamma_t[\chi_{\ddagger}A^1, \phi])$$

$$(7.52)$$

looks to be formally of the same size (with regards to eigenvalues) as Proposition 7.2's endomorphism, it's norm is far, far smaller. This is said formally by the next lemma.

**Lemma 7.6**:  *There exists* $\kappa > 1$ *with the following significance: Fix* $m < \frac{1}{\kappa}$ *and then fix* $r = \frac{n}{m} > \kappa$. *The norm of the corresponding version of the endomorphism that is depicted in* (7.52) *is at most* $\kappa \, e^{-r/\kappa}$.

It follows from this lemma that the addition of (7.52)'s endomorphism to Proposition 7.2's endomorphism does not change the conclusions of Proposition 7.2 if the latter's version of $\kappa$ is increased by a $c_0$ factor (see Lemmas 7.4 and 7.5). This is to say that the conclusions of Proposition 7.2 hold for the endomorphism

$$\phi \to \Pi\big(c^1\rho_t[\sigma_3, \phi] + \gamma_t[\chi_{\ddagger}A^1, \phi] + (\phi)_c\big)$$

$$(7.53)$$

when the stated version of $\kappa$ is increased by a $c_0$ multiplicative factor.

*Proof of Lemma 7.6*:  The proof has five steps.

Step 1:  The first thing to keep in mind is that $\gamma_t[\chi_{\ddagger}A^1, \cdot\,]$ is supported only on the holomorphic coordinate chart disks around the points in $\Theta$. Having chosen such a point, then an element from the $\mathbb{S}_{SL}{}^{\perp}$ kernel of $\mathfrak{D}_0$ with components $(\mathfrak{b}_1 e^1 + \mathfrak{b}_2 e^2, c = c_1 e^1 + c_2 e^2)$ (when written using an orthonormal frame for $T^*\Sigma|_D$) has the following form:

- $\mathfrak{b}_1{}' = -\chi_{\ddagger}[A^1, c_2]$ ,
- $\mathfrak{b}_2{}' = \chi_{\ddagger}[A^1, c_1]$ ,
- $c_1{}' = -\chi_{\ddagger}[A^1, \mathfrak{b}_2]$ ,
- $c_2{}' = +\chi_{\ddagger}[A^1, \mathfrak{b}_1]$ .

$$(7.54)$$

As was the case with Proposition 7.2's endomorphism, it is convenient to write the endomorphism in terms of $\mathfrak{b} = \mathfrak{b}_1 + i\mathfrak{b}_2$ and $c = \mathfrak{c}_1 - i\mathfrak{c}_2$. This what is depicted below (with $\mathfrak{b}^*$ denoting $\mathfrak{b}_1 - i\mathfrak{b}_2$ and $c^*$ denoting $\mathfrak{c}_1 + i\mathfrak{c}_2$):

- $\mathfrak{b}' = i\,\chi_{\ddagger}[A^1, c^*]$ ,
- $c' = -i\,\chi_{\ddagger}[A^1, \mathfrak{b}^*]$ .

(7.55)

The associated quadratic form $(\phi, \phi') \rightarrow \langle \phi, \gamma_t[\chi_{\ddagger}A^1, \phi']\rangle$ when $\phi$ is determined by $(\mathfrak{b}, c)$ and $\phi'$ by $(\mathfrak{b}', c')$ is as follows:

$$i\chi_{\ddagger}\big(\langle A^1[c^*, \mathfrak{b}'^*]\rangle + \langle A^1[c'^*, \mathfrak{b}^*]\rangle\big) .$$

(7.56)

The real part of the integral of (7.56) over $\Sigma$ gives the pairing between $\phi$ and $\Pi(\gamma_t[\chi_{\ddagger}A^1, \phi'])$. The assertion of the lemma follows from the assertion that the bound

$$\Big|\int_{\Sigma} i\chi_{\ddagger}\big(\langle A^1[c^*, \mathfrak{b}'^*]\rangle + \langle A^1[c'^*, \mathfrak{b}^*]\rangle\big)\Big| \leq c_0 e^{-r/c_0}\,\|\phi\|_{\mathbb{L}}\|\phi'\|_{\mathbb{L}}$$

(7.57)

holds whenever $\phi$ and $\phi'$ are from the $\mathbb{S}_{SL}^{\perp}$ kernel of $\mathfrak{D}_0$.

With regards to this integral: The appearance in the integrand of the function $\chi_{\ddagger}$ implies that the integrand is support in the union of the holomorphic coordinate chart disks around the points in $\Theta$. As a consequence, it is sufficient for purposes of the proof to verify that the norm of the contribution of each such disk to the left hand side integral in (7.57) is no greater than $c_0 e^{-r/c_0}$.

Step 2: To see what (7.56) looks like on the disk D about a given point in $\Theta$ (this is where the coordinate z obeys $|z| < r_0$), write $\mathfrak{b}$ and $c$ on that disk using anti-meromorphic Laurent series $P$ and $Q$ with an order 1 pole at z = 0 (if any) as done previously:

$$\mathfrak{b} = P\,\phi^* - Q\,\Xi\,[\tfrac{i}{2}\sigma_3, \phi^*] \ \ and \ \ c = -Q\,(\tfrac{i}{2}\tfrac{\bar{z}}{s}\tfrac{d}{ds}\Xi)\,\sigma_3 ;$$

(7.58)

Write $\mathfrak{b}'$ and $c'$ analogously using anti-meromorphic $P'$ and $Q'$. (As per Lemma 5.2, replacing $\mu_-$ by $Q\Xi$ and then $\varsigma_3$ by $-\tfrac{i}{2}Q\,(\tfrac{\bar{z}}{s}\tfrac{d}{ds}\Xi)$ incurrs an error of size at most $c_0 e^{-r/c_0}\,\|\phi\|_{\mathbb{L}}\|\phi'\|_{\mathbb{L}}$ with regards to the integral on the left hand side of (7.57).) Also, write $A^1$ on this disk as

$$A^1 = -\tfrac{1}{3}\,(z\phi + \bar{z}\phi^*) .$$

(7.59)

Use the preceding expressions to write $\chi_{\ddagger}\langle A^1[c^*, \mathfrak{b}'^*]\rangle$ as follows:

$$\chi_\ddagger \langle A^1[c^*, \not{b}'^*] \rangle = -\tfrac{i}{6} \chi_\ddagger \, \overline{Q} \, \tfrac{z}{s} \, \tfrac{d}{ds} \Xi \, \langle (z\varphi + \overline{z}\varphi^*), [\sigma_3, (\overline{P}' \, \varphi + \overline{Q}' \Xi \, [\tfrac{i}{2}\sigma_3, \varphi])] \rangle \, ;$$

(7.60)

and then expand the expression above to obtain the identity

$$\chi_\ddagger \langle A^1[c^*, \not{b}'^*] \rangle = -\tfrac{1}{3} \chi_\ddagger \, \overline{Q} \, \overline{Q}' \, \tfrac{1}{s} \Xi \tfrac{d}{ds} \Xi \; z^2 \langle \varphi^2 \rangle + \tfrac{1}{6} \chi_\ddagger s \tfrac{d}{ds} \Xi \, (\overline{Q} \, \overline{P}' \langle i\sigma_3[\varphi^*, \varphi] \rangle - \overline{Q} \, \overline{Q}' 2\Xi |\varphi|^2)$$

(7.61)

The next steps analyze the integral of the right hand side of (7.61).

Step 3: As explained directly, the left most term on the right hand side of (7.60) integrates to zero on the constant $|z|$ circles in D (this is the term with $\langle \varphi^2 \rangle$). Thus, its integral over the $|z| < r_0$ part of $\Sigma$ is identically zero. The reason is this: Because $\langle \varphi^2 \rangle$ is proportional to z and because $\overline{Q}$ and $\overline{Q}'$ are meromorphic with a pole of order at most 1 at $z = 0$, that left most term is the product of a holomorphic function of z vanishing at $z = 0$ with the radially symmetric function $\chi_\ddagger \tfrac{1}{s} \Xi \tfrac{d}{ds} \Xi$.

With regards to the other integrals: Replacing $\langle i\sigma_3[\varphi^*, \varphi] \rangle$ and $|\varphi|^2$ by the expressions in (7.7) incurs an error of size at most $c_0 e^{-r/c_0} \|\phi\|_{\mathbb{L}} \|\phi'\|_{\mathbb{L}}$ when computing the integral on the left in (7.57). With those replacements made here and henceforth, then the term with the factor $(\overline{Q} \, \overline{P}' \langle i\sigma_3[\varphi^*, \varphi] \rangle - \overline{Q} \, \overline{Q}' 2\Xi |\varphi|^2)$ in (7.60) integrates to zero on the constant $|z|$ circles except for contributions to that integral from the product of the Laurent series terms described below:

- *The product of the constant terms in the Laurent series for* $Q$ *and those for* $P'$ *and* $Q'$.
- *The product of the first order pole term in the Laurent series for* $Q$ *and the order z terms in the Laurent series for* $P'$ *and* $Q'$; *and also the product of the respective first order pole term in the Laurent series for* $P'$ *and* $Q'$ *with the order z term in the Laurent series for* $Q$.

(7.62)

(Any other product of respective Laurent series terms will have either a positive power of z or a negative power of z.)

The case described by the top bullet in (7.62) is analyzed in Step 4 of the proof and the cases described in the second bullet of (7.62) are analyzed in Step 5 of the proof.

Step 4: To consider the case in the top bullet of (7.62), let a denote the zero order term (it is a complex number) in the Laurent expansion for $P$ (and thus for $Q$) and let $a'$ denote the zero order term in the Laurent expansion for $P'$ (and thus for $Q'$). The contribution of the product of these respective zero'th order terms to

$$\tfrac{1}{6} \chi_\ddagger s \tfrac{d}{ds} \Xi \, (\overline{Q} \, \overline{P}' \langle i\sigma_3[\varphi^*, \varphi] \rangle - \overline{Q} \, \overline{Q}' 2\Xi |\varphi|^2)$$

(7.63)

(which is from (7.61)) is this:

$$\frac{1}{6} \chi_{\ddagger} s \, \overline{a} \, \overline{a}' \, \frac{d}{ds} \Xi \left( \langle i\sigma_3[\varphi^*, \varphi] \rangle - 2\Xi |\varphi|^2 \right) .$$

(7.64)

To see about the $|z| < r_0$ integral of the preceding function: Use the fact that $\Xi$ obeys (7.14) on $\mathbb{C}$ (which is (5.12)) to write the function $\langle i\sigma_3[\varphi^*, \varphi] \rangle - 2\Xi |\varphi|^2$ as

$$\langle i\sigma_3[\varphi^*, \varphi] \rangle - 2\Xi |\varphi|^2 = -\frac{1}{2} \frac{1}{s} \frac{d}{ds} (s \frac{d}{ds} \Xi) .$$

(7.65)

Doing this identifies the function depicted in (7.64) with the function depicted below:

$$-\frac{1}{24} \chi_{\ddagger} \, \overline{a} \, \overline{a}' \, \frac{1}{s} \frac{d}{ds} \left( (s \frac{d}{ds} \Xi)^2 \right) .$$

(7.66)

As explained directly, the norm of the $|z| < r_0$ integral of the function in (7.63) is bounded by $c_0 e^{-r/c_0} |a||a'|$. Indeed, the small size of this integral follows from three facts: First, the integral of the function in (7.66) is $2\pi$ times its integral with respect to s on $[0, \infty)$ using the measure s ds. Second, integration by parts (note that the function has compact support where $|z| < r_0$) identifies the latter integral with the integral on the domain $[0, \infty)$ using the Euclidean measure ds of the function depicted below:

$$\frac{1}{24} \, \overline{a} \, \overline{a}' \, \frac{d}{ds} \chi_{\ddagger} \left( s \frac{d}{ds} \Xi \right)^2 .$$

(7.67)

Third, the function $|s \frac{d}{ds} \Xi|$ whose square appears in (7.66) is bounded by $c_0 \, e^{-r/c_0}$ on the support of the function $|\frac{d}{ds} \chi_{\ddagger}|$ whose integral is bounded by $c_0$.

<u>Step 5</u>: This step considers the contribution to the $|z| < r_0$ integral of the function depicted in (7.56) from the product of the Laurent series terms that are described by the second bullet of (7.62). To this end, write the Laurent series for P and Q on a given $p \in \Theta$ version of D as

$$P = a_{-} \frac{1}{\overline{z}} + a_0 + a_1 \, \overline{z} + \cdots \quad \textit{and} \quad Q = a_{-} \frac{3}{\overline{z}} + a_0 + a_1 \frac{3}{5} \, \overline{z} + \cdots$$

(7.68)

Write P′ and Q′ analogously with primed coefficients. The relevant contribution to the function $\chi_{\ddagger} \langle A^1[c^*, \hat{b}'^*] \rangle$ is this:

$$-\frac{i}{2} \, \overline{a}_{-} \, \overline{a}_1' \, \chi_{\ddagger} \frac{z}{s} \frac{d}{ds} \Xi \langle (z\varphi + \overline{z}\varphi^*), [\sigma_3, (\varphi + \frac{3}{5} \Xi [\frac{i}{2}\sigma_3, \varphi])] \rangle$$
$$-\frac{i}{10} \, \overline{a}_{-}' \, \overline{a}_1 \, \chi_{\ddagger} \frac{z}{s} \frac{d}{ds} \Xi \langle (z\varphi + \overline{z}\varphi^*), [\sigma_3, (\varphi + 3 \, \Xi [\frac{i}{2}\sigma_3, \varphi])] \rangle .$$

(7.69)

Add to this the corresponding function $\chi_{\ddagger}\langle \Lambda^1[c'^*, \hat{b}^*]\rangle$ to obtain

$$-\frac{i}{2}\left(\bar{a}_- \,\bar{a}_1{}' + \bar{a}_-{}' \,\bar{a}_1\right)\chi_{\ddagger}\frac{z}{s}\frac{d}{ds}\Xi\left\langle(z\varphi + \bar{z}\varphi^*), [\sigma_3, (\tfrac{6}{5}\varphi + \tfrac{6}{5}\Xi\,[\tfrac{i}{2}\sigma_3, \varphi])]\right\rangle$$

(7.70)

which integrates to within $c_0\,e^{-r/c_0}$ zero over the $|z| < r_0$ part of $\mathbb{C}$ (it has compact support there) since it is equal to

$$\frac{3i}{10}\left(\bar{a}_- \,\bar{a}_1{}' + \bar{a}_-{}' \,\bar{a}_1\right)\chi_{\ddagger}\frac{z}{s}\frac{d}{ds}\Xi\left\langle(z\varphi + \bar{z}\varphi^*), [\sigma_3, (\varphi + \Xi\,[\tfrac{i}{2}\sigma_3, \varphi])]\right\rangle$$

(7.71)

which is a constant multiple of the $a = a' = 1$ version of the expresion in (7.64).

*Part 2*: This part finishes the proof of Proposition 7.1. The last key input to the proof is Lemma 7.7 below which supplies a priori bounds for the norm of various components of the endomorphism $\phi \to \Pi(\mathfrak{l}(\phi))$ of the $\mathbb{S}_{SL}^{\perp}$ kernel of $\mathfrak{D}_0$. With regards to the statement of the lemma: The lemma writes this endomorphism as block diagonal matrix with respect to the $\mathbb{C} \oplus K$ decomposition of the $\mathbb{S}_{SL}^{\perp}$ kernel of $\mathfrak{D}_0$, thus as

$$\begin{pmatrix} s & t \\ t^{\dagger} & S \end{pmatrix}$$

(7.72)

with $s$ mapping the $\mathbb{C}$ summand to itself and $S$ mapping the $K$ summand to itself.

**Lemma 7.7**: *There exists $\kappa > 1$ with the following significance: Fix $m < \frac{1}{\kappa}$ and then fix $r = \frac{n}{m} > \kappa$. Use (6.58) to define the corresponding version of $\mathfrak{l}(\cdot)$ and then write the resulting endomorphism $\phi \to \Pi(\phi)$ of the $\mathbb{S}_{SL}^{\perp}$ kernel of $\mathfrak{D}_0$ in block diagonal form as done in (7.69). The norms of the component matrices $s$, $S$ and $t$ obey the bounds below:*

- $|s| \leq \kappa(m^2 r^{-2/3} + e^{-r/\kappa})\,r^{-2}$.
- $|S| \leq \kappa(m^2 r^{-2/3} + e^{-r/\kappa})\,r^{-2/3}$.
- $|t| \leq \kappa(m^2 r^{-2/3} + e^{-r/\kappa})\,r^{-4/3}$

This lemma is proved momentarily.

To finish the proof of Proposition 7.1: Add Lemma 7.6's contribution and the matrix in (7.72) to $m$ times the matrix in (7.50) to obtain a block diagonal representation of Proposition 7.1's endomorphism $\Pi\mathfrak{k}\Pi$. This will have the form below

$$\begin{pmatrix} m\mathfrak{h} + s' & m\ell + t' \\ m\ell^{\dagger} + t'^{\dagger} & m\mathrm{H} + S' \end{pmatrix}$$

(7.73)

with $s'$, $t'$ and $S'$ denoting the sum of their unprimed counterparts in (7.72) and the corresponding entries from Lemma 7.6's matrix. In particular, the primed and unprimed versions differ by at most $c_0 e^{-r/c_0}$. With (7.73) understood and with Lemma 7.7's bounds in hand, then the conclusions of Proposition 7.1 follow directly with an appeal to Lemma 7.4 (for the effect of adding $s'$ to $m\mathfrak{h}$ and $S'$ to $mH$) and then Lemma 7.5 for the effect of adding the off-diagonal blocks with $m\ell + t'$ and its adjoint. (And keep in mind that the norm of $\ell$ is bounded by $c_0 r^{-2}$.)

*Proof of Lemma 7.7:* The task at hand is to bound the norm of the bilinear form

$$(\phi', \phi) \to \int_\Sigma \langle \phi', \mathfrak{l}(\phi) \rangle .$$

(7.74)

To this end: The depiction of $\mathfrak{l}(\cdot)$ in (6.59) and what is said surrounding (5.32)-(5.35) imply the pointwise bound below for the norm of $\langle \phi', \mathfrak{l}(\phi) \rangle$:

$$|\langle \phi', \mathfrak{l}(\phi) \rangle| \le c_0 \, (m^2 r^{-2/3} \, e^{-r\mathfrak{d}^{3/2}/c_0} + e^{-r/c_0}) |\phi||\phi'| \; + c_0 m^2 r^{-2/3} \, |[\hat{\tau}, \phi]||[\hat{\tau}, \phi']|$$

(7.75)

It follows from this inequality and from what is said by Lemmas 5.2 and 3.3 that the contribution to the $\Sigma$ integral of (7.75) from the $\mathfrak{d} > r_0$ part of $\Sigma$ is no greater than $c_0 e^{-r/c_0} \|\phi\|_{\mathbb{L}} \|\phi'\|_{\mathbb{L}}$. This follows in particular from the third bullet of Lemma 5.2 and the implication from Lemma 3.3 that

$$|[\hat{\tau}, \varphi]| \le r^{2/3} \, e^{-r\mathfrak{d}^{3/2}/c_0} .$$

(7.76)

With the preceding understood, fix $p \in \Theta$ to see about bounding the integral of the right hand side of (7.75) on the holomorphic coordinate chart disk D centered at p. At the expense of an error of size at most $c_0 e^{-r/c_0} \|\phi\|_{\mathbb{L}} \|\phi'\|_{\mathbb{L}}$ when computing the $\Sigma$ integral of the right hand side of (7.75), the complex endomorphism valued 1-forms $(\acute{b}, c)$ that determine $\phi$ can be written as in (7.58) in terms of the Laurent series P and Q of the variable $\bar{z}$, and likewise $\phi'$ can be written using primed versions of P and Q.

Three cases arise with the first being the case when both $\phi$ and $\phi'$ are from the $\mathbb{C}$ summand of the $\mathbb{S}_{SL}{}^\perp$ kernel of $\mathfrak{D}_0$. In this case, both P and P′ are constant. Denote their constant value by a and a′. In this case, Lemma 5.2 with (7.76) lead to the bounds below for the norms of $\phi$ and $[\hat{\tau}, \phi]$ on the disk D:

$$|\phi| \le c_0 \max(r \, \mathfrak{d}^{1/2}, r^{2/3}) \, |a| \; \textit{and} \; \; |[\hat{\tau}, \phi]| \le c_0 r^{2/3} \, e^{-r|z|^{3/2}/c_0} |a| .$$

(7.77)

The analogous bounds are obeyed in this case for $|\phi'|$. The preceding bounds with what is said about $\Xi$ in Lemma 5.2 lead directly to a

$$c_0(m^2 r^{-2/3} + e^{-r/c_0})|a||a'|$$

(7.78)

bound for the D integral of the expression on the right in (7.71).

The next case has $\phi$ from the K summand of the $\mathbb{S}_{SL}{}^{\perp}$ kernel of $\mathfrak{D}_0$ and $\phi'$ from the $\mathbb{C}$ summand (or vice-versa). In the case when $\phi$ is from the K summand, write $P$ in an abbreviated version of (7.68) as done below:

$$P = a_- \frac{1}{\bar{z}} + \mathfrak{w}$$

(7.79)

where $\mathfrak{w}$ denotes here an anti-holomorphic function on the $|z| < 2r_0$ disk. The norms of $\phi$ and $[\hat{\tau}, \phi]$ on the disk D obey these bounds:

- $|\phi| \leq c_0 \min(r\,\mathfrak{d}^{-1/2}, r^{4/3})\,|a_-| + c_0 \max(r\,\mathfrak{d}^{1/2}, r^{2/3}) \sup_{|z|<r_0}|\mathfrak{w}|$.
- $|[\hat{\tau}, \phi]| \leq c_0 e^{-r|z|^{3/2}/c_0}(r^{4/3}|a_-| + \sup_{|z|<r_0}|\mathfrak{w}|)$.

(7.80)

The preceding bounds with those in (7.77) lead to the following bound for the integral of the right hand side of (7.71) over the disk D:

$$c_0(m^2 r^{-2/3} + e^{-r/c_0})(r^{2/3}|a_-| + \sup_{|z|<r_0}|\mathfrak{w}|)|a'|.$$

(7.81)

An analogous bound holds when $\phi$ is from $\mathbb{C}$ and $\phi'$ is from K (take the prime off of $a'$ and put one on $a_-$ and another on $\mathfrak{w}$).

The final case to consider has both $\phi$ and $\phi'$ being from K. In this case, an appeal to (7.80) and its primed analog leads to the bound below for the integral of the right hand side of (7.76) over D:

$$c_0(m^2 r^{-2/3} + e^{-r/c_0})\big(r^{4/3}|a_-||a_-'| + r^{2/3}\,(\sup_{|z|<r_0}|\mathfrak{w}|\,|a_-'| + \sup_{|z|<r_0}|\mathfrak{w}'||a_-|)$$
$$+ \sup_{|z|<r_0}|\mathfrak{w}|\cdot\sup_{|z|<r_0}|\mathfrak{w}'|\big).$$

(7.82)

These three cases now lead to Lemma 7.7's bound on the norm of the block diagonal componens $s$, $S$ and $\ell$ of $\Pi(\cdot)\Pi$ because, as per (7.45), the norms $|a|$ and $|a'|$ are at most $c_0 r^{-1}$ times the respective $\mathbb{L}$ norms of $\phi$ and $\phi'$, as are $|a_-| + \sup_D|\mathfrak{w}|$ and $|a_-'| + \sup_D|\mathfrak{w}'|$.

**c) The endomorphism** $\phi \to \Pi\mathfrak{k}(\phi + \mathfrak{T}_\lambda(\phi))$ **of the** $\mathbb{S}_{SL}^\perp$ **kernel of** $\mathfrak{D}_0$ **when** $\lambda$ **is small**

For any given real number $\lambda$, the what is depicted in (6.64) with $\mathfrak{T}_\lambda(\cdot)$ defined in (6.62) defines a symmetric endomorphism of the $\mathbb{S}_{SL}^\perp$ kernel of $\mathfrak{D}_0$. A real number $\eta$ is said below to be a *linear eigenvalue* of the endomorphism $\phi \to \Pi\mathfrak{k}(\phi + \mathfrak{T}_\lambda(\phi))$ when there exists a non-zero element (denoted below by $\phi$) from the $\mathbb{S}_{SL}^\perp$ kernel of $\mathfrak{D}_0$ such that

$$\Pi\mathfrak{k}(\phi + \mathfrak{T}_\lambda(\phi)) = \eta\phi \ .$$

(7.83)

(Of ultimate interest are the norms for the numbers $\lambda$ for which (7.83) has $\lambda$ itself as an linear eigenvalue.) In any event, the following lemma says what is needed about the linear eigenvalues when $|\lambda|$ is small.

**Lemma 7.8**: *There exists* $\kappa > 1$ *with the following significance: Fix* $m < \frac{1}{\kappa}$ *and then fix* $r = \frac{n}{m} > \kappa$. *Also fix* $\lambda \in \mathbb{R}$ *with* $|\lambda| < \frac{1}{\kappa}$. *Then the corresponding* $(m, r, \lambda)$ *version of the endomorphism*

$$\phi \to \Pi\mathfrak{k}(\phi + \mathfrak{T}_\lambda(\phi))$$

*of the* $\mathbb{S}_{SL}^\perp$ *kernel of* $\mathfrak{D}_0$ *has two linear eigenvalues (counted with multiplicity) with norm between* $\frac{1}{\kappa} mr^{-2}$ *and* $\kappa mr^{-2}$, *and* $6g$ - $8$ *eigenvalues (counted with multiplicity) with norm between* $\frac{1}{\kappa} mr^{-2/3}$ *and* $\kappa\, mr^{-2//3}$.

The proof follows momentarily. First comes the following immediate corollary:

**Corollary 7.9**: *There exists* $\kappa > 1$ *with the following significance: Fix* $m < \frac{1}{\kappa}$ *and then fix* $r = \frac{n}{m} > \kappa$. *If* $\lambda \in \mathbb{R}$ *and if* $|\lambda| < \frac{1}{\kappa} mr^{-2}$, *then* $\lambda$ *can not be a linear eigenvalue of the corresponding* $(m, r, \lambda)$ *version of the endomorphism* $\phi \to \Pi\mathfrak{k}(\phi + \mathfrak{T}_\lambda(\phi))$.

*Proof of Lemma 7.8*: Fix $\varepsilon > 0$ for the moment and then a real number (denoted by $\lambda$) with norm at most $\varepsilon$. The chosen number $\lambda$ determines the bilinear form on the $\mathbb{S}_{SL}^\perp$ kernel of $\mathfrak{D}_0$ that is defined by the rule

$$(\phi', \phi) \to \int_\Sigma \langle \phi', \mathfrak{k}(\mathfrak{T}_\lambda(\phi)) \rangle$$

(7.84)

This bilinear form defines a matrix operator on the $\mathbb{S}_{SL}^\perp$ of $\mathfrak{D}_0$ which will be written in block diagonal form with respect to the $\mathbb{C} \oplus K$ splitting of this kernel. The block diagonal description uses the notation in (7.69) with $s$ here denoting the $2 \times 2$ part mapping the $\mathbb{C}$ summand to itself, with $S$ here denoting the (6g - 8) × (6g - 8) summand mapping the K part to itself, and with $\ell$ now denoting the part that maps the K summand to the $\mathbb{C}$ summand.

The inequality below in (7.85) will be used to derive a priori norms for $s$, $S$ and $\ell$:

$$\left| \int_{\Sigma} \langle \phi', \mathfrak{k}(\mathfrak{T}_\lambda(\phi)) \rangle \right| \leq c_0 m \, \|\mathfrak{k}\phi'\|_{\mathbb{L}} \, \|\mathfrak{k}\phi\|_{\mathbb{L}} \,.$$

(7.85)

This follows from (6.63) and the fact that $\mathfrak{k}$ is a symmetric endomorphism with respect to the fiber inner product on $\mathbb{S}_{\text{SL}}{}^{\perp}$ (which is to say that $\langle (\cdot), \mathfrak{k}(\cdot) \rangle = \langle \mathfrak{k}(\cdot), (\cdot) \rangle$). Bounds for the norm of $s$ are obtained by taking $\phi'$ and $\phi$ in (7.80) from the $\mathbb{C}$ summand, those for $S$ are obtained by taking $\phi'$ and $\phi$ from the K summand, and those for $\ell$ are obtained by taking $\phi'$ from the $\mathbb{C}$ summand and $\phi$ from the K summand.

With the preceding understood, suppose first that $\phi$ is from the $\mathbb{C}$ summand of the $\mathbb{S}_{\text{SL}}{}^{\perp}$ kernel of $\mathfrak{D}_0$. The corresponding versions of $\mu_+$ is a constant which will be denoted below by a. In this case, the norms of $\phi$ and $[\hat{\tau}, \phi]$ are described in (7.73). These bounds imply (via (5.32)-(5.34)) and via (7.73)) that

$$\|\mathfrak{k}\phi\|_{\mathbb{L}} \leq c_0(m + e^{-r/c_0}) \, |a|$$

(7.86)

Supposing instead that $\phi$ is from the K summand of the $\mathbb{S}_{\text{SL}}{}^{\perp}$ kernel of $\mathfrak{D}_0$, then the corresponding version of $\mu_+$ has the form depicted in (7.79), and then (7.80) holds. Those bounds imply in turn that

$$\|\mathfrak{k}\phi\|_{\mathbb{L}} \leq c_0(m + e^{-r/c_0})(r^{2/3}|a_-| + \sup_{|z| < r_0} |\mathfrak{w}|) \,.$$

(7.87)

With regards to the norms of a and $a_-$ and $\mathfrak{w}$: In the case when $\phi$ is from the $\mathbb{C}$ summand, then $|a| \leq c_0 r^{-1}\|\phi\|_{\mathbb{L}}$; and when $\phi$ is from the K summand, then $|a_-| + \sup_{\mathbb{D}}|\mathfrak{w}| \leq c_0 r^{-1}\|\phi\|_{\mathbb{L}}$. (See (7.45).)

Given what was said in the preceding paragraph, it now follows directly that when $\Pi(\mathfrak{k}\mathfrak{T}_\lambda(\Pi(\cdot)))$ is written with respect to the $\mathbb{C} \oplus \text{K}$ splitting of the $\mathbb{S}_{\text{SL}}{}^{\perp}$ kernel of $\mathfrak{D}_0$ to obtain a new version of (7.72), then norms of the new versions of the matrices $s$, $S$ and $t$ obey

- $|s| \leq c_0(m^3 + e^{-r/c_0}) \, r^{-2}$,
- $|S| \leq c_0(m^3 + e^{-r/c_0}) \, r^{-2/3}$,
- $|t| \leq c_0(m^3 + e^{-r/c_0}) \, r^{-4/3}$.

(7.88)

Keeping in mind what was said previously from Proposition 7.1 and Lemmas 7.6 and 7.7 about the entries in (7.73): The assertion of Lemma 7.8 follow immediately from (7.88) using Lemmas 7.4 and 7.5. Indeed, the addition of Lemma 7.8's endomorphism to $\Pi \mathfrak{k} \Pi$ replaces $s'$ and $t'$ and $S'$ (7.73) with versions whose norms obey

$$|s'| + r^{2/3}|t'| + r^{4/3}|S'| \leq c_0\, c_0(m^3 + e^{-r/c_0})\,.$$

$$(7.89)$$

If $m < c_0^{-1}$ and $r > c_0$, then Lemmas 7.4 and 7.5 can be brought to bear to draw the required conclusions (Lemma 7.4 to deal with $s'$ and $S'$ and Lemma 7.5 to deal with $t'$).

### d) Proof of Proposition 6.2

Assume in what follows that $m < c_0^{-1}$ and $r > c_0$ so as to invoke Lemmas 6.4-6.9 and the lemmas and propositions in Sections 7a-c. Also, a reminder regarding notation: The operator $\mathfrak{L}$ on $\Sigma$ is defined in (5.56) as a perturbation of the operator $\mathfrak{D}_0$ from (5.1) and the endomorphism $\mathfrak{l}(\cdot)$ is defined in (6.58).

To start the proof: If $\lambda$ is an eigenvalue of $\mathcal{L}$ with norm at most $c_0^{-1}$, then Lemmas 6.4 and 6.6 say in effect that $\lambda$ differs from an eigenvalue of $\mathfrak{L} + \mathfrak{l}(\cdot)$ by at most $c_0 e^{-r/c_0}$. Because both $\mathfrak{L}$ and $\mathfrak{l}(\cdot)$ preserve the orthogonal splitting of $\mathbb{S}$ as $\mathbb{S}_{SL} \oplus \mathbb{S}_{SL}{}^{\perp}$, there is a basis of eigenvectors for $\mathfrak{L} + \mathfrak{l}(\cdot)$ whose elements are either sections of $\mathbb{S}_{SL}$ or sections of $\mathbb{S}_{SL}{}^{\perp}$. As noted at the start of Section 6e, no eigenvalue of the operator $\mathfrak{L} + \mathfrak{l}(\cdot)$ on $C^{\infty}(\Sigma; \mathbb{S}_{SL})$ has norm less than $c_0^{-1}m$. Meanwhile, this operator on $C^{\infty}(\Sigma; \mathbb{S}_{SL}{}^{\perp})$ has 6g - 6 linearly independent eigenvectors with norm no greater than $c_0 m$.

The manipulations in Section 6e identify each eigenvector of $\mathfrak{L} + \mathfrak{l}(\cdot)$ on $C^{\infty}(\Sigma; \mathbb{S}_{SL}{}^{\perp})$ with small normed eigenvalue with an element in the $\mathbb{S}_{SL}{}^{\perp}$ kernel of the operator $\mathfrak{D}_0$ that obeys a *non-linear* eigenvalue equation as defined by that same eigenvalue. To elaborate: Let $\lambda$ denotes the eigenvalue of the given eigenvector in $C^{\infty}(\Sigma; \mathbb{S}_{SL}{}^{\perp})$. Then there is a corresponding (non-zero) element in the $\mathbb{S}_{SL}{}^{\perp}$ kernel of $\mathfrak{D}_0$ (denoted by $\phi$) that obeys

$$\Pi \mathfrak{k}(\phi + \mathfrak{T}_{\lambda}(\phi)) = \lambda \phi$$

$$(7.90)$$

with $\mathfrak{k}$ depicted in (7.1) and with $\mathfrak{T}_{\lambda}$ given by (6.62). With this last point regarding (7.90) and $\mathcal{L}$ well understood: Corollary 7.9 says that (7.90) has no solution if $|\lambda| \leq c_0^{-1} m r^{-2}$.

## 8. Spectral flow

Let P denote a principle SU(2) or SO(3) bundle over either $S^1 \times S^1 \times \Sigma$ or just $S^1 \times \Sigma$. The associated equations in (1.1) on $S^1 \times S^1 \times \Sigma$ and those in (1.7) on $S^1 \times \Sigma$ are the variational equations of a functional on the space of pairs of connection on P and ad(P)-valued self-dual 2 form or 1-form as the case may be. The functional in question in these two cases are

- $(A, \omega) \rightarrow \int_{S^1 \times S^1 \times \Sigma} \left( \langle \omega, F_A^+ \rangle - \frac{1}{6} \langle \omega, [\omega; \omega] \rangle - \frac{m}{2} |\omega|^2 \right).$

- $(A, \mathfrak{a}) \rightarrow \int_{S^1 \times \Sigma} \left( \langle \mathfrak{a} \wedge F_A \rangle - \frac{1}{3} \langle \mathfrak{a} \wedge \mathfrak{a} \wedge \mathfrak{a} \rangle - \frac{m}{2} |\mathfrak{a}|^2 \right).$

(8.1)

In both instances, there is a 'relative' Morse index for any given critical point of the functional in question which assigns an integer to the critical point. This relative Morse index is (by definition) the spectral flow along a continuous, 1-parameter family of Fredholm, self-adjoint operators that interpolates between specific operators that are defined by the solution pair $(A, \omega)$ or $(A, \mathfrak{a})$ and a fiducial pair $(A_0, \omega_0)$ or $(A_0, \mathfrak{a}_0)$ as the case may be.

To be sure: The spectral flow along a continuous path of Fredholm, self-adjoint operators is a certain weighted count of the number of points along the path where the operator in question has 0 is an eigenvalue (thus, a non-trivial kernel). In the simplest case (which can be achieved via a suitably generic compact perturbation of the operators along the path), eigenvalues will cross zero transversally as the path is traversed. In this case, an eigenvalue that crosses zero from negative to positive contributes +1 to the count, and an eigenvalue that crosses zero from positive to negative contributes -1 to the count. Note that there will be some ambiguity regarding the definition of the spectral flow if either the initial operator or the final operator has 0 as an eigenvalue. The ambiguity in this case is the sum of the dimensions of the corresponding eigenvalue 0 eigenspaces. (See [T1] for more about the Vafa-Witten spectral flow question.)

In the case of a solution to (1.7), the initial operator is the that pair's version of the operator on the space of sections of $\mathbb{S}$ depicted below:

$$\psi \rightarrow \mathcal{D}\psi + m\psi_\epsilon.$$

(8.2)

(Remember that for any given pair $(A, \mathfrak{a})$, there is a corresponding operator $\mathcal{D}$ which is $\gamma_k \nabla_{Ak} + \rho_k[\mathfrak{a}_k, \cdot]$.) In the case of a solution to (1.1) that comes from a solution to (1.7) on $S^1 \times \Sigma$ via pullback, the initial operator is the operator on the space of sections of the pullback of $\mathbb{S}$ to $S^1 \times S^1 \times \Sigma$ given by the rule

$$\psi \rightarrow \gamma_s \frac{\partial}{\partial s} \psi + \mathcal{D}\psi + m\psi_\epsilon$$

(8.3)

where the notation is as follows: First, s is used for the Euclidean $\mathbb{R}/2\pi\mathbb{Z}$ coordinate on the left-most $S^1$ factor of $S^1 \times S^1 \times \Sigma$. Meanwhile, $\gamma_s$ is used to denote a specific, covariantly constant anti-symmetric endomorphism of $\mathbb{S}$ with the following properties with respect to the $\gamma$ and $\rho$ endomorphisms that are depicted in (4.6)

- $\gamma_s{}^2 = -1$ ,
- $\gamma_s\gamma_i + \gamma_i\gamma_s = 0$ ,
- $\gamma_s\rho_i + \rho_i\gamma_s = 0$ .

$$(8.4)$$

As for the final operator for the path of operators: That will be the operator in either (8.2) or (8.3) with $m = 0$ as defined using a fiducial pair $(A_0, \mathfrak{a}_0)$ on $S^1 \times \Sigma$ in the case of (8.2) or a fiducial pair $(A_0, \omega_0)$ that is obtained via pullback from that same $(A_0, \mathfrak{a}_0)$ pair on $S^1 \times \Sigma$. With regards to this fiducial pair: For the present purposes the fiducial connection $A_0$ will be the product connection (denoted below by $\theta$) from a given product structure for the principle bundle over $S^1 \times \Sigma$. As for $\mathfrak{a}_0$ (and thus $\omega_0$): This is the identically zero version of $\mathfrak{a}$ or $\omega$ as the case may be.

With regards to the product structure: The spectral flow is invariant with respect to the action of the automorphism group of the principle bundle (the gauge group) in the sense that it does not depend on the chosen product structure. (Said differently: Any solution pair that is obtained from any given solution pair by the action of bundle's automorphism group has the same spectral flow as the given pair. See [T1] for more about this.) The final operator along the Fredholm operator path is the $(A_0 = \theta, \mathfrak{a}_0 \equiv 0)$ version of the operator that is depicted (8.2) or (8.3).

An instance of the upcoming Proposition 8.1 is the assertion that the spectral flow in the context of both (1.7) and (1.1) is bounded along the diverging sequence of solutions to (1.7) and (1.1) that is obtained by fixing some small, positive value for $m$ and then taking $r$ ever larger in the constructions of the previous sections.

The content of the upcoming proposition is quite general: Let Y denote a compact, oriented Riemannian 3-manifold and let $P \to Y$ denote a principle SU(2) or SO(3) bundle. Any given pair $(A, \mathfrak{a})$ of connection on P and ad(P) valued 1-form defines its version of the operator $\mathcal{D} \equiv \gamma_k\nabla_{Ak} + \rho_k[\mathfrak{a}_k, \cdot]$ on the space of sections of Y's version of the bundle $\mathbb{S}$ (which is the bundle $((\oplus_2 T^*Y) \oplus (\oplus_2 \underline{\mathbb{R}})) \otimes \mathrm{ad}(P))$. Let $\theta$ denote a flat connection on P. The proposition concerns the spectral flow for perturbations of the $(A, \mathfrak{a})$ version of $\mathcal{D}$ along paths of self adjoint to the version of $\mathcal{D}$ that is defined by $(A_0 = \theta, \mathfrak{a}_0 = 0)$.

The given pair $(A, \mathfrak{a})$ also defines a self-adjoint operator on the space of sections over $S^1 \times Y$ of the pull-back of $\mathbb{S}$ to $S^1 \times Y$, this being $\gamma_s\frac{\partial}{\partial s} + \mathcal{D}$. This is also the case for the pair $(\theta, 0)$. Proposition 8.1 also remarks on the spectral flow along paths between a perturbation of the $(A, \mathfrak{a})$ version of $\gamma_s\frac{\partial}{\partial s} + \mathcal{D}$ and the $(\theta, 0)$ version with no perturbation.

**Proposition 8.1**: *Given a positive integer* N*, there exists* $\kappa > 0$ *with the following significance: Let* $(A, \mathfrak{a})$ *denote pair of connection on* P *and* ad(P) *valued 1-form on* Y *whose version of* $\mathcal{D}$ *has* N *or less eigenvalues counted with multiplicity with norm less than* $\frac{1}{\kappa}$. *Let* $\mathfrak{e}$ *denote an endomorphism of* $\mathbb{S}$ *with norm* $|\mathfrak{e}| < \frac{1}{\kappa^2}$. *The norm of the spectral flow along a* *path from the* $(A, \mathfrak{a})$ *version of* $\mathcal{D} + \mathfrak{e}$ *to the* $(\theta, 0)$ *version of* $\mathcal{D}$ *is then bounded by* $\kappa$. *This is also the case for the norm of the spectral flow along paths from the* $(A, \mathfrak{a})$ *version of* $\gamma_s \frac{\partial}{\partial s} + \mathcal{D} + \mathfrak{e}$ *to the* $(\theta, 0)$ *version of* $\gamma_s \frac{\partial}{\partial s} + \mathcal{D}$.

This proposition is proved momentary. Some remarks follow directly.

The versions of $(A, \mathfrak{a})$ in Proposition 8.1 for the purposes of this paper are the solutions to (1.7) that are constructed in the previous sections for a fixed small value of $m$ using ever larger values of $r$ (which is to say the integer n when $r$ is written as $\frac{n}{m}$). In the context of the large $r$ solutions to (1.7) constructed in Sections 3-6: The initial operator in (8.2), has trivial kernel because it has the form $\mathcal{L} + \mathfrak{e}$ where $\mathcal{L}$ is the operator from Proposition 6.2 as defined by the pair from Proposition 5.3 and where $\mathfrak{e}$ is an endomorphism of $\mathbb{S}$ with norm bounded by $c_0 e^{-r/c_0}$. And, as per Corollary 7.9 and Lemma 6.6, $\mathcal{L}$ has no eigenvalue with norm less than $c_0^{-1} m\, r^{-2}$ which implies the same for $\mathcal{L} + \mathfrak{e}$ when $r$ is sufficiently large. On the other hand, the final operator on the path will have a kernel. But even so, the kernel dimension is a priori fixed which implies that the spectral flow ambiguity is independent of $r$ when $r$ is large.

The relevant version of $\mathfrak{e}$ for this paper is $m(\cdot)_\mathfrak{c}$. As for the integer N: If $m$ is small, less than $c_0^{-1}$, then it follows from Corollary 7.9 and Lemma 6.6 that $\mathcal{D} + m(\cdot)_\mathfrak{c}$ has at most 12g - 12 eigenvalues counted with multiplicities with norm less than $c_0^{-1}$.

The observations about the integer N in the preceding paragraph apply directly to the operator in (8.3) because any (hypothetical) element in its kernel is independent of the coordinate s for the $S^1$ factor (it is annihilated by $\frac{\partial}{\partial s}$) when $m < c_0^{-1}$. (This can be seen by depicting any hypothetical kernel element as a Fourier series with respect to the s-coordinate. More is said in this regard in Part 2 of the proof of the proposition.)

**Proof of Proposition 8.1**: The proof has two parts. The first part proves the proposition's assertion for the $\mathcal{D} + \mathfrak{e}$. The second part explains why the arguments in the first part of the proof also prove the proposition's assertion for the operator $\gamma_s \frac{\partial}{\partial s} + \mathcal{D} + \mathfrak{e}$.

*Part 1*: The path of operators that will be used to bound the spectral flow is has five stages (the path is continuous and piecewise differentiable).

STAGE 1: This stage consists of a [0, 1] parametrized path with the parameter x version of the operator along the path being $\mathcal{D} + (1-x)\mathfrak{e}$. If $|\mathfrak{e}| < c_0^{-1}$, then, given the assumptions, the spectral flow along this path is no greater than N. Note in particular that the end member of this path is the operator $\mathcal{D} = \gamma_k \nabla_{Ak} + \rho_k[\mathfrak{a}_k, \cdot\,]$ that is defined by the given pair $(A, \mathfrak{a})$.

STAGE 2: To see about the spectral flow from the $(A, \mathfrak{a})$ version of $\mathcal{D}$ to the $(\theta, 0)$ version, digress for a moment to introducea new endomorphism of $\mathbb{S}$ which is denoted here and subsequently by $\Gamma$. To define its action, write any given element from $\mathbb{S}$ as a data set $((\flat, \mathfrak{c}), (\flat_0, \mathfrak{c}_0))$ with $\flat$ and $\mathfrak{c}$ being sections of $T^*Y \otimes ad(P)$ and with $\flat_0$ and $\mathfrak{c}_0$ being sections of $ad(P)$. The endomorphism $\Gamma$ sends this element to the element defined by the data set $((-\flat, \mathfrak{c}), (-\flat_0, \mathfrak{c}_0))$. Thus, it reverses the signs of $\flat$ and $\flat_0$. This endomorphism is symmetric, its square is 1, it commutes with the $\rho$ endomorphisms in (4.6) and it anti-commutes with the $\gamma$ endomorphisms. It also anti-commutes with $\gamma_s$.

Consider now the family of operators parametrized by $[0, \infty)$ according to the rule whereby the parameter $R \in [0, \infty)$ version is

$$\mathcal{D} + R\Gamma \, .$$

$$(8.5)$$

The Bochner-Weitzenboch formula for that operator says that

$$(\mathcal{D} + r\Gamma)^2 = (\gamma_k \nabla_{Ak})^2 + \gamma_k \rho_j [\nabla_{Ak} \mathfrak{a}_j, \cdot \,] + (\rho_k [\mathfrak{a}_k, \cdot \,] + R\Gamma)^2 \, .$$

$$(8.6)$$

This implies in particular that if $R >> \sup_Y(|\mathfrak{a}| + |\nabla_A \mathfrak{a}|)$, then (8.5) has no kernel since in this case,

$$(\mathcal{D} + r\Gamma)^2 \geq (\gamma_k \nabla_{Ak})^2 + \tfrac{1}{2} R^2 \, .$$

$$(8.7)$$

As explained after the description of STAGES 3-5, the spectral flow is bounded by $c_0$ along the family $R \to \mathcal{D} + R\Gamma$ and it is identically zero if $\mathcal{D}$ has trivial kernel.

STAGE 3: Supposing that $R$ is sufficiently large so that (8.7) holds, consider the $[0, 1]$ parametrized family $x \to \gamma_k \nabla_{Ak} + (1 - x)\rho_k [\mathfrak{a}_k, \cdot \,] + R\Gamma$. This family interpolates between the operators $\mathcal{D} + R\Gamma$ and $\gamma_k \nabla_{Ak} + R\Gamma$. There is no spectral flow along this family because, for any given $x$, the same calculation that leads to (8.7) leads to the inequality below:

$$(\gamma_k \nabla_{Ak} + (1 - x)\rho_k [\mathfrak{a}_k, \cdot \,] + R\Gamma)^2 \geq (\gamma_k \nabla_{Ak})^2 + \tfrac{1}{2} R^2 \, .$$

$$(8.8)$$

STAGE 4: Now let $x \to A(x)$ denote any given $[0, 1]$-parametrized family of connections on P that starts at A and ends at $\theta$. There is the corresponding path of operators $\{\gamma_k \nabla_{A(x)k} + R\Gamma\}_{x \in [0, 1]}$; and there is zero spectral flow along this path because

$$(\gamma_k \nabla_{A(x)k} + R\Gamma)^2 = (\gamma_k \nabla_{A(x)k})^2 + R^2 \, .$$

$$(8.9)$$

Stage 5: Finally, consider the path of operators parametized by [0, 1] whose parameter x member is $\gamma_k \nabla_{\theta k} + (1 - x)R\Gamma$. (This path decreases R so that at x = 1, the operator is the fiducial ($\theta$, 0) version of $\mathcal{D}$.) Of particular note: There is no spectral flow on this path until x = 1 because the square of any parameter x version of $\gamma_k \nabla_{\theta k} + (1 - x)R\Gamma$

$$(\gamma_k \nabla_{\theta k})^2 + (1 - x)^2 R^2 .$$

(8.10)

What was said above about Stages 1 and 3-5 imply the assertion of Proposition 8.1 for the operator in $\mathcal{D} + \mathfrak{e}$ provided there is zero spectral flow (or, in any event, spectral flow with norm bounded in terms of N) along the Stage 2 path $R \to \mathcal{D} + R\Gamma$. And, that claim is true because the endomorphism $\gamma_s$ anti-commutes with each $R \in [0, \infty)$ version of $\mathcal{D} + R\Gamma$. Indeed, this anti-commutation implies in particular that if $\lambda$ is an eigenvalue of $\mathcal{D} + R\Gamma$, then so is -$\lambda$; and that $\lambda$ and -$\lambda$ have the same multiplicities (if $\psi$ denotes an eigenvector with eigenvalue $\lambda$, then $\gamma_s\psi$ is an eigenvector with eigenvalue -$\lambda$). As a consequence, any eigenvalue that crosses zero from below as R is increased is offset by an eigenvalue that crosses zero from above. This implies that the spectral flow is zero along the R > 0 part of the path. This also holds for the whole [0, $\infty$) path if the kernel of $\mathcal{D}$ is zero. If not, then the norm of the spectral flow is in any event bounded by N.

*Part 2*: This part of the proof explains why the calculation in Part 1 for the spectral flow imply Propositions 8.1's assertions about the operator $\gamma_s\frac{\partial}{\partial s} + \mathcal{D} + \mathfrak{e}$ acting on sections over $S^1 \times Y$ of the pullback of $\mathbb{S}$. The relationship between the spectral flows for the two operators depends on the fact that the operators $\mathcal{D} + \mathfrak{e}$, $\mathcal{D}$ and the all of the other operators along the paths in Stages 1-5 in Part 1 commutes with the action of $\frac{\partial}{\partial s}$ on the space of sections of the pullback of $\mathbb{S}$. In this regard, a key observation is that any given (A, $\mathfrak{a}$) and R version of the operator $(\gamma_s\frac{\partial}{\partial s} + \mathcal{D} + R\Gamma)$ has square given by the formula below:

$$(\gamma_s\frac{\partial}{\partial s} + \mathcal{D} + R\Gamma)^2 = -\frac{\partial^2}{\partial s^2} + (\mathcal{D} + R\Gamma)^2 .$$

(8.11)

This implies in particular that the large R version of $(\gamma_s\frac{\partial}{\partial s} + \mathcal{D} + R\Gamma)$ has zero kernel and that any kernel element (call it $\psi$) the occurs in Stages 2-5 when (A, $\mathfrak{a}$) and/or R are changed obeys $\frac{\partial}{\partial s}\psi = 0$ and thus is a kernel element for just the $\mathcal{D} + R\Gamma$ part. Hence, the spectral flow along the Stages 2-5 segments of the path for $\gamma_s\frac{\partial}{\partial s} + \mathcal{D} + R\Gamma$ is the same as for the Stages 2-4 segments of the path for $\mathcal{D} + R\Gamma$. As for the Stage 1 segment, the analog of (8.11) for that is the following identity:

$$(\gamma_s\tfrac{\partial}{\partial s} + \mathcal{D} + (1\text{-}x)\mathfrak{e})^2 = -\tfrac{\partial^2}{\partial s^2} + (\mathcal{D} + (1\text{-}x)\mathfrak{e})^2 + (\gamma_s\mathfrak{e} + \mathfrak{e}\gamma_s)\tfrac{\partial}{\partial s}.$$

(8.12)

In particular, if $|\mathfrak{e}| < c_0^{-1}$, then this identity implies that any kernel element (call it $\psi$ again) along this segment obeys $\tfrac{\partial}{\partial s}\psi = 0$ and thus is a kernel element for the STAGE 1 path for the operator $\mathcal{D} + \mathfrak{e}$.

## Appendix

The purpose of this appendix is to supply the proofs for Lemmas 3.1, 3.2, 3.3 and 5.2. The order of proof will be as follows:  Lemma 3.2 first, then Lemmas 3.1 and 3.3 together, and finally, Lemma 5.2.

*Proof of Lemma 3.2:*  Let u denote an initial choice of holomorphic coordinate for a neighborhood of the given point p such that the norm of du at p is equal to $\sqrt{2}$.  The holomorphic quadratic differential $q_1$ (and remember that $q = r^2 q_1$) when written using u as the coordinate has the form

$$q_1 = \alpha_1 u(1 + \vartheta(u))(du)^2$$

(A.1)

with $\vartheta$ denoting a holomorphic function on a disk in $\mathbb{C}$ centered at the origin that vanishes at the origin.  Fix $r_1 > 0$ so that $|\vartheta| < \tfrac{1}{100}$ on the $|u| \le 2r_1$ disk.  Now suppose that z is a different holomorphic coordinate on a neighborhood of the u = 0 point disk with z = 0 being the u = 0 origin and with dz = du at u = 0.  The original coordinate u can be written on the $|u| < 2r_1$ part of the domain of z as a holomorphic function of z, and when this is done, then $q_1$ when written using z has the form

$$q_1 = \alpha_1 u(z)\big(1 + \vartheta(u(z))\big)\big(\tfrac{du}{dz}\big)^2(dz)^2.$$

(A.2)

Thus, the task is to find a holomorphic function z such that $z = u + \mathcal{O}(u^2)$ and

$$z\,\big(\tfrac{dz}{du}\big)^2 = u\,(1 + \vartheta(u))\,.$$

(A.3)

The strategy for finding the desired function z is to look instead for a holomorphic function of u, to be denoted by $\beta$, which vanishes at u = 0 and such that $z \equiv u(1 + \beta(u))$ obeys (A.3).  In particular (A.3) is satisfied if the function $\beta$ obeys the equation

$$u \frac{d\beta}{du} + \beta = (1 + \vartheta)^{1/2}(1 + \beta)^{-1/2} - 1$$

(A.4)

To see about solving (A.4), fix $r \in (0, r_1]$ for the moment and then use it to construct a Banach space of holomorphic functions on the $|u| < r$ disk by completing the space of holomorphic functions on open neighborhoods in $\mathbb{C}$ of the the $|u| \leq r$ disk which vanish at the origin using the norm

$$\beta \to \|\beta\| = \sup_{|u|=r} |\beta|.$$

(A.5)

Denote this Banach space by $\mathcal{H}_r$. Given $\delta \in (0, \frac{1}{100})$, let $\mathcal{B}_{r,\delta}$ denote the radius $\delta$ ball in $\mathcal{H}_r$ as defined using the norm $\| \cdot \|$. (Thus $\sup_{|z| \leq r} |\beta| \leq \delta$ when $\beta \in \mathcal{B}_{r,\delta}$.)

Now supposing that $\varepsilon > 0$ but small (so that $\delta + 2\varepsilon < \frac{1}{100}$), let $\beta$ denote for the moment a holomorphic function defined on a neighborhood of the $|u| \leq r$ disk with norm bounded by $\delta + \varepsilon$ on the circle where $|u| = r$. Define a new holomorphic function on a neighborhood of the origin in $\mathbb{C}$ according to the rule

$$\mathbb{T}[\beta] \equiv \int_0^1 ((1 + \vartheta(xu))^{1/2}(1 + \beta(xu))^{-1/2} - 1) \, dx \, .$$

(A.6)

Some things to note about $\mathbb{T}[\beta]$: First, this function is holomorphic for $|u| \leq r$ since there exists $r' > r$ such that $|\beta| \leq \delta + 2\varepsilon$ on the $|u| < r'$ disk which implies that the integrand in (A.6) on this $|u| \leq r'$ disk is bounded by

$$\tfrac{1}{2} \left( c_0 \sup_{|u|<2r} |\vartheta| + (\delta + 2\varepsilon) + c_0((\delta + 2\varepsilon))^2 \right) .$$

(A.7)

In particular, $\mathbb{T}[\beta]$ obeys

$$|\mathbb{T}[\beta]| \leq \tfrac{1}{2} \left( c_0 \sup_{|u|<2r} |\vartheta| + (\delta + 2\varepsilon) + c_0((\delta + 2\varepsilon))^2 \right)$$

(A.8)

where $|u| < r'$. The second point to note is that $\mathbb{T}[\beta]|_{z=0} = 0$. And the third point is that $\mathbb{T}$ obeys the equation

$$u \frac{d\mathbb{T}}{du} + \mathbb{T} = (1 + \vartheta)^{1/2}(1 + \beta)^{-1/2} - 1$$

(A.9)

Note in particular that if $\mathbb{T}[\beta] = \beta$, then $\beta$ is a solution to (A.4). This is to say that a fixed point of the map $\mathbb{T}$ is a solution to (A.4)

These points have the following consequences: First, (A.8) and the fact that $\mathbb{T}[\beta]$ vanishes at $z = 0$ implies that $\mathbb{T}$ maps the ball $\mathcal{B}_{r,\delta}$ into $\mathcal{H}_r$. Moreover, because $|\vartheta| < c_0 r$ on the $|u| < r$ disk (remember that $\vartheta$ vanishes at the origin), the bound in (A.7) implies that $\|\mathbb{T}(\beta)\| < \frac{3}{4}\delta$ if $r < c_0^{-1}\delta$. Thus, if this bound on $r$ is assumed (which it is henceforth), then $\mathbb{T}$ maps $\mathcal{B}_{r,\delta}$ to itself. Granted that, and supposing that $\mathbb{T}$ is a contraction mapping, it must have a fixed point in $\mathcal{B}_{r,\delta}$; and that fixed point is the sought for solution to (A.4). (See (A.9.)

To see if $\mathbb{T}$ is a contraction mapping: Supposing that $\beta_1$ and $\beta_2$ are in $\mathcal{B}_{r,\delta}$, then the definition of $\mathbb{T}$ in (A.6) can be used to see that

$$|\mathbb{T}[\beta_1] - \mathbb{T}[\beta_2]| \leq \tfrac{1}{2}\int_0^1 ((1 + \vartheta(xu))^{1/2}|\beta_1(xu) - \beta_2(xu)|(1 + c_0\delta)\, dx\,.$$

(A.10)

In particular, if $\delta < c_0^{-1}$ and $r < c_0^{-1}\delta$, then (A.10) implies directly that

$$|\mathbb{T}[\beta_1] - \mathbb{T}[\beta_2]| \leq \tfrac{3}{4}|\beta_1 - \beta_2|\,,$$

(A.11)

and the latter inequality implies that $\mathbb{T}$ is indeed a contraction mapping.

***Proof of Lemmas 3.1 and 3.3***: The proof of these lemmas has six parts.

*Part 1*: This first part of the proof establishes the claim in Lemma 3.1 that there does indeed exist a solution to the $q = r^2 q_1$ version of (3.4). This part also addresses the uniqueness assertion of Lemma 3.1.

To prove existence, define a functional on $C^\infty(\Sigma)$ to be denoted by $\mathcal{E}$ by the rule whereby

$$\mathcal{E}(w) = \tfrac{1}{2}\int_\Sigma \left(|dw|^2 + \tfrac{1}{2}r^2(e^{2w}|\varpi_{1+}|^2 + e^{-2w}|\varpi_{1-}|^2)\right).$$

(A.12)

This function $\mathcal{E}$ is relevant by virtue of the fact that a critical point of $\mathcal{E}$ is formally a solution to the $q = r^2 q_1$ version of (3.4). To prove that $\mathcal{E}$ has a critical point, note that $\mathcal{E}$ is a sum of non-negative terms so it is bounded from below. Therefore, it has an infimum which will be denoted by $\mathcal{E}_0$. (Take $w = 0$ to see that $\mathcal{E}_0 \leq \tfrac{1}{4}r^2$ times the sum of the $\Sigma$ integrals of $|\varpi_+|^2$ and $|\varpi_-|^2$.) To see that this minimum is achieved, it proves useful to extend $\mathcal{E}$'s domain; and to this end, introduce the Hilbert space $L^2_1$ which is the completion of $C^\infty(\Sigma)$ using the norm whose square is given by the rule

$$w \to \int_\Sigma \left(|dw|^2 + w^2\right).$$

(A.13)

With regards to $\mathcal{E}$'s extension to $L^2_1$: If $c$ is any real number, then the map $w \to e^{cw}$ defines a smooth map from $L^2_1$ into the Hilbert space $L^2$ (which is the completion of $C^\infty(\Sigma)$ using the norm whose square assigns to any given function $w$ the value of the $\Sigma$ integral of $w^2$). See for example Theorem 2.46 in [A] for a proof of this fact. With the preceding point understood, it then follows that $\mathcal{E}$ extends to $L^2_1$ as a smooth function.

The extension to $L^2_1$ is also coercive in the sense that the inequality below holds for any given w from $L^2_1$:

$$\mathcal{E}(w) \geq c_0^{-1} \int_\Sigma \left( |dw|^2 + r^{2/3} w^2 \right).$$

(A.14)

To see why this is: Note first that the following inequality

$$r^2 \left( e^{2w} |\varpi_+|^2 + e^{-2w} |\varpi_-|^2 \right) \geq r^2 |w|^2 \, \mathfrak{d}^2$$

(A.15)

holds for any function w because $e^{2w} > 2|w|^2$ if $w \geq 0$ and $e^{-2w} > 2w^2$ if $w < 0$. Meanwhile: For any $r \in (0, r_0)$, a standard Dirichlet eigenvalue bound applied to the $\mathfrak{d} \leq r$ disks about the points in $\Theta$ leads to this bound:

$$\frac{1}{r} \int_{\mathfrak{d}<r} |w|^2 \leq c_0 r \int_{\mathfrak{d}<2r} |dw|^2 + c_0 \frac{1}{r} \int_{r<\mathfrak{d}<2r} |w|^2$$

(A.16)

The latter bound with $r = c_0 r^{-2/3}$ and (A.14) together lead to (A.14).

With (A.14) understood: Let $\{w_k\}_{k=1,2,\dots}$ denote a sequence in $C^\infty(\Sigma)$ with the corresponding sequence $\{\mathcal{E}(w_k)\}_{k=1,2,\dots}$ being a non-increasing sequence with limit value $\mathcal{E}_0$. Because the $\mathcal{E}(w_k)$ sequence is non-increasing, the sequence $\{w_k\}_{k=1,2,\dots}$ is necessarily a norm-bounded sequence in $L^2_1$ (this is the point of (A.14)). And, being that $\{w_k\}_{k=1,2,\dots}$ is a norm bounded sequence in a (separable) Hilbert space, it has a weakly convergent subsequence. Let $\mu$ denote a weak limit. Weak $L^2_1$ convergence implies that

$$\int_\Sigma |d\mu|^2 \leq \lim_{k\to\infty} \int_\Sigma |dw_k|^2$$

(A.17)

and that

$$\int_\Sigma \left( e^{2\mu} |\varpi_{1+}|^2 + e^{-2\mu} |\varpi_{1-}|^2 \right) \leq \lim_{k\to\infty} \int_\Sigma \left( e^{2w_k} |\varpi_{1+}|^2 + e^{-2w_k} |\varpi_{1-}|^2 \right).$$

(A.18)

(The latter inequality is a priori an equality because the map $w \to e^{cw}$ is a compact mapping from $L^2_1$ to $L^2$ (see Theorem 2.46 in [Aubin]).) These inequalities in (A.17) and (A.18) imply that $\mathcal{E}(\mu) \leq \mathcal{E}_0$ which requires that $\mathcal{E}(\mu) = \mathcal{E}_0$ since $\mathcal{E}_0$ is the infimum of $\mathcal{E}$. Thus, $\mathcal{E}$ achieves its infimum at $\mu$. That observation implies in turn that $\mu$ is an $L^2_1$ solution to the equation in (3.4). As such, standard elliptic regularity theorems can be brought to bear to see that $\mu$ is smooth (see for example Chapter 5 of [M] or Chapter 8 of [GT]).

The last task for this part of Lemma 3.1's proof addresses the uniqueness assertion. To see that, suppose for the moment that $\mu$ and $\mu'$ are solutions to (3.4). If this is so, then

$$d^\dagger d(\mu - \mu') + r^2(e^{2\mu} - e^{2\mu'})|\varpi_+|^2 - r^2(e^{-2\mu} - e^{-2\mu'})|\varpi_-|^2 = 0 \ .$$

$$(A.19)$$

Appeal to the maximum principle to see that $\mu - \mu'$ can not have either a positive maximum or negative minimum. (Check the signs of the three terms above at any such point.) The only recourse is that $\mu$ and $\mu'$ are equal.

*Part 2*: This second part of the proof constructs the function $\mathfrak{f}$ that appears in Lemma 3.3. There are three steps to doing this. By way of a look ahead: The plan is to minimize a functional to obtain the desired $\mathfrak{f}$, much like what was done in Part 1's proof for (3.4). But, in this regard, a verbatim copying of Part 1's argument can't be made because the analog of the function $e^{2w}|\varpi_+|^2 + e^{-2w}|\varpi_-|^2$ would be $e^{2w}|z|^2 + e^{-2w}$ which does not have a finite integral on $\mathbb{C}$ for any choice of w.

By way of notation in the subsequent steps: What is denoted by s signifies the radial coordinate $|z|$ on $\mathbb{C}$.

Step 1: Introduce the radially symmetric function ô on $\mathbb{C}$ by the rule

$$s \rightarrow \hat{o}(s) \equiv -\tfrac{1}{2}(1 - \chi(s))\ln(s).$$

$$(A.20)$$

To be sure, the function ô is zero where $s < \tfrac{1}{4}$ and it is equal to $-\tfrac{1}{2}\ln(s)$ where $s > \tfrac{3}{4}$. The function ô is used directly to define a functional on the $L^2{}_1$ Hilbert space to be denoted by $\mathfrak{E}$ which will play the role played by the functional $\mathcal{E}$ in Part 1. This functional $\mathfrak{E}$ is defined by setting its value on any given compactly supported radial function w as follows:

$$\mathfrak{E}(w) = \int_{[0,\infty)} \tfrac{1}{2}|\tfrac{d}{ds}w|^2 \, sds \ - \ \tfrac{1}{2}\int_{[0,\infty)} ((\tfrac{d^2}{ds^2}\chi + \tfrac{1}{s}\tfrac{d}{ds}\chi)\ln(s) + 2\tfrac{d}{ds}\chi)w) \, sds$$
$$+ \tfrac{1}{4}\int_{[0,\infty)} (e^{2w} s^{1+\chi} + e^{-2w}s^{1-\chi} - 2s)sds \ .$$

$$(A.21)$$

Note in particular that all three integrals in (A.21) are finite when w has compact support (this is because $\chi$ and $d\chi$ have compact support).

As explained in the next step, this functional extends to a complete, locally convex function space where it is bounded from below and has a global minimizer. The relevance of the preceding observation stems from the fact that if $w_0$ is a critical point of $\mathfrak{E}$ on that space (for example, its global minimizer), then the function $\mathfrak{f} \equiv w_0 + \hat{o}$ satisfies the equation in the top bullet of Lemma 3.3.

Step 2:  Let $L^2{}_{1\diamond}$ denote the Hilbert space of radially symmetric functions (functions of s) that is obtained by completing the space of compactly supported, radially symmetric functions using the norm whose square assigns the value

$$\int_{[0,\infty)}\ (|\tfrac{d}{ds}\,w|^2 + s\,w^2)\ sds$$

(A.22)

to any given compactly supported function.  The norm whose square is depicted above is denoted here by $\|\cdot\|_{\diamond}$.

The first task is show that $\mathfrak{E}$ extends to this Hilbert space $L^2{}_{1\diamond}$ as a smooth functional. To do this, note first that if $w \in L^2{}_{1\diamond}$, then there exists $R \geq 2$ such that $|w| < \frac{1}{1000}$ where $s > R$ and that R can be bounded by w's Hilbert space norm.  Indeed, such is the case because of the following facts:

- *If $\rho > 2$ and $|\rho - \rho'| \leq 1$, then $|w(\rho) - w(\rho')| \leq c_0 \frac{1}{\rho} (\int_{\rho-1}^{\rho+1}\ |\tfrac{d}{ds}\,w|^2 sds)^{1/2}$ .*
- *Thus, if $|w(\rho)| > \frac{1}{1000}$, then $w > \frac{1}{2000}$ on $[\rho - 1, \rho + 1]$ when $\int_{\rho-1}^{\rho+1} |\tfrac{d}{ds}\,w|^2 sds \leq c_0 \rho^2$.*
- *Thus, if $|w(\rho)| > \frac{1}{1000}$ and $\|w\|_{\diamond}^2 \leq c_0 \rho^2$, then $\int_{[0,\infty)}\ w^2\,s^2 ds \geq c_0^{-1} \rho^2$.*

(A.23)

Indeed, the third bullet above implies that

$$|w(s)| \leq \tfrac{1}{1000}\ \ where\ \ s \geq c_0 \|w\|_{\diamond}.$$

(A.24)

Meanwhile, if $|w| < \frac{1}{1000}$ where $s \geq R$ for any given $R > 2$, then

$$\tfrac{1}{4}\ \int_{[R,\infty)}\ (e^{2w}\,s^{1+\chi} + e^{-2w} s^{1-\chi} - 2s) sds\ \leq c_0 \int_{[R,\infty)}\ w^2 s^2 ds\ .$$

(A.25)

The preceding observation can be used to see that the contribution to the right most integral in (A.21) from the $s \geq c_0 \|w\|_{\diamond}$ part of the integration domain $[0, \infty)$ extends to $L^2{}_{1\diamond}$ as smooth functional.  On the flip side, what is said by Theorem 2.46 in [A] implies that the contribution from the $s \leq c_0 \|w\|_{\diamond}$ part of the right most integral in (A. 21) also extends as a smooth function on $L^2{}_{1\diamond}$.

Step 3:  With it understood that $\mathfrak{E}$ extends to the Hilbert space $L^2{}_{1\diamond}$, note next that this extension is bounded from below and coercive.  This is because $\mathfrak{E}(w)$ obeys the bound

$$\mathfrak{E}(w) \geq \int_{[0,\infty)}\ \tfrac{1}{2}|\tfrac{d}{ds}\,w|^2\,sds\ + \int_{[2,\infty)}\ \sinh^2(w)\,s^2 ds - \tfrac{1}{2}\int_{[0,\infty)}\ \mathfrak{g}(s)\,w\,sds$$

(A.26)

with $\mathfrak{g}$ denoting $(\frac{d^2}{ds^2}\chi + \frac{1}{s}\frac{d}{ds}\chi)\ln(s) + 2\frac{d}{ds}\chi$ which is a smooth function with compact support on the interval $[\frac{1}{4}, 1]$. (And keep in mind that $\sinh^2 w \geq w^2$).

It follows from the preceding observations that any minimizing sequence of $\mathfrak{E}$ in $L^2{}_{1\Diamond}$ has a weakly convergent subsequence in $L^2{}_{1\Diamond}$; and as before, any limit of that subsequence is an element in $L^2{}_{1\Diamond}$ where $\mathfrak{E}$ achieves its minimal value. (This is because the respective values of the left most two integrals in (A.26) for the limit function are no greater than the limit of their values on the minimizing sequence; and the right most integral for the limit is equal to the limit of its values on the minimizing sequence.) It follows as a consequence that the limit function (denoted by $w_{\ddagger}$) is a critical point of $\mathfrak{E}$ and thus $\mathfrak{f} = w_{\ddagger} + \hat{o}$ obeys the equation in the first bullet of Lemma 3.3.

By way of a parenthetical remark: The functional $\mathfrak{E}$ has only one critical point in the Hilbert space $L^2{}_{1\Diamond}$. (Reuse the uniqueness argument from Part 1 to prove this).

*Part 3*: This part of the proof establishes the claims about the function $\mathfrak{f}$ from the third bullet of Lemma 3.3. The next part of the proof deals with the claims in the second bullet. As a preliminary step, note that $\mathfrak{f}$ is a smooth function on $[0, \infty)$ by virtue of it satisfying the equation in the second bullet of Lemma 3.3. (This is a standard result from the theory of ordinary differential equations.) Therefore, the derivatives to any given order are a priori bounded on the interval $[0, 2]$. Thus, there is some positive $c_1 < c_0$ such that

$$|\mathfrak{f}(s)| + |\frac{d}{ds}\mathfrak{f}| + |\frac{d^2}{ds^2}\mathfrak{f}| \leq c_1 \ \textit{where} \ s \leq 2.$$

(A.27)

With the preceding as some background, define now a function on $(0, \infty)$ to be denoted by u by writing $\mathfrak{f}$ as $u - \frac{1}{2}\ln(s)$. (The function u where $s \geq 1$ is the function $w_{\ddagger}$ from Part 2 of this proof.) The equation in Lemma 3.3's second bullet for $\mathfrak{f}$ when written in terms of u says that

$$-\frac{1}{s}\frac{d}{ds}(s\frac{d}{ds}u) + s\sinh(2u) = 0 .$$

(A.28)

Since $\mathfrak{f}$ is bounded where $s < 1$ (see (A.27), and $-\frac{1}{2}\ln(s)$ is unbounded and grows ever larger as $s \to 0$, the function u is negative for small s. That fact with the maximum principle applied to (A.28) implies that u is negative for all s. (This is one claim from the Lemma 3.3's third bullet.)

Since $u < 0$, so is $\sinh(2u)$ and thus so is $s\sinh(2u)$ which is $e^{2\mathfrak{f}}s^2 - e^{-2\mathfrak{f}}$. (This is another claim from that third bullet.)

Since $\sinh(2u)$ is negative, the equation in (A.28) says in part that $s\frac{d}{ds}u$ is a decreasing function of s which implies in turn that $\frac{d}{ds}u$ must be strictly positive for all s. Indeed, if it is ever zero at some finite s, it would turn negative for larger s and get ever more negative in which case $w_{\ddagger}$ would not limit to zero as $s \to \infty$. (Another claim from that third bullet.) By way of a parenthetical remark: Because $\frac{d}{ds}u > 0$, the function u lacks critical points.

Since $\sinh(2u)$ is negative and $\frac{d}{ds}u$ is positive, the equation in (A.28) implies that $\frac{d^2}{ds^2}u$ is strictly negative (this is the final claim from the third bullet).

*Part 4*: This part proves the bounds that are asserted by the second bullet of Lemma 3.3. In this regard: These bounds can be written using u,

$$|u(s)| + s|\frac{d}{ds}u| + s^2|\frac{d^2}{ds^2}u| \leq c_0 e^{-s^{3/2}/c_0} \quad where \, s \geq 1 \, .$$

(A.29)

With u being negative, the inequality $\sinh(2u) < 2u$ holds and this inequality with (A.29) implies the inequality below:

$$-\frac{d^2}{ds^2}u - \frac{1}{s}\frac{d}{ds}u + 2su > 0 \, .$$

(A.30)

To see the implications of (A.30), fix $c > \frac{9}{8}$ for the moment and let $q(s) \equiv e^{-s^{3/2}/c}$. This function obeys the same differential inequality as u:

$$-\frac{d^2}{ds^2}q - \frac{1}{s}\frac{d}{ds}q + 2sq > 0$$

(A.31)

on $[1, \infty)$. As a consequence of this and (A.30), any $C > 0$ version of $Cq + u$ obeys (A.30)/(A.31) also on $[1, \infty)$ which implies that $Cq + u$ can not have a negative, local minimum on $[1, \infty)$. If $C > -e^{1/c}u|_{s=1}$, then $Cq+u$ will be positive at $s = 1$ and since it limits to zero as $s \to \infty$, it must be positive on $[1, \infty)$. This is to say that u on $[1, \infty)$ obeys

$$|u| < c_0 e^{-s^{3/2}/c}$$

(A.32)

as long as $c > \frac{9}{8}$.

It follows from (A.31) that $\sinh(2u) > -c_0 e^{-s^{3/2}/c}$ on $[1, \infty)$. As a consequence

$$\frac{d}{ds}(s\frac{d}{ds}u) \geq -c_0 \, s^2 e^{-s^{3/2}/c}$$

(A.33)

Since the right hand side of this last inequality is greater than $-c_0 s^{1/2} e^{-s^{3/2}/2c}$ where $s \geq 1$, integration of this inequality on the interval $[s, \infty)$ leads to the bound below

$$- s\frac{d}{ds}u > - \frac{4}{3}\, c_0 c\, e^{-s^{3/2}/2c}$$

$$(A.34)$$

which is the desired bound for $\frac{d}{ds}u$ because the latter function is positive.

Finally, consider that $s^2\frac{d^2}{ds^2}u$ can be written as

$$s^2\frac{d^2}{ds^2}u = s^3 \sinh(2u) - s\frac{d}{ds}u \ ,$$

$$(A.35)$$

which implies that it is negative and then, what with (A.33) and (A.34), that it obeys

$$s^2|\frac{d^2}{ds^2}u| \leq c_0 c e^{-s^{3/2}/4c}$$

$$(A.36)$$

where $s \geq 1$.

*Part 5*:  This part of the proof addresses a portion of the assertion in the first bullet of Lemma 3.1 regarding the pointwise bound for the norm of its function $\mu_r$ and its derivatives to second order.   To be precise, what follows establishes the assertion below

$$|\mu_r - \mu_\diamond| + |d(\mu_r - \mu_\diamond)| + |\nabla d(\mu_r - \mu_\diamond)| \leq c_0 e^{-r/c_0} \ \ where\, \mathfrak{d} > \frac{1}{100}\, r_0 \ .$$

$$(A.37)$$

This inequality implies directly that the inequality asserted by the first bullet of Lemma 3.1 holds on the $\mathfrak{d} > \frac{1}{100}\, r_0$ part of $\Sigma$.

The key step towards proving (A.37) is to prove that $|\mu_r - \mu_\diamond|$ is bounded by $c_0 e^{-r/c_0}$ on a slightly larger domain:

$$|\mu_r - \mu_\diamond| \leq c_0 e^{-r/c_0} \ \ where\, \mathfrak{d} > \frac{1}{150}\, r_0 \ .$$

$$(A.38)$$

Granted the preceding bound (which is proved momentarily), and granted that $\mu_r - \mu_\diamond$ obeys

$$d^\dagger d(\mu_r - \mu_\diamond) + r^2|q_1| \sinh(2(\mu_r - \mu_\diamond)) = 0 \ ,$$

$$(A.39)$$

standard elliptic regularity techniques applied to (A.39) give the first and second derivative bounds that are asserted in (A.37).  A brief elaboration is given momentarily (see for example Chapter 5 in [M] or Chapter 8 in [GP]).

As for proving (A.38):  Given (A.39), and given that the function $x \to \frac{1}{x}\sinh(x)$ is no smaller than 1 for any given real number $x$, the identity in (A.39) implies in turn the inequality below for the norm of $\mu_r - \mu_\diamond$:

$$d^\dagger d|\mu_r - \mu_\diamond| + c_0^{-1}r^2\,|\mu_r - \mu_\diamond| \le 0 \quad \text{where}\, \eth > \tfrac{1}{400}\,r_0\,.$$

(A.40)

Let $c_\ddagger$ denote the version of $c_0$ that appears in (A.39). Also, let $\Sigma_0$ denote the part of $\Sigma$ where $\eth$ is between $\frac{1}{400}\,r_0$ and $\frac{1}{200}\,r_0$. The Green's function for the operator $d^\dagger d + c_\ddagger^{-1}r^2$ can be used with (A.39) to see that

$$|\mu_r - \mu_\diamond| \le c_0 e^{-r/c_0 c_\ddagger} \int_{\Sigma_0}\;|\mu_r - \mu_\diamond| \quad \text{where}\, \eth > \tfrac{1}{150}\,r_0\,.$$

(A.41)

In this regard, keep in mind that Green's function for that operator with pole at any given point p is bounded by $c_0(1 + |\ln(r\,\mathrm{dist}(p, \cdot)|)\;e^{-r\,\mathrm{dist}(\cdot,p)/c_0 c_\ddagger}$.

To obtain (A.38) from (A.41), return to a remark in Part 1 to the effect that the minimum of the functional $\mathcal{E}(\cdot)$ is less than $c_0\,r^2$. Thus, $\mathcal{E}(\mu_r) < c_0\,r^2$. That implies in turn the bound below:

$$\int_{\Sigma_0}\;r^2|q_1|\,\cosh(2(\mu_r - \mu_\diamond)) \le c_0 r^2\,,$$

(A.42)

and thus (because $\cosh(x) \ge 1 + \tfrac{1}{2}x^2$ for any real number x) that

$$\int_{\Sigma_0}\;|\mu_r - \mu_\diamond|^2 \le c_0\,.$$

(A.43)

This last bound with the bound in (A.41) implies directly the desired bound in (A.38).

As for the derivative bounds in (A.37): To avoid looking at an analysis book, here is a brief description to the derivation: Differentiate (A.39) to obtain an elliptic equation for the 1-form $d(\mu_r - \mu_\diamond)$ which implies the differential inequality below on the $\eth > \frac{1}{150}\,r_0$ part of $\Sigma$ (this assumes $r > c_0$ and (A.38)):

$$d^\dagger d|d(\mu_r - \mu_\diamond)| + c_0^{-1}r^2|d(\mu_r - \mu_\diamond)| \le c_0 e^{-r/c_0}\,.$$

(A.44)

Let $c_{\ddagger 1}$ denote the version of $c_0$ in (A.44). Mimic what was done for (A.40) with this inequality using the Green's for the operator $d^\dagger d + c_\ddagger^{-1}r^2$ to obtain the desired bound for the pointwise norm of $d(\mu_r - \mu_\diamond)$. To obtain the second derivative bounds in (A.37), differentiate (A.39) twice to derive a similar inequality for $d^\dagger d|\nabla d(\mu_r - \mu_\diamond)|$. But in the case of the second derivative, note that the equation obtained from (A.38) for $d(\mu_r - \mu_\diamond)$ by differentiating once must be used first to obtain a $c_0 e^{-r/c_0}$ bound for the integral of the function $|\nabla d(\mu_r - \mu_\diamond)|^2$ over the $\eth > \frac{1}{125}\,r_0$ part of $\Sigma$. That bound is needed to obtain the final pointwise bound for $|\nabla d(\mu_r - \mu_\diamond)|$ from the $|\nabla d(\mu_r - \mu_\diamond)|$ analog of (A.41).

*Part 6*:  This last part of the proof of Lemmas 3.1 and 3.3 proves the assertions in Lemma 3.3 regarding the function $y$.  With regards to completing the proof of Lemma 3.3: Lemma 3.3's assertions about $y$ and the assertions in Lemma 3.3's second bullet about $f$ complete the proof of Lemma 3.1 because they directly imply what is asserted by Lemma 3.1's second bullet and they directly imply that Lemma 3.1's first bullet inequalities holds on the $\partial < r_0$ part of $\Sigma$.

To start the analysis of $y$:  Fix a point $p \in \Theta_+$ and then define the function $v$ on the $\partial < r_0$ disk D centered at p using the $\mu = \mu_r$ version of (3.7).  Meanwhile: Use the function $f$ from Lemma 3.3 to define the function $v_f$ using the rule in (7.6).  Note that the function $y$ from Lemma 3.3 is $v$ - $v_f$.  Note also that both $v$ and $v_f$ obey the equation in Lemma 3.2 where $|z| < r_0$ and so their difference obeys the equation below where $|z| < r_0$:

$$-4\partial\overline{\partial}(v - v_f) + \tfrac{1}{2}\,\alpha\,((e^{2v} - e^{2v_f})\,|z|^2 - (e^{-2v} - e^{-2v_f})) = 0\;.$$

(A.45)

Of particular relevance:  This equation implies that $v$ - $v_f$ (which is $y$) has neither a positive local maximum or a negative local minimum where $|z| < r_0$.  Therefore, $|v - v_f|$ takes on its maximum on the $|z| \leq \tfrac{3}{4} r_0$ disk on the boundary of this disk.  But note that $|v - v_f|$ is no greater than $c_0 e^{-r/c_0}$ there, this being a consequence of (A.38) and what is said in the second bullet of Lemma 3.3 about $f$.  Thus $|y| < c_0 e^{-r/c_0}$ where $\partial < \tfrac{3}{4} r_0$ also.  This $c_0 e^{-r/c_0}$ for $|y|$ is the norm bound claim about $y$ in Lemma 3.3.

With the preceding bound on $|y|$ settled, consider next Lemma 3.3's asserted bound for the derivatives of $y$.  To obtain that bound, take the derivative of (A.44) to obtain this:

$$- 4\partial\overline{\partial}(d(v - v_f)) + \alpha\,((e^{2v} - e^{2v_f})\,|z|d|z| + \alpha\,((e^{2v} - e^{2v_f})\,|z|^2 + (e^{-2v} - e^{-2v_f}))(d(v - v_f)) = 0\;.$$

(A.46)

It follows as a consequence of this and the bound on $|y|$ that

$$- 4\partial\overline{\partial}|dy| \leq c_0 e^{-r/c_0}\;\; where\;\; |z| < \tfrac{3}{4}\,r_0\;.$$

(A.47)

Given the bounds in (A.37) where $\partial > \tfrac{1}{100}\, r_0$, and given the first derivative bounds on $f$ in the second bullet of Lemma 3.3, this last inequality can be satisfied only if $|dy| \leq c_0 e^{-r/c_0}$ where $\partial < \tfrac{3}{4}\, r_0$.  (Use the comparison principle with $c_0 e^{-r/c_0}(1 - |z|^2)$ being the comparison function.)  This $c_0 e^{-r/c_0}$ bound for $|dy|$ is the norm bound for $|dy|$ in Lemma 3.3.

The norm bound from Lemma 3.3 on $\nabla dy$ is obtained by differentiating (A.46) again and using the bounds now established for the norms of $y$ and $dy$ to see that $- 4\partial\overline{\partial}|\nabla dy|$ is also bounded by $c_0 e^{-r/c_0}$ *where* $|z| < \frac{3}{4} r_0$. This being the case, then the bounds from (A.37) and Lemma 3.3's bounds for the second derivative of $f$ can be used with the comparison principle to justify the claimed $c_0 e^{-r/c_0}$ bound in Lemma 3.3 for $|\nabla dy|$.

***Proof of Lemma 5.2:*** By way of a look ahead, the proof is much like the proof of Lemma 3.2. In any event, there are four parts.

*Part 1*: This part of the proof addresses the bounds asserted by the third bullet of the lemma. Borrow the function $\chi_{\ddagger}$ from (5.37) (it is the bump function $\chi(\frac{8\mathfrak{d}}{r_0} - 1)$ which is equal to 1 where $\mathfrak{d} < \frac{1}{8} r_0$ and equal to zero where $\mathfrak{d} > \frac{1}{4} r_0$) and use it to define a function (denoted by $\chi_{\ddagger 1}$) with compact support on the $\mathfrak{d} < r_0$ part of $\Sigma$ by the rule whereby $\chi_{\ddagger 1} = \pm \chi_{\ddagger}$ with the + sign used near the points in $\Theta_+$ and with the - sign used near the points in $\Theta_-$. With this definition understood, use (5.12) with (5.14) and (5.15) to see that the $\mathbb{C}$-valued function $\eta \equiv \varsigma_- + \chi_{\ddagger 1}\varsigma_+$ is a smooth function on $\Sigma$ that obeys the equation below:

$$d^{\dagger}d\eta + 4|\hat{\mathfrak{a}}|^2\eta = 4(\chi_{\ddagger 1}|\hat{\mathfrak{a}}|^2 + *_{\Sigma}\langle\sigma_3\hat{\mathfrak{a}} \wedge \hat{\mathfrak{a}}\rangle)\varsigma_+ + d^{\dagger}d\chi_{\ddagger 1}\varsigma_+ - 2\langle d\chi_{\ddagger 1}, d\varsigma_+\rangle$$

$$(A.48)$$

An essential point is that the right hand side of this identity is a smooth function which is pointwise bounded by $c_0 r^{4/3}$ (see Lemmas 3.1 and 3.3). Keeping this in mind, multiply both sides of (A.48) by the complex conjugate of $\eta$ and then integrate both sides of the resulting identity over $\Sigma$. Doing that leads directly to the inequality below:

$$\int_{\Sigma} (|d\eta|^2 + 4|\hat{\mathfrak{a}}|^2|\eta|^2) \leq c_0 r^{8/3} .$$

$$(A.49)$$

Now let $\Sigma_1$ denote the part of $\Sigma$ where $\mathfrak{d} > \frac{1}{4} r_0$. The inequality below implies this next one:

$$\int_{\Sigma_1} (|d\varsigma_-|^2 + r^2|\varsigma_-|^2) \leq c_0 r^{8/3} .$$

$$(A.50)$$

This is a poor bound but it suffices for the purposes at hand.

To continue, use (A.48) on $\Sigma_1$ with what is said in Lemma 3.3 about the function $y$ and what is said in that lemma's final bullet to see that the inequality below holds on $\Sigma_1$:

$$d^{\dagger}d|\varsigma_-| + c_0^{-1}r^2|\varsigma_-| \leq c_0 e^{-r/c_0} .$$

$$(A.51)$$

Now use the Green's function for $d^{\dagger}d + c_0^{-1}r^2$ with (A.51) in the manner of (A.38)'s derivation to see that $|\varsigma_-| \leq c_0 e^{-r/c_0}$ on the $\mathfrak{d} > \frac{3}{8} r_0$ part of $\Sigma$.

To obtain the desired bounds for $|d\varsigma_-|$ on the $\mathfrak{d} > \frac{1}{2}\, r_0$ part of $\Sigma$: Differentiate both sides of (A.48) on $\Sigma_1$ and use the bounds just proved for $|\varsigma_-|$ to see that $|d\varsigma_-|$ on the $\mathfrak{d} > \frac{3}{8}\, r_0$ part of $\Sigma_1$ obeys its own version of (A.51) (with larger versions of $c_0$). The integral bound in (A.50) can then be used in conjuction with the $d^{\dagger}d + c_0^{-1} r^2$ Green's function to bound $|d\varsigma_-|$ by $c_0 e^{-r/c_0}$ on the $\mathfrak{d} > \frac{1}{2}\, r_0$ part of $\Sigma$.

*Part 2*: This part of the proof establishes the assertions in the fourth bullet of Lemma 5.2 about the function $\Xi$. To this end, it is sufficient to prove these assertions for the case when $\alpha = 1$ because the $\alpha = r^2\alpha_1$ version is obtained from the $\alpha = 1$ version by rescaling the coordinate z. To establish those assertions, note first that the $\alpha = 1$ version of $\Xi$ is $\mp \frac{2}{3}\, (s\frac{d}{ds}\mathfrak{f} + \frac{1}{2})$ with $\mathfrak{f}$ from Lemma 3.3. This being the case, the assertion in Item a) of the fourth bullet of Lemma 5.2 follows directly from the bounds in Lemma 3.3's second bullet. With regards to Item b) of Lemma 5.2's fourth bullet: These observations are direct consequence of what is said in the third bullet of Lemma 3.3.

*Part 3*: This part of the proof addresses the claim made by Lemma 5.2's first bullet. To this end, fix a point in $\Theta_+$ to consider the disk D around that point where $\mathfrak{d} < r_0$. (The arguments below apply to the points in $\Theta_+$ since the arguments for the points in $\Theta_-$ are identical but for some sign changes.) The equation in (5.12) in the case $\varsigma_+ = 1$ when written on the $\mathfrak{d} < r_0$ disk centered at the given point using the corresponding holomorphic coordinate z is the equation below:

$$-4\partial\bar{\partial}\varsigma_- + \alpha\, (e^{2v}|z|^2 + e^{-2v})\varsigma_- \,=\, \alpha\, (e^{2v}|z|^2 - e^{-2v})\ .$$

(A.52)

Here, v is the function in (3.7). Meanwhile, Lemma 5.2's function $\Xi$ obeys (7.14) which is the analogous equation with v replaced by Lemma 3.3's function $v_{\mathfrak{f}}$. Noting that v and $v_{\mathfrak{f}}$ differ by at most $c_0 e^{-r/c_0}$, the equation in (7.15) and the equation in (A.52) together imply a differential equation for $\varsigma_- - \Xi$ which has this form:

$$-4\partial\bar{\partial}(\varsigma_- - \Xi) \,+ \alpha\, (e^{2v}|z|^2 + e^{-2v})(\varsigma_- - \Xi) \,=\, \mathfrak{e}\ ,$$

(A.53)

where $\mathfrak{e}$ is a function on the $|z| < r_0$ part of $\mathbb{C}$ whose norm and first derivative norm obey

$$|\mathfrak{e}| + |d\mathfrak{e}| \leq c_0 e^{-r/c_0}.$$

(A.54)

Because of this, the norm of $|\varsigma_- - \Xi|$ obeys the differential inequality

$$-4\partial\bar{\partial}|\varsigma_- - \Xi| \leq c_0 e^{-r/c_0}\ .$$

(A.55)

Because the norms of both $\varsigma_-$ and $\Xi$ are bounded by $c_0\,e^{-r/c_0}$ where $|z| > \frac{1}{2}r_0$ (the third and fourth bullets of Lemma 5.2), the Dirichlet Green's function on the $|z| \le \frac{3}{4}$ disk can be used in conjunction with (A.55) in the manner as previous Green's function applications to obtain a $c_0e^{-r/c_0}$ bound for $|\varsigma_- - \Xi|$ where $\mathfrak{d} \le \frac{5}{8}r_0$.

A $c_0e^{-r/c_0}$ bound for $|d(\varsigma_- - \Xi)|$ can be derived by first differentating both sides of (A.53) and then using (A.54) with the $c_0e^{-r/c_0}$ bound for $|\varsigma_- - \Xi|$ to see that $-4\partial\bar\partial|d(\varsigma_- - \Xi)|$ is no greater than $c_0e^{-r/c_0}$ (the analog of (A.55) with some larger version of $c_0$). This $|d(\varsigma_- - \Xi)|$ version of (A.55) with the fact that both $|d\varsigma_-|$ and $|d\Xi|$ are bounded by $c_0e^{-r/c_0}$ where $\mathfrak{d} > \frac{1}{2}r_0$ leads directly (use the Dirichlet Green's function for $-4\partial\bar\partial$) to a $c_0e^{-r/c_0}$ bound for the norm of the derivative of $(\varsigma_- - \Xi)$ where $|z| \le \frac{1}{2}r_0$.

*Part 4*: This last part of the proof addresses the claim made by Lemma 5.2's second bullet. The discussion that follows focuses on the case when $p \in \Theta_+$ because the $p \in \Theta_-$ case differs only cosmetically. With $p$ chosen, then the equation in (5.12) for the function $\varsigma_-$ can be written using the associated complex coordinate on the $|z| < 2r_0$ disk for $p$ is this:

$$-4\partial\bar\partial(\varsigma_- + \varsigma_+) + \alpha\,(e^{2v}|z|^2 + e^{-2v})(\varsigma_- + \varsigma_+) = 2\alpha\,e^{2v}|z|^2\,\varsigma_+ \ .$$

(A.56)

Hold onto this equation for the moment for a brief digression.

To start the digression: Supposing that $\varsigma_+$ has the Laurent series as depicted in the second bullet of Lemma 5.2 (which is convergent where $0 < |z| < 2r_0$), denote by $Q$ the Laurent series below

$$Q \equiv \frac{3a_{-1}}{z} + \sum_{k \ge 0} \frac{3}{3+2k}a_k\,z^k \ .$$

(A.57)

This Laurent series is such that $Q\Xi$ obeys the $v_f$ analog of (A.56) on the $|z| < 2r_0$ disk:

$$-4\partial\bar\partial(Q\Xi + \varsigma_+) + \alpha\,(e^{2v_f}\,|z|^2 + e^{-2v_f})(Q\Xi + \varsigma_+) = 2\alpha\,e^{2v_f}\,|z|^2\,\varsigma_+ \ .$$

(A.58)

(To prove this, use (7.9) to rewrite the contribution to $\partial\bar\partial Q\Xi$ from $\partial Q\bar\partial\Xi$.) Since $|v - v_f|$ and $|d(v - v_f)|$ are both bounded by $c_0e^{-r/c_0}$, the equations in (A.56) and (A.58) lead to the equation below for $\varsigma_- - Q\Xi$:

$$-4\partial\bar\partial(\varsigma_- - Q\Xi) + \alpha\,(e^{2v}\,|z|^2 + e^{-2v})(\varsigma_- - Q\Xi) = \epsilon' \ ,$$

(A.59)

where $\epsilon'$ denotes a complex function on the $|z| < r_0$ disk whose norm and that if its derivative obey $|\epsilon'| + |d\epsilon'| \le c_0e^{-r/c_0}$.

Granted the equation in (A.59), and granted these $c_0e^{-r/c_0}$ bounds for $|\epsilon'|$ and $|d\epsilon'|$, and granted the $c_0e^{-r/c_0}$ bounds in the third and fourth bullets of Lemma 5.2, then an almost verbatim repetition of the arguments from Part 3 will establish the desired $c_0 e^{-r/c_0}$ bounds for the norms of both $\varsigma_- - Q\Xi$ and $d(\varsigma_- - Q\Xi)$.